\documentclass[a4paper,11pt,english]{report}

\usepackage[latin1]{inputenc}

\usepackage{graphicx,palatino,verbatim}

\usepackage{amsmath}

\usepackage[T1]{fontenc}
\usepackage[english]{babel}
\usepackage{newlfont}

\usepackage{amsthm}

\usepackage{amssymb}

\usepackage{amsmath}

\usepackage{eufrak}

\usepackage{latexsym}

\usepackage{young}
\usepackage[mathscr]{eucal}

\usepackage{amsfonts}
\usepackage{syntonly}
\usepackage{graphicx}

\input xy
\xyoption{all}
\newtheorem{teo}{Theorem}

\newtheorem{teorema}{Theorem}[section]

\newtheorem{oss}[teorema]{Remark}

\newtheorem{lemma}[teorema]{Lemma}
\newtheorem{definizione}[teorema]{Definition}
\newtheorem{esempi}[teorema]{Example}

\newtheorem{proposizione}[teorema]{Proposition}

\newtheorem{cor}[teorema]{Corollary}
\newtheorem{con}[teorema]{Conjecture}

\begin{document}

\begin{center}
\Huge{\textbf{Semi-invariants of Symmetric Quivers}}
\end{center}

\begin{center}
Riccardo Aragona
\end{center}
\begin{center}
Università degli Studi di Roma \textquotedblleft
Tor Vergata\textquotedblright\\
Dipartimento di Matematica\\
Via della Ricerca Scientifica 1, 00133 Rome (Italy)
\end{center}

\begin{center}
E-mail: $\quad$ ric$\_$aragona@yahoo.it
\end{center}\begin{center}
\textbf{Abstract}
\end{center}
A symmetric quiver $(Q,\sigma)$ is a 
 finite quiver without oriented cycles $Q=(Q_0,Q_1)$ equipped with a contravariant involution $\sigma$ on $Q_0\sqcup Q_1$. The
involution allows us to define a nondegenerate bilinear form $<,>$
on a representation $V$ of $Q$. We shall say that $V$ is orthogonal
if $<,>$ is symmetric and symplectic if $<,>$ is
skew-symmetric. Moreover we define an action of products of
classical groups on the space of orthogonal representations and on
the space of symplectic representations. So we prove that if $(Q,\sigma)$ is a symmetric quiver of finite type or of tame type then the rings of semi-invariants for this
action are spanned by the semi-invariants of determinantal type
$c^V$ and, in the case when matrix defining $c^V$ is
skew-symmetric, by the Pfaffians $pf^V$.

\tableofcontents

\newpage
\begin{center}
\Large{\textbf{Acknowledgements}}
\end{center}
I wish to express deep thanks to Prof. Jerzy Weyman for the
patience with which he followed me during thesis draft and for
giving me the opportunity to work with him, having me as a guest
at
Northeastern University of Boston for two years.\\
I am grateful to Prof. Elisabetta Strickland for the great
helpfulness
she showed me during the three years of PhD and for her precious advise. \\
I wish also thank Prof. Fabio Gavarini for everytime that I have
gone to his office if I had a problem of mathematics or any
other problems.\\
Besides I am grateful to Prof. Alessandro D'Andrea, Dott. Giuseppe
Marchei, Dott. Giovanni Cerulli Irelli and Dott. Cristina Di
Trapano for the numerous and fruitful discussions we had together
and to Prof. Corrado De Concini for giving me important
suggestions before I began to work on my thesis draft.\\
I would also thank Prof. Marialuisa J. de Resmini for being so
close
to me during all years that I spent to University.

\chapter*{Introduction}\addcontentsline{toc}{chapter}{Introduction}
The representations of quivers can be viewed as a formalization of
some linear algebra problems. Symmetric quivers have been
introduced by Derksen and
Weyman in [DW2] to provide similar formalization for other classical groups.\\
In the recent years the quiver representations were used to prove
interesting results related to general linear groups.\\
Derksen and Weyman in [DW1] gave a proof of saturation property
for Littlewood-Richardson coefficients.\\
Magyar, Weyman and Zelevinsky in [MWZ1] classified products of
flag varieties with finitely many orbits under the diagonal action
of general linear groups. We hope that the representations of
symmetric quivers are a tool to solve similar problems for classical groups.\\
Another interesting aspect and direction for future research is
the connection with Cluster algebras (see [FZ1]). Igusa, Orr,
Todorov and Weyman in [IOTW] generalized the semi-invariants of
quivers to virtual representations of quivers. They associated,
via virtual semi-invariants of quivers, a simplicial complex
$\mathcal{T}(Q)$ with each quiver $Q$. In particular, if $Q$ is of
finite type, then the simplices of $\mathcal{T}(Q)$ correspond to
tilting objects in a corresponding Cluster category (defined in
[BMRRT]). It would be interesting to carry out a similar
construction for symmetric quivers of finite type and to relate it
to Cluster algebras (see [FZ2]).\\
The results of this thesis are first steps in this direction. We
describe the ring of semi-invariants for symmetric quivers of
finite and tame type.\\
A symmetric quiver is a pair $(Q,\sigma)$ where $Q$ is a quiver
(called \textit{underlying quiver of} $(Q,\sigma)$) and $\sigma$
is a contravariant involution on the union of the set of arrows
and the set of vertices of $Q$. The involution allows us to define
a nondegenerate bilinear form $<,>$ on a representation $V$ of
$Q$. We call the pair $(V,<,>)$ orthogonal representation
(respectively symplectic) of $(Q,\sigma)$ if $<,>$ is symmetric
(respectively skew-symmetric). We define $SpRep(Q,\beta)$ and
$ORep(Q,\beta)$ to be respectively the space of symplectic
$\beta$-dimensional representations and the space of orthogonal
$\beta$-dimensional representations  of $(Q,\sigma)$. Moreover we
can define an action of a product of classical groups, which we
call $SSp(Q,\beta)$ in the symplectic case and $SO(Q,\beta)$ in
the orthogonal case, on these space. We describe a set of
generators of the ring of semi-invariants of $ORep(Q,\beta)$
   $$
  OSI(Q,\beta)=\mathbb{K}[ORep(Q,\beta)]^{SO(Q,\beta)}=
$$
$$
  \{f\in\mathbb{K}[ORep(Q,\beta)]|g\cdot
  f=f\;\forall g\in SO(Q,\beta)\}
  $$
  and of the ring of semi-invariants of $SpRep(Q,\alpha)$
   $$
  SpSI(Q,\beta)=\mathbb{K}[SpRep(Q,\beta)]^{SSp(Q,\beta)}=
  $$
 $$
  \{f\in\mathbb{K}[SpRep(Q,\beta)]|g\cdot
  f=f\;\forall g\in SSp(Q,\beta)\},
   $$
where  $\mathbb{K}[ORep(Q,\beta)]$ is the ring of polynomial
functions on
  $ORep(Q,\beta)$ and  $\mathbb{K}[SpRep(Q,\beta)]$ is the ring of polynomial functions on
  $SpRep(Q,\beta)$.\\
  Let $(Q,\sigma)$ be a symmetric quiver and $V$ a representation of the underlying quiver $Q$ such that
  $\langle\underline{dim}\,V,\beta\rangle=0$, where $\langle\cdot,\cdot\rangle$ is the Euler form of $Q$. Let
$$
0\longrightarrow P_1\stackrel{d^V}{\longrightarrow}
P_0\longrightarrow V\longrightarrow 0
$$
be the canonical projective resolution of $V$ (see [R1]). We
define the semi-invariant $c^V:=det(Hom_Q(d^V,\cdot))$ of $
OSI(Q,\beta)$ and $ SpSI(Q,\beta)$ (see
  [DW1] and [S]).\\
   Let $\tau$ be the Auslander-Reiten translation functor and let $\nabla$ be the duality functor.
    We will
  prove in the symmetric case the following
  \begin{teo}\label{tp1}
Let $(Q,\sigma)$ be a symmetric quiver of finite type or of tame
type such that the underlying quiver $Q$ is without oriented
cycles and let $\beta$ be a symmetric dimension vector. The ring
$SpSI(Q,\beta)$ is generated by semi-invariants
\begin{itemize}
\item[(i)] $c^V$ if $V\in Rep(Q)$ is such that
$\langle\underline{dim}\,V,\beta\rangle=0$,
\item[(ii)] $pf^V:=\sqrt{c^V}$ if $V\in Rep(Q)$ is such that
$\langle\underline{dim}\,V,\beta\rangle=0$, $\tau V=\nabla V$ and
the almost split sequence $0\rightarrow\nabla V\rightarrow
Z\rightarrow V\rightarrow 0$ has the middle term $Z$ in $ORep(Q)$.
\end{itemize}
\end{teo}
\begin{teo}\label{tp2}
Let $(Q,\sigma)$ be a symmetric quiver of finite type or of tame
type such that the underlying quiver $Q$ is without oriented
cycles and let $\beta$ be a symmetric dimension vector. The ring
$OSI(Q,\beta)$ is generated by semi-invariants
\begin{itemize}
\item[(i)] $c^V$ if $V\in Rep(Q)$ is such that
$\langle\underline{dim}\,V,\beta\rangle=0$,
\item[(ii)] $pf^V:=\sqrt{c^V}$ if $V\in Rep(Q)$ is such that
$\langle\underline{dim}\,V,\beta\rangle=0$, $\tau V=\nabla V$ and
the almost split sequence $0\rightarrow\nabla V\rightarrow
Z\rightarrow V\rightarrow 0$ has the middle term $Z$ in
$SpRep(Q)$.
\end{itemize}
\end{teo}
A similar result has been obtained by Lopatin in [Lo], in which the author considers a different form of the generators of the ring of semi-invariants of a quiver $Q$, using ideas from [LoZ].\\
The strategy of the proofs is the following. First we set the
technique of reflection functors on the symmetric quivers. Then we
prove that we can reduce theorems \ref{tp1} and \ref{tp2}, by this
technique, to particular orientations of the symmetric quivers.
Finally, we check theorems \ref{tp1} and \ref{tp2} for these
orientations.\\
In the first chapter we give general notions and results about
symmetric quivers and their representations. First, we state main
results \ref{tp1} and \ref{tp2}. Next, we adjust to symmetric
quivers the technique of reflection functors  and we describe
particular orientations for every symmetric quiver of finite type
and tame type. Finally, we prove general results about
semi-invariants of symmetric quivers and we check that we can
reduce theorems \ref{tp1} and \ref{tp2} to these particular
orientations.\\
In the second chapter, using classical invariant theory and the
technique of Schur functors, we prove case by case theorems
\ref{tp1} and
\ref{tp2} for symmetric quivers of finite type with the orientations described in chapter 1.\\
In the third chapter we prove theorems \ref{tp1} and \ref{tp2} for
symmetric quivers of tame type with the orientations described in
chapter 1. First, we deal with symplectic and orthogonal
representations of dimension $\beta=ph$, where $p\in\mathbb{N}$
and $h$ is the homogeneous simple regular dimension vector. We
give a proof of theorems \ref{tp1} and \ref{tp2} case by case.
Next, we adjust to symmetric quivers some general results of Dlab
and Ringel about regular representations of tame quivers (see
[DR]) and we describe generic decomposition of dimension vectors
of symplectic and orthogonal representations (see [K1] and [K2]).
Finally, by these results, we describe case by case the ring of
semi-invariants of symmetric quivers of tame type for any regular
dimension vectors.\\
At last, in appendix \textit{A} we recall some results of
representations of general linear group and of invariant theory.
In appendix \textit{B} we recall general definitions and results
about quiver representations and semi-invariants of quivers.

\chapter{Main results}
\section{Symmetric quivers}
Throughout all this section, we use the notation of section
\ref{B1}.
\begin{definizione}
 A symmetric quiver is a pair $(Q,\sigma)$ where $Q$ is a quiver (called the underlying quiver of $(Q,\sigma)$) and
  $\sigma$ is an involution from the disjoint union $Q_0\coprod Q_1$ to itself, such that
\begin{itemize}
   \item[(i)] $\sigma(Q_0)=Q_0$ and $\sigma(Q_1)=Q_1$,
   \item[(ii)] $t\sigma(a)=\sigma(ha)$ and $h\sigma(a)=\sigma(ta)$ for all $a\in
   Q_1$,
   \item[(iii)] $\sigma(a)=a$ whenever $a\in Q_1$ and $\sigma(ta)=ha$.
   \end{itemize}
 \end{definizione}
 \begin{definizione}
 Let $(Q,\sigma)$ be a symmetric quiver and
 $$V=\big\{\{V(x)\}_{x\in Q_0},\{V(a)\}_{a\in Q_1}\big\}$$ be a
 representation of the underlying quiver $Q$. We define the
 duality functor $\nabla:V\rightarrow V^*$ with $V^*=\big\{\{V^*(x)\}_{x\in Q_0},\{V^*(a)\}_{a\in
 Q_1}\big\}$ where $V^*(x):=V(\sigma(x))^*$ for every $x\in Q_0$ and
 $V^*(a):=-V(\sigma(a))^*$ for every $a\in Q_1$. Moreover if $W$ is another representation of $Q$ and $f:V\rightarrow W$ is a
 morphism, then
 $\nabla f:\nabla W\rightarrow \nabla V$ is defined by
 $(\nabla f)(x):=f(\sigma(x))^*:W^*(x)\rightarrow V^*(x)$, for every $x\in Q_0$. We shall call $V$ selfdual if $\nabla V=V$.
 \end{definizione}
 \begin{definizione}\label{defortsp}
 An orthogonal (resp. symplectic) representation of a symmetric quiver $(Q,\sigma)$ is a pair $(V,<\cdot,\cdot >)$, where $V$ is
  a representation of the underlying quiver $Q$ with a nondegenerate symmetric (resp. skew-symmetric) scalar product $<\cdot,\cdot >$ on $\bigoplus_{x\in Q_0}V(x)$
   such that
   \begin{itemize}
   \item[(i)] the restriction of $<\cdot,\cdot>$ to $V(x)\times V(y)$ is 0 if $y\neq\sigma(x)$,
   \item[(ii)] $<V(a)(v),w>+<v,V(\sigma(a))(w)>=0$ for all $v\in V(ta)$ and all $w\in V(\sigma(ha))$.
   \end{itemize}
  \end{definizione}
  By properties \textit{(i)} and \textit{(ii)} of definition \ref{defortsp}, an orthogonal or symplectic representation $(V,<\cdot,\cdot>)$ of a symmetric quiver is
  selfdual.
\begin{definizione}\label{defindspo}
An orthogonal (respectively symplectic) representation is called
indecomposable orthogonal (respectively indecomposable symplectic)
if it cannot be expressed as a direct sum of orthogonal
(respectively symplectic) representations.
\end{definizione}
We denote
 $ Q^{\sigma}_0$
 (respectively $Q^{\sigma}_1$) the set of vertices (respectively
arrows) fixed by $\sigma$. Thus we have partitions
$$
Q_0 = Q^+_ 0 \cup Q^{\sigma}_ 0 \cup Q^-_0
$$
$$
 Q_1 = Q^+_ 1 \cup Q^{\sigma}_ 1 \cup
Q^-_1
$$
 such that $Q^-_0 = \sigma(Q^+_ 0 )$ and $Q^-_1 = \sigma(Q^+_ 1 )$, satisfying:
 \begin{itemize}
\item[i)] $\forall a \in Q^+_ 1$ , either $\{ta, ha\} \subset Q^+_ 0$ or one of the elements in $\{ta, ha\}$ is
in $Q^+_ 0$ while the other is in $Q^{\sigma}_ 0$;
\item[ii)] $\forall x\in Q^+_0$, if $a\in Q_1$ with $ta = x$ or $ha = x$, then $a\in Q^+_ 1 \cup
Q^{\sigma}_1$.
\end{itemize}
  \begin{definizione}\label{delta}
 Let $(Q,\sigma)$ be a symmetric quiver. We define a linear map $\delta:\mathbb{Z}^{Q_0}_{\geq 0}\rightarrow\mathbb{Z}^{Q_0}_{\geq
  0}$ by setting $\{\delta\alpha(i)\}_{i\in Q_0}=\{\alpha(\sigma(i))\}_{i\in Q_0}$ for every dimension vector
  $\alpha$.
   \end{definizione}
    \begin{oss}
  Since $\sigma$ is an involution, also $\delta$ is one.
  \end{oss}
   \begin{oss} If $V$ is a representation of dimension $\alpha$ then
  $\delta\alpha=\underline{dim}(\nabla V)$. In particular if $V$
  is an orthogonal or symplectic representation of $(Q,\sigma)$ of
  dimension $\alpha$, then $\delta\alpha=\alpha$. Such $\alpha$ is
  called \textit{symmetric dimension vector}.
  \end{oss}
\begin{proposizione}
  Let $\delta:\mathbb{Z}^{Q_0}_{\geq 0}\rightarrow\mathbb{Z}^{Q_0}_{\geq
  0}$ as in definition \ref{delta}. If $\alpha$ and $\beta$ are
  dimension vectors, then
 \begin{eqnarray}
  \langle\alpha,\beta\rangle=\langle\delta\beta,\delta\alpha\rangle.
  \end{eqnarray}
  \end{proposizione}
  \textit{Proof.}
  \begin{eqnarray}\label{langle}
  &\langle\alpha,\beta\rangle=
  \sum_{i\in Q_0^+\cup Q_0^{\sigma}}\alpha(i)\beta(i)+\sum_{i\in Q_0^+}\alpha(\sigma(i))\beta(\sigma(i))\nonumber&\\
  &+\sum_{a\in Q_1^+\cup Q_1^{\sigma}}\alpha(ta)\beta(ha)+\sum_{a\in
  Q_1^+}\alpha(t\sigma(a))\beta(h\sigma(a))&.
  \end{eqnarray}
  By definition of $\sigma$, we have
  $$
  \langle\delta\beta,\delta\alpha\rangle=
  $$
  $$
  \sum_{i\in Q_0^+\cup Q_0^{\sigma}}\beta(\sigma(i))\alpha(\sigma(i))+\sum_{i\in Q_0^+}\beta(\sigma(\sigma(i)))\alpha(\sigma(\sigma(i)))+
  $$
  $$
  \sum_{a\in Q_1^+\cup Q_1^{\sigma}}\beta(\sigma(ta))\alpha(\sigma(ha))+\sum_{a\in Q_1^+}\beta(\sigma(t\sigma(a)))\alpha(\sigma(h\sigma(a)))=
  $$
  $$
  \sum_{i\in Q_0^+}\beta(\sigma(i))\alpha(\sigma(i))+\sum_{i\in Q_0^{\sigma}}\beta(i)\alpha(i)+
  $$
  $$
  \sum_{i\in Q_0^+}\beta(i)\alpha(i)+\sum_{a\in Q_1^+}\beta(h\sigma(a))\alpha(t\sigma(a))+
  $$
  $$
  \sum_{a\in  Q_1^{\sigma}}\beta(h\sigma(a))\alpha(t\sigma(a))+\sum_{a\in Q_1^+}\beta(\sigma^2(ha))\alpha(\sigma^2(t\sigma(a))
  $$
  which is the right hand side of
  (\ref{langle}), recalling that $\sigma$ is an involution. $\Box$\\

The space of orthogonal $\alpha$-dimensional
  representations of a symmetric quiver $(Q,\sigma)$ can be
  identified with
  \begin{equation}
   ORep(Q,\alpha)=\bigoplus_{a\in Q^{+}_1 }Hom(\mathbb{K}^{\alpha(ta)},\mathbb{K}^{\alpha(ha)})\oplus\bigoplus_{a\in
Q^{\sigma}_1}\bigwedge^2(\mathbb{K}^{\alpha(ta)})^*.
  \end{equation}
  The space of symplectic $\alpha$-dimensional representations
  can be identified with
 \begin{equation}
 SpRep(Q,\alpha)=\bigoplus_{a\in Q^+_1}Hom(\mathbb{K}^{\alpha(ta)},\mathbb{K}^{\alpha(ha)})\oplus\bigoplus_{a\in Q^{\sigma}_1 }S_2(\mathbb{K}^{\alpha(ta)})^*.
   \end{equation}
  We define the group
 \begin{equation}
   O(Q,\alpha)=\prod_{x\in Q^+_0}GL(\mathbb{K},\alpha(x))\times\prod_{x\in Q^{\sigma}_0}O(\mathbb{K},\alpha(x))
   \end{equation}
  and the subgroup
 \begin{equation}
 SO(Q,\alpha)=\prod_{x\in
Q^+_0}SL(\mathbb{K},\alpha(x))\times\prod_{x\in
Q^{\sigma}_0}SO(\mathbb{K},\alpha(x)).
   \end{equation}
  Here $O(\mathbb{K},\alpha(x))$ is the group of orthogonal
  transformations for the symmetric form
  $<\cdot,\cdot>$ restricted to $V(x)$.\\
  \\
  Assuming that $\alpha(x)$ is even
  for every $x\in Q_0^{\sigma}$, we define the group
\begin{equation}
Sp(Q,\alpha)=\prod_{x\in Q^+_0}GL(\mathbb{K},\alpha(x))\times\prod_{x\in Q^{\sigma}_0}Sp(\mathbb{K},\alpha(x))
   \end{equation}
  and the subgroup
 \begin{equation}
   SSp(Q,\alpha)=\prod_{x\in Q^+_0}SL(\mathbb{K},\alpha(x))\times\prod_{x\in
   Q^{\sigma}_0}Sp(\mathbb{K},\alpha(x)).
   \end{equation}
Here $Sp(\mathbb{K},\alpha(x))$ is the group of isometric
  transformations for the skew-symmetric form
  $<\cdot,\cdot>$ restricted to $V(x)$.\\
  \\
The action of these groups is defined by
$$
g\cdot
V=\{g_{ha}V(a)g_{ta^{-1}}\}_{a\in Q_1^+\cup Q_1^{\delta}}
$$
where
$g=(g_x)_{x\in Q_0}\in O(Q,\alpha)$ (respectively
  $g\in Sp(Q,\alpha)$) and $V\in ORep(Q,\alpha)$ (respectively in $SpRep(Q,\alpha)$). In particular we can suppose $g_{\sigma(x)}=(g_x^{-1})^t$
  for every $x\in Q_0$.
  \begin{esempi}
  \begin{itemize}
  \item[(1)] Consider the symmetric quiver $(Q,\sigma)$
  $$
  \begin{array}{ccc}
  \circ\rightarrow\bullet\rightarrow\circ
  \end{array}
  $$
  where $\sigma$ interchanges the antipodal nodes and fixes the
  closed node. An orthogonal representation of $(Q,\sigma)$ is a
  quadruple $(V_1,V_2,\phi,\langle\cdot,\cdot\rangle)$ where $V_1$ and $V_2$ are
   vector spaces, $\phi:V_1\rightarrow V_2$ is a linear map and $\langle\cdot,\cdot\rangle$ is a non-degenerate symmetric bilinear form on $V_2$. We also have the dual map
  $-\phi^*:V_2^*\cong V_2\rightarrow V_1^*$ and so we have the
  following diagram:
  $$
  \begin{array}{ccc}
  V_1\stackrel{\phi}{\rightarrow} V_2\stackrel{-\phi^*}{\rightarrow}
  V_1^*.
  \end{array}
  $$
  Hence the isomorphism classes of orthogonal representations of $(Q,\sigma)$ are the
  $GL(V_1)\times O(V_2)$-orbits in $Hom(V_1,V_2)$.
  \item[(2)] Consider the symmetric quiver $(Q,\sigma)$
 $$
  \begin{array}{ccc}
  \circ\rightarrow\circ\rightarrow\circ\rightarrow\circ
  \end{array}
  $$
  where $\sigma$ sends the first vertex to the last one and the
  second one to the third one. A symplectic representation of $(Q,\sigma)$
  is a quadriple $(V_1,V_2,\phi,\psi)$ where $V_1$ and $V_2$ are vector
  spaces, $\phi:V_1\rightarrow V_2$ is linear map and $\psi\in S_2V_2^*$. We also
  have the dual map $-\phi^*: V_2^*\rightarrow V_1^*$. We consider the following
  diagram:
   $$
  \begin{array}{ccc}
  V_1\stackrel{\phi}{\rightarrow}V_2\stackrel{\psi}{\rightarrow} V_2^*\stackrel{-\phi^*}{\rightarrow}
  V_1^*.
  \end{array}
  $$
  Hence the isomorphism classes of symplectic representations of $(Q,\sigma)$ are the
  $GL(V_1)\times GL(V_2)$-orbits in $Hom(V_1,V_2)\oplus S_2V_2^*$.

  \end{itemize}
  \end{esempi}
  \begin{definizione}
  \begin{itemize}
\item[(i)] Let $\mathbb{K}[ORep(Q,\alpha)]$ be the ring of polynomial
functions on $ORep(Q,\alpha)$.
   $$
  OSI(Q,\alpha)=\mathbb{K}[ORep(Q,\alpha)]^{SO(Q,\alpha)}=
$$
\begin{equation}
  \{f\in\mathbb{K}[ORep(Q,\alpha)]|g\cdot
  f=f\;\forall g\in SO(Q,\alpha)\}
   \end{equation}
  is the ring of orthogonal semi-invariants of $(Q,\alpha)$.
  \item[(ii)] Let
  $\mathbb{K}[SpRep(Q,\alpha)]$ be the ring of polynomial functions on
  $SpRep(Q,\alpha)$,
   $$
  SpSI(Q,\alpha)=\mathbb{K}[SpRep(Q,\alpha)]^{SSp(Q,\alpha)}=
  $$
  \begin{equation}
  \{f\in\mathbb{K}[SpRep(Q,\alpha)]|g\cdot
  f=f\;\forall g\in SSp(Q,\alpha)\}
   \end{equation}
  is the ring of symplectic semi-invariants of $(Q,\alpha)$.
  \end{itemize}
  \end{definizione}
\subsection{Symmetric quivers of finite type}
  \begin{definizione}
A symmetric quiver is said to be of finite representation type if
it has only finitely many indecomposable orthogonal (resp.
symplectic) representations up to isomorphisms.
\end{definizione}
 We recall the
following theorem proved by Derksen and Weyman in [DW2]
\begin{teorema}\label{ctf}
A symmetric quiver $(Q,\sigma)$ is  of finite type if and only if
the underlying quiver $Q$ is of type $A_n$.
\end{teorema}
\textit{Proof.} See [DW2, theorem 3.1 and proposition 3.3]
$\Box$\\
\subsection{Symmetric quivers of tame type}
\begin{definizione}
A symmetric quiver is said to be of tame representation type if is
not of finite representation type, but in every dimension vector
the indecomposable orthogonal (symplectic) representations occur
in families of dimension $\leq 1$.
\end{definizione}
\begin{teorema}
A symmetric quiver $(Q,\sigma)$ with $Q$ connected is tame if and
only if the underlying quiver $Q$ is an extended Dynkin quiver.
\end{teorema}
\textit{Proof.} See [DW2, theorem 4.1]. $\Box$\\
\\
One can classify the symmetric tame quivers with connected
underlying quiver.
\begin{proposizione}\label{ctt}
Let $(Q,\sigma)$ be a symmetric tame quiver with $Q$ connected.
Then $(Q,\sigma)$ is one of the following symmetric quivers.
\begin{itemize}
\item[(1)] Of type $\widetilde{A}^{2,0,1}_n$:
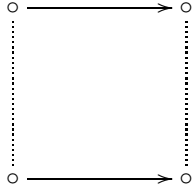
$$
\xymatrix{
\circ\ar[rr]\ar@{.}[dd]&&\circ\ar@{.}[dd]\\
&&&\\
 \circ\ar[rr]&&\circ}
$$
with arbitrary orientation reversed under $\sigma$ if
$Q=\widetilde{A}_{2n+1}$ ($\geq 1$). Here $\sigma$ is a reflection
with respect to a central vertical line (so $\sigma$ fixes two
arrows and no vertices).
\item[(2)] Of type $\widetilde{A}^{2,0,2}_n$:
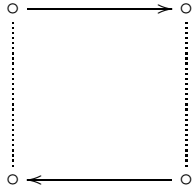
$$
\xymatrix{
\circ\ar[rr]\ar@{.}[dd]&&\circ\ar@{.}[dd]\\
&&&\\
 \circ&&\ar[ll]\circ}
$$
with arbitrary orientation reversed under $\sigma$ if
$Q=\widetilde{A}_{2n+1}$ ($\geq 1$). Here $\sigma$ is a reflection
with respect to a central vertical line (so $\sigma$ fixes two
arrows and no vertices).

\item[(3)] Of type $\widetilde{A}^{0,2}_n$:
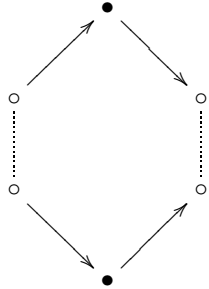
$$
\xymatrix{&\bullet\ar[dr]&\\
\circ\ar[ur]\ar@{.}[d]&&\circ\ar@{.}[d]\\
\circ\ar[dr]&&\circ\\
&\bullet\ar[ur]&}
$$
with arbitrary orientation reversed under $\sigma$ if
$Q=\widetilde{A}_{2n-1}$ ($n\leq 1$). Here $\sigma$ is a
reflection with respect to a central vertical line (so $\sigma$
fixes two vertices and no arrows).
\item[(4)] Of type $\widetilde{A}^{1,1}_n$:
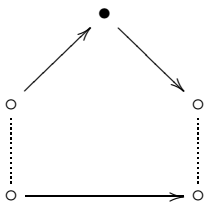
$$
\xymatrix{&\bullet\ar[dr]&\\
\circ\ar[ur]\ar@{.}[d]&&\circ\ar@{.}[d]\\
\circ\ar[rr]&&\circ}
$$
with arbitrary orientation reversed under $\sigma$ if
$Q=\widetilde{A}_{2n}$ ($n\geq 1$). Here $\sigma$ is a reflection
with respect to a central vertical line (so $\sigma$ fixes one
arrow and one vertex).
\item[(5)] Of type $\widetilde{A}^{0,0}_n$:
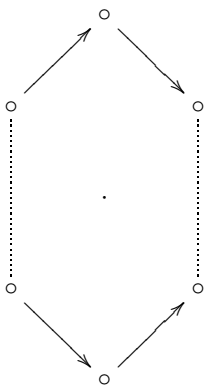
$$
\xymatrix{&\circ\ar[dr]&\\
\circ\ar[ur]\ar@{.}[dd]&&\circ\ar@{.}[dd]\\
&\cdot&\\
\circ\ar[dr]&&\circ\\
&\circ\ar[ur]&}
$$
with arbitrary orientation reversed under $\sigma$ if
$Q=\widetilde{A}_{2n+1}$ ($n\geq 1$). Here $\sigma$ is a central
symmetry (so $\sigma$ fixes neither arrows nor vertices).
\item[(6)] Of type $\widetilde{D}^{1,0}_n$
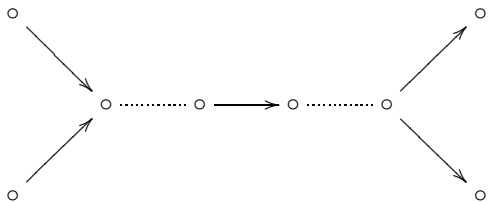
$$
\xymatrix{\circ\ar[dr]&&&&&\circ\\
&\circ\ar@{.}[r]&\circ\ar[r]&\circ\ar@{.}[r]&\circ\ar[ur]\ar[dr]&\\
\circ\ar[ur]&&&&&\circ}
$$
 with arbitrary orientation reversed
under $\sigma$ if $Q=\widetilde{D}_{2n}$ ($n\geq 2$). Here
$\sigma$ is a reflection with respect to a central vertical line
(so $\sigma$ fixes one arrow and no vertices).
\item[(7)] Of type $\widetilde{D}^{0,1}_n$
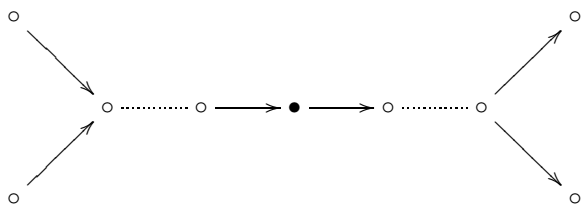
$$
\xymatrix{\circ\ar[dr]&&&&&&\circ\\
&\circ\ar@{.}[r]&\circ\ar[r]&\bullet\ar[r]&\circ\ar@{.}[r]&\circ\ar[ur]\ar[dr]&\\
\circ\ar[ur]&&&&&&\circ}
$$
 with arbitrary orientation reversed
under $\sigma$ if $Q=\widetilde{D}_{2n-1}$ ($n\geq 2$). Here
$\sigma$ is a reflection with respect to a central vertical line
(so $\sigma$ fixes one vertex and no arrows).

\end{itemize}
\end{proposizione}
\textit{Proof.} See [DW2, proposition 4.3]. $\Box$
\section{The main results}
In this thesis we describe the rings of semi-invariants of
symmetric quivers in the finite type and in the tame cases. We
also conjecture in general the following results. Below we use the
notations of section \ref{B3} and we conjecture the following
theorems
\begin{con}\label{mt1}
Let $(Q,\sigma)$ a symmetric quiver such that the underlying
quiver $Q$ is without oriented cycles and let $\beta$ be a
symmetric dimension vector. The ring $SpSI(Q,\beta)$ is generated
by semi-invariants
\begin{itemize}
\item[(i)] $c^V$ if $V\in Rep(Q)$ is such that
$\langle\underline{dim}\,V,\beta\rangle=0$,
\item[(ii)] $pf^V:=\sqrt{c^V}$ if $V\in Rep(Q)$ is such that
$\langle\underline{dim}\,V,\beta\rangle=0$, $V=\tau^-\nabla V$ and
the almost split sequence $0\rightarrow\nabla V\rightarrow
Z\rightarrow V\rightarrow 0$ has the middle term $Z$ in $ORep(Q)$.
\end{itemize}
\end{con}
\begin{con}\label{mt2}
Let $(Q,\sigma)$ a symmetric quiver such that the underlying
quiver $Q$ is without oriented cycles and let $\beta$ be a
symmetric dimension vector. The ring $OSI(Q,\beta)$ is generated
by semi-invariants
\begin{itemize}
\item[(i)] $c^V$ if $V\in Rep(Q)$ is such that
$\langle\underline{dim}\,V,\beta\rangle=0$,
\item[(ii)] $pf^V:=\sqrt{c^V}$ if $V\in Rep(Q)$ is such that
$\langle\underline{dim}\,V,\beta\rangle=0$, $V=\tau^-\nabla V$ and
the almost split sequence $0\rightarrow\nabla V\rightarrow
Z\rightarrow V\rightarrow 0$ has the middle term $Z$ in
$SpRep(Q)$.
\end{itemize}
\end{con}
We prove these conjectures for symmetric quivers of finite type
(chapter 2) and for symmetric quivers of tame type and regular
dimension vectors $\beta$ (chapter 3).\\
We use the following strategy. First we adjust to symmetric
quivers the technique of reflection functors. Next we prove with
this technique that we can reduce the conjectures \ref{mt1} and
\ref{mt2} to a particular orientation of the quiver. Then we state
and prove conjectures \ref{mt1} and \ref{mt2} for these
orientations.
\begin{definizione}\label{popssp}
We will say that $V\in Rep(Q)$ satisfies property \textit{(Op)} if
\begin{itemize}
\item[(i)] $V=\tau^-\nabla V$
\item[(ii)] the almost split sequence $0\rightarrow\nabla V\rightarrow
Z\rightarrow V\rightarrow 0$ has the middle term $Z$ in $ORep(Q)$.
\end{itemize}
Similarly we will say that $V\in Rep(Q)$ satisfies property
\textit{(Spp)} if
\begin{itemize}
\item[(i)] $V=\tau^-\nabla V$
\item[(ii)] the almost split sequence $0\rightarrow\nabla V\rightarrow
Z\rightarrow V\rightarrow 0$ has the middle term $Z$ in
$SpRep(Q)$.
\end{itemize}
\end{definizione}

\section{Reflection functors for symmetric
quivers} In this section we describe the technique of reflection
functors for the symmetric quivers.
\subsection{Admissible sink-source pairs}
We use the notation of section \ref{B2}.
  \begin{definizione}
  Let $(Q,\sigma)$ be a symmetric quiver. A sink (respectively
  source) $x\in Q_0$ is called admissible if there are no arrows connecting $x$ and $\sigma(x)$.
  \end{definizione}
  By definition of $\sigma$, $x$ is a sink (respectively a
  source) if and only if $\sigma(x)$ is a source, so we can define
  the quiver $c_{\sigma(x)}c_x(Q)$. We shall call
$(x,\sigma(x))$ \textit{the admissible sink-source pair}. The
corresponding reflection is denoted by
$c_{(x,\sigma(x))}:=c_{\sigma(x)}c_x$.
\begin{lemma}
  If $(Q,\sigma)$ is a symmetric quiver and $x$ is an admissible sink or source, then $(c_{(x,\sigma(x))}(Q),\sigma)$
   is symmetric.
  \end{lemma}
  \textit{Proof.}
  Let $x\in Q_0$ be an admissible sink of $(Q,\sigma)$. When we apply $c_{(x,\sigma(x))}$  to $Q$, the only arrows which
we reverse are the arrows connecting to $x$ and those connecting
to $\sigma(x)$. Now in $c_{(x,\sigma(x))}(Q)$, $x$ becomes a
source and $\sigma(x)$ becomes a sink. So if $a$ is an arrow
connecting to $x$ or to $\sigma(x)$ we have
$\sigma(tc_{(x,\sigma(x))}(a))=\sigma(ha)=t\sigma(a)=h\sigma(c_{(x,\sigma(x))}(a))$
and
$\sigma(hc_{(x,\sigma(x))}(a))=\sigma(ta)=h\sigma(a)=t\sigma(c_{(x,\sigma(x))}(a))$.
Hence $c_{(x,\sigma(x))}(Q)$ is a symmetric quiver. One proves
similarly if $x$
   is a source. $\Box$
 \begin{definizione}
 Let $(Q,\sigma)$ be a symmetric quiver. A sequence $x_1,\ldots,x_m$ of vertices of $Q$ is an admissible sequence of sinks (or sources) for admissible sink-source
  pairs if $x_{i+1}$ is a sink such that there are no arrows linking $x_{i+1}$ and $\sigma(x_{i+1})$ in $c_{(x_i,\sigma(x_i))}\cdots  c_{(x_1,\sigma(x_1))}(Q)$ for $i=1,\ldots,m-1$.
 \end{definizione}
 \begin{proposizione}\label{Q=Q'}
 Let $(Q,\sigma)$ and $(Q',\sigma)$ be two symmetric connected quivers, without cycles, with the same underlying graph and such that $Q'$ differs from $Q$
 only by changing the orientation of some arrows. Then there exists a sequence $x_1,\ldots,x_m\in Q_0$ which is an admissible sequence of sinks (or sources)
 for admissible sink-source pairs such that
 $$
 Q'=c_{(x_m,\sigma(x_m))}\cdots  c_{(x_1,\sigma(x_1))}(Q).
 $$
   \end{proposizione}
For the proof of proposition \ref{Q=Q'}, we need a lemma.
\begin{lemma}\label{l}
If $(Q,\sigma)$ is a symmetric quiver with
$|\{x\rightarrow\sigma(x)|x\in Q_0\}|>1$, then $(Q,\sigma)$ has
cycles or it is not connected.
\end{lemma}
\textit{Proof of lemma \ref{l}.} If there are more than one arrow
$x\rightarrow\sigma(x)$ for the same $x$ in $Q$ then $Q$ has
cycles. Otherwise we suppose that $Q$ is connected and that there
are two arrows  $x\stackrel{a}{\rightarrow}\sigma(x)$ and
$y\stackrel{b}{\rightarrow}\sigma(y)$, with $x\ne y$ in $Q$. Since
$Q$ is connected, these two arrows have to be linked with a
sequence of other arrows (this regarding their orientation). If
there exists a sequence of arrows $a_1,\ldots,a_t$ from $x$ to $y$
then, by definition of $\sigma$, there exists a sequence of arrows
$\sigma(a_1),\ldots,\sigma(a_t)$ from $\sigma(y)$ to $\sigma(x)$,
reversed respect to $a_1,\ldots,a_t$. So $a_1\cdots a_t
a\sigma(a_t)\cdots\sigma(a_1)b$ is a cycle.
 By a similar reasoning for the other possible three links between  $x\rightarrow\sigma(x)$ and   $y\rightarrow\sigma(y)$ (from $x$ to $\sigma(y)$, from $y$ to $\sigma(x)$ and from $\sigma(x)$ to $\sigma(y)$), we obtain the$\;$ same$\;$ conclusion. $\Box$\\
\\
\textit{Proof of proposition \ref{Q=Q'}.}  By lemma \ref{l} we can
suppose that $Q$ has at most one arrow $x\rightarrow\sigma(x)$ for
some $x\in Q_0$.  First of all we notice that the underlying graph
of $Q$ and $Q'$, being a connected graph without cycles, is a
tree, i.e. a graph where every vertex $x$ has one parent and a
several of children each connected by one edge to the vertex $x$.
We define ancestor and descendants in obvious way and we call
$x\in Q_0$ a vertex without children if there is only one edge
connected to $x$. Let $S$ be a set of vertices without children in $Q$.\\
We observe, by definition of $\sigma$, that if $Q\neq A_2$, in
that case there are no admissible sink or source, and if $S$
contains $x\in Q_0$ then it contains $\sigma(x)$.\\
Observe that, using reflection of the admissible sink-source pair
at the vertex without children $x$, we can change arbitrarily
orientation of arrow connected to $x$ and so of the arrow
connected to
$\sigma(x)$.\\
We proceed by induction on the number $m$ of generations in the
tree. If the number of generations is one, each vertex but one is
without children, applying reflection at the admissible
sink-source pairs we can pass from orientation of $Q$ to
orientation of $Q'$, by which we
observed before.\\
Assume proposition true for the trees with $m-1$ generations. We
remove all vertices without children from $Q$ and $Q'$, so the
resulting quivers $\tilde{Q}$ and $\tilde{Q}'$ have $m-1$
generations and are symmetric. By inductive assumption, we can go
from $\tilde{Q}$ to $\tilde{Q}'$ by a sequence of
reflections at admissible sink-source pairs.\\
To pass from $Q$ to $Q'$ we use the same sequence of reflections
at each point, adjusting the orientations of arrows incident to
$S$, to get the next admissible sink-source pair if necessarily. $\Box$\\
\\
We prove some results on orientations of symmetric quivers of tame
type. The underlying graph of $\widetilde{D}$ is a tree, so by
proposition \ref{Q=Q'}, we will consider a particular orientation
of $\widetilde{D}$
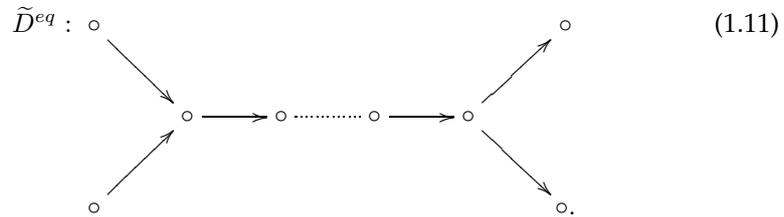
\begin{equation}\label{Dequi}
\widetilde{D}^{eq}:\xymatrix {\circ\ar[dr] & & & & & \circ\\
& \circ\ar[r]& \circ\ar@{.}[r]&\circ\ar[r]&\circ\ar[ur]\ar[dr] &\\
\circ\ar[ur] & & & & & \circ.}
\end{equation}
Applying a compositions of reflections at admissible sink-source
pairs we can get
any orientation of $\widetilde{D}$ from $\widetilde{D}^{eq}$.\\
Now we deal with orientation of symmetric quivers with underlying
quiver of type $\widetilde{A}$. First we prove lemma about
possible exchange of orientation of a quiver $Q$ of type $A_n$,
that does not involve reflections at the end points of $Q$. We
denote vertices of $Q$ with $\{1,\ldots,n\}$ from left to right.
\begin{lemma}\label{Afix}
Let
$$
Q: \def\objectstyle{\scriptstyle}\xymatrix @-1pc{& & \circ \ar[dr]
\ar[dl] & & & & & &
\circ \ar[dr] \ar[dl] &  &  &  &  &   \\
&  \circ\ar@{.}[dl] & & \circ\ar@{.}[dr] & & & & \circ\ar@{.}[dl] & &\circ\ar@{.}[dr] & & & & \circ\ar@{.}[dl]\\
 \circ & & & & \circ\ar[dr] & & \circ  \ar[dl] & & & & \circ\ar[dr]& & \circ  \ar[dl]& \\
 & & & & & \circ & & & & & & \circ & &  ,}
 $$
with $k$ south-west arrows and $h$ south-east arrows. Then there
exists a sequence of admissible sinks $x_1,\ldots,x_l$ with
$x_i\neq 1,n$ for every $i\in\{1,\ldots,l\}$, such that
$c_{x_1}\cdots c_{x_l}Q$ is
$$
Q':\def\objectstyle{\scriptstyle}
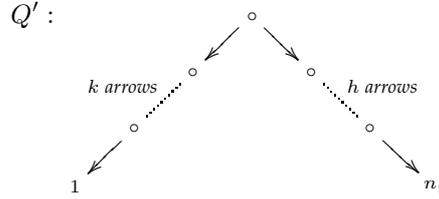
\xymatrix @-1pc{ &  & & \circ
\ar[dr] \ar[dl] & & &   \\
& & \circ\ar@{.}[dl]_{k\;\textrm{arrows}} &  &\circ \ar@{.}[dr]^{h\;\textrm{arrows}} & & \\
 & \circ\ar[dl] & & & &  \circ  \ar[dr] &  \\
 1 & & & & & &  n ,}
$$
i.e. $Q'$ has $1,n$ as only sinks, with $k$ south-west arrows and
with $h$ south-east arrows.
\end{lemma}
\textit{Proof.} Let $x$ and $y$ be two sinks closest to 1.
$$
Q:\def\objectstyle{\scriptstyle}
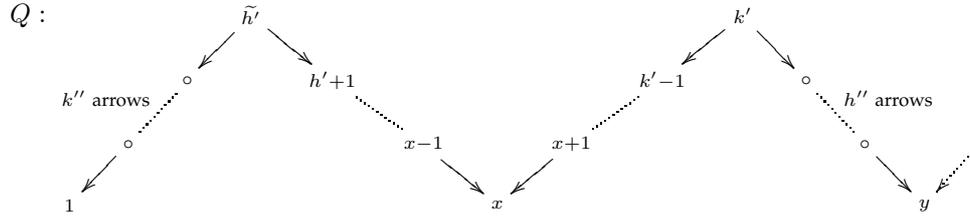
\xymatrix @-1pc{
& & & \widetilde{h'}\ar[dl]\ar[dr] & & & & & & k'\ar[dl]\ar[dr] & & &&\\
& & \circ\ar@{.}[dl]_{k''\;\textrm{arrows}} & & h'+1\ar@{.}[dr] & & & & k'-1\ar@{.}[dl] & & \circ\ar@{.}[dr]^{h''\;\textrm{arrows}} & &&\\
&\circ\ar[dl] & & & & x-1\ar[dr] & & x+1\ar[dl] & & & & \circ\ar[dr]&&\ar@{.>}[dl]\\
1 & & & & & & x & & & & & &  y&}
$$
From 1 to $y$, $Q$ has $k'+k''$ south-west arrows and $h'+h''$
south-east arrows. We remove $x$ by applying only reflections at
vertices with number smaller than $y$, as follows. We suppose
$k'\geq h'$ (the other case is similar). Applying $c_x$ we get
$$
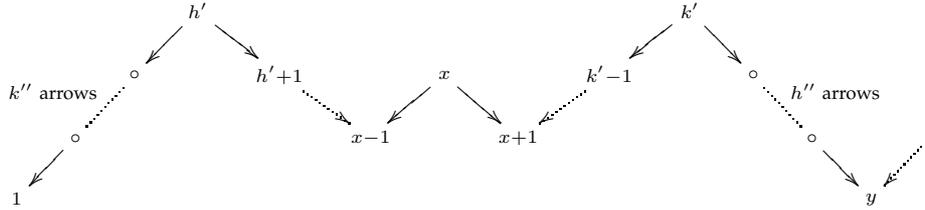
\def\objectstyle{\scriptstyle}\xymatrix @-1pc{
& & & h'\ar[dl]\ar[dr] & & & & & & k'\ar[dl]\ar[dr] & & &&\\
& & \circ\ar@{.}[dl]_{k''\;\textrm{arrows}} & & h'+1\ar@{.>}[dr] && x \ar[dl]\ar[dr]&  & k'-1\ar@{.>}[dl] & & \circ\ar@{.}[dr]^{h''\;\textrm{arrows}} & &&\\
&\circ\ar[dl] & & & & x-1 &  & x+1 &  & & & \circ\ar[dr]&&\ar@{.>}[dl]\\
1 & & & & & &  & & & & & &  y&.}
$$
Now we can apply $c_{x-1}c_{x+1}$ and so on we obtain
$$
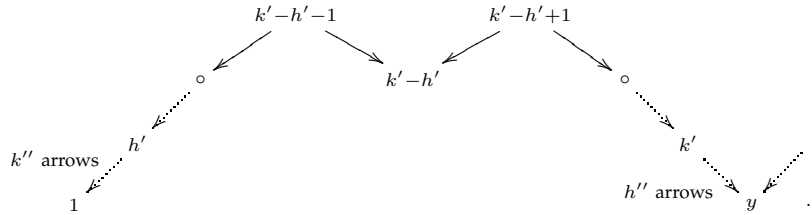
\def\objectstyle{\scriptstyle}\xymatrix @-1pc{
& & & k'-h'-1\ar[dl]\ar[dr]& & k'-h'+1\ar[dl]\ar[dr]& & &&\\
& & \circ\ar@{.>}[dl]& & k'-h' & & \circ\ar@{.>}[dr]& & &\\
& h'\ar@{.>}[dl]_{k''\;\textrm{arrows}}& & & & & &
k'\ar@{.>}[dr]_{h''\;\textrm{arrows}} &&\ar@{.>}[dl]\\
1 & & & & & & & & y&.}
$$
Finally, applying $c_{k'-h'}$ we get
$$
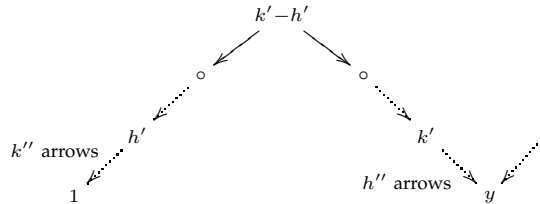
\def\objectstyle{\scriptstyle}\xymatrix @-1pc{ &  & & k'-h'
\ar[dr] \ar[dl] & & & &  \\
& & \circ\ar@{.>}[dl] &  &\circ \ar@{.>}[dr] & && \\
 & h'\ar@{.>}[dl]_{k''\;\textrm{arrows}} & & & &  k' \ar@{.>}[dr]_{h''\;\textrm{arrows}} & &\ar@{.>}[dl] \\
 1 & & & & & &  y&}
$$
in which there are $(k'-h')+h'+k''=k'+k''$ south-west arrows and
$k'-(k'-h')+h''=h'+h''$ south-east arrows. Removing internal sinks
in this way proves lemma. $\Box$
\begin{definizione}
We will say that a symmetric quiver is of type $(s,t,k,l)$ if
\begin{itemize}
\item[(i)] it is of type $\widetilde{A}$,
\item[(ii)] $|Q_1^{\sigma}|=s$ and $|Q_0^{\sigma}|=t$,
\item[(iii)] it has $k$ counterclockwise arrows and $l$ clockwise arrows in $Q_1^+\sqcup
Q_1^-$.
\end{itemize}
\end{definizione}
By proposition \ref{ctt}, $s,t\in\{0,1,2\}$ and if either $s$ or
$t$ are not zero, then $s+t=2$. Moreover, by symmetry, we note
that $k$ and $l$ have to be even.
\begin{proposizione}\label{oA}
Let $(Q,\sigma)$ be a symmetric quiver of type $\widetilde{A}$
such that $Q$ is without oriented cycles. Then there is an
admissible sequence of sinks $x_1,\ldots,x_s$ of $Q$ for
admissible sink-source pairs such that
$c_{(x_1,\sigma(x_1))}\cdots c_{(x_s,\sigma(x_s))}Q$ is one of the
quivers:
\begin{itemize}
\item[(1)]
$$
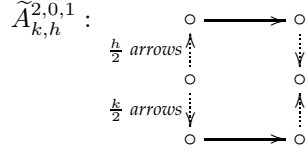
\widetilde{A}^{2,0,1}_{k,h}:\xymatrix@-1pc{
\circ\ar[rr]&&\circ\ar@{.>}[d]\\
\circ\ar@{.>}[d]_{\frac{k}{2}\;\textrm{arrows}}\ar@{.>}[u]^{\frac{h}{2}\;\textrm{arrows}}&&\circ\\
 \circ\ar[rr]&&\circ\ar@{.>}[u],}
$$
and
\item[(2)]
$$
\widetilde{A}^{2,0,2}_{k,h}:
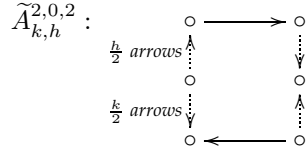
\xymatrix@-1pc{
\circ\ar[rr]&&\circ\ar@{.>}[d]\\
\circ\ar@{.>}[d]_{\frac{k}{2}\;\textrm{arrows}}\ar@{.>}[u]^{\frac{h}{2}\;\textrm{arrows}}&&\circ\\
 \circ&&\ar[ll]\circ\ar@{.>}[u],}
$$
if $(Q,\sigma)$ is of type $(2,0,k,l)$;
\item[(3)]
$$
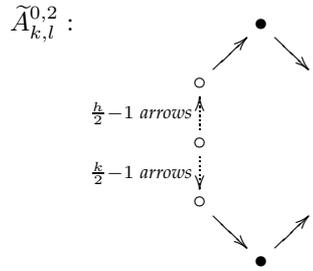
\widetilde{A}^{0,2}_{k,l}:\xymatrix@-1pc{&\bullet\ar[dr]&\\
\circ\ar[ur]&&\circ\ar@{.>}[d]\\
\circ\ar@{.>}[d]_{\frac{k}{2}-1\;\textrm{arrows}}\ar@{.>}[u]^{\frac{h}{2}-1\;\textrm{arrows}}&&\circ\\
\circ\ar[dr]&&\circ\ar@{.>}[u]\\
&\bullet\ar[ur]&,}
$$
if $(Q,\sigma)$ is of type $(0,2,k,l)$;
\item[(4)]
$$
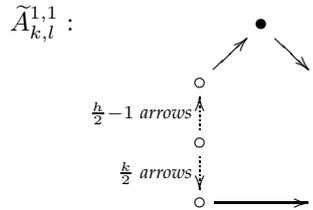
\widetilde{A}^{1,1}_{k,l}:\xymatrix@-1pc{&\bullet\ar[dr]&\\
\circ\ar[ur]&&\circ\ar@{.>}[d]\\
\circ\ar@{.>}[u]^{\frac{h}{2}-1\;\textrm{arrows}}\ar@{.>}[d]_{\frac{k}{2}\;\textrm{arrows}}&&\circ\\
 \circ\ar[rr]&&\circ\ar@{.>}[u],}
$$
if $(Q,\sigma)$ is of type $(1,1,k,l)$;
\item[(5)]
$$
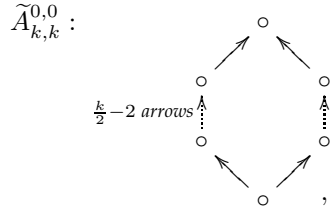
\widetilde{A}^{0,0}_{k,k}:\xymatrix@-1pc{&\circ&\\
\circ\ar[ur]&&\circ\ar[ul]\\
\circ\ar@{.>}[u]^{\frac{k}{2}-2\;\textrm{arrows}}&&\circ\ar@{.>}[u]\\
&\circ\ar[ul]\ar[ur]&,}
$$
if $(Q,\sigma)$ if of type $(0,0,k,k)$.
\end{itemize}
\end{proposizione}
\textit{Proof.} For $(Q,\sigma)$ of types $(2,0,k,l)$, $(0,2,k,l)$
and $(1,1,k,l)$ we apply lemma \ref{Afix} respectively to the
subquivers whose the underlying graphs are
$$
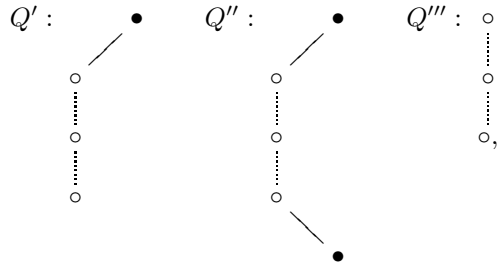
\begin{array}{ccccc}
Q':\xymatrix@-1pc{&\bullet\\
\circ\ar@{-}[ur]&\\
\circ\ar@{.}[u]\ar@{.}[d]&\\
 \circ&} && Q'':\xymatrix@-1pc{&\bullet\\
\circ\ar@{-}[ur]&\\
\circ\ar@{.}[d]\ar@{.}[u]&\\
\circ\ar@{-}[dr]&\\
&\bullet} && Q''':\xymatrix@-1pc{
\circ\\
\circ\ar@{.}[d]\ar@{.}[u]\\
 \circ,}
\end{array}
$$
i.e. the subquivers which have as first and last vertex
respectively: the $\sigma$-fixed vertex and $ta$, where $a$ is the
$\sigma$-fixed arrow, for $Q'$; the $\sigma$-fixed vertices for
$Q''$; $ta$ and $tb$, where $a$ and $b$ are the $\sigma$-fixed
arrows, for $Q'''$. We note that these three quivers have
$\frac{k}{2}$ counterclockwise arrows and $\frac{l}{2}$ clockwise
arrows. So for each one of $Q'$, $Q''$ and $Q'''$ there exists a
sequence of sinks $x_1,\ldots,x_s$ such that $c_{x_1}\cdots
c_{x_s}Q'$, $c_{x_1}\cdots c_{x_s}Q''$ and $c_{x_1}\cdots
c_{x_s}Q'''$ are respectively
$$
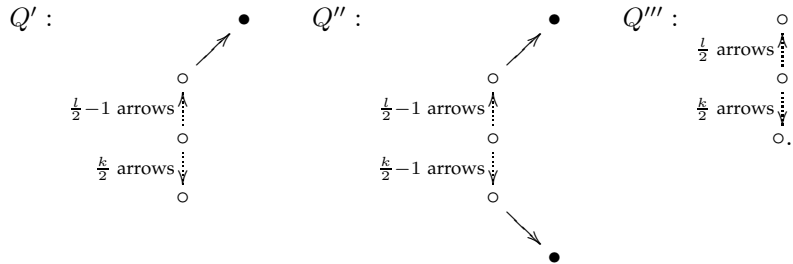
\begin{array}{ccccc}
Q':\xymatrix@-1pc{&\bullet\\
\circ\ar[ur]&\\
\circ\ar@{.>}[u]^{\frac{l}{2}-1\;\textrm{arrows}}\ar@{.>}[d]_{\frac{k}{2}\;\textrm{arrows}}&\\
 \circ&} && Q'':\xymatrix@-1pc{&\bullet\\
\circ\ar[ur]&\\
\circ\ar@{.>}[d]_{\frac{k}{2}-1\;\textrm{arrows}}\ar@{.>}[u]^{\frac{l}{2}-1\;\textrm{arrows}}&\\
\circ\ar[dr]&\\
&\bullet} && Q''':\xymatrix@-1pc{
\circ\\
\circ\ar@{.>}[d]_{\frac{k}{2}\;\textrm{arrows}}\ar@{.>}[u]^{\frac{l}{2}\;\textrm{arrows}}\\
 \circ.}
\end{array}
$$
Hence, by symmetry, applying $c_{(x_1,\sigma(x_1))}\cdots
c_{(x_s,\sigma(x_s))}$, we
obtain the desired orientations.\\
For $(Q,\sigma)$ of type $(0,0,k,k)$ we consider a sink $x$ of $Q$
and we apply lemma \ref{Afix} to the subquiver $Q'$ which has as
first and last vertex respectively $x$ and $\sigma(x)$. So there
exists a sequence of sinks $x_1,\ldots,x_s$ such that
$c_{x_1}\cdots c_{x_s}Q'$ is
$$
\xymatrix{x&&\circ\ar@{.>}[ll]_{k'\;\textrm{arrows}}\ar@{.>}[rr]^{k''\;\textrm{arrows}}&&\sigma(x).}
$$
Hence, by symmetry, applying $c_{(x_1,\sigma(x_1))}\cdots
c_{(x_s,\sigma(x_s))}$ we obtain
$$
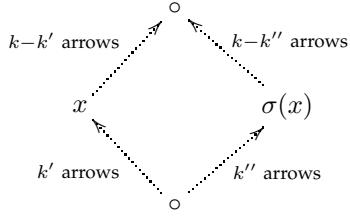
\xymatrix{&\circ&\\
x\ar@{.>}[ur]^{k-k'\;\textrm{arrows}}&&\sigma(x)\ar@{.>}[ul]_{k-k''\;\textrm{arrows}}\\
&\circ\ar@{.>}[ul]^{k'\;\textrm{arrows}}\ar@{.>}[ur]_{k''\;\textrm{arrows}}&}
$$
i.e. the desired orientation. $\Box$
\subsection{Reflection functors for symmetric quivers}
Let $(Q,\sigma)$ be a symmetric quiver, $(x,\sigma(x))$ a
sink-source admissible pair. For every $V\in Rep(Q)$, we define
the reflection functors
$$C_{(x,\sigma(x))}^+V:=C_{\sigma(x)}^-C_x^+V$$
and
$$C_{(\sigma(x),x)}^-V:=C_x^-C_{\sigma(x)}^+V.$$
We note that $C_{\sigma(x)}^-C_x^+V=C_x^+C_{\sigma(x)}^-V$
(respectively $C_x^-C_{\sigma(x)}^+V=C_{\sigma(x)}^-C_x^+V$) since
there are no arrows connecting $x$ and $\sigma(x)$.
\begin{proposizione}\label{C+nabla}
Let $(Q,\sigma)$ be a symmetric quiver and $V$ be a representation
of the underlying quiver.
\begin{itemize}
\item[(i)] If $x$ is an admissible sink, then $\nabla
C^+_{(x,\sigma(x))}V\cong C^+_{(x,\sigma(x))}\nabla V$.
\item[(ii)] If $x$ is an admissible source, then $\nabla
C^-_{(x,\sigma(x))}V\cong C^-_{(x,\sigma(x))}\nabla V$.
\end{itemize}
In particular for every $x$ admissible sink and $y$ admissible
source we have
$$
V=\nabla V \Leftrightarrow C^+_{(x,\sigma(x))}V=\nabla
C^+_{(x,\sigma(x))} V\Leftrightarrow C^-_{(y,\sigma(y))}V=\nabla
C^-_{(y,\sigma(y))} V.
$$
\end{proposizione}
\textit{Proof.} We prove (i) (the proof of (ii) is similar).
Recall that $x\neq\sigma(x)$, otherwise $x$ is not a sink. Let
$\{a_1,\ldots,a_k\}$ be the set of arrows whose head is $x$.
$$
(\nabla C^+_{(x,\sigma(x))}
V)_y=(C^+_xC^-_{\sigma(x)}V)^*_{\sigma(y)}=
$$
$$
\left\{\begin{array}{ll}
(V_{\sigma(y)})^* & \sigma(y)\neq\sigma(x),x\\
(Coker(V_{\sigma(x)}\stackrel{\tilde{h}}{\longrightarrow}\bigoplus_{i=1}^k
V_{h\sigma(a_i)}))^* & \sigma(y)=\sigma(x)\\
(Ker(\bigoplus_{i=1}^k V_{ta_i}\stackrel{h'}{\longrightarrow}
V_x))^* & \sigma(y)=x,
\end{array}\right.
$$
where $\tilde{h}(v)=(V(\sigma(a_1))(v),\ldots,V(\sigma(a_k))(v))$
with $v\in V_{\sigma(x)}$ and
$h'(v_1,\ldots,v_k)=V(a_1)(v_1)+\cdots+V(a_k)(v_k)$ with
$(v_1,\ldots,v_k)\in \bigoplus_{i=1}^k V_{ta_i}$.
$$
(C^+_{(x,\sigma(x))}\nabla V)_y=
$$
$$
\left\{\begin{array}{ll}
(\nabla V_y)_y & y\neq\sigma(x),x\\
Coker((\nabla
V)_{\sigma(x)}\stackrel{\tilde{h'}}{\longrightarrow}\bigoplus_{i=1}^k
(\nabla V)_{h\sigma(a_i)}) & y=\sigma(x)\\
Ker(\bigoplus_{i=1}^k (\nabla
V)_{ta_i}\stackrel{h}{\longrightarrow}(\nabla V)_x) & y=x,
\end{array}\right.
$$
where $\tilde{h'}(v)=(\nabla V(\sigma(a_1))(v),\ldots,\nabla
V(\sigma(a_k))(v))$ with $v\in(\nabla V)_{\sigma(x)}$ and
$h(v_1,\ldots,v_k)=\nabla V(a_1)(v_1)+\cdots+\nabla V(a_k)(v_k)$
with
$(v_1,\ldots,v_k)\in \bigoplus_{i=1}^k(\nabla V)_{ta_i}$.\\
Since $(\nabla V)_y=(V_{\sigma(y)})^*$ for every $y\in Q_0$ and
$\nabla V(a)=-V(\sigma(a))^*$, we have $h=-\tilde{h}^*$ and
$h'=-\tilde{h'}^*$; moreover if $\varphi$ is a linear map, in
general we have $(Ker(\varphi))^*\cong Coker(\varphi^*)$ and
$(Coker(\varphi))^*\cong Ker(\varphi^*)$, so $(\nabla
C^+_{(x,\sigma(x))}V)_y\cong
(C^+_{(x,\sigma(x))}\nabla V)_y$ for every $y\in Q_0$.\\
We note that
$$
\sigma(c_{(x,\sigma(x))}a)=\sigma(c_x
c_{\sigma(x)}a)=\left\{\begin{array}{ll} c_{\sigma(x)}\sigma(a_i)
&
a=a_i\;\textrm{with}\;i\in\{1,\ldots,k\}\\
c_x a_i & a=\sigma(a_i)\;\textrm{with}\;i\in\{1,\ldots,k\}\\
\sigma(a) & a\neq
a_i,\sigma(a_i)\;\textrm{with}\;i\in\{1,\ldots,k\}.
\end{array}\right.
$$
So we have
$$
(\nabla
C^+_{(x,\sigma(x))}V)(c_{(x,\sigma(x))}a)=-((C^+_xC^-_{\sigma(x)}V)(\sigma(c_x
c_{\sigma(x)}a)))^*=
$$
$$
\left\{\begin{array}{ll} -V(\sigma(a))^* & a\neq
a_j,\sigma(a_j)\;\textrm{with}\;j\in\{1,\ldots,k\}\\
-(V_x\hookrightarrow\bigoplus_{i=1}^k V_{ta_i}\twoheadrightarrow
V_{ta_j})^* & a=\sigma(a_j)\;\textrm{with}\;j\in\{1,\ldots,k\}\\
-(V_{h\sigma(a_j)}\hookrightarrow\bigoplus_{i=1}^k
V_{h\sigma(a_i)}\twoheadrightarrow V_{\sigma(x)})^* &
a=a_j\;\textrm{with}\;j\in\{1,\ldots,k\}
\end{array}\right.
$$
and
$$
( C^+_{(x,\sigma(x))}\nabla V)(c_{(x,\sigma(x))}a)=
$$
$$
\left\{\begin{array}{ll} \nabla V(a) & a\neq
a_j,\sigma(a_j)\;\textrm{with}\;j\in\{1,\ldots,k\}\\
(\nabla V)_x\hookrightarrow\bigoplus_{i=1}^k (\nabla
V)_{ta_i}\twoheadrightarrow
(\nabla V)_{ta_j} & a=a_j\;\textrm{with}\;j\in\{1,\ldots,k\}\\
(\nabla V)_{h\sigma(a_j)}\hookrightarrow\bigoplus_{i=1}^k (\nabla
V)_{h\sigma(a_i)}\twoheadrightarrow(\nabla V)_{\sigma(x)} &
a=a_j\;\textrm{with}\;j\in\{1,\ldots,k\}.
\end{array}\right.
$$
Hence $\nabla C^+_{(x,\sigma(x))}V\cong C^+_{(x,\sigma(x))}\nabla
V$. $\Box$
\begin{cor}\label{C+taunabla}
Let $(Q,\sigma)$ and $(Q',\sigma)$ be two symmetric quivers with
the same underlying graph. We suppose there exists a sequence
$x_1,\ldots,x_m$ of admissible sinks for admissible sink-source
pairs such that $Q'=c_{(x_m,\sigma(x_m))}\cdots
c_{(x_1,\sigma(x_1))}Q$. Let $V\in Rep(Q)$ and
$V'=C^+_{(x_m,\sigma(x_m))}\cdots C^+_{(x_1,\sigma(x_1))}V\in Rep(Q')$. Then\\
$$V=\tau^-\nabla V\Longleftrightarrow V'=\tau^-\nabla V'.$$
\end{cor}
\textit{Proof.} By proposition \ref{C+nabla}, we have
$$\tau^-\nabla V'=\tau^-\nabla C^+_{(x_m,\sigma(x_m))}\!\!\!\cdots
C^+_{(x_1,\sigma(x_1))}V=\tau^-C^+_{(x_m,\sigma(x_m))}\!\!\!\cdots
C^+_{(x_1,\sigma(x_1))}\nabla V=$$
$$\tau^-C^-_{\sigma(x_m)}C^+_{x_m} \cdots
C^-_{\sigma(x_1)}C^+_{x_1}\tau^+ V=C^-_{\sigma(x_m)}\tau^-
C^+_{x_m}\!\!\! \cdots C^-_{\sigma(x_1)}\tau^+ C^+_{x_1}
V=\cdots=$$ $$C^-_{\sigma(x_m)}\!\!\! \cdots
C^-_{\sigma(x_1)}\tau^-\tau^+ C^+_{x_m} \!\!\!\cdots C^+_{x_1}
V=C^+_{(x_m,\sigma(x_m))}\cdots
C^+_{(x_1,\sigma(x_1))}V=V'.\quad\Box$$
\begin{proposizione}\label{sptosp}
Let $(Q,\sigma)$ be a symmetric quiver and let $x$ be an
admissible sink. Then
\begin{itemize}
\item[(i)] $V$ is a symplectic
representation of $(Q,\sigma)$ if and only if
$C^+_{(x,\sigma(x))}V$ is a symplectic representation;
\item[(ii)] $V$ is a orthogonal
representation of $(Q,\sigma)$ if and only if
$C^+_{(x,\sigma(x))}V$ is a orthogonal representation.
\end{itemize}
Similarly if $x$ is an admissible source then
$C^-_{(x,\sigma(x))}$ sends symplectic representations to
symplectic representations and orthogonal representations to
orthogonal representations.
\end{proposizione}
\textit{Proof.} By proposition \ref{C+nabla} we have $V=\nabla V$
if and only if $C^+_{(x,\sigma(x))}V=\nabla C^+_{(x,\sigma(x))}V$.
To define an orthogonal (respectively symplectic) structure on
$C^+_{(x,\sigma(x))}V$ the only problem could occur at the
vertices fixed by $\sigma$. But, by definition of admissible sink
and of the involution $\sigma$, fixed vertices and fixed arrows
don't change under our reflection. The proof is similar for
$C^-_{(x,\sigma(x))}$ with
$x$ an admissible source. $\Box$\\
\\
Next we prove that the reflection functors for symmetric quivers
preserve the rings of orthogonal and symplectic semi-invariants.
We need some basic property of Grasmannians.
\begin{definizione} Let $W$ be a vector
space of dimension $n$. Consider the set of all decomposable
tensor $w_1\wedge\ldots\wedge w_r$, with $w_1,\ldots,w_r\in W$,
inside $\bigwedge^r W$. This set is an affine subvariety of the
space vector $\bigwedge^r W$, called \textit{affine cone over the
Grasmannian}. It will be denoted by $\widetilde{Gr}(r,W)$.
\end{definizione}
\begin{definizione}
The Grasmannian $Gr(r,W)$ is the projective subvariety of
$\mathbb{P}(\bigwedge^r W)$ corresponding to
$\widetilde{Gr}(r,W)$.
\end{definizione}
This variety can be thought as the set of $r$-dimensional
subspaces of $W$. The identification between $\bigwedge^rW$ and
$\bigwedge^{n-r}W^*$ induces an identification between
$\widetilde{Gr}(r,W)$ and $\widetilde{Gr}(n-r,W^*)$ and so between
$Gr(r,W)$ and $Gr(n-r,W^*)$. By the first fundamental theorem
(FFT) for $SL\,V$ (see [P, chapter 11 section 1.2 ]), it follows
that
$$
\mathbb{K}[V\otimes
W]^{SL\,V}\cong\mathbb{K}[\widetilde{Gr}(r,W)],
$$
where $r=dim(V)$.
\begin{lemma}\label{kac}
If $x$ is an admissible sink or source for a symmetric quiver
$(Q,\sigma)$ and $\alpha$ is a dimension vector such that
$c_{(x,\sigma(x))}\alpha(x)\geq 0$, then
\begin{itemize}
\item[i)] if $c_{(x,\sigma(x))}\alpha(x)> 0$ there exist isomorphisms
$$
SpSI(Q,\alpha)\stackrel{\varphi^{Sp}_{x,\alpha}}{\longrightarrow}
SpSI(c_{(x,\sigma(x))}Q,c_{(x,\sigma(x))}\alpha)
$$
and
$$
OSI(Q,\alpha)\stackrel{\varphi^{O}_{x,\alpha}}{\longrightarrow}
OSI(c_{(x,\sigma(x))}Q,c_{(x,\sigma(x))}\alpha),
$$
\item[ii)] if $c_{(x,\sigma(x))}\alpha(x)= 0$ there exist isomorphisms
$$
SpSI(Q,\alpha)\stackrel{\varphi^{Sp}_{x,\alpha}}{\longrightarrow}
SpSI(c_{(x,\sigma(x))}Q,c_{(x,\sigma(x))}\alpha)[y]
$$
and
$$
OSI(Q,\alpha)\stackrel{\varphi^{O}_{x,\alpha}}{\longrightarrow}
OSI(c_{(x,\sigma(x))}Q,c_{(x,\sigma(x))}\alpha)[y]
$$
\end{itemize}
where $A[y]$ denotes a polynomial ring with coefficients in $A$.
\end{lemma}
\textit{Proof.} We will prove the lemma for the symplectic case
because the orthogonal case is similar. Let $x\in Q_0$ be an
admissible sink. Put $r=\alpha(x)$ and $n=\sum_{ha=x}\alpha(ta)$.
We note that $c_{(x,\sigma(x))}\alpha(x)=n-r$. Put
$V=\mathbb{K}^r$, $V'=\mathbb{K}^{n-r}$ and
$W=\bigoplus_{ha=x}\mathbb{K}^{\alpha(ta)}\cong\mathbb{K}^n$. We
define $$Z=\bigoplus_{{a\in Q_1^+\atop ha\neq
x}}Hom(\mathbb{K}^{\alpha(ta)},\mathbb{K}^{\alpha(ha)})\oplus\bigoplus_{a\in
Q_1^{\sigma}}S^2(\mathbb{K}^{\alpha(ta)})^*$$
and
$$G=\prod_{{y\in
Q_0^+\atop y\neq x}}SL(\alpha(y))\times\prod_{y\in
Q_0^{\sigma}}Sp(\alpha(y)).
$$
\textit{Proof of i).} If $c_{(x,\sigma(x))}\alpha(x)> 0$ we have
$$
SpSI(Q,\alpha)=\mathbb{K}[SpRep(Q,\alpha)]^{SSp(Q,\alpha)}=
$$
$$
\mathbb{K}[Z\times Hom(W,V)]^{G\times
SL\,V}=(\mathbb{K}[Z]\otimes\mathbb{K}[ Hom(W,V)]^{SL\,V}]^G=
$$
$$
(\mathbb{K}[Z]\otimes\mathbb{K}[\widetilde{Gr}(r,W^*)])^G
$$
and
$$
SpSI(c_{(x,\sigma(x))}Q,c_{(x,\sigma(x))}\alpha)=
$$
$$
\mathbb{K}[SpRep(c_{(x,\sigma(x))}Q,c_{(x,\sigma(x))}\alpha)]^{SSp(c_{(x,\sigma(x))}Q,c_{(x,\sigma(x))}\alpha)}=
$$
$$
\mathbb{K}[Z\times Hom(V',W)]^{G\times
SL\,V'}=(\mathbb{K}[Z]\otimes\mathbb{K}[ Hom(V',W)]^{SL\,V'}]^G=
$$
$$
(\mathbb{K}[Z]\otimes\mathbb{K}[\widetilde{Gr}(n-r,W)])^G.
$$
Since $\widetilde{Gr}(r,W^*)$ and $\widetilde{Gr}(n-r,W)$ are
isomorphic as $G$-varieties, it follows that $SpSI(Q,\alpha)$ and
$SpSI(c_{(x,\sigma(x))}Q,c_{(x,\sigma(x))}\alpha)$ are
isomorphic.\\
\textit{Proof ii).} If $c_{(x,\sigma(x))}\alpha(x)= 0$, then $n=r$
and $V'=0$. So $\widetilde{Gr}(0,W)$ is a point and hence
\begin{equation}\label{era}
SpSI(Q,\alpha)=(\mathbb{K}[Z]\otimes
\mathbb{K}[Hom(W,V)])^{G\times SL\,V}
\end{equation}
is isomorphic to
$$
SpSI(c_{(x,\sigma(x))}Q,c_{(x,\sigma(x))}\alpha)=(\mathbb{K}[Z]\otimes
\mathbb{K}[Hom(V',W)])^{G\times SL\,V}=\mathbb{K}[Z]^{G\times
SL(V)}.
$$
Now let $A=\{a\in Q^+_1|\,ha=x\}$. Using theorem \ref{c}, each
summand of (\ref{era}) contains $(\bigotimes_{a\in
A}S_{\lambda(a)}V)^{SL\,V}$ as factor. By proposition \ref{i2}
each $\lambda(a)$, with $a\in A$, has to contain a column of
height $\alpha(ta)$, hence $\lambda(a)=\mu(a)+(1^{\alpha(ta)})$,
for some $\mu(a)$ in the set of partitions $\Lambda$. So as factor
we have
$$
\bigotimes_{a\in
A}(S_{(1^{\alpha(ta)})}\mathbb{K}^{\alpha(ta)})^{SL\,V_{ta}}\otimes\left(\bigotimes_{a\in
A}S_{(1^{\alpha(ta)})}V\right)^{SL\,V}
$$
which is generated by
$det(\bigoplus_{ha=x}\mathbb{K}^{\alpha(ta)}\rightarrow\mathbb{K}^{\alpha(x)})$.
On the other hand we have $\mathbb{K}[Hom(W,V)])^{G\times
SL\,V}=\mathbb{K}[det(\bigoplus_{ha=x}\mathbb{K}^{\alpha(ta)}\rightarrow\mathbb{K}^{\alpha(x)})]$
and so we have the statement \textit{ii)}, with
$y=det(\bigoplus_{ha=x}\mathbb{K}^{\alpha(ta)}\rightarrow\mathbb{K}^{\alpha(x)})$.
$\Box$

\section{Semi-invariants of symmetric quivers}
In this section we prove some general results about
semi-invariants of symmetric quivers with underlying
  quiver without oriented cycles.\\
  We assume that $(Q,\sigma)$ is a symmetric quiver with underlying quiver $Q$ without oriented
  cycles for rest of the thesis.\\
  We recall that, by definition, symplectic groups or orthogonal
groups act
  on the spaces which are defined on the $\sigma$-fixed
  vertices, so we have
  \begin{definizione}\label{wsQ}
  Let $V$ be a
  representation of the underlying quiver $Q$ with $\underline{dim}V=\alpha$ such that
  $\langle\alpha,\beta\rangle=0$ for some symmetric
  dimension vector $\beta$. The weight of $c^V$ on
  $SpRep(Q,\beta)$ (respectively on $ORep(Q,\beta)$) is
  $\langle\alpha,\cdot\rangle-\sum_{x\in
  Q_0^{\sigma}}\varepsilon_{x,\alpha}$, where
  \begin{equation}\label{epsilon}
  \varepsilon_{x,\alpha}(y)=\left\{\begin{array}{ll}
  \langle\alpha,\cdot\rangle(x) & y=x\\
  0 & \textrm{otherwise}.\end{array}\right.
  \end{equation}
  \end{definizione}
  In general we define an involution $\gamma$ on the space of
  weights $\langle\alpha,\cdot\rangle$ with $\alpha$ dimension vector.
  \begin{definizione} Let $\alpha$ be the dimension vector of a representation $V$ of the underlying quiver $Q$ and let
   $\langle\alpha,\cdot\rangle=\chi=\{\chi(i)\}_{i\in Q_0}$ be the
  weight of $c^V$. We define
  $\gamma \chi=\{\gamma \chi(i)\}_{i\in Q_0}$ where $\gamma \chi(i)=-\chi(\sigma(i))$ for every $i\in Q_0$
\end{definizione}
We number vertices in such way that $ta<ha$ for every $a\in
  Q_1$. We note that $\chi=\langle\alpha,\cdot\rangle=(\alpha(j)-\sum_{i<j}b_{i,j}\alpha(i))_{j\in Q_0}$,\\
  where
  $b_{i,j}:=|\{a\in Q_1|ta=i,ha=j\}|=|\{a\in Q_1|ta=\sigma(j),ha=\sigma(i)\}|=:b_{\sigma(j),\sigma(i)}$.
  \begin{lemma}\label{gammachi}
  \begin{eqnarray}
  \gamma \chi=\langle\tau^-\delta\alpha,\cdot\rangle=\langle\underline{dim}(\tau^-\nabla
  V),\cdot\rangle,
  \end{eqnarray}
  i.e. $\gamma \chi$ is the weight of $c^{\tau^-\nabla
  V}$. Moreover $\gamma$ is an involution.
  \end{lemma}
  \textit{Proof.} By definition of $\gamma$, $\gamma
  \chi(j)=-\alpha(\sigma(j))+\sum_{i<j}b_{i,j}\alpha(\sigma(i))$.
  Now it follows
  by theorem \ref{ard} that
  $\langle\tau^-\delta\alpha,\cdot\rangle=-\langle\cdot,\delta\alpha\rangle$,
  thus, for every $j\in Q_0$,
  $\langle\tau^-\delta\alpha,\cdot\rangle(j)=-\langle\cdot,\delta\alpha\rangle(j)=-\delta\alpha(j)+\sum_{i<j}b_{i,j}\delta\alpha(i)=\gamma
  \chi(j)$. Hence $\gamma
  \chi=\langle\tau^-\delta\alpha,\cdot\rangle$. $\Box$\\
  \\
  Moreover, since $\gamma\gamma
  \chi(i)=\gamma(-\chi(\sigma(i)))=\chi(\sigma\sigma(i))=\chi(i)$ for every
  $i\in Q_0$, $\gamma$ is an involution.\\
  If $\beta$ is the dimension vector of a
  representation W of the underlying quiver $Q$, we have
  \begin{eqnarray}\label{alphabeta0}
  \langle\alpha,\beta\rangle=0\Leftrightarrow\langle\tau^-\delta\alpha,\delta\beta\rangle=0.
  \end{eqnarray}
  Indeed, by theorem \ref{ard},
 \begin{eqnarray}\label{alphabetataunabla}
  \langle\alpha,\beta\rangle=\langle\delta\beta,\delta\alpha\rangle=-\langle\tau^-\delta\alpha,\delta\beta\rangle.
  \end{eqnarray}
  Since $\beta$ is the dimension vector of an
  orthogonal or symplectic representation $W$,
  we have that $\beta$ is a symmetric dimension vector and so
  \begin{eqnarray}
  \langle\alpha,\beta\rangle=0\Leftrightarrow\langle\tau^-\delta\alpha,\beta\rangle=0.
  \end{eqnarray}

\begin{lemma}\label{cV=cVnabla}
Let $(Q,\sigma)$ be a symmetric quiver. For every representation
  $V$ of the underlying quiver $Q$ and for every orthogonal or
  symplectic representation $W$ such that
  $\langle\underline{dim}(V),\underline{dim}(W)\rangle=0$, we have
  $$
  c^{V}(W)=c^{\tau^-\nabla V}(W).
  $$
\end{lemma}
\textit{Proof.} It follows directly from lemma \ref{cV=cVnabla1}. $\Box$\\
\\
Now we prove in general a crucial lemma which will be useful
later. Let $(Q,\sigma)$ be a symmetric quiver. If $V$ is a
  representation of the underlying quiver $Q$ such that
  $V=\tau^-\nabla V$ then, by the theorem \ref{ar}, there exists an
  almost split sequence $0\rightarrow\nabla V\rightarrow
  Z\rightarrow V\rightarrow 0$ with $Z\in Rep(Q)$. Moreover for such $V\in Rep(Q)$ with $\underline{dim}V=\alpha$ we have
  $\alpha=\tau^-\delta\alpha$ and
$\gamma \chi=\chi$, where $\chi=\langle\alpha,\cdot\rangle$. So
$\chi(i)=\gamma \chi(i)=-\chi(\sigma(i))$ for every
$i\in\{1,\ldots,n\}$, in particular $\chi(i)=0$ if $\sigma(i)=i$.
\begin{definizione}
A weight $\chi$ such that $\gamma \chi=\chi$ is called a symmetric
weight.
\end{definizione}

\begin{lemma}\label{ss1}
  Let $(Q,\sigma)$ be a symmetric quiver of finite type or of tame type. Let $d_{min}^V$
  be the matrix of the minimal projective presentation of $V\in Rep(Q,\alpha)$ and
  let $\beta$ be a symmetric dimension vector such that
  $\langle\alpha,\beta\rangle=0$. Then
  \begin{itemize}
  \item[(1)] $Hom_Q(d_{min}^V,\cdot)$ is skew-symmetric on $SpRep(Q,\beta)$  if and only if $V$ satisfies property \textit{(Op)};
  \item[(2)] $Hom_Q(d_{min}^V,\cdot)$ is skew-symmetric on $ORep(Q,\beta)$ if and only if $V$ satisfies property \textit{(Spp)}.
  \end{itemize}
\end{lemma}
\textit{Proof.} We use notation of section \ref{Qtt}. We call
$(Q',\sigma)$ the symmetric quiver with the same underlying graph
of $(Q,\sigma)$ such that
\begin{itemize}
\item[(i)] if $Q$ is of type $A$, then $Q'$ has all the arrows
with the same orientations;
\item[(ii)] if $Q$ is of type $\widetilde{A}$, then $Q'$ is one
of the quiver as in proposition \ref{oA} (it depends on which kind
of quiver is $Q$);
\item[(iii)] if $Q$ is of type $\widetilde{D}$, then $Q'$ is
$\widetilde{D}^{eq}$ (see picture (\ref{Dequi})).
\end{itemize}
By propositions \ref{Q=Q'} and \ref{oA}, there exists a sequence
$x_1,\ldots,x_m$ of admissible sink for admissible sink-source
pairs such that $c_{(x_m,\sigma(x_m))}\cdots
c_{(x_1,\sigma(x_1))}Q=Q'$. We call
$V':=C^+_{(x_m,\sigma(x_m))}\cdots C^+_{(x_1,\sigma(x_1))}V$ for
every $V\in Rep(Q)$ and if $\alpha=\underline{dim}\,V$, then
$\alpha':=c_{(x_m,\sigma(x_m))}\cdots
c_{(x_1,\sigma(x_1))}\alpha$. We note that, by corollary
\ref{C+taunabla} and proposition \ref{sptosp}, $V$ satisfies
property \textit{(Op)} (respectively property \textit{(Spp)}) if
and only if $V'$ satisfies property \textit{(Op)} (respectively
property \textit{(Spp)}). We prove only \textit{(1)}, because the
proof of \textit{(2)} is similar.\\
\\
\textbf{Type $\boldsymbol{A}$.} Let $(A_n,\sigma)$ be a symmetric
quiver of type $A$. We enumerate vertices with $1,\ldots,n$ from
left to right and we call $a_i$ the arrow with $i$ on the left and
$i+1$ on the right. We define $\sigma$ by $\sigma(i)=n-i+1$ for
every $i\in Q_0$ and $\sigma(a_i)$ for every
$i\in\{1,\ldots,n-1\}$. Let $V'=V_{i,\sigma(i)-1}$, i.e. is the
indecomposable of $A_n$ such that
$$
(\underline{dim} V_{i,\sigma(i)-1})_j=\left\{\begin{array}{lr}1 &
j\in\{i,\ldots,\sigma(i)-1\}\\
0 & \textrm{otherwise}.\end{array}\right.
$$
We note that $\nabla V'=V_{i+1,\sigma(i)}=\tau^+ V'$ and
$Z'=V_{i,\sigma(i)}\oplus V_{i+1,\sigma(i)-1}$. So, by definition
\ref{defortsp}, on $Z'$ we can define a structure of orthogonal
representation if $n$ is odd and a structure of symplectic
representation if $n$ is even. So it's enough to check when
$Hom_Q(d_{min}^V,\cdot)$ is skew symmetric and, for type $A$, we
do it
explicitly.\\
Let $\chi=\langle\alpha,\cdot\rangle-\sum_{x\in
Q_0^{\sigma}}\varepsilon_{x,\alpha}$ be the symmetric weight
associated to $\alpha$. If $m_1$ is the first vertex such that
$\chi(m_1)\neq 0$, in particular we suppose $\chi(m_1)=1$, then
the last vertex $m_s$ such that $\chi(m_s)\neq 0$ is
$m_s=\sigma(m_1)$ and $\chi(m_s)=-1$. Between $m_1$ and $m_s$, -1
and 1 alternate in correspondence respectively of sinks and of
sources. Moreover, by definition of symmetric weight, we have
$s=2l$ for some $l\in\mathbb{N}$. We call $i_2,\ldots,i_l$ the
sources, $j_1,\ldots,j_{l-1}$ the sinks, $i_1=m_1$ and $j_l=m_s$.
Hence we have $\sigma(i_t)=j_{l-t+1}$ and
$i_1<j_1<\ldots<i_l<j_l$. Now the minimal projective resolution
for V is
\begin{equation}
0\longrightarrow\bigoplus_{j=j_1}^{j_l}
P_{j}\stackrel{d_{min}^V}{\longrightarrow}\bigoplus_{i=i_1}^{i_l}P_{i}\longrightarrow
V\longrightarrow 0
\end{equation}
and for the remark above we have
\begin{equation}
0\longrightarrow\bigoplus_{j=j_1}^{j_l}
P_{j}\stackrel{d_{min}^V}{\longrightarrow}\bigoplus_{j=j_1}^{j_l}P_{\sigma(j)}\longrightarrow
V\longrightarrow 0,
\end{equation}
with
\begin{equation}
(d_{min}^V)_{hk}=\left\{\begin{array}{ll} -a_{i_{k+1},j_k} &
\textrm{if}\quad
h=l-k\\
a_{i_k,j_k} & \textrm{if}\quad h=l-k+1\\
0 & \textrm{otherwise},
\end{array}\right.
\end{equation}
where $a_{i,j}$ is the oriented path from $i$ to $j$.\\
Hence
\begin{equation}
Hom(d_{min}^V,W):\bigoplus_{j=j_1}^{j_l}W(\sigma(j))=\bigoplus_{j=j_1}^{j_l}W(j)^*\longrightarrow\bigoplus_{j=j_1}^{j_l}W(j)
\end{equation}
where
\begin{equation}
(Hom(d_{min}^V,W))_{hk}=\left\{\begin{array}{ll}
-W(a_{i_{h+1},j_{h}}) &
\textrm{if}\quad k=l-h\\
W(a_{i_h,j_h}) &
\textrm{if}\quad k=l-h+1\\
0 & \textrm{otherwise}.
\end{array}\right.
\end{equation}
Now $W$ is orthogonal or symplectic, so for $k\neq h$, if
$k=l-h+1$ we have
$$
(Hom(d_{min}^V,W))_{hk}=W(a_{i_h,j_{h}})=W(a_{\sigma(j_{l-h+1}),j_{h}})=
-W(a_{\sigma(j_{h}),j_{l-h+1}})^t=
$$
$$
-W(a_{i_{l-h+1},j_{l-h+1}})^t=
-W(a_{i_k,j_k})^t=-((Hom(d_{min}^V,W))_{kh})^t.
$$
In a similar way it proves that if $k=l-h$ then
$(Hom(d_{min}^V,W))_{hk}=-((Hom(d_{min}^V,W))_{kh})^t$.\\
Finally the only cases for which $(Hom(d_{min}^V,W))_{hh}\neq 0$
are when $h=l-h+1$ and $h=l-h$. In the first case (the second one
is similar) we have
$(Hom(d_{min}^V,W))_{hh}=W(a_{i_h,j_{h}})=W(a_{\sigma(j_{h}),j_{h}})$
and $-((Hom(d_{min}^V,W))_{hh})^t=-W(a_{i_{h},j_{h}})^t=
-W(a_{\sigma(j_{h}),j_{h}})^t$. But
$W(a_{\sigma(j_{h}),j_{h}})=-W(a_{\sigma(j_{h}),j_{h}})^t$ for $n$
even if and only if $W\in ORep(Q)$, for $n$ odd if and only if
$W\in SpRep(Q)$.\\
\\
We consider the tame case. First we note, by Auslander-Reiten
quiver of $Q$, that if $(Q,\sigma)$ is a symmetric quiver of tame
type, then the only representations $V\in Rep(Q)$ such that
$\tau^-\nabla V=V$ are regular ones.\\
\\
\textbf{Type $\boldsymbol{\widetilde{A}}$}. We prove lemma only
for $Q$ of type $(1,1,k,l)$ because for the other cases it
proceeds similarly. We consider the following labelling for
$Q'=\widetilde{A}^{1,1}_{k,l}$:
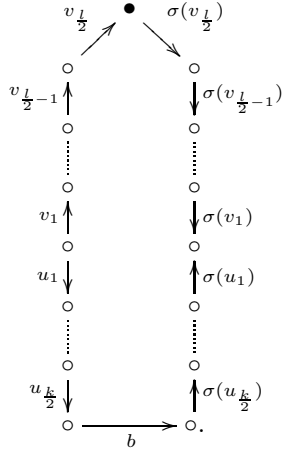
$$
\xymatrix@-1pc{ &\bullet\ar[dr]^{\sigma(v_{\frac{l}{2}})}&\\
\circ\ar[ur]^{v_{\frac{l}{2}}}&&\circ\ar[d]^{\sigma(v_{\frac{l}{2}-1})}\\
\circ\ar[u]^{v_{\frac{l}{2}-1}}\ar@{.}[d]&&\circ\\
\circ&&\circ\ar[d]^{\sigma(v_{1})}\ar@{.}[u]\\
\circ\ar[d]_{u_1}\ar[u]^{v_1}&&\circ\\
\circ&&\circ\ar[u]_{\sigma(u_{1})}\ar@{.}[d]\\
\circ\ar[d]_{u_{\frac{k}{2}}}\ar@{.}[u]&&\circ\\
 \circ\ar[rr]_{b}&&\circ\ar[u]_{\sigma(u_{\frac{k}{2}})}.}
$$
The following indecomposable representations $V'\in Rep(Q')$
satisfy property \textit{(Op)}. The other regular indecomposable
representations of $Rep(Q')$ satisfying property \textit{(Op)} are
extensions of these.
\begin{itemize}
\item[(a)] $V_{(0,1)}$; in this case $Z'=E_{h}^1\oplus E_{2,0}$
where $E_{h}^1$ is the regular indecomposable representation of
dimension $e_1+h$ with socle $E_1$.
\item[(b)] $E_{i,j-1}$, with $1\leq i<j\leq l+1$, such that
$\nabla E_{i,j-1}=E_{i+1,j}$; in this case we have
$Z'=E_{i+1,j-1}\oplus E_{i,j}$.
\item[(c)] $E'_{i,j-1}$, with $2\leq j<i-1\leq k+1$, such that
$\nabla E'_{i,j-1}=E'_{i+1,j}$; in this case we have
$Z'=E'_{i+1,j-1}\oplus E'_{i,j}$.
\end{itemize}
Let $\chi$ be the symmetric weight associated to $\alpha$. We
order vertices of $Q$ clockwise from $tb=1$ to $hb=k+l+1$. We use
the same notation of type $A$ for vertices on which the
components of $\chi$ are not zero.\\
Let $W$ be a symplectic representation. We prove that
$Hom_Q(d_{min}^V,W)$ is skew-symmetric for every regular
indecomposable representation $V$ of type (a), (b) and (c). First
we observe that the associated to $V$ symmetric weight $\chi$ have
components equal to 0, 1 and -1. In particular, $\chi(m_1)=\pm
1=-\chi(m_s)$ and $\chi(m_i)=1,-1$, for every
$i\in\{2,\ldots,s-1\}$, respectively if $m_i$ is a source or a
sink. We note that, for every $Hom_Q(d_{min}^V,W)$ with $V$ one
representation of type (a), (b) and (c), we can restrict to the
symmetric subquiver of type $A$ which has first vertex $m_1$ and
last vertex $m_s$ and passing through the $\sigma$-fixed vertex of
$Q$. Hence it proceeds as done for type
$A$.\\
Finally, if $V$ is the middle term of a short exact sequence
$0\rightarrow V^1\rightarrow V\rightarrow V^2\rightarrow 0$, with
$V^1$ and $V^2$ one of the representations of type (a), (b) or
(c), we have the blocks matrix
$$
Hom_Q(d_{min}^V,\cdot)=\left(\begin{array}{cc}Hom_Q(d_{min}^{V^1},\cdot) & 0 \\
Hom_Q(B,\cdot) & Hom_Q(d^{V^2},\cdot)\end{array}\right).
$$
where $d_{min}^{V^1}:P^1_1\rightarrow P^1_0$ is the minimal
projective presentation of $V^1$, $d_{min}^{V^2}:P^2_1\rightarrow
P^2_0$ is the minimal projective presentation of $V^2$ and for
some $B\in Hom_Q(P^2_1,P^1_0)$. In general for every blocks matrix
we have $\left(\begin{array}{cc}A & 0\\0 &
C\end{array}\right)=\left(\begin{array}{cc}Id & 0\\-BA^{-1}&
Id\end{array}\right)\cdot \left(\begin{array}{cc}A & 0\\B &
C\end{array}\right)$ if $A$ is invertible. Hence using rows
operations on $Hom_Q(d_{min}^V,\cdot)$, we obtain
$$
Hom_Q(d_{min}^V,\cdot)\approx\left(\begin{array}{cc}Hom_Q(d_{min}^{V^1},\cdot) & 0 \\
0 & Hom_Q(d_{min}^{V^2},\cdot)\end{array}\right).
$$
So it's enough to prove the skew-symmetry of
$Hom_Q(d_{min}^V,\cdot)$ for $V$ one of
representations of type (a), (b) and (c).\\
\\
\textbf{Type $\boldsymbol{\widetilde{D}}$}. We prove lemma only
for $Q=\widetilde{D}^{0,1}_{n}$ because for the case
$\widetilde{D}^{1,0}_{n}$ it proceeds similarly. We consider the
following labelling for $(\widetilde{D}^{0,1}_{n})^{eq}$:
$$
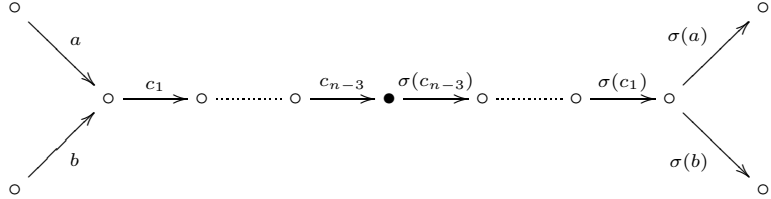
\xymatrix{\circ\ar[dr]^a&&&&&&&&\circ\\
&\circ\ar[r]^{c_1}&\circ\ar@{.}[r]&\circ\ar[r]^{c_{n-3}}&\bullet\ar[r]^{\sigma(c_{n-3})}&\circ\ar@{.}[r]&\circ\ar[r]^{\sigma(c_1)}&\circ\ar[ur]^{\sigma(a)}\ar[dr]_{\sigma(b)}&\\
\circ\ar[ur]_b&&&&&&&&\circ}
$$
We consider again indecomposable representations $V'\in Rep(Q')$
satisfying property \textit{(Op)}. The other regular
indecomposable representations of $Rep(Q')$ satisfying property
\textit{(Op)}. are extensions of these.
\begin{itemize}
\item[(a)] $E_{i,j-1}$, with $1\leq i<j\leq 2n-3$ or $2\leq j<i-1\leq 2n-4$, such that
$\nabla E_{i,j-1}=E_{i+1,j}$; in this case we have
$Z'=E_{i+1,j-1}\oplus E_{i,j}$.
\item[(b)] $E''_0$ and $E''_1$. We note that $\nabla
E''_0=E''_1=\tau^+ E''_0$, $\nabla E''_1=E''_0=\tau^+ E''_1$ and
the respective $Z'$ are
$$
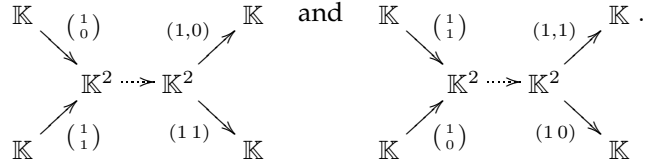
\begin{array}{ccccc}
\xymatrix@-1pc{\mathbb{K}\ar[dr]^{\left({1\atop 0}\right)}&&&\mathbb{K}\\
&\mathbb{K}^2\ar@{.>}[r]&\mathbb{K}^2\ar[dr]_{(1\,1)}\ar[ur]^{(1,0)}&\\
\mathbb{K}\ar[ur]_{\left({1\atop 1}\right)}&&&\mathbb{K}} &\textrm{and}& \xymatrix@-1pc{\mathbb{K}\ar[dr]^{\left({1\atop 1}\right)}&&&\mathbb{K}\\
&\mathbb{K}^2\ar@{.>}[r]&\mathbb{K}^2\ar[dr]_{(1\,0)}\ar[ur]^{(1,1)}&\\
\mathbb{K}\ar[ur]_{\left({1\atop 0}\right)}&&&\mathbb{K}}.
\end{array}
$$
where linear maps defined on $c_i$, with $1\leq i\leq n-3$, are
identity maps.
\item[(c)] $V_{(0,1)}$ and $V_{(1,1)}$; respectively $Z'=E_{h}^{n-1}\oplus E_{0,2n-6}$ and $Z'=E_{h}^1\oplus E_{2n-6,0}$
where $E_{h}^1$ and $E_{h}^{n-1}$ are the regular indecomposable
representations respectively of dimension $e_1+h$ and $e_{n-1}+h$.
\end{itemize}
We consider the following labelling of vertices end arrows for
$\widetilde{D}^{0,1}_{n}$:
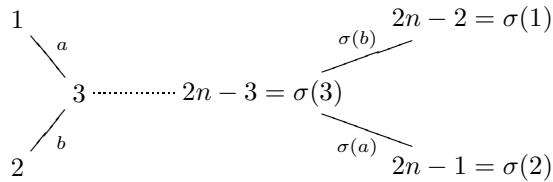
$$
\xymatrix@-1pc{1\ar@{-}[dr]^{a}&&&& 2n-2=\sigma(1)\\
&3\ar@{.}[rr]&& 2n-3=\sigma(3)\ar@{-}[dr]_{\sigma(a)}\ar@{-}[ur]^{\sigma(b)}&\\
2\ar@{-}[ur]_{b}&&&& 2n-1=\sigma(2)}
$$
and we call $c_{i-2}$ the arrow such that $tc_{i-2}=i$.\\
Let $\chi$ be the symmetric weight associated to $V$. We use the
same notation of type $A$ for vertices from 3 to $2n-3$ on which
the components of $\chi$ are not zero. Suppose that 1 and 2 are
source (the other cases are similar). We check when
$Hom_Q(d^V_{min},\cdot)$ is skew-symmetric, for $V$ of type (a),
(b) and (c)
\begin{itemize}
\item[(a)] Let $V$ be one of representation of type (a). We note
that either $\chi(1)=0=\chi(2)$ or $\chi(1)\neq 0\neq\chi(2)$. If
$\chi(1)=0=\chi(2)$, then we have $\chi(m_1)=\pm 1=-\chi(m_s)$ and
$\chi(m_i)=1,-1$, for every $i\in\{2,\ldots,s-1\}$, respectively
if $m_i$ is a source or a sink. Hence it proceeds as in type $A$.
If $\chi(1)\neq 0\neq\chi(2)$ then
$-\chi(2n-2)=\chi(1)=1=\chi(2)=-\chi(2n-1$ and we have
$\chi(m_i)=1,-1$, for every $i\in\{1,\ldots,s\}$, respectively if
$m_i$ is a source or a sink. Let $i_1<\ldots<i_t$ be the sources
from $3$ to $2n-3$ and let $j_1<\ldots<j_t$ be the sinks from $3$
to $2n-3$. We also note that $j_1<i_1<\ldots<j_t<i_t$.\\
$d^V_{min}$ is a matrix $(t+2)\times (t+2)$ whose entries are
$$
(d^V_{min})_{h,k}=\left\{\begin{array}{ll} -a_{i_{k+1},j_k} &
h=t-k\,\textrm{and}\,1\leq k\leq t-1\\
a_{i_k,j_k} & h=t-k+1\,\textrm{and}\,1\leq k\leq t\\
-a_{i,j_1} & h=t+i\,\textrm{and}\,k=1\,\textrm{for}\,i=1,2\\
-\sigma(a_{i,j_1}) &
h=1\,\textrm{and}\,k=t+i\,\textrm{for}\,i=1,2\\
0 & \textrm{otherwise}
\end{array}\right.
$$
where $a_{i,j}$ is oriented path from $i$ to $j$.\\
Finally, as for the type $A$, we note that $Hom_Q(d^V_{min},W)$ is
skew-symmetric if and only if $W\in
SpRep(\widetilde{D}^{0,1}_n,\beta)$.
\item[(b)] Let $V$ be a representation of type (b). We note that
if $\chi$ if the weight associated to $E''_0$, then
$-\chi(2n-2)=\chi(1)=1$ and $\chi(m_i)=1,-1$, for every
$i\in\{1,\ldots,s\}$, respectively if $m_i$ is a source or a sink.
So we can proceed as in type $A$.
\item[(c)] Let $V$ be a representation of type (c). We use the same
notation of part (a) of type $\widetilde{D}$. We note that
$-\chi(2n-2)=\chi(1)=1=\chi(2)=-\chi(2n-1$ and we have
$\chi(m_i)=2,-2$, for every $i\in\{1,\ldots,s\}$, respectively if
$m_i$ is a source or a sink.\\
In the remainder of the proof, we use notation of section
\ref{dVWrf}. In this case, $d^V_{min}$ is a blocks $(2t+2)\times
(2t+2)$-matrix $\left(\begin{array}{cc} A & C\\B &
0\end{array}\right)$. Here
\begin{itemize}
\item[(i)] $A$ is a $2t\times 2t$-matrix with $2\times 2$-blocks $A_{h,k}$,
defined as follows
$$
A_{h,k}=\left\{\begin{array}{ll} (-a_{i_{k+1},j_k})_{Id_2} &
h=t-k\,\textrm{and}\,1\leq k\leq t-1\\
(a_{i_k,j_k})_{Id_2} & h=t-k+1\,\textrm{and}\,1\leq k\leq t\\
0 & \textrm{otherwise}.\end{array}\right.
$$
\item[(ii)] $B$ is a $2\times 2t$-matrix, whose entries $b_{h,k}$
are
$$
\left\{\begin{array}{cc} (-1)^{h+k+1}a_{h,j_1} &
h=1,2\,\textrm{and}\,k=1,2\\
0 & \textrm{otherwise}.\end{array}\right.
$$
\item[(iii)] $C$ is a $2t\times 2$-matrix, whose entries $c_{h,k}$
are
$$
\left\{\begin{array}{cc} (-1)^{h+k+1}\sigma(a_{k,j_1}) &
h=1,2\,\textrm{and}\,k=1,2\\
0 & \textrm{otherwise}.\end{array}\right.
$$
\end{itemize}
Finally, as for the type $A$, we note that $Hom_Q(d^V_{min},W)$ is
skew-symmetric if and only if $W\in
SpRep(\widetilde{D}^{0,1}_n,\beta)$.
\end{itemize}
At last it remains to prove that lemma is true also for every $V$
decomposable representation. But we note that if $V=V^1\oplus
V^2$, then
\begin{itemize}
\item[(i)] V satisfies property \textit{(Op)} if and only if $V^1$
and $V^2$ satisfy property \textit{(Op)};
\item[(ii)] $d^V_{min}=\left(\begin{array}{cc} d^{V^1}_{min} & 0\\0 & d_{min}^{V^2}
\end{array}\right)$.
\end{itemize}
This concludes the proof. $\Box$

\section{Relations between semi-invariants of $(Q,\sigma)$ and of\\
 $(c_{(x,\sigma(x))}(Q),\sigma)$}
Let $(Q,\sigma)$ be a symmetric quiver and let $x$ be an
admissible sink of $(Q,\sigma)$. First we consider the action of
$c_{(x,\sigma(x))}$ on the weights of semi-invariants
\begin{lemma}\label{lecxchi}
Let $(Q,\sigma)$ be a symmetric quiver and let $x$ be an
admissible sink-source of $Q$. If
$\chi=\langle\alpha,\cdot\rangle-\sum_{x\in
  Q_0^{\sigma}}\varepsilon_{x,\alpha}$ is a weight for some dimension vector $\alpha$ (see definition \ref{wsQ}), then
\begin{equation}\label{cxchi}
(c_{(x,\sigma(x))}\chi)(y)=\left\{\begin{array}{ll} -\chi(x) & y=x\\
-\chi(\sigma(x)) & y=\sigma(x)\\
\chi(y)+b_{x,y}\chi(x) &  y\not\in Q_0^{\sigma}\cup\{x\}\\
\chi(y)+b_{\sigma(x),y}\chi(x) &  y\not\in
Q_0^{\sigma}\cup\{\sigma(x)\}\\
0 & \textrm{otherwise},
\end{array}\right.
\end{equation}
where $b_{x,y}$ is the number of arrows linking $x$ and $y$.
\end{lemma}
\textit{Proof.} First we note that, by definition, $\chi(y)=0$ for
every $y\in Q_0^{\sigma}$.
\begin{itemize}
\item[(i)] If $y=x$, then $y\not\in Q_0^{\sigma}$ and
$$
(c_{(x,\sigma(x))}\chi)(x)=(c_{(x,\sigma(x))}\alpha)(x)=\sum_{{a\in
Q_1:\atop ha=x}}\alpha(ta)-\alpha(x)=-\chi(x).
$$
Similarly one proves the case $y=\sigma(x)$.
\item[(ii)] If $y=ta\not\in Q_0^{\sigma}\cup\{x\}$ such that $ha=x$ in $Q$, then $y=hc_{(x,\sigma(x))}a$ such that
$tc_{(x,\sigma(x))}a=x$ in $c_{(x,\sigma(x))}Q$ and
$$
(c_{(x,\sigma(x))}\chi)(y)=
$$
$$
(c_{(x,\sigma(x))}\alpha)(y)-\sum_{{a\in
c_{(x,\sigma(x))}Q_1:\atop ha=y\;\textrm{and}\;ta\neq
x}}(c_{(x,\sigma(x))}\alpha)(ta)-\sum_{{a\in
c_{(x,\sigma(x))}Q_1:\atop ha=y\;\textrm{and}\;ta=
x}}(c_{(x,\sigma(x))}\alpha)(x)=
$$
$$
\alpha(y)-\sum_{{a\in Q_1:\atop ha=y}}\alpha(ta)+\sum_{{a\in
Q_1:\atop ha= x}}(\alpha(x)-\sum_{{a\in Q_1:\atop ha=
x}}\alpha(ta))=
$$
$$
\chi(y)+b_{x,y}\chi(x).
$$
Similarly one proves the case $y=h\sigma(a)\not\in
Q_0^{\sigma}\cup\{\sigma(x)\}$ such that $t\sigma(a)=x$ in $Q$.
\item[(iii)] Finally we have to consider $y$ such that there are no arrows linking $y$ and $x$
 (i.e. $b_{x,y}=0$) and no arrows linking $y$ and $\sigma(x)$. In
 this case
$$
(c_{(x,\sigma(x))}\chi)(y)=
$$
$$
(c_{(x,\sigma(x))}\alpha)(y)-\sum_{{a\in
c_{(x,\sigma(x))}Q_1:\atop ha=y}}(c_{(x,\sigma(x))}\alpha)(ta)=
$$
$$
\alpha(y)-\sum_{{a\in Q_1:\atop ha=y}}\alpha(ta)=
$$
$$
\chi(y).
$$
Similarly one proves for $\sigma(x)$. $\Box$
\end{itemize}
Next we study the relation between $SpSI(Q,\alpha)$ and
$SpSI(c_{(x,\sigma(x))}Q,c_{(x,\sigma(x))}\alpha)$ (respectively
between $OSI(Q,\alpha)$ and\\
$OSI(c_{(x,\sigma(x))}Q,c_{(x,\sigma(x))}\alpha)$) with the
following lemmas
\begin{lemma}\label{SpcxQ}
Let $(Q,\sigma)$ be a symmetric quiver, let $x$ be a sink and let
$\alpha$ be the dimension vector of a symplectic representation.
\begin{itemize}
\item[(i)] If $V\in Rep(Q)$ is indecomposable, not projective, such that
$C_{(x,\sigma(x))}^+V$ is not projective and
$\langle\underline{dim}V,\alpha\rangle=0$, then $c^V\in
SpSI(Q,\alpha)$ and $c^{C_{(x,\sigma(x))}^+V}\in
SpSI(c_{(x,\sigma(x))}Q,c_{(x,\sigma(x))}\alpha)$.
\item[(ii)] If $V=S_x$ and
$\langle\underline{dim}S_x,c_{(x,\sigma(x))}\alpha\rangle=0$, then
$c^{S_x}$ and $c^{C^-S_{\sigma(x)}}$ in\\
$SpSI(c_{(x,\sigma(x))}Q,c_{(x,\sigma(x))}\alpha)$, where $S_x$
and $S_{\sigma(x)}$ are considered as representation of
$c_{(x,\sigma(x))}Q$, but $c^{S_x}$ and $c^{C^-S_{\sigma(x)}}$ are
zero for $Q$. Moreover $c^{S_x}=c^{C^-S_{\sigma(x)}}$.
\item[(iii)] If $V=C^-S_x$ and
$\langle\underline{dim}C^-S_x,\alpha\rangle=0$, then we have
$c^{C^-S_x},c^{S_{\sigma(x)}}\in SpSI(Q,\alpha)$ but they are zero
for $c_{(x,\sigma(x))}Q$. Moreover $c^{S_{\sigma(x)}}=c^{C^-S_x}$.
\end{itemize}
\end{lemma}
\begin{lemma}\label{OcxQ}
Let $(Q,\sigma)$ be a symmetric quiver, let $x$ be a sink and let
$\alpha$ be the vector dimension of an orthogonal representation.
\begin{itemize}
\item[(i)] If $V\in Rep(Q)$ is indecomposable, not projective and such that
$C_{(x,\sigma(x))}^+V$ is not projective and
$\langle\underline{dim}V,\alpha\rangle=0$, then $c^V\in
OSI(Q,\alpha)$ and $c^{C_{(x,\sigma(x))}^+V}\in
OSI(c_{(x,\sigma(x))}Q,c_{(x,\sigma(x))}\alpha)$.
\item[(ii)] If $V=S_x$ and
$\langle\underline{dim}S_x,c_{(x,\sigma(x))}\alpha\rangle=0$, then
we have $c^{S_x}$ and $c^{C^-S_{\sigma(x)}}$ in
$OSI(c_{(x,\sigma(x))}Q,c_{(x,\sigma(x))}\alpha)$, where $S_x$ and
$S_{\sigma(x)}$ are considered as representation of
$c_{(x,\sigma(x))}Q$, but $c^{S_x}$ and $c^{C^-S_{\sigma(x)}}$ are
zero for $Q$. Moreover $c^{S_x}=c^{C^-S_{\sigma(x)}}$.
\item[(iii)] If $V=C^-S_x$ and
$\langle\underline{dim}C^-S_x,\alpha\rangle=0$, then we have
$c^{C^-S_x},c^{S_{\sigma(x)}}\in OSI(Q,\alpha)$ but they are zero
for $c_{(x,\sigma(x))}Q$. Moreover $c^{S_{\sigma(x)}}=c^{C^-S_x}$.
\end{itemize}
\end{lemma}
We prove only lemma \ref{SpcxQ} because the proof of lemma
\ref{OcxQ} is similar.\\
 \textit{Proof.} First of all we note that if $x$ is an admissible sink, then
$S_{\sigma(x)}\neq\tau^-\nabla S_{\sigma(x)}$ and
$C^-S_{x}\neq\tau^-\nabla C^-S_{x}$ and so, by lemma \ref{ss1}, we
can not define both $pf^{S_{\sigma(x)}}$ and $pf^{C^-S_{x}}$. It's
enough to prove the first one because, by lemma \ref{C=tau},
$\tau^-\nabla C^-S_x=\tau^-\nabla
\tau^-S_x=\nabla\tau^+\tau^-S_x=\nabla S_x=S_{\sigma(x)}$. If
$S_{\sigma(x)}=\tau^-\nabla S_{\sigma(x)}$, by theorem \ref{ar}
there exists an almost split sequence
\begin{equation}\label{assSx}
0\longrightarrow\nabla S_{\sigma(x)}=S_x\longrightarrow
Z\longrightarrow S_{\sigma(x)}\longrightarrow 0.
\end{equation}
Hence $(\underline{dim}\,Z)_y=\left\{\begin{array}{ll} 1 &
\textrm{if}\;y=x,\sigma(x)\\
0 & \textrm{otherwise}\end{array}\right.$ and so either
$Z=S_x\oplus S_{\sigma(x)}$ which is an absurd because
(\ref{assSx}) would be a split sequence, or $Z$ is indecomposable
and thus there is an arrow $\sigma(x)\rightarrow x$ which is not
possible since
$x$ is an admissible sink.\\
We recall that
$(\underline{dim}S_{\sigma(x)})_y=\left\{\begin{array}{ll} 1 &
\textrm{if}\;\sigma(x)=y\\
0 & \textrm{otherwise}\end{array}\right.$, that
$\alpha_x=\alpha_{\sigma(x)}$ for every $x\in Q_0$ and, by theorem
\ref{ard}, that $\langle
\underline{dim}C^-S_{x},\alpha\rangle=-\langle\alpha,\underline{dim}S_{x}\rangle$.
So, for a dimension vector $\alpha$ of a symplectic (respectively
orthogonal) representation, $\langle
\underline{dim}S_{\sigma(x)},\alpha\rangle=\alpha_{\sigma(x)}-\sum_{a\in
Q_1:ha=x}\alpha_{\sigma(ta)} = \alpha_x-\sum_{a\in
Q_1:ha=x}\alpha_{ta}=\langle\alpha,\underline{dim}S_{x}\rangle=-\langle\underline{dim}C^-S_{x},\alpha\rangle$.
Similarly we have $\langle
\underline{dim}S_x,c_{(x,\sigma(x))}\alpha\rangle=-\langle\underline{dim}C^-S_{\sigma(x)},c_{(x,\sigma(x))}\alpha\rangle$.
Hence, since $x$ is a sink of $Q$ and $\sigma(x)$ is a sink of
$c_{(x,\sigma(x))}Q$, it's enough to apply lemma \ref{SIcxQ} to
both $Q$ and $c_{(x,\sigma(x))}Q$. Finally $\tau^-\nabla
S_{\sigma(x)}=\tau^-S_x=C^-S_x$ and $\tau^-\nabla C^-S_{\sigma(x)}
=\tau^-\tau^+\nabla S_{\sigma(x)}=S_x$, so, by lemma
\ref{cV=cVnabla}, $c^{S_x}=c^{C^-S_{\sigma(x)}}$
and $c^{S_{\sigma(x)}}=c^{C^-S_x}$. $\Box$\\
\\
We observe that, by proposition \ref{C+nabla}, $\tau^-\nabla V=V$
if and only if $\tau^-\nabla
C^+_{(x,\sigma(x))}V=C^+_{(x,\sigma(x))}V$. Let $\alpha$ be a
symmetric dimension vector. We recall that
$\alpha_{y}=c_{(x,\sigma(x))}\alpha_{y}$ for every $y\neq
x,\sigma(x)$ and $(c_{(x,\sigma(x))}\alpha)_{x}=\sum_{a\in
Q_1:ha=x}\alpha_{ta}-\alpha_{x}= \sum_{a\in
Q_1:ha=x}\alpha_{\sigma(ta)}-\alpha_{\sigma(x)}
=(c_{(x,\sigma(x))}\alpha)_{\sigma(x)}$, so we consider three
cases.
\begin{itemize}
\item[(i)] $0\neq\alpha_x\neq\sum_{a\in
Q_1:ha=x}\alpha_{ta}$, i.e. $\langle
\underline{dim}S_{\sigma(x)},\alpha\rangle\neq 0$ and $\langle
\underline{dim}S_x,c_{(x,\sigma(x))}\alpha\rangle\neq 0$.
\item[(ii)] $0=\alpha_x\neq\sum_{a\in
Q_1:ha=x}\alpha_{ta}$, i.e. $\langle
\underline{dim}S_{\sigma(x)},\alpha\rangle\neq 0$ and $\langle
\underline{dim}S_x,c_{(x,\sigma(x))}\alpha\rangle=0$.
\item[(iii)] $0\neq\alpha_x=\sum_{a\in
Q_1:ha=x}\alpha_{ta}$, i.e. $\langle
\underline{dim}S_{\sigma(x)},\alpha\rangle= 0$ and $\langle
\underline{dim}S_x,c_{(x,\sigma(x))}\alpha\rangle\neq 0$.
\end{itemize}
We note that $0=\alpha_x=\sum_{a\in Q_1:ha=x}\alpha_{ta}$ is not
possible, unless $\alpha_{ta}=0$ for every $a$ such that $ha=x$.
\begin{proposizione}\label{SIcx1}
Let $(Q,\sigma)$ be a symmetric quiver. Let $\alpha$ be a
symmetric dimension vector, $x$ be an admissible
sink and $\varphi^{Sp}_{x,\alpha}$ be as defined in lemma \ref{kac}.\\
Then $\varphi^{Sp}_{x,\alpha}(c^V)=c^{C_{(x,\sigma(x))}V}$ and
$\varphi^{Sp}_{x,\alpha}(pf^W)=pf^{C_{(x,\sigma(x))}W}$, where $V$
and $W$ are indecomposables of $Q$ such that $\langle
\underline{dim}\,V,\alpha
\rangle=0=\langle\underline{dim}\,W,\alpha \rangle$ and $W$
satisfies property \textit{(Op)}. In particular
\begin{itemize}
\item[(i)] if $0=\alpha_x\neq\sum_{a\in Q_1:ha=x}\alpha_{ta}$, then
$(\varphi^{Sp}_{x,\alpha})^{-1}(c^{S_x})=0$;
\item[(ii)] if $0\neq\alpha_x=\sum_{a\in Q_1:ha=x}\alpha_{ta}$, then
$\varphi^{Sp}_{x,\alpha}(c^{S_{\sigma(x)}})=0$.
\end{itemize}
\end{proposizione}
\textit{Proof.} We consider the same notation of proof of lemma
\ref{kac}. If $x$ is an admissible sink of $(Q,\sigma)$, then we
have
$$
C_{(x,\sigma(x))}^-(Z\times
Hom(V',W))=C_{(x,\sigma(x))}^-(SpRep(c_{(x,\sigma(x))}Q,c_{(x,\sigma(x))}\alpha))=
$$
$$
SpRep(Q,\alpha)=Z\times Hom(W,V).
$$
So, by definition,
$$C_{(x,\sigma(x))}^-|_{Z}(SpRep(c_{(x,\sigma(x))}Q,c_{(x,\sigma(x))}\alpha))=Z$$
and
$$C_{(x,\sigma(x))}^-|_{Hom(V',W)}(SpRep(c_{(x,\sigma(x))}Q,c_{(x,\sigma(x))}\alpha))=Hom(W,V).$$
Now $C_{(x,\sigma(x))}^-$ induces a ring morphism
$$
\begin{array}{rcl}
\phi_{x,\alpha}^{Sp}:\mathbb{K}[SpRep(Q,\alpha)]&\longrightarrow&\mathbb{K}[SpRep(c_{(x,\sigma(x))}Q,c_{(x,\sigma(x))}\alpha)]\\
f&\longmapsto& f\circ C_{(x,\sigma(x))}^-\end{array}
$$
By proof of lemma \ref{kac}, we note that
$$
\mathbb{K}[C_{(x,\sigma(x))}^-Z\times
C_{(x,\sigma(x))}^-Hom(V',W)]^{SSp(Q,\alpha)}=\mathbb{K}[Z\times
Hom(W,V)]^{SSp(Q,\alpha)}
$$
is isomorphic by $\varphi_{x,\alpha}^{Sp}$ to $\mathbb{K}[Z\times
Hom(V',W)]^{SSp(c_{(x,\sigma(x))}Q,c_{(x,\sigma(x))}\alpha)}$.
Hence
$\varphi_{x,\alpha}^{Sp}=\phi_{x,\alpha}^{Sp}|_{SpSI(Q,\alpha)}$
and so for every representation $Z$ of dimension vector $\alpha$
of $(Q,\sigma)$ we have
\begin{equation}\label{varphi}
\varphi_{x,\alpha}^{Sp}(c^V)(C_{(x,\sigma(x))}^+Z)=(c^V\circ
C_{(x,\sigma(x))}^-)(C_{(x,\sigma(x))}^+Z)=c^V(Z)
\end{equation}
and
\begin{equation}\label{varphipf}\varphi_{x,\alpha}^{Sp}(pf^W)(C_{(x,\sigma(x))}^+Z)=(pf^W\circ
C_{(x,\sigma(x))}^-)(C_{(x,\sigma(x))}^+Z)=pf^W(Z).\end{equation}
By lemma \ref{cVW=c+xVW} and \ref{cVW=c-xVW} we have
$c^V(Z)=\lambda\cdot
c^{C_{(x,\sigma(x))}^+V}(C_{(x,\sigma(x))}^+Z)$, for some
$\lambda\in\mathbb{K}$. So, by (\ref{varphi}),
$\varphi_{x,\alpha}^{Sp}$ sends $c^V$ to
$c^{C_{(x,\sigma(x))}^+V}$ up to a constant in $\mathbb{K}$.
Similarly for $pf^W$. Finally \textit{(i)} and \textit{(ii)}
follow by lemma \ref{SpcxQ}. $\Box$
\begin{proposizione}\label{SIcx2}
Let $(Q,\sigma)$ be a symmetric quiver. Let $\alpha$ be a
symmetric dimension vector, $x$ be an admissible
sink and $\varphi^{O}_{x,\alpha}$ be as defined in lemma \ref{kac}.\\
Then $\varphi^{O}_{x,\alpha}(c^V)=c^{C_{(x,\sigma(x))}V}$ and
$\varphi^{O}_{x,\alpha}(pf^W)=pf^{C_{(x,\sigma(x))}W}$, where $V$
and $W$ are indecomposables of $Q$ such that $\langle
\underline{dim}\,V,\alpha
\rangle=0=\langle\underline{dim}\,W,\alpha \rangle$ and $W$
satisfies property \textit{(Spp)}. In particular
\begin{itemize}
\item[(i)] if $0=\alpha_x\neq\sum_{a\in Q_1:ha=x}\alpha_{ta}$, then
$(\varphi^{O}_{x,\alpha})^{-1}(c^{S_x})=0$;
\item[(ii)] if $0\neq\alpha_x=\sum_{a\in Q_1:ha=x}\alpha_{ta}$, then
$\varphi^{O}_{x,\alpha}(c^{S_{\sigma(x)}})=0$.
\end{itemize}
\end{proposizione}
\textit{Proof.} It
is similar to that one of proposition \ref{SIcx1}. $\Box$\\
\\
By previous propositions and by lemma \ref{kac} it follows that if
the conjectures \ref{mt1} and \ref{mt2} are true for a symmetric
quiver $(Q,\sigma)$, then they are true for
$(c_{(x,\sigma(x))}Q,\sigma)$.

\section{Composition lemmas}
We conclude this chapter with general lemmas which will be useful
in our proofs.
\begin{lemma}\label{cl}
Let
$$
(Q,\sigma):\cdots
y\stackrel{a}{\longrightarrow}x\stackrel{b}{\longrightarrow}z\cdots\sigma(z)\stackrel{\sigma(b)}{\longrightarrow}\sigma(x)\stackrel{\sigma(a)}{\longrightarrow}\sigma(y)\cdots
$$
be a symmetric quiver. Assume the underlying quiver with $n$
vertices. Also assume there exist only two arrows in $Q_1^+$
incident to $x\in Q_0^+$, $a:y\rightarrow x$ and $b:x\rightarrow
z$ with $y,z\in Q_0^+\cup Q_0^{\sigma}$. Let $V$
be an orthogonal or symplectic representation with symmetric dimension vector $(\alpha_1,\ldots,\alpha_n)=\alpha$ such that $\alpha_x\geq max\{\alpha_y,\alpha_z\}$.\\
We define the symmetric quiver $Q'=((Q_0',Q_1'),\sigma)$ with
$n-2$ vertices such that $Q_0'=Q_0\setminus\{x,\sigma(x)\}$ and
$Q_1'=Q_1\setminus\{a,b,\sigma(a),\sigma(b)\}\cup\{ba,\sigma(a)\sigma(b)\}$,
i.e.
$$
Q':\cdots y\stackrel{ba}{\longrightarrow}
z\cdots\sigma(z)\stackrel{\sigma(a)\sigma(b)}{\longrightarrow}
\sigma(y)\cdots,
$$
and let $\alpha'$ be the dimension of $V$ restricted to $Q'$. \\
We have:
\begin{itemize}
\item[(Sp)] Assume $V$ symplectic. Then
\begin{itemize}
\item[(a)] if $\alpha_x> max\{\alpha_y,\alpha_z\}$ then
$SpSI(Q,\alpha)=SpSI(Q',\alpha')$,
\item[(b)] if $\alpha_x=\alpha_y>\alpha_z$ then
$SpSI(Q,\alpha)=SpSI(Q',\alpha')[detV(a)]$,
\item[(b')] if $\alpha_x=\alpha_z>\alpha_y$ then
$SpSI(Q,\alpha)=SpSI(Q',\alpha')[detV(b)]$,
\item[(c)] if $\alpha_x=\alpha_y=\alpha_z$ then
 $SpSI(Q,\alpha)=SpSI(Q',\alpha')[detV(a),detV(b)]$.
\end{itemize}
\item[(O)] Assume $V$ orthogonal. Then
\begin{itemize}
\item[(a)] if $\alpha_x> max\{\alpha_y,\alpha_z\}$ then
$OSI(Q,\alpha)=OSI(Q',\alpha')$,
\item[(b)] if $\alpha_x=\alpha_y>\alpha_z$ then
$OSI(Q,\alpha)=OSI(Q',\alpha')[detV(a)]$,
\item[(b')] if $\alpha_x=\alpha_z>\alpha_y$ then
$OSI(Q,\alpha)=OSI(Q',\alpha')[detV(b)]$,
\item[(c)] if $\alpha_x=\alpha_y=\alpha_z$ then
 $OSI(Q,\alpha)=OSI(Q',\alpha')[detV(a),detV(b)]$.
\end{itemize}
\end{itemize}
\end{lemma}
\textit{Proof.} We use the notation of section \ref{A1}.\\
\textit{(Sp)}
Using Cauchy formula (theorem \ref{c}) we have
$$
SpSI(Q,\alpha)=\left(\bigoplus_{{\lambda:Q^+_1\rightarrow\Lambda
\atop \mu:Q^{\sigma}_1\rightarrow ER\Lambda} }\;\bigotimes_{c\in
Q^+_1}(S_{\lambda(c)}V_{tc}\otimes
S_{\lambda_(c)}V_{hc}^*)\otimes\left(\bigotimes_{d\in
Q_1^{\sigma}}S_{\mu(d)}V_{td}\right)\right)^{SSp(Q,\alpha)}
$$
where $\Lambda$ is the set of all partitions and $ER\Lambda$ is the set of the partitions with even rows.\\
\textit{(a)} If $\alpha_x> max\{\alpha_y,\alpha_z\}$, by theorem
\ref{pS},
$$
S_{\lambda(a)}V_{x}^*=S_{(\underbrace{\scriptstyle{0,\ldots,0}}_{\alpha_x-\alpha_y},\underbrace{\scriptstyle{-\lambda(a)_{\alpha_y},\ldots,-\lambda(a)_{1}}}_{\alpha_y})}V_{x},
$$
where $\lambda(a)=(\lambda(a)_1,\ldots,\lambda(a)_{\alpha_y})$. By
proposition \ref{i2}, $\lambda(a)$ and $\lambda(b)$ have to
satisfy the following equations
\begin{equation}
\left\{\begin{array}{lr}
\lambda(b)_i-\lambda(b)_{i+1}=0, & i\in\{\alpha_y+1,\ldots,\alpha_x-1\}\\
\lambda(b)_{\alpha_y}-\lambda(b)_{\alpha_y+1}=\lambda(a)_{\alpha_y}
& \\
\lambda(b)_{\alpha_y-i}-\lambda(b)_{\alpha_y-i+1}=\lambda(a)_{\alpha_y-i}-\lambda(a)_{\alpha_y-i+1},
& i\in\{1,\ldots,\alpha_y-1\}.
\end{array}\right.
\end{equation}
We call $\lambda(b)_i=k\geq 0$ for every
$i\in\{\alpha_y+1,\ldots,\alpha_x\}$ and so
$$
\lambda(b)=(\lambda(b)_1,\ldots,\lambda(b)_{\alpha_x})=(\underbrace{\lambda(a)_1+k,\ldots,\lambda(a)_{\alpha_y}+k}_{\alpha_y},\underbrace{k,\ldots,k}_{\alpha_x-\alpha_y}).
$$
Now, by theorem \ref{pS}, $S_{\lambda(b)}V_z^*=0$ unless
$ht(\lambda(b))\leq\alpha_z$. If $\alpha_y\leq\alpha_z$, then
$S_{\lambda(b)}V_z^*=0$ unless
$\lambda(b)_{\alpha_z+1}=\ldots=\lambda(b)_{\alpha_x}=0$, i.e.
$k=0$, so
$\lambda(b)=(\lambda(a)_1,\ldots,\lambda(a)_{\alpha_y},\underbrace{0,\ldots,0}_{\alpha_x-\alpha_y})=\lambda(a)$.
If $\alpha_z<\alpha_y$, then $S_{\lambda(b)}V_z^*=0$ unless
$\lambda(b)_{\alpha_z+1}=\ldots=\lambda(b)_{\alpha_x}=0$, i.e.
$k=0$ and
$\lambda(a)_{\alpha_z+1}=\ldots=\lambda(a)_{\alpha_y}=0$, so
$\lambda(b)=\lambda(a)$ again.\\
So let $\lambda(a)=\lambda(b)=\bar{\lambda}$. By proposition
\ref{i2}, $S_{\bar{\lambda}}V_{x}^*\otimes S_{\bar{\lambda}}V_x$
contains a semi-invariant of weight zero, which is hence a
$GL(V_x)$-invariant. Since $V_{y}^*\otimes V_x \oplus V_x^*\otimes
V_{z}=V_{x}^{\alpha_{y}}\oplus(V_x^*)^{\alpha_{z}}$ and since
$S_{\bar{\lambda}}V_{x}^*\otimes S_{\bar{\lambda}}V_x$ is a
summand in the Cauchy formula of
$\mathbb{K}[V_{x}^{\alpha_{y}}\oplus(V_x^*)^{\alpha_{z}}]$, using
FFT for $GL$ (theorem \ref{fft}) we obtain $SL(V)$ acts trivially
on $S_{\bar{\lambda}}V_x^*\otimes S_{\bar{\lambda}}V_{x}$ and so
$(S_{\lambda(a)}V_{x}^*\otimes
S_{\lambda(b)}V_x)^{SL\,V_x}=\mathbb{K}$. So we have
$$
SpSI(Q,\alpha)\cong SpSI(Q',\alpha').
$$
\textit{(b)} If $\alpha_{x}=\alpha_y>\alpha(z)$, by theorem
\ref{pS},
$$
S_{\lambda(a)}V_{x}^*=S_{(-\lambda(a)_{\alpha_y=\alpha_x},\ldots,-\lambda(a)_{1})}V_{x}.
$$
 By proposition \ref{i2}, $\lambda(a)$ and $\lambda(b)$ have to satisfy $\lambda(a)_i-\lambda(a_{i+1})=\lambda(b)_i-\lambda(b)_{i+1}$ for every
 $i\in\{1,\ldots,\alpha_x\}$ and moreover $S_{\lambda(b)}V_{z}^*=0$ unless
 $ht(\lambda(b))\leq\alpha_z<\alpha_x$. Hence we have
 \begin{equation}
 \left\{\begin{array}{ll}
 \lambda(b)_i=0 & i\in\{\alpha_{z+1},\ldots,\alpha_x\}\\
 \lambda(a)_i-\lambda(a)_{i+1}=\lambda(b)_i-\lambda(b)_{i+1} & i\in\{1,\ldots,\alpha_x-1\}
 \end{array}\right.
 \end{equation}
 and thus
\begin{equation}
 \left\{\begin{array}{lll}
 \lambda(a)_i-\lambda(a)_{i+1}=\lambda(b)_i-\lambda(b)_{i+1} &
 i\in\{1,\ldots,\alpha_z-1\}\\
 \lambda(a)_{\alpha_z}=\lambda(a)_{\alpha_z+1}+\lambda(b)_{\alpha_z} &\\
 \lambda(a)_i=\lambda(a)_{i+1} &
 i\in\{\alpha_z+1,\ldots,\alpha_x-1\}.
\end{array}\right.
\end{equation}
 Hence $\lambda(a)$ contains a column of length $\alpha_x=\alpha_y$  for some $k\in\mathbb{N}$, so
 we have
 $\lambda(a)=(\lambda(b)_1+k,\ldots,\lambda(b)_{\alpha_z}+k,k,\ldots,k)$
 then $S_{\lambda(a)}V_y\otimes
 S_{\lambda(a)}V_x^*=S_{\lambda(b)}V_y\otimes(\bigwedge^{\alpha_y} V_y)^k\otimes(\bigwedge^{\alpha_x}
 V_x^*)^k\otimes S_{\lambda(b)}V_x^*$. Now $(\bigwedge^{\alpha_y}V_y)^k\otimes(\bigwedge^{\alpha_x}
 V_x^*)^k$ is spanned by $(det\,V(a))^k$. So we have a semi-invariant $f$ of the form $(detV(a))^k
 f'$ where $f'$ is of weight zero, hence using theorem FFT for $GL$ (\ref{fft}) as before and by lemma
 \ref{VxW}, we have
$$
SpSI(Q,\alpha)= SpSI(Q',\alpha')[detV(a)].
$$
In the similar way we prove \textit{(b')}.\\
\textit{(c)} If $\alpha(x)=\alpha(y)=\alpha(z)$, by theorem
\ref{pS},
$$
S_{\lambda(a)}V_x^*=S_{(-\lambda(a)_{\alpha_y=\alpha_x},\ldots,-\lambda(a)_{1})}V_x
$$
and
$$
S_{\lambda(b)}V_z^*=S_{(-\lambda(b)_{\alpha_x=\alpha_z},\ldots,-\lambda(b)_{1})}V_z,
$$
where $\lambda(a)=(\lambda(a)_{1},\ldots,\lambda(a)_{\alpha_y})$
and $\lambda(b)=(\lambda(b)_{1},\ldots,\lambda(b)_{\alpha_x})$. By
proposition \ref{i2}, $\lambda(a)$ and $\lambda(b)$ have to
satisfy the following equations
\begin{eqnarray}\label{eq}
\lambda(a)_{i-1}-\lambda(a)_i=\lambda(b)_{i-1}-\lambda(b)_i
\end{eqnarray}
for every $i\in\{2,\ldots,\alpha_x=\alpha_y\}$. Thus
$\lambda(a)_i=\lambda(b)_{i}-\lambda(b)_{\alpha_x}+\lambda(a)_{\alpha_x}$
for every $i\in\{1,\ldots,\alpha_x\}$. Hence if we set
$\lambda(b)_{\alpha_x}=h$ and $\lambda(a)_{\alpha_x}=k$ we have
\begin{equation}\label{lambda(a)}
\lambda(a)_i=\lambda(b)_{i}-h+k
\end{equation}
for every $i\in\{1,\ldots,\alpha_x\}$. So in our case
$\lambda(a)=(\lambda(b)-(h^{\alpha_x}))+(k^{\alpha_x})$ and
$\lambda(b)=(\lambda(a)-(k^{\alpha_x}))+(h^{\alpha_x})$. We call
$\lambda(b)-(h^{\alpha_x})=\lambda(b)'$ and
$\lambda(a)-(k^{\alpha_x})=\lambda(a)'$ and we note that
$\lambda(a)'=\lambda(b)'$ by the system (\ref{lambda(a)}). Then
$S_{\lambda(a)}V_y\otimes S_{\lambda(a)}V_x^*\otimes
S_{\lambda(b)}V_x\otimes
S_{\lambda(b)}V_z^*=S_{\lambda(b)'}V_y\otimes<(detV(a))^k>\otimes
S_{\lambda(b)'}V_x^*\otimes
S_{\lambda(a)'}V_x\otimes<(detV(b))^h>\otimes
S_{\lambda(a)'}V_z^*$. So we have a semi-invariant $f$ of the form
$(detV(a))^k (detV(b))^h
 f'$ where $f'$ is of weight zero, hence using theorem FFT for $GL$ (\ref{fft}) as before and by lemma
 \ref{VxW}, we have
$$
SpSI(Q,\alpha)\cong SpSI(Q',\alpha')[detV(a),detV(b)].
$$
\textit{(O)} Using Cauchy formula we have
$$
OSI(Q,\alpha)=\left(\bigoplus_{{\lambda:Q^+_1\rightarrow\Lambda
\atop \mu:Q^{\sigma}_1\rightarrow EC\Lambda} }\;\bigotimes_{c\in
Q^+_1}(S_{\lambda(c)}V_{tc}\otimes
S_{\lambda_(c)}V_{hc}^*)\otimes\left(\bigotimes_{d\in
Q_1^{\sigma}}S_{\mu(d)}V_{td}\right)\right)^{SO(Q,\alpha)}
$$
where $\Lambda$ is the set of all partitions and $EC\Lambda$ is
the set of the partitions with even columns. The rest of the proof
is similar of the symplectic case. \quad $\Box$
\begin{lemma}\label{cls}
Let $(Q,\sigma)$ be a symmetric quiver with $n$ vertices such that
there exist only two arrows $a$ and $b$ incident to the vertex $x$
in $Q_0$ and $b$ is fixed by $\sigma$, i.e.
$$
Q:\cdots y\stackrel{a}{\longrightarrow}
x\stackrel{b}{\longrightarrow}
\sigma(x)\stackrel{\sigma(a)}{\longrightarrow}\sigma(y)\cdots
$$
Let
$$
V:\cdots V_y\stackrel{V(a)}{\longrightarrow}
V_x\stackrel{V(b)}{\longrightarrow}
V_x^*\stackrel{-V(a)^t}{\longrightarrow} V_y^*\cdots
$$
be an orthogonal or symplectic representation of $(Q,\sigma)$ with
$\underline{dim}V=\alpha$ such that $\alpha_x\geq\alpha_y$.
Moreover define the symmetric quiver
$(Q',\sigma)=((Q_0',Q_1'),\sigma)$ with $n-2$ vertices such that
$Q_0'=Q_0\setminus\{x,\sigma(x)\}$ and
$Q_1'=Q_1\setminus\{a,b,\sigma(a)\}\cup\{\sigma(a)ba\}$, i.e
$$
Q':\cdots y\stackrel{\sigma(a)ba}{\longrightarrow}\sigma(y)\cdots.
$$
Let $\alpha'$ be the dimension of $V$ restricted to $Q'$.
\begin{itemize}
\item[(Sp)] If $V$ is symplectic, then
\begin{itemize}
\item[(i)]  $\alpha_x>\alpha_y\Longrightarrow
SpSI(Q,\alpha)=SpSI(Q',\alpha')[detV(b)]$
\item[(ii)]  $\alpha_x=\alpha_y\Longrightarrow
SpSI(Q,\alpha)=SpSI(Q',\alpha')[detV(a)]$.
\end{itemize}
\item[(O)] If $V$ is orthogonal, then
\begin{itemize}
\item[(i)]  $\alpha_x>\alpha_y\;\textrm{and}\;\alpha_x\;\textrm{is even}\Longrightarrow
OSI(Q,\alpha)=OSI(Q',\alpha')[pfV(b)]$
\item[(ii)]  $\alpha_x=\alpha_y\Longrightarrow
OSI(Q,\alpha)=OSI(Q',\alpha')[detV(a)]$.
\end{itemize}
\end{itemize}
\end{lemma}
\textit{Proof}. We consider again the Cauchy
formulas.\\
\textit{(Sp)} If $\alpha_y\leq\alpha_x$, by proposition \ref{i2},
$\lambda(a)$ and $\lambda(b)$ have to satisfy
$\lambda(a)_{i-1}-\lambda(a)_i=\lambda(b)_{i-1}-\lambda(b)_i$ for
every $i\in\{2,\ldots,\alpha_y\}$.\\
\textit{(i)} Let $\alpha_y<\alpha_x$, we have
$$
S_{\lambda(a)}V^*_x=S_{(0,\ldots,0,-\lambda(a)_{\alpha_y},\ldots,-\lambda(a)_1)}V_x
$$
and so
$$
\lambda(b)=(\overbrace{\lambda(a)_1,\ldots,\lambda(a)_{\alpha_y},0,\ldots,0}^{\alpha_x})+(\overbrace{2k,\ldots,2k}^{\alpha_x}),
$$
for some $k\in \mathbb{Z}_{\geq 0}$ and with $\lambda(a)_i$ even
for every $i$. Then $S_{\lambda(a)}V_x^*\otimes
 S_{\lambda(b)}V_x=S_{\lambda(a)}V_x^*\otimes S_{\lambda(a)}V_x \otimes(\bigwedge^{\alpha_x} V_x)^{2k}$. Now $(\bigwedge^{\alpha_x}
 V_x)^{2k}$ is spanned by $(det\,V(b))^k$. So we have a semi-invariant $f$ of the form $(detV(b))^k
 f'$ where $f'$ is of weight zero, hence using theorem FFT for $GL$ (\ref{fft}) as before and by lemma
 \ref{VxW}, we have
$$
SpSI(Q,\alpha)\cong SpSI(Q',\alpha')[det\;V(b)].
$$
\textit{(ii)} If $\alpha_x=\alpha_y$, the proof is similar to the part \textit{(b)} of lemma \ref{cl}.\\
\\
\textit{(O)} If $\alpha_y\leq\alpha_x$, by proposition \ref{i2},
$(S_{\lambda(a)}V^*_x\otimes S_{\lambda(b)}V_x)^{SL(V_x)}\ne 0$ if
and only if
$\lambda(a)_{i-1}-\lambda(a)_i=\lambda(b)_{i-1}-\lambda(b)_i$ for
every $i\in\{2,\ldots,\alpha_y\}$.\\
Now the proof is similar to the symplectic case, recalling that
$V(b)$, in this case, is skew-symmetric, so we can define
$pf\;V(b)$. $\Box$

\chapter{Semi-invariants of symmetric quivers of finite type}
In this chapter we prove conjectures \ref{mt1} and \ref{mt2} for
the symmetric quivers of finite type. We recall that, by theorem
\ref{ctf}, a symmetric quiver of finite type has the underlying
quiver of type $A_n$. Throughout this chapter we enumerate
vertices with $1,\ldots,n$ from left to right and we call $a_i$
the arrow with $i$ on the left and $i+1$ on the right; moreover we
define $\sigma$ by $\sigma(i)=n-i+1$, for every
$i\in\{1,\ldots,n\}$, and $\sigma(a_i)=a_{n-i}$, for every
$i\in\{1,\ldots,n-1\}$.\\
First we prove a lemma valid for $Q=A_n$, which is a particular
case of lemma \ref{ss1}.
\begin{lemma}\label{ss2}
Let $(A_n,\sigma)$ be a symmetric quiver of type $A$. Let $V\in
Rep(Q)$ such that $V=\tau^-\nabla V$ and let $W$ a selfdual
representation such that
$\langle\underline{dim}V,\underline{dim}W\rangle$=0, then we have
the following.
\begin{itemize}
\item[(i)] If $n$ is even, $d^V_W$ is skew-symmetric if and only if $W\in
ORep(Q,\underline{dim}W)$.
\item[(ii)] If $n$ is odd $d^V_W$ is skew-symmetric if and only if $W\in
SpRep(Q,\underline{dim}W)$.
\end{itemize}
\end{lemma}
\textit{Proof.} It checked in the proof of lemma \ref{ss1}. $\Box$ \\
\\
By proof of lemma \ref{ss1} we noted also that an indecomposable
representation $V$ of $A_n$ satisfies property \textit{(Spp)} if
$n$ is even and it satisfies property
\textit{(Op)} if $n$ is odd.\\
The conjectures \ref{mt1} and \ref{mt2} for symmetric quivers of
finite type become
\begin{teorema}\label{tfs}
Let $(Q,\sigma)$ be a symmetric quiver of finite type. Let
$\alpha$ be the dimension vector of a symplectic representation.
Then $SpSI(Q,\alpha)$ is generated by the following
semi-invariants.
\begin{itemize}
\item[($n$ even)] $c^V$ with $V$ indecomposable in
$Rep(Q)$ such that $\langle\underline{dim}V,\alpha\rangle=0$.
\item[($n$ odd)] $(i)$ $c^V$ with $V$ indecomposable in
$Rep(Q)$ such that $\langle\underline{dim}V,\alpha\rangle=0$;\\
$(ii)$ $pf^V$ with $V\in Rep(Q)$ such that $V=\tau^-\nabla V$.
\end{itemize}
\end{teorema}
\begin{teorema}\label{tfo}
Let $(Q,\sigma)$ be a symmetric quiver of finite type. Let
$\alpha$ be the dimension vector of an orthogonal representation.
Then $OSI(Q,\alpha)$ is generated by the following
semi-invariants.
\begin{itemize}
\item[($n$ odd)] $c^V$ with $V$ indecomposable in
$Rep(Q)$ such that $\langle\underline{dim}V,\alpha\rangle=0$.
\item[($n$ even)] $(i)$ $c^V$ with $V$ indecomposable in
$Rep(Q)$ such that $\langle\underline{dim}V,\alpha\rangle=0$;\\
$(ii)$ $pf^V$ with $V\in Rep(Q)$ such that $V=\tau^-\nabla V$.
\end{itemize}
\end{teorema}
By proposition \ref{Q=Q'} and by propositions \ref{SIcx1} and
\ref{SIcx2}, it's enough to study the equioriented case, i.e. the
case in which all the arrows have orientation from left to right.
\begin{lemma}\label{spsiosipr}
Let $(Q,\sigma)$ be a symmetric quiver of finite type. Then
$SpSI(Q,\beta)$ and $OSI(Q,\beta)$ are polynomial rings, for every
symmetric dimension vector $\beta$.
\end{lemma}
\textit{Proof.} Since the isomorphism classes of
$\beta$-dimensional symplectic (resp. orthogonal) representations
of $(Q,\sigma)$ correspond to the orbits of the action of
$Sp(Q,\beta)$ (resp. of $O(Q,\beta)$) on $SpRep(Q,\beta)$ (resp.
on $ORep(Q,\beta)$), then lemma follows by definition of symmetric
quiver finite type and by lemma \ref{sk}. $\Box$

\section{Equioriented symmetric quivers of finite type}
In this section we state and prove case by case theorems \ref{tfs}
and \ref{tfo} for equioriented case. Throughout this section we
call $V_{j,i}$ the indecomposable of $A_n$ with dimension vector
$$(v_{j,i})_{k}=\left\{\begin{array}{ll} 1 & \textrm{if}\;\;j\leq
k\leq i \\
0 & \textrm{otherwise}.
\end{array}\right.$$
\subsection{The symplectic case for $A_{2n}$}
We rewrite theorem \ref{tfs} in the following way
\begin{teorema}\label{tfsij1}
Let $(Q,\sigma)$ be an equioriented symmetric quiver of type
$A_{2n}$ and let $\alpha$ be the dimension vector of a symplectic representation of $(Q,\sigma)$.\\
Then $SpSI(Q,\alpha)$ is generated by the following indecomposable
semi-invariants:
\begin{itemize}
\item[(i)] $c^{V_{j,i}}$ of weight $\langle\underline{dim}V_{j,i},\cdot\rangle$
for every $1\leq j\leq i\leq n-1$ such that
$\langle\underline{dim}\,V_{j,i},\alpha\rangle=0$,
\item[(ii)] $c^{V_{i,2n-i}}$ of weight
$\langle\underline{dim}V_{i,2n-i},\cdot\rangle$ for every
$i\in\{1,\ldots,n\}$.
\end{itemize}
\end{teorema}
The result follows from the following statement
\begin{teorema}\label{tfse}
Let $(Q,\sigma)$ be an equioriented symmetric quiver of type
$A_{2n}$, where
$$
Q=A_n^{eq}:1\stackrel{a_1}{\longrightarrow}2\cdots
n\stackrel{a_n}{\longrightarrow}n+1\cdots
2n-1\stackrel{a_{2n-1}}{\longrightarrow}2n,
$$
$\sigma(i)=2n-i+1$ and $\sigma(a_i)=a_{2n-i}$ for every
$i\in\{1,\ldots,n\}$ and let $V$
be a symplectic representation, $\underline{dim}(V)=(\alpha_1,\ldots,\alpha_n)=\alpha$.\\
Then $SpSI(Q,\alpha)$ is generated by the following indecomposable
semi-invariants:
\begin{itemize}
\item[(i)] $det(V(a_i)\cdots V(a_j))$ with $j\leq i\in\{1,\ldots,n-1\}$ if $min\{\alpha_{j+1},\ldots,\alpha_{i}\}>\alpha_j=\alpha_{i+1}$;
\item[(ii)]  $det(V(a_{2n-i})\cdots V(a_i))$ with $i\in\{1,\ldots,n\}$ if $min\{\alpha_{i+1},\ldots,\alpha_{n}\}>\alpha_i.$
\end{itemize}
\end{teorema}
\textit{Proof.} First we recall that if $V$ is a symplectic
representation of dimension $\alpha=(\alpha_1,\ldots,\alpha_n)$ of
a symmetric quiver of type $A_{2n}$, then we have
$$
 SpRep(Q,\alpha)=\bigoplus_{i=1}^{n-1}V(ta_i)^*\otimes V(ha_i)\oplus S_2V_n^*.
$$
We proceed by induction on $n$.  For $n=1$
we have the symplectic representation
$$
V_1\stackrel{V(a)}{\longrightarrow} V_1^*
$$
where $V_1$ is a vector space of dimension $\alpha$ and $V(a)$ is
a linear map such that $V(a)=V(a)^t$. So
$$
SpRep(Q,\alpha)=S^2V_1^*
$$
and by theorem \ref{c}
$$
SpSI(Q,\alpha)=\bigoplus_{\lambda\in
ER\Lambda}(S_{\lambda}V_1)^{SL(V_1)},
$$
where $ER\Lambda$ is the set of the partitions with even rows. By
proposition \ref{i1} and since $\lambda\in ER\Lambda$,
$SpSI(Q,\alpha)\neq 0$ if and only if
$\lambda=(\overbrace{2k,\ldots,2k}^{\alpha})$ for some
$k\in\mathbb{Z}_{\geq 0}$ and we have that
$(S_{\lambda}V_1)^{SL(V_1)}$ is generated by a semi-invariant of
weight $2k$. Since $g^k\cdot
detV(a)=det((g^t)^kV(a)g^k)=(det\,g)^{2k}detV(a)$ for every $g\in
GL(V)$, we note that $V(a)\in S_2V_1^*\mapsto(detV(a))^{k}$ is a
semi-invariant of weight $2k$. So
$(detV(a))^{k}$ is a generator of $(S_{\lambda}V_1)^{SL(V_1)}$ and thus $SpSI(Q,\alpha)=\mathbb{K}[detV(a)]$.\\
 Now we prove the induction step. By theorem \ref{c} we obtain
$$
SpSI(Q,\alpha)=\big(\mathbb{K}[X]\big)^{SL(V)}=
$$

$$
\bigoplus_{{\lambda(a_1),\ldots,\lambda(a_{n-1})\; and  \atop
\lambda(a_n)\in ER\Lambda}}(S_{\lambda(a_1)}V_1)^{SL(V_1)}\otimes
(S_{\lambda(a_1)}V_2^*\otimes
S_{\lambda(a_2)}V_2)^{SL(V_2)}\otimes
$$
$$
\cdots\otimes (S_{\lambda(a_{n-1})}V_n^*\otimes
S_{\lambda(a_n)}V_n)^{SL(V_n)}.
$$
where $SL(V)=SL(V_1)\times\cdots\times SL(V_n)$. We suppose that
there exists $i\in \{1,\dots,n-2\}$ such that
$\alpha_1\leq\cdots\leq\alpha_{i}>\alpha_{i+i}$. By lemma
\ref{cl},
$$
SpSI(Q,\alpha)=SpSI(Q^1,\alpha^1)
$$
where $Q^1$ is the smaller quiver $1\longrightarrow 2\cdots
i-1\longrightarrow i+1\cdots 2n-i+1\longrightarrow 2n-i+3\cdots
2n-1\longrightarrow 2n$ and $\alpha^1$ is the restriction of
$\alpha$ in $Q^1$.\\
If $i$ does't exist, we have $\alpha_1\leq\cdots\leq\alpha_{n-1}$.
So, by lemma \ref{cl}, we have the generators
$det\,V(a_i)=det\,V(\sigma(a_i))$ if
$\alpha_i=\alpha_{i+1}$, $1\leq i \leq n-2$.\\
We note that, by proposition \ref{i1},
$$
\lambda(a_1)=(\overbrace{k_1,\ldots,k_1}^{\alpha_1})
$$
is a rectangle with $k_1$ columns of height $\alpha_1$, for some
$k_1\in\mathbb{Z}_{\geq 0}$. Since
$\alpha_1\leq\cdots\leq\alpha_{n-1}$, by proposition \ref{i2}, we
obtain that there exist $k_1,\ldots,k_{n-1}\in\mathbb{Z}_{\geq 0}$
such that
$$
\lambda(a_{i})=(\overbrace{k_{i}+\cdots+k_1,\ldots,k_{i}+\cdots+k_1}^{\alpha_1},\ldots,\overbrace{k_{i},\ldots,k_{i}}^{\alpha_{i}-\alpha_{i-1}}),
$$
for every $i\in\{1,\ldots,n-1\}$.
 We also know that $\lambda_n$ must have even rows. If
 $\alpha_n=\alpha_j\leq \alpha_{j+1}\leq\cdots\leq \alpha_{n-1}$
 for some $j\in\{1,\ldots,n-1\}$ then $S_{\lambda_{n-1}}V_n^*=0$
 unless $ k_{n-1}+\cdots+k_{j+1}=0$, so
 $\lambda(a_{n-1})=\cdots=\lambda(a_{j+1})=\lambda(a_j)$. By proposition
 \ref{i2}, $(S_{\lambda(a_{n-1})}V_n^*\otimes S_{\lambda(a_n)}V_n)^{SL(V_n)}=(S_{\lambda(a_{j})}V_n^*\otimes
 S_{\lambda(a_n)}V_n)^{SL(V_n)}$ contains a semi-invariant if and
 only if
 $$
 \lambda(a_n)=(\overbrace{k_{n}+k_{j-1}+\cdots+k_1,\ldots,k_{n}+k_{j-1}+\cdots+k_1}^{\alpha_1},\ldots,\overbrace{k_{n},\ldots,k_{n}}^{\alpha_{n}-\alpha_{j-1}}),
 $$
 but
 $k_{n}+k_{j-1}+\cdots+k_1,k_{n}+k_{j-1}+\cdots+k_2,\ldots,k_n$
 have to be even and then $k_n,k_{j-1},\ldots,k_1$ have to be
 even. As before, by lemma \ref{cl}, we can consider the smaller
 quiver $Q^2:1\longrightarrow 2\cdots j\longrightarrow n\longrightarrow n+1 \longrightarrow 2n-j+1\cdots 2n-1\longrightarrow 2n$  and then
 $$
SpSI(Q,\alpha)\cong SpSI(Q^2,\alpha^2)=
$$
$$
(S_{\lambda(a_1)}V_1)^{SL(V_1)}\otimes\cdots\otimes(S_{\lambda(a_{j-1})}V_j^*\otimes
S_{\lambda(a_j)}V_j)^{SL(V_j)}\otimes(S_{\lambda(a_{j})}V_n^*\otimes
S_{\lambda(a_n)}V_n)^{SL(V_n)}.
$$
Now to complete the proof it's enough to find the generators of
$SpSI(Q^2,\alpha^2)$ for $\alpha_n=\alpha_j\leq
\alpha_{j+1}\leq\cdots\leq \alpha_{n-1}$.
\begin{itemize}
\item[(a)] By proposition \ref{i2}, for every $l\in \{1,\ldots,j\},\;(S_{\lambda(a_{l-1})}V_l^*\otimes
S_{\lambda(a_l)}V_l)^{SL(V_l)}$ is generated by a semi-invariant
of weight $(0,\ldots,0,k_l,0,\dots,0)$ where $k_l=2h$ with
$h\in\mathbb{Z}_{\geq 0}$, is $l$-th component. Since $g^h\cdot
det(V(a_{2n-l})\cdots V(a_l))=det((g_{\sigma(l)}^{-1})^h
V(a_{2n-l})\cdots V(a_l)(g_{l})^h)=det((g_{l}^t)^h
V(a_{2n-l})\cdots V(a_l)(g_{l})^h)=(det\,g_{l})^{2h}
det(V(a_{2n-l})\cdots V(a_l))$ for every $g=\{g_i\}_{i\in Q_0}\in
GL(V)$, we note that $V(a_{2n-l})\cdots V(a_l)\in
SpSI(Q,\alpha)\mapsto (det(V(a_{2n-l})\cdots V(a_l)))^h$ is a
semi-invariant of weight $(0,\ldots,0,k_l,0,\dots,0)$, so it
generates\\ $(S_{\lambda(a_{l-1})}V_l^*\otimes
S_{\lambda(a_l)}V_l)^{SL(V_l)}$. Now
$\lambda(a_l)=\lambda(a_{l-1})+(k_l^{\alpha_l})$ hence, using
lemma \ref{VxW}, $det(V(a_{2n-l})\cdots V(a_l))$ is a generator of
$SpSI(Q,\alpha)$.
\item[(b)] In the summand of $SpSI(Q,\alpha)$ indexed by the
families of partitions in which
$\lambda(a_j)=(\overbrace{k_j,\ldots,k_j}^{\alpha_j=\alpha_n})$,
with $k_j\in\mathbb{Z}_{\geq 0}$, we have that
$(S_{\lambda(a_{j})}V_j)^{SL(V_j)}\otimes
 (S_{\lambda(a_j)}V_n^*)^{SL(V_n)}$ is generated by
 a semi-invariant of weight\\
 $(0,\ldots,0,k_j,0,\dots,0,-k_j)$ where $k_j$ and $-k_j$ are
 respectively the $j$-th and the $n$-th component and we note, as
 before, that $(det(V(a_{n-1})\cdots V(a_j)))^{k_j}$ is a
 semi-invariant of weight $(0,\ldots,0,k_j,0,\dots,0,-k_j)$. Since
 $\lambda(a_j)=\lambda(a_{j-1})+(k_j^{\alpha_j=\alpha_n})$,
 $det(V(a_{n-1})\cdots V(a_j))$ is a generator of $SpSI(Q,\alpha)$;
 \item[(c)] in the summand of $SpSI(Q,\alpha)$ indexed by the
families of partitions in which
$\lambda(a_n)=(\overbrace{k_n,\ldots,k_n}^{\alpha_n})$ with
$k_n\in 2 \mathbb{Z}_{\geq 0}$, we note again that
 $(S_{\lambda(a_n)}V_n)^{SL(V_n)}$ is generated by
 $(det(V(a_n)))^{k_n}$ of weight $(0,\ldots,0,k_n)$ where $n$-th component $k_n$ is
 even. Since $\lambda(a_n)=\lambda(a_{j-1})+(k_n^{\alpha_n})$, $det(V(a_n))$ is a generator of
 $SpSI(Q,\alpha)$. $\Box$
 \end{itemize}
\textit{Proof theorem \ref{tfsij1}.} First we note that
$det(V(a_i)\cdots V(a_j))=det(V_j\rightarrow
V_{i+1})=c^{V_{j,i}}(V)$ and $\alpha_j=\alpha_{i+1}$ is equivalent
to $\langle\underline{dim}\,V_{j,i},\underline{dim}\, V\rangle=0$.
We recall, in fact, that the definition of $c^{V_{j,i}}$ doesn't
depend to the choose of projective resolution of $V_{j,i}$. If we
consider the minimal projective resolution of $V_{j,i}$, we have
$$
0\longrightarrow P_{i+1}\stackrel{a_i\cdots a_j}{\longrightarrow}
P_j\longrightarrow V_{j,i}\longrightarrow 0
$$
and applying the $Hom$-functor we have
$$
Hom(a_i\cdots a_j,V):Hom(P_j,V)=V_j\stackrel{V(a_i\cdots a_j)}{
\longrightarrow} V_{i+1}=Hom(P_{i+1},V).
$$
In the same way one proves that $det(V(a_{2n-i})\cdots
V(a_i))=det(V_i\longrightarrow
V_{2n-i+1}=V_i^*)=c^{V_{i,2n-i}}(V)$, but in this case, since
$\underline{dim}V=\underline{dim}\nabla V$, we have
$\alpha_i=\alpha_{2n-i+1}$ and so
$\langle\underline{dim}V_{i,2n-i},\underline{dim}V\rangle=0$ for
every $i\in\{1,\ldots,n\}$. Moreover we note that
\begin{itemize}
\item[(i)] $c^{V_{2n-i,2n-j}}(V)=c^{V_{j,i}}(V)$, by lemma \ref{cV=cVnabla}, since
$\tau^-\nabla V_{j\,i}=V_{2n-i\,2n-j}$;\\
\item[(ii)] for every $j\in\{1,\ldots,n-1\}$ and for every $i\in\{n+1,\ldots,2n-1\}\setminus\{2n-j\}$ there exists $j<k\in\{1,\ldots,n-1\}$ such that $2n-k=i$ and so
$c^{V_{j,i}}(V)=c^{V_{j,k-1}}(V)\cdot c^{V_{k,2n-k}}(V)$.
\end{itemize}
Now, using theorem \ref{tfse}, we obtain the statement of the
theorem. $\Box$

\subsection{The orthogonal case for $A_{2n}$}
We rewrite theorem \ref{tfo} in the following way
\begin{teorema}\label{tfoij1}
Let $(Q,\sigma)$ be an equioriented symmetric quiver of type
$A_{2n}$ and let $\alpha$ be the dimension vector of an orthogonal representation of $(Q,\sigma)$.\\
Then $OSI(Q,\alpha)$ is generated by the following indecomposable
semi-invariants:
\begin{itemize}
\item[(i)] $c^{V_{j,i}}$ of weight $\langle\underline{dim}V_{j,i},\cdot\rangle$
for every $1\leq j\leq i\leq n-1$ such that
$\langle\underline{dim}\,V_{j,i},\alpha\rangle=0$,
\item[(ii)] $pf^{V_{i,2n-i}}$ of weight
$\frac{\langle\underline{dim}V_{i,2n-i},\cdot\rangle}{2}$ for
every $i\in\{1,\ldots,n\}$ such that $\alpha_i$ is even.
\end{itemize}
\end{teorema}
The result follows from the following statement
\begin{teorema}\label{tfoe}
Let $(Q,\sigma)$ be an equioriented symmetric quiver of type
$A_{2n}$, where
$$
Q=A_n^{eq}:1\stackrel{a_1}{\longrightarrow}2\cdots
n\stackrel{a_n}{\longrightarrow}n+1\cdots
2n-1\stackrel{a_{2n-1}}{\longrightarrow}2n,
$$
$\sigma(i)=2n-i+1$ and $\sigma(a_i)=a_{2n-i}$ for every
$i\in\{1,\ldots,n\}$ and let $V$
be an orthogonal representation, $\underline{dim}(V)=(\alpha_1,\ldots,\alpha_n)=\alpha$.\\
Then $OSI(Q,\alpha)$ is generated by the following indecomposable
semi-invariants:
\begin{itemize}
\item[(i)] $det(V(a_i)\cdots V(a_j))$ with $j\leq i\in\{1,\ldots,n-1\}$ if $min(\alpha_{j+1},\ldots,\alpha_{i})>\alpha_j=\alpha_{i+1}$;
\item[(ii)]  $pf(V(a_{2n-i})\cdots V(a_i))$ with $i\in\{1,\ldots,n\}$ if
$min(\alpha_{i+1},\ldots,\alpha_{n})>\alpha_i$and $\alpha_i$ is
even.
\end{itemize}
\end{teorema}
\textit{Proof.} First we recall that if $V$ is a orthogonal
representation of dimension $\alpha=(\alpha_1,\ldots,\alpha_n)$ of
a symmetric quiver of type $A_{2n}$, then
$$
 ORep(Q,\alpha)=\bigoplus_{i=1}^{n-1}V(ta_i)^*\otimes V(ha_i)\oplus \bigwedge^2V_n^*.
$$
We proceed by induction on $n$.  For $n=1$ we have
the orthogonal representation
$$
V_1\stackrel{V(a)}{\longrightarrow} V_1^*
$$
where $V_1$ is a vector space of dimension $\alpha$ and $V(a)$ is
a linear map such that $V(a)=-V(a)^t$.
$$
ORep(Q,\alpha)=\bigwedge^2V_1^*
$$
and by theorem \ref{c}
$$
OSI(Q,\alpha)=\bigoplus_{\lambda\in
EC\Lambda}(S_{\lambda}V_1)^{SL(V_1)}
$$ where with $EC\Lambda$ we
denote the set of partitions with even columns. By proposition
\ref{i1} since $\lambda\in EC\Lambda$, $OSI(Q,\alpha)\neq 0$ if
and only if $\lambda=(\overbrace{k,\ldots,k}^{\alpha})$ with
$\alpha$ even, for some $k$. Since for every $g\in GL(V_1)$,
$g^k\cdot
pfV(a)=g^k\cdot\sqrt{detV(a)}=\sqrt{det((g^t)^{\frac{k}{2}}V(a)g^{\frac{k}{2}})}=(det\,g)^k
pfV(a)$, we note that $V(a)\in\bigwedge^2V^*\mapsto (pfV(a))^k$ is
a semi-invariant of weight $k$ so $(S_{\lambda}V_1)^{SL(V_1)}$
is generated by the semi-invariant $(pf\,V(a))^k$ if $\alpha$ is even and $OSI(Q,\alpha)=\mathbb{K}[pfV(a)]$.\\
Now we prove the induction step. Let $ X=ORep(Q,\alpha)$
 and by theorem \ref{c} we obtain
$$
OSI(Q,\alpha)=\big(\mathbb{K}[X]\big)^{SL(V)}=
$$

$$
\bigoplus_{{\lambda(a_1),\ldots,\lambda(a_{n-1})\; and  \atop
\lambda(a_n)\in EC\Lambda}}(S_{\lambda(a_1)}V_1)^{SL(V_1)}\otimes
(S_{\lambda(a_1)}V_2^*\otimes
S_{\lambda(a_2)}V_2)^{SL(V_2)}\otimes
$$
$$
\cdots\otimes (S_{\lambda(a_{n-1})}V_n^*\otimes
S_{\lambda(a_n)}V_n)^{SL(V_n)},
$$
where $SL(V)=SL(V_1)\times\cdots\times SL(V_n)$.\\
The proof of this theorem is the same of the proof of the theorem
\ref{tfse} up to when we have to consider $\alpha_n$. As in the
previous proof we can suppose $\alpha_1\leq\cdots\leq
\alpha_{n-1}$, otherwise, by induction, we can reduce to a smaller
quiver.\\
By lemma \ref{cl}, we have the generators
$det\,V(a_i)=det\,V(\sigma(a_i))$ if
$\alpha_i=\alpha_{i+1}$, $1\leq i\leq n-2$.\\
By proposition \ref{i2}, we obtain that there exist
$k_1,\ldots,k_{n-1}\in\mathbb{Z}_{\geq 0}$ such that
$$
\lambda(a_{i})=(\overbrace{k_{i}+\cdots+k_1,\ldots,k_{i}+\cdots+k_1}^{\alpha_1},\ldots,\overbrace{k_{i},\ldots,k_{i}}^{\alpha_{i}-\alpha_{i-1}}),
$$
for every $i\in\{1,\ldots,n-1\}$.\\
Now we consider the hypothesis on $\lambda(a_n)$ by which it must
have even columns. If
 $\alpha_n=\alpha_j\leq \alpha_{j+1}\leq\cdots\leq \alpha_{n-1}$
 for some $j\in\{1,\ldots,n-1\}$ then $S_{\lambda(a_{n-1})}V_n^*=0$
 unless $ k_{n-1}+\cdots+k_{j+1}=0$, so
 $\lambda(a_{n-1})=\cdots=\lambda(a_{j+1})=\lambda(a_j)$. By proposition
 \ref{i2}, $(S_{\lambda(a_{n-1})}V_n^*\otimes S_{\lambda(a_n)}V_n)^{SL(V_n)}=(S_{\lambda(a_{j})}V_n^*\otimes
 S_{\lambda(a_n)}V_n)^{SL(V_n)}$ contains a semi-invariant if and
 only if
 $$
 \lambda(a_n)=(\overbrace{k_{n}+k_{j-1}+\cdots+k_1,\ldots,k_{n}+k_{j-1}+\cdots+k_1}^{\alpha_1},\ldots,\overbrace{k_{n},\ldots,k_{n}}^{\alpha_{n}-\alpha_{j-1}}),
 $$
 but
 $\alpha_1,\alpha_2-\alpha_1,\ldots,\alpha_n-\alpha_{j-1}$
 have to be even and then $\alpha_1,\ldots,\alpha_{j-1},\alpha_n$ have to be
 even. As before, by lemma \ref{cl}, we can consider the smaller
 quiver $Q^1:1\longrightarrow 2\cdots j\longrightarrow n\longrightarrow n+1 \longrightarrow 2n-j+1\cdots 2n-1\longrightarrow 2n$  and then
 $$
OSI(Q,\alpha)\cong
(S_{\lambda(a_1)}V_1)^{SL(V_1)}\otimes\cdots\otimes(S_{\lambda(a_{j-1})}V_j^*\otimes
S_{\lambda(a_j)}V_j)^{SL(V_j)}\otimes
$$
$$
(S_{\lambda(a_{j})}V_n^*\otimes S_{\lambda(a_n)}V_n)^{SL(V_n)}.
$$
Now to complete the proof it's enough to find the generator of
this algebra for $\alpha_n=\alpha_j\leq \alpha_{j+1}\leq\cdots\leq
\alpha_{n-1}$.
\begin{itemize}
\item[(a)] By proposition \ref{i2}, for every $l\in \{1,\ldots,j\}$ such that $\alpha_l$ is even, $(S_{\lambda(a_{l-1})}V_l^*\otimes
S_{\lambda(a_l)}V_l)^{SL(V_l)}$ is generated by a semi-invariant
of weight $(0,\ldots,0,k_l,0,\dots,0)$ where
$k_l\in\mathbb{Z}_{\geq 0}$, is $l$-th component. Since $g^k\cdot
pf(V(a_{2n-l})\cdots
V(a_l))=\sqrt{det((g_{\sigma(l)}^{-1})^{\frac{k}{2}}
V(a_{2n-l})\cdots V(a_l)(g_{l})^{\frac{k}{2}})}=\\
(det\,g_{l})^{k} pf(V(a_{2n-l})\cdots V(a_l))$ for every
$g=\{g_i\}_{i\in Q_0}\in GL(V)$, we note that $V(a_{2n-l})\cdots
V(a_l)\in OSI(Q,\alpha)\mapsto (pf(V(a_{2n-l})\cdots
V(a_l)))^{k_l}$ is a semi-invariant of weight
 $(0,\ldots,0,k_l,0,\dots,0)$, so it
generates\\ $(S_{\lambda(a_{l-1})}V_l^*\otimes
S_{\lambda(a_l)}V_l)^{SL(V_l)}$. Since
$\lambda(a_l)=\lambda(a_{l-1})+(k_l^{\alpha_l})$,\\
$pf(V(a_{2n-l})\cdots V(a_l))$ is a generator of $OSI(Q,\alpha)$.
\item[(b)] In the summand of $OSI(Q,\alpha)$ indexed by the
families of partitions in which
$\lambda(a_j)=(\overbrace{k_j,\ldots,k_j}^{\alpha_j=\alpha_n})$,
with $k_j\in\mathbb{Z}_{\geq 0}$, we have that
$(S_{\lambda(a_{j})}V_j)^{SL(V_j)}\otimes
 (S_{\lambda(a_j)}V_n^*)^{SL(V_n)}$ is generated by
 a semi-invariant of weight\\
 $(0,\ldots,0,k_j,0,\dots,0,-k_j)$ where $k_j$ and $-k_j$ are
 respectively the $j$-th and the $n$-th component and we note, as
 before, that $(det(V(a_{n-1})\cdots V(a_j)))^{k_j}$ is a
 semi-invariant of weight
 $(0,\ldots,0,k_j,0,\dots,0,-k_j)$. Since
 $\lambda(a_j)=\lambda(a_{j-1})+(k_j^{\alpha_j=\alpha_n})$,
 $det(V(a_{n-1})\cdots V(a_j))$ is a generator of $OSI(Q,\alpha)$;
 \item[(c)] in the summand of $OSI(Q,\alpha)$ indexed by the
families of partitions in which
$\lambda(a_n)=(\overbrace{k_n,\ldots,k_n}^{\alpha_n})$ with
$k_n\in \mathbb{Z}_{\geq 0}$, we note again that if $\alpha_n$ is
even
 $(S_{\lambda(a_n)}V_n)^{SL(V_n)}$ is generated by
 $(pf(V(a_n)))^{k_n}$ of weight $(0,\ldots,0,k_n)$. Since $\lambda(a_n)=\lambda(a_{j-1})+(k_n^{\alpha_n})$, $pf(V(a_n))$ is a generator of
 $SpSI(Q,\alpha)$. $\Box$
 \end{itemize}
\textit{Proof of theorem \ref{tfoij1}.} By lemma \ref{ss2}, we can
define $pf^V$ if $V=\tau^-\nabla V$, since we are dealing with
orthogonal case. Moreover we note that $V_{i, 2n-i}=\tau^-\nabla
V_{i, 2n-i}$. Hence using the theorem \ref{tfoe}, the  proof is
similar to the proof of theorem \ref{tfsij1}. $\Box$

\subsection{The symplectic case for $A_{2n+1}$}
We rewrite theorem \ref{tfs} in the following way
\begin{teorema}\label{tfsij2}
Let $(Q,\sigma)$ be an equioriented symmetric quiver of type
$A_{2n+1}$ and let $\alpha$ be the dimension vector of an symplectic representation of $(Q,\sigma)$.\\
Then $SpSI(Q,\alpha)$ is generated by the following indecomposable
semi-invariants:
\begin{itemize}
\item[(i)] $c^{V_{j,i}}$ of weight $\langle\underline{dim}V_{j,i},\cdot\rangle-\varepsilon_{n+1,v_{j,i}}$
for every $1\leq j\leq i\leq n$ such that
$\langle\underline{dim}\,V_{j,i},\alpha\rangle=0$, where\\
$\varepsilon_{n+1,v_{j,i}}(h)=\left\{\begin{array}{ll}
\langle\underline{dim}V_{j,i},\cdot\rangle(n+1) & \textrm{if}\quad h=n+1\\
0 & \textrm{otherwise},\end{array}\right.$
\item[(ii)] $pf^{V_{i,2n+1-i}}$ of weight
$\frac{\langle\underline{dim}V_{i,2n+1-i},\cdot\rangle}{2}$ for
every $i\in\{1,\ldots,n\}$ such that $\alpha_i$ is even.
\end{itemize}
\end{teorema}
The result follows from the following statement
\begin{teorema}\label{tfse1}
Let $(Q,\sigma)$ be an equioriented symmetric quiver of type
$A_{2n+1}$, where
$$
Q:1\stackrel{a_1}{\longrightarrow}2\cdots
n\stackrel{a_n}{\longrightarrow}n+1\stackrel{a_{n+1}}{\longrightarrow}n+2\cdots
2n\stackrel{a_{2n}}{\longrightarrow}2n+1,
$$
$\sigma(i)=2n-i+2$ and $\sigma(a_i)=a_{2n-i+1}$ for every
$i\in\{1,\ldots,n+1\}$ and let $V$
be an symplectic representation, $\underline{dim}(V)=(\alpha_1,\ldots,\alpha_{n+1})=\alpha$.\\
Then $SpSI(Q,\alpha)$ is generated by the following indecomposable
semi-invariants:
\begin{itemize}
\item[(i)] $det(V(a_i)\cdots V(a_j))$ with $j\leq i\in\{1,\ldots,n+1\}$ if $min(\alpha_{j+1},\ldots,\alpha_{i})>\alpha_j=\alpha_{i+1}$;
\item[(ii)]  $pf(V(a_{2n-i+1})\cdots V(a_i))$ with $i\in\{1,\ldots,n\}$ if
$min(\alpha_{i+1},\ldots,\alpha_{n+1})>\alpha_i$ and $\alpha_i$ is
even.
\end{itemize}
\end{teorema}
\textit{Proof.} First we recall that if $V$ is a symplectic
representation of dimension
$\alpha=(\alpha_1,\ldots,\alpha_{n+1})$ of a symmetric quiver of
type $A_{2n+1}$, in the symplectic case,
  $V_{n+1}=V_{n+1}^*$ is a symplectic space, so if $V_{n+1}\neq 0$ then $dim\,V_{n+1}$ has to be
  even. We proceed by induction on $n$.  For $n=1$ we have
the symplectic representation
$$
V_1\stackrel{V(a)}{\longrightarrow}
V_2=V_2^*\stackrel{-V(a)^t}{\longrightarrow} V_1^*.
$$
By theorem \ref{c}
$$
 SpSI(Q,\alpha)=\bigoplus_{\lambda\in
\Lambda}(S_{\lambda}V_1)^{SL(V_1)}\otimes(S_{\lambda}V_2)^{Sp(V_2)}.
$$
By proposition \ref{i1} and proposition \ref{i3},
$SpSI(Q,\alpha)\neq 0$ if and only if
$\lambda=(\overbrace{k,\ldots,k}^{\alpha_1})$, for some $k$, and
$ht(\lambda)$ has to be even. If $\alpha_1>\alpha_2$ then
$S_{\lambda}V_2=0$ unless $\lambda=0$ and in this case
$SpSI(Q,\alpha)=\mathbb{K}$. If $\alpha_1=\alpha_2$ then
$ht(\lambda)=\alpha_1=\alpha_2$. For every $(g_1,g_2)\in
GL(V_1)\times Sp(V_2)$, $(g_1,g_2)^k\cdot detV(a)=det(g_1)^k
det(g_2^{-1})^k detV(a)=det(g_1)^k
 detV(a)$, because $g_2\in Sp(V_2)$ so we note that
$detV(a)^k$ is a semi-invariant of weight $(k,0)$. Hence
$(S_{\lambda}V_1)^{SL(V_1)}\otimes(S_{\lambda}V_2)^{Sp(V_2)}$ is
generated by the semi-invariant $detV(a)^k$, so
$SpSI(Q,\alpha)=\mathbb{K}[detV(a)]$. Finally if
$\alpha_1<\alpha_2$ then $ht(\lambda)=\alpha_1$ has to be even. We
recall that in the symplectic case $-V(a)^tV(a)$ is
skew-symmetric. Since for every $(g_1,g_2)\in GL(V_1)\times
Sp(V_2)$, $(g_1,g_2)^k\cdot
pf(-V(a)^tV(a))=(g_1,g_2)^k\cdot\sqrt{det(-V(a)^tV(a))}=$\\
$\sqrt{det((g_1^t)^{\frac{k}{2}}(-V(a)^t)(g_2)^{\frac{k}{2}}(g_2^{-1})^{\frac{k}{2}}(V(a))g_1^{\frac{k}{2}})}=
(det\,g_1)^k pf(-V(a)^tV(a))$, we note that $pf(-V(a)^tV(a))^k$ is
a semi-invariant of weight $(k,0)$ so
$(S_{\lambda}V_1)^{SL(V_1)}\otimes(S_{\lambda}V_2)^{Sp(V_2)}$ is
generated by the semi-invariant $pf(-V(a)^tV(a))^k$ if $\alpha_1$
is even and thus $SpSI(Q,\alpha)=\mathbb{K}[pf(-V(a)^tV(a))].$\\
 Now we prove the
induction step. Let $X=SpRep(Q,\alpha)$
 and by theorem \ref{c} we obtain
$$
SpSI(Q,\alpha)=\big(\mathbb{K}[X]\big)^{SSp(V)}=
$$

$$
\bigoplus_{{\lambda(a_1),\ldots,\lambda(a_{n})\in
\Lambda}}(S_{\lambda(a_1)}V_1)^{SL(V_1)}\otimes
(S_{\lambda(a_1)}V_2^*\otimes
S_{\lambda(a_2)}V_2)^{SL(V_2)}\otimes
$$
$$
\cdots\otimes (S_{\lambda(a_{n-1})}V_n^*\otimes
S_{\lambda(a_n)}V_n)^{SL(V_n)}\otimes
(S_{\lambda(a_n)}V_{n+1})^{Sp(V_{n+1})} ,
$$
where $SSp(V)=SL(V_1)\times\cdots\times SL(V_n)\times Sp(V_{n+1})$.\\
The proof of this theorem is the same of the proof of the theorem
\ref{tfse} up to when we have to consider $\alpha_{n+1}$. As in
the proof of theorem \ref{tfse} we can suppose
$\alpha_1\leq\cdots\leq \alpha_{n}$, otherwise, by induction, we
can reduce to a smaller quiver.\\
By lemma \ref{cl}, we have the generators
$det\,V(a_i)=det\,V(\sigma(a_i))$ if
$\alpha_i=\alpha_{i+1}$, $1\leq i\leq n-1$.\\
By proposition \ref{i2}, we obtain that there exist
$k_1,\ldots,k_{n}\in\mathbb{Z}_{\geq 0}$ such that
$$
\lambda(a_{i})=(\overbrace{k_{i}+\cdots+k_1,\ldots,k_{i}+\cdots+k_1}^{\alpha_1},\ldots,\overbrace{k_{i},\ldots,k_{i}}^{\alpha_{i}-\alpha_{i-1}}),
$$
for every $i\in\{1,\ldots,n\}$.\\
Now, by proposition \ref{i3}, $\lambda(a_n)$ must have even
columns. If
 $\alpha_{n+1}=\alpha_j\leq \alpha_{j+1}\leq\cdots\leq \alpha_{n}$
 for some $j\in\{1,\ldots,n\}$ then $S_{\lambda(a_{n})}V_{n+1}^*=0$
 unless $ k_{n}+\cdots+k_{j+1}=0$, so
 $\lambda(a_{n})=\cdots=\lambda(a_{j+1})=\lambda(a_j)$. As before, by lemma \ref{cl}, we can consider the smaller
 quiver $Q^1:1\longrightarrow 2\cdots j\longrightarrow n+1\longrightarrow 2n-j+2\cdots 2n\longrightarrow 2n+1$  and then
 $$
SpSI(Q,\alpha)\cong
(S_{\lambda(a_1)}V_1)^{SL(V_1)}\otimes\cdots\otimes(S_{\lambda(a_{j-1})}V_j^*\otimes
S_{\lambda(a_j)}V_j)^{SL(V_j)}\otimes
$$
\begin{equation}\label{01}
(S_{\lambda(a_{j-1})}V_n^*\otimes
S_{\lambda(a_j)}V_n)^{SL(V_n)}\otimes(S_{\lambda(a_j)}V_{n+1})^{Sp(V_{n+1})},
\end{equation}
where
$$
\lambda(a_{j})=(\overbrace{k_{j}+\cdots+k_1,\ldots,k_{j}+\cdots+k_1}^{\alpha_1},\ldots,\overbrace{k_{j},\ldots,k_{j}}^{\alpha_{n+1}-\alpha_{j-1}}),
$$
and $\alpha_1,\alpha_2-\alpha_1,\dots,\alpha_{n+1}-\alpha_{j-1}$
have to be even otherwise, by proposition \ref{i3},
$(S_{\lambda(a_j)}V_{n+1})^{Sp(V_{n+1})}=0$.
 Now to complete the proof it's enough to find the generators of the algebra (\ref{01}) for $\alpha_{n+1}=\alpha_j\leq \alpha_{j+1}\leq\cdots\leq \alpha_{n}$.
\begin{itemize}
\item[(a)] By proposition \ref{i2}, for every $l\in \{1,\ldots,j\}$ such that $\alpha_l$ is even, $(S_{\lambda(a_{l-1})}V_l^*\otimes
S_{\lambda(a_l)}V_l)^{SL(V_l)}$ is generated by a semi-invariant
of weight $(0,\ldots,0,k_l,0,\dots,0)$ where
$k_l\in\mathbb{Z}_{\geq 0}$, is $l$-th component. Since $g^k\cdot
pf(V(a_{2n-l+1})\cdots
V(a_l))=\sqrt{det((g_{\sigma(l)}^{-1})^{\frac{k}{2}}
V(a_{2n-l+1})\cdots V(a_l)(g_{l})^{\frac{k}{2}})}=\\
(det\,g_{l})^{k} pf(V(a_{2n-l+1})\cdots V(a_l))$ for every
$g=\{g_i\}_{i\in Q_0}\in Sp(V)$, we note that $V(a_{2n-l+1})\cdots
V(a_l)\in SpSI(Q,\alpha)\mapsto (pf(V(a_{2n-l+1})\cdots
V(a_l)))^{k_l}$ is a semi-invariant of weight
 $(0,\ldots,0,k_l,0,\dots,0)$, so it
generates\\ $(S_{\lambda(a_{l-1})}V_l^*\otimes
S_{\lambda(a_l)}V_l)^{SL(V_l)}$. Since
$\lambda(a_l)=\lambda(a_{l-1})+(k_l)^{\alpha_l}$, then\\
$pf(V(a_{2n-l+1})\cdots V(a_l))$ is a generator of
$SpSI(Q,\alpha)$.
\item[(b)] In the summand of $SpSI(Q,\alpha)$ indexed by the
families of partitions in which
$\lambda(a_j)=(\overbrace{k_j,\ldots,k_j}^{\alpha_j=\alpha_{n+1}})$,
with $k_j\in\mathbb{Z}_{\geq 0}$, we have that
$(S_{\lambda(a_{j})}V_j)^{SL(V_j)}\otimes
 (S_{\lambda(a_j)}V_{n+1})^{Sp(V_{n+1})}$ is generated by
 a semi-invariant of weight\\
 $(0,\ldots,0,k_j,0,\dots,0,0)$ where $k_j$ is
  the $j$-th component and we note, as
 before, that $(det(V(a_{n})\cdots V(a_j)))^{k_j}$ is a
 semi-invariant of weight
 $(0,\ldots,0,k_j,0,\dots,0,0)$. Since
 $\lambda(a_j)=\lambda(a_{j-1})+(k_j)^{\alpha_j=\alpha_{n+1}}$,\\
 $det(V(a_{n})\cdots V(a_j))$ is a generator of $SpSI(Q,\alpha)$.
 $\Box$
 \end{itemize}
\textit{Proof of theorem \ref{tfsij2}.} By lemma \ref{ss2}, we can
define $pf^V$ if $V=\tau^-\nabla V$, since we are dealing with
symplectic case. Moreover we note that $V_{i, 2n+1-i}=\tau^-\nabla
V_{i, 2n+1-i}$, for every $i\in\{1,\ldots,n\}$ . Hence using the
theorem \ref{tfse1}, the proof is similar to the proof of theorem
\ref{tfsij1}. $\Box$

\subsection{The orthogonal case for $A_{2n+1}$}
We rewrite the theorem \ref{tfo} in the following way
\begin{teorema}\label{tfoij2}
Let $(Q,\sigma)$ be an equioriented symmetric
quiver of type
$A_{2n+1}$ and let $\alpha$ be the dimension vector of an orthogonal representation of $(Q,\sigma)$.\\
Then $OSI(Q,\alpha)$ is generated by the following indecomposable
semi-invariants:
\begin{itemize}
\item[(i)] $c^{V_{j,i}}$ of weight
$\langle\underline{dim}V_{j,i},\cdot\rangle-\varepsilon_{n+1,v_{j,i}}$
for every $1\leq j\leq i\leq n$ such that
$\langle\underline{dim}\,V_{j,i},\alpha\rangle=0$,
where\\
$\varepsilon_{n+1,v_{j,i}}(h)=\left\{\begin{array}{ll}
\langle\underline{dim}V_{j,i},\cdot\rangle(n+1) & \textrm{if}\quad h=n+1\\
0 & \textrm{otherwise},\end{array}\right.$
\item[(ii)] $c^{V_{i,2n+1-i}}$ of weight
$\langle\underline{dim}V_{i,2n+1-i},\cdot\rangle$ for every
$i\in\{1,\ldots,n\}$.
\end{itemize}
\end{teorema}
The result follows from the following statement
\begin{teorema}\label{tfoe1}
Let $(Q,\sigma)$ be an equioriented symmetric quiver of type
$A_{2n+1}$, where
$$
Q:1\stackrel{a_1}{\longrightarrow}2\cdots
n\stackrel{a_n}{\longrightarrow}n+1\stackrel{a_{n+1}}{\longrightarrow}n+2\cdots
2n\stackrel{a_{2n}}{\longrightarrow}2n+1,
$$
$\sigma(i)=2n-i+2$ and $\sigma(a_i)=a_{2n-i+1}$ for every
$i\in\{1,\ldots,n+1\}$ and let $V$
be an orthogonal representation, $\underline{dim}(V)=(\alpha_1,\ldots,\alpha_{n+1})=\alpha$.\\
Then $OSI(Q,\alpha)$ is generated by the following indecomposable
semi-invariants:
\begin{itemize}
\item[(i)] $det(V(a_i)\cdots V(a_j))$ with $j\leq i\in\{1,\ldots,n+1\}$ if $min(\alpha_{j+1},\ldots,\alpha_{i})>\alpha_j=\alpha_{i+1}$;
\item[(ii)]  $det(V(a_{2n-i+1})\cdots V(a_i))$ with $i\in\{1,\ldots,n\}$ if
$min(\alpha_{i+1},\ldots,\alpha_{n+1})>\alpha_i$.
\end{itemize}
\end{teorema}
\textit{Proof.} First we recall that if $V$ is a orthogonal
representation of dimension
$\alpha=(\alpha_1,\ldots,\alpha_{n+1})$ of a symmetric quiver of
type $A_{2n+1}$, in the orthogonal case,
  $V_{n+1}=V_{n+1}^*$ is a orthogonal space. We proceed by induction on $n$.  For $n=1$ we have
the orthogonal representation
$$
V_1\stackrel{V(a)}{\longrightarrow}
V_2=V_2^*\stackrel{-V(a)^t}{\longrightarrow} V_1^*
$$
where $V_1$ is a vector space of dimension $\alpha_1$, $V_2$ is a
orthogonal space of dimension $\alpha_2$ and $V(a)$ is a linear
map. By theorem \ref{c}
$$
 OSI(Q,\alpha)=\bigoplus_{\lambda\in
\Lambda}(S_{\lambda}V_1)^{SL(V_1)}\otimes(S_{\lambda}V_2)^{SO(V_2)}.
$$
By proposition \ref{i1} and proposition \ref{i3},
$OSI(Q,\alpha)\neq 0$ if and only if
$\lambda=(\overbrace{k,\ldots,k}^{\alpha_1})$, for some $k\in
2\mathbb{Z}$. If $\alpha_1>\alpha_2$ then $S_{\lambda}V_2=0$
unless $\lambda=0$ and in this case $OSI(Q,\alpha)=\mathbb{K}$. If
$\alpha_1=\alpha_2$ then $ht(\lambda)=\alpha_1=\alpha_2$. For
every $(g_1,g_2)\in GL(V_1)\times SO(V_2)$, $(g_1,g_2)^k\cdot
detV(a)=det(g_1)^k det(g_2^{-1})^k detV(a)=det(g_1)^k
 detV(a)$, because $g_2\in SO(V_2)$ so we note that
$detV(a)^k$ is a semi-invariant of weight $(k,0)$. Hence
$(S_{\lambda}V_1)^{SL(V_1)}\otimes(S_{\lambda}V_2)^{SO(V_2)}$ is
generated by the semi-invariant $detV(a)^k$, so
$OSI(Q,\alpha)=\mathbb{K}[detV(a)]$. Finally if
$\alpha_1<\alpha_2$ for every $(g_1,g_2)\in GL(V_1)\times
SO(V_2)$, $(g_1,g_2)^k\cdot
det(-V(a)^tV(a))=(g_1,g_2)^k\cdot det(-V(a)^tV(a))=$\\
$det((g_1^t)^{k}(-V(a)^t)(g_2)^{k}(g_2^{-1})^{k}(V(a))g_1^k)=
(det\,g_1)^k det(-V(a)^tV(a))$, we note that $det(-V(a)^tV(a))^k$
is a semi-invariant of weight $(k,0)$ so
$(S_{\lambda}V_1)^{SL(V_1)}\otimes(S_{\lambda}V_2)^{SO(V_2)}$ is
generated by the semi-invariant $det(-V(a)^tV(a))^k$ and thus $OSI(Q,\alpha)=\mathbb{K}[det(-V(a)^tV(a))].$\\
 Now we prove the
induction step. Let $X=ORep(Q,\alpha$
 and by theorem \ref{c} we obtain
$$
OSI(Q,\alpha)=\big(\mathbb{K}[X]\big)^{SO(V)}=
$$

$$
\bigoplus_{{\lambda(a_1),\ldots,\lambda(a_{n})\in
\Lambda}}(S_{\lambda(a_1)}V_1)^{SL(V_1)}\otimes
(S_{\lambda(a_1)}V_2^*\otimes
S_{\lambda(a_2)}V_2)^{SL(V_2)}\otimes
$$
$$
\cdots\otimes (S_{\lambda(a_{n-1})}V_n^*\otimes
S_{\lambda(a_n)}V_n)^{SL(V_n)}\otimes
(S_{\lambda(a_n)}V_{n+1})^{SO(V_{n+1})} ,
$$
where $SO(V)=SL(V_1)\times\cdots\times SL(V_n)\times SO(V_{n+1})$.\\
The proof of this theorem is the same of the proof of the theorem
\ref{tfse} up to when we have to consider $\alpha_{n+1}$. As in
the proof of theorem \ref{tfse} we can suppose
$\alpha_1\leq\cdots\leq \alpha_{n}$, otherwise, by induction, we
can reduce to a smaller quiver.\\
By lemma \ref{cl}, we have the generators
$det\,V(a_i)=det\,V(\sigma(a_i))$ if
$\alpha_i=\alpha_{i+1}$, $1\leq i\leq n-2$.\\
By proposition \ref{i2}, we obtain that there exist
$k_1,\ldots,k_{n}\in\mathbb{Z}_{\geq 0}$ such that
$$
\lambda(a_{i})=(\overbrace{k_{i}+\cdots+k_1,\ldots,k_{i}+\cdots+k_1}^{\alpha_1},\ldots,\overbrace{k_{i},\ldots,k_{i}}^{\alpha_{i}-\alpha_{i-1}}),
$$
for every $i\in\{1,\ldots,n\}$.\\
Now, by proposition \ref{i3}, $\lambda(a_n)$ must have even rows.
If
 $\alpha_{n+1}=\alpha_j\leq \alpha_{j+1}\leq\cdots\leq \alpha_{n}$
 for some $j\in\{1,\ldots,n\}$ then $S_{\lambda(a_{n})}V_{n+1}^*=0$
 unless $ k_{n}+\cdots+k_{j+1}=0$, so
 $\lambda(a_{n})=\cdots=\lambda(a_{j+1})=\lambda(a_j)$. As before, by lemma \ref{cl}, we can consider the smaller
 quiver $Q^1:1\longrightarrow 2\cdots j\longrightarrow n+1\longrightarrow 2n-j+2\cdots 2n\longrightarrow 2n+1$  and then
 $$
OSI(Q,\alpha)\cong
(S_{\lambda(a_1)}V_1)^{SL(V_1)}\otimes\cdots\otimes(S_{\lambda(a_{j-1})}V_j^*\otimes
S_{\lambda(a_j)}V_j)^{SL(V_j)}\otimes
$$
\begin{equation}\label{02}
(S_{\lambda(a_{j-1})}V_n^*\otimes
S_{\lambda(a_j)}V_n)^{SL(V_n)}\otimes(S_{\lambda(a_j)}V_{n+1})^{SO(V_{n+1})},
\end{equation}
where
$$
\lambda(a_{j})=(\overbrace{k_{j}+\cdots+k_1,\ldots,k_{j}+\cdots+k_1}^{\alpha_1},\ldots,\overbrace{k_{j},\ldots,k_{j}}^{\alpha_{n+1}-\alpha_{j-1}}),
$$
and $k_{j}+\cdots+k_1,\ldots,k_j$ have to be even otherwise, by
proposition \ref{i3}, $(S_{\lambda(a_j)}V_{n+1})^{SO(V_{n+1})}=0$.
Hence $k_l$ has to be even for every $l\in\{1,\ldots,j\}$.
 Now to complete the proof it's enough to find the generators of the algebra (\ref{02}) for $\alpha_{n+1}=\alpha_j\leq \alpha_{j+1}\leq\cdots\leq \alpha_{n}$.
\begin{itemize}
\item[(a)] By proposition \ref{i2}, for every $l\in \{1,\ldots,j\}$, $(S_{\lambda(a_{l-1})}V_l^*\otimes
S_{\lambda(a_l)}V_l)^{SL(V_l)}$ is generated by a semi-invariant
of weight $(0,\ldots,0,k_l,0,\dots,0)$ where $k_l\in 2
\mathbb{Z}_{\geq 0}$, is $l$-th component. Since $g^k\cdot
det(V(a_{2n-l+1})\cdots V(a_l))=det((g_{\sigma(l)}^{-1})^{k}
V(a_{2n-l+1})\cdots V(a_l)(g_{l})^k)=\\
(det\,g_{l})^{k} det(V(a_{2n-l+1})\cdots V(a_l))$ for every
$g=\{g_i\}_{i\in Q_0}\in SO(V)$, we note that
$$V(a_{2n-l+1})\cdots V(a_l)\in OSI(Q,\alpha)\mapsto
(det(V(a_{2n-l+1})\cdots V(a_l)))^{\frac{k_l}{2}}$$ is a
semi-invariant of weight
 $(0,\ldots,0,k_l,0,\dots,0)$, so it
generates\\ $(S_{\lambda(a_{l-1})}V_l^*\otimes
S_{\lambda(a_l)}V_l)^{SL(V_l)}$. Since
$\lambda(a_l)=\lambda(a_{l-1})+(k_l)^{\alpha_l}$, then\\
$det(V(a_{2n-l+1})\cdots V(a_l))$ is a generator of
$OSI(Q,\alpha)$.
\item[(b)] In the summand of $OSI(Q,\alpha)$ indexed by the
families of partitions in which
$\lambda(a_j)=(\overbrace{k_j,\ldots,k_j}^{\alpha_j=\alpha_{n+1}})$,
with $k_j\in 2\mathbb{Z}_{\geq 0}$, we have that
$(S_{\lambda(a_{j})}V_j)^{SL(V_j)}\otimes
 (S_{\lambda(a_j)}V_{n+1})^{SO(V_{n+1})}$ is generated by
 a semi-invariant of weight\\
 $(0,\ldots,0,k_j,0,\dots,0,0)$ where $k_j$ is
  the $j$-th component and we note, as
 before, that $(det(V(a_{n})\cdots V(a_j)))^{k_j}$ is a
 semi-invariant of weight
 $(0,\ldots,0,k_j,0,\dots,0,0)$. Since
 $\lambda(a_j)=\lambda(a_{j-1})+(k_j)^{\alpha_j=\alpha_{n+1}}$,\\
 $det(V(a_{n})\cdots V(a_j))$ is a generator of $OSI(Q,\alpha)$.
 $\Box$
 \end{itemize}
\textit{Proof of theorem \ref{tfoij2}.} Using the theorem
\ref{tfoe1}, the proof is similar to the proof of theorem
\ref{tfsij1}. $\Box$

\chapter{Semi-invariants of symmetric quivers of tame type}

In this chapter we prove conjectures \ref{mt1} and \ref{mt2} for
the symmetric quivers of tame type. We recall that the underlying
quiver of a symmetric quiver of tame type is either
$\widetilde{A}$ or $\widetilde{D}$ as in proposition \ref{ctt}. As
done for the finite case we again reduce the proof to particular
orientations (orientations in proposition \ref{oA} for
$\widetilde{A}$ and orientation of $\widetilde{D}^{eq}$ for
$\widetilde{D}$). In section 3.1, we prove the conjectures for
dimension vector $ph$ (for definition, see proposition \ref{h}).
In section 3.2, we treat the other regular dimension vectors.

\section{Semi-invariants of symmetric
quivers of tame type for dimension vector $ph$}

In this section we deal with dimension vector $ph$. By lemma
\ref{kac} and proposition \ref{SIcx1} and \ref{SIcx2}, it's enough
to consider particular orientations of symmetric quivers of type
$\widetilde{A}$ in proposition \ref{oA} and orientation of
symmetric quiver $\widetilde{D}^{eq}$. First we prove case by case
some theorems by which conjectures \ref{mt1} and \ref{mt2} follow.
Finally, in section \ref{fineph}, we conclude proofs of
conjectures \ref{mt1} and \ref{mt2}. We note that $h$ is preserved
under reflection functor.

\subsection{$\widetilde{A}^{2,0,1}_{k,l}$ for dimension vector $ph$}
\begin{teorema}
Let $(Q,\sigma)$ be a symmetric quiver of type $(2,0,k,l)$ of
orientation
$$
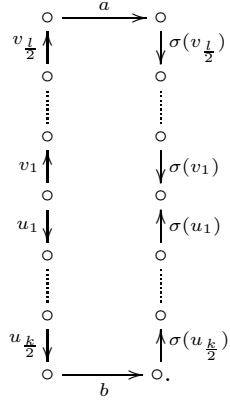
\xymatrix@-1pc{
\circ\ar[rr]^a&&\circ\ar[d]^{\sigma(v_{\frac{l}{2}})}\\
\circ \ar[u]^{v_{\frac{l}{2}}}\ar@{.}[d]&&\circ\\
\circ&&\circ\ar[d]^{\sigma(v_{1})}\ar@{.}[u]\\
\circ\ar[d]_{u_1}\ar[u]^{v_1}&&\circ\\
\circ&&\circ\ar[u]_{\sigma(u_{1})}\ar@{.}[d]\\
\circ\ar[d]_{u_{\frac{k}{2}}}\ar@{.}[u]&&\circ\\
\circ\ar[rr]_b&&\circ\ar[u]_{\sigma(u_{\frac{k}{2}})}.}
$$
Then\\
Sp) $SpSI(Q,ph)$ is generated by the following indecomposable
semi-invariants:
\begin{itemize}
\item[a)] $det\,V(u_j)$ with $j\in\{1,\ldots,\frac{k}{2}\}$;
\item[b)] $det\,V(v_j)$ with $j\in\{1,\ldots,\frac{l}{2}\}$;
\item[c)] $det\,V(a)$ and $det\,V(b)$;
\item[d)] the coefficients $c_i$ of $\varphi^{p-i}\psi^i$, $0\leq i\leq p$, in $det(\psi V(\bar{a})+\varphi
V(\bar{b}))$, where $\bar{a}=\sigma(v_{1})\cdots
\sigma(v_{\frac{l}{2}})a v_{\frac{l}{2}}\cdots v_1$ and
$\bar{b}=\sigma(u_{1})\cdots \sigma(u_{\frac{k}{2}})b
u_{\frac{k}{2}}\cdots u_1$.
\end{itemize}
O) $OSI(Q,ph)$ is generated by the following indecomposable
semi-invariants:\\
if $p$ is even,
\begin{itemize}
\item[a)] $det\,V(u_j)$ with $j\in\{1,\ldots,\frac{k}{2}\}$;
\item[b)] $det\,V(v_j)$ with $j\in\{1,\ldots,\frac{l}{2}\}$;
\item[c)] $pf\,V(a)$ and $pf\,V(b)$;
\item[d)] the coefficients $c_i$ of $\varphi^{p-2i}\psi^{2i}$, $0\leq i\leq \frac{p}{2}$, in $pf(\psi V(\bar{a})+\varphi
V(\bar{b}))$, where $\bar{a}=\sigma(v_{1})\cdots
\sigma(v_{\frac{l}{2}})a v_{\frac{l}{2}}\cdots v_1$ and
$\bar{b}=\sigma(u_{1})\cdots \sigma(u_{\frac{k}{2}})b
u_{\frac{k}{2}}\cdots u_1$;
\end{itemize}
if $p$ is odd,
\begin{itemize}
\item[a)] $det\,V(u_j)$ with $j\in\{1,\ldots,\frac{k}{2}\}$;
\item[b)] $det\,V(v_j)$ with $j\in\{1,\ldots,\frac{l}{2}\}$.
\end{itemize}
\end{teorema}
\textit{Proof.} We proceed by induction on
$\frac{k}{2}+\frac{l}{2}$. The smallest case is
$\widetilde{A}^{2,0,1}_{0,0}$
$$
\xymatrix@-1pc{ 1\ar@{=>}[rr]^a_b&&\sigma(1).}
$$
The induction step follows by lemma \ref{cls}, so it's enough to prove the theorem for $\widetilde{A}^{2,0,1}_{0,0}$.\\
Let $V$ be a representation of $\widetilde{A}^{2,0,1}_{0,0}$ of
dimension $ph$ for some $p\in\mathbb{Z}_{\geq 0}$, in this case $h=1$.\\
\textbf{Sp)} The ring of symplectic semi-invariants is
$$
SpSI(\widetilde{A}^{2,0,1}_{0,0},ph)=\bigoplus_{\lambda(a),\lambda(b)\in
ER\Lambda}(S_{\lambda(a)}V\otimes S_{\lambda(b)}V)^{SL\,V}.
$$
By proposition \ref{i2} we have
\begin{equation}\label{1}
\lambda(a)_j+\lambda(b)_{p+j-1}=t
\end{equation}
for some $t\in\mathbb{Z}_{\geq 0}$ and for every $0\leq j\leq
p$.\\
We consider the summand in which $t=2$ because the other ones are
generated by products of powers of the generators of this summand.
The solutions of (\ref{1}) are $\lambda(a)=(2^i)$ and
$\lambda(b)=(2^{p-i})$ for every $0\leq i \leq p$. So the
considered summand $\bigoplus_{i=0}^p(S_{(2^i)}V\otimes
S_{(2^{p-i})}V)^{SL\,V}$ is generated by semi-invariants of weight
2, i.e. the coefficients $c_i$ of $\varphi^{p-i}\psi^{i}$ in
$det(\psi V(a)+\varphi V(b))$ (see [R2]). In particular we have
$c_0=det\,V(b)$ and $c_p=det\,V(a)$.\\
\textbf{O)} The ring of orthogonal semi-invariants is
$$
OSI(\widetilde{A}^{2,0,1}_{0,0},ph)=\bigoplus_{\lambda(a),\lambda(b)\in
EC\Lambda}(S_{\lambda(a)}V\otimes S_{\lambda(b)}V)^{SL\,V}.
$$
By proposition \ref{i2} we have
\begin{equation}\label{2}
\lambda(a)_j+\lambda(b)_{p+j-1}=t
\end{equation}
for some $t\in\mathbb{Z}_{\geq 0}$ and for every $0\leq j\leq
p$.\\
We consider the summand in which $t=1$ because the other ones are
generated by products of powers of the generators of this summand.
Let $p$ be odd. $(\lambda(a)')_1$ and $(\lambda(b)')_p$ have to be
even but $(\lambda(a)')_1+(\lambda(b)')_p=p$ is odd, this is an
absurd, and so $OSI(\widetilde{A}^{2,0,1}_{0,0},ph)=\mathbb{K}$.\\
Let $p$ be even. The solutions of (\ref{2}) are
$\lambda(a)=(1^{2i})$ and $\lambda(b)=(1^{p-2i})$ for every $0\leq
i \leq \frac{p}{2}$. So the considered summand
$\bigoplus_{i=0}^p(S_{(1^{2i})}V\otimes S_{(1^{p-2i})}V)^{SL\,V}$
is generated by semi-invariants of weight 1, i.e. the coefficients
$c_i$ of $\varphi^{p-2i}\psi^{2i}$ in $pf(\psi V(a)+\varphi
V(b))$. In particular we have $c_0=pf\,V(b)$ and
$c_{\frac{p}{2}}=pf\,V(a)$. $\Box$
\subsection{$\widetilde{A}^{2,0,2}_{k,l}$ for dimension vector $ph$}
\begin{teorema}
Let $(Q,\sigma)$ be a symmetric quiver of type $(2,0,k,l)$ with
orientation
$$
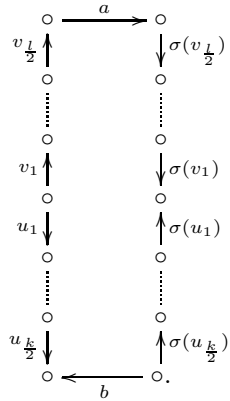
\xymatrix@-1pc{
\circ\ar[rr]^a&&\circ\ar[d]^{\sigma(v_{\frac{l}{2}})}\\
\circ\ar[u]^{v_{\frac{l}{2}}}\ar@{.}[d]&&\circ\\
\circ&&\circ\ar[d]^{\sigma(v_{1})}\ar@{.}[u]\\
\circ\ar[d]_{u_1}\ar[u]^{v_1}&&\circ\\
\circ&&\circ\ar[u]_{\sigma(u_{1})}\ar@{.}[d]\\
\circ\ar[d]_{u_{\frac{k}{2}}}\ar@{.}[u]&&\circ\\
 \circ&&\circ\ar[ll]^b\ar[u]_{\sigma(u_{\frac{k}{2}})}.}
$$
Then\\
Sp) $SpSI(Q,ph)$ is generated by the following indecomposable
semi-invariants:
\begin{itemize}
\item[a)] $det\,V(u_j)$ with $j\in\{1,\ldots,\frac{k}{2}\}$;
\item[b)] $det\,V(v_j)$ with $j\in\{1,\ldots,\frac{l}{2}\}$;
\item[c)] $det\,V(a)$ and $det\,V(b)$;
\item[d)] the $c_i$ coefficients of $\varphi^{i}\psi^i$, $0\leq i\leq p$, in
$det\left(\begin{array}{cc}\varphi V(\bar{a})&V(c)\\
V(\sigma(c))&\psi V(b))\end{array}\right)$, where
$\bar{a}=\sigma(v_{1})\cdots \sigma(v_{\frac{l}{2}})a
v_{\frac{l}{2}}\cdots v_1$ and $c=u_{\frac{k}{2}}\cdots u_1$.
\end{itemize}
O) $OSI(Q,ph)$ is generated by the following indecomposable
semi-invariants:\\
if $p$ is even,
\begin{itemize}
\item[a)] $det\,V(u_j)$ with $j\in\{1,\ldots,\frac{k}{2}\}$;
\item[b)] $det\,V(v_j)$ with $j\in\{1,\ldots,\frac{l}{2}\}$;
\item[c)] $pf\,V(a)$ and $pf\,V(b)$;
\item[d)] the coefficients $c_i$ of $\varphi^{i}\psi^{i}$, $0\leq i\leq \frac{p-1}{2}$, in
$pf\left(\begin{array}{cc}\varphi V(\bar{a})&V(c)\\
V(\sigma(c))&\psi V(b))\end{array}\right)$, where
$\bar{a}=\sigma(v_{1})\cdots \sigma(v_{\frac{l}{2}})a
v_{\frac{l}{2}}\cdots v_1$ and $c=u_{\frac{k}{2}}\cdots u_1$.
\end{itemize}
if $p$ is odd,
\begin{itemize}
\item[a)] $det\,V(u_j)$ with $j\in\{1,\ldots,\frac{k}{2}\}$;
\item[b)] $det\,V(v_j)$ with $j\in\{1,\ldots,\frac{l}{2}\}$;
\item[c)] the coefficients $c_i$ of $\varphi^{i}\psi^{i}$, $0\leq i\leq \frac{p-1}{2}$, in
$pf\left(\begin{array}{cc}\varphi V(\bar{a})&V(c)\\
V(\sigma(c))&\psi V(b))\end{array}\right)$, where
$\bar{a}=\sigma(v_{1})\cdots \sigma(v_{\frac{l}{2}})a
v_{\frac{l}{2}}\cdots v_1$ and $c=u_{\frac{k}{2}}\cdots u_1$.
\end{itemize}
\end{teorema}
\textit{Proof.} We proceed by induction on
$\frac{k}{2}+\frac{l}{2}$. The smallest case is
$\widetilde{A}^{2,0,2}_{2,0}$
$$
\xymatrix@-1pc{ 1\ar[dd]_c\ar[rr]^a&&\sigma(1)\\
&&\\
 2 &&\sigma(2)\ar[ll]^b\ar[uu]_{\sigma(c)}}
$$
and so it's enough to study the semi-invariants of
$\widetilde{A}^{2,0,2}_{2,0}$.\\
The induction step follows by lemma \ref{cls} and by lemma \ref{cl}, so it's enough to prove the theorem for $\widetilde{A}^{2,0,2}_{2,0}$.\\
\textbf{Sp)} The ring of symplectic semi-invariants is
$$
SpSI(\widetilde{A}^{2,0,2}_{2,0},ph)=\bigoplus_{{\lambda(a),\lambda(b)\in
ER\Lambda\atop \lambda(c)\in\Lambda}}(S_{\lambda(a)}V_1\otimes
S_{\lambda(c)}V_1)^{SL\,V_1}\otimes (S_{\lambda(b)}V_2^*\otimes
S_{\lambda(c)}V_2^*)^{SL\,V_2}.
$$
By proposition \ref{i2} we have
\begin{equation}\label{3}
\left\{\begin{array}{c}
\lambda(a)_j+\lambda(c)_{p+j-1}=k_1\\
\lambda(b)_j+\lambda(c)_{p+j-1}=k_2
\end{array}\right.
\end{equation}
for some $k_1,k_2\in\mathbb{Z}_{\geq 0}$ and for every $0\leq
j\leq
p$.\\
We consider the summands in which $k_1=0,1,2$ and $k_2=0,1,2$
because the other ones are generated by products of powers of the
generators of this summands. If $k_1=2$ and $k_2=0$ we have
$\lambda(b)=0=\lambda(c)$ and so the summand is
$(S_{(2^p)}V_1)^{SL\,V_1}$ which is generated by a semi-invariant
of weight $(2,0)$, i.e. $det\,V(a)$. If $k_1=0$ and $k_2=2$ as
before we obtain the generator of ring of semi-invariant
$det\,V(b)$ of weight $(0,-2)$. The summand in which $k_1=1$ and
$k_2=0$ (respectively $k_1=0$ and $k_2=1$) doesn't exist because
otherwise we have $\lambda(a)$ (respectively $\lambda(b)$) with
odd columns. If $k_1=1=k_2$ we have $\lambda(a)=0=\lambda(b)$ and
$\lambda(c)=(1^p)$ and so the summand is
$(S_{(1^p)}V_1)^{SL\,V_1}\otimes (S_{(1^p)}V_2^*)^{SL\,V_2}$ which
is generated by a semi-invariant of weight $(1,-1)$ which is
$det\,V(c)=det\,V(\sigma(c))$. If $k_1=2=k_2$, the solutions of
(\ref{3}) are $\lambda(a)=(2^i)=\lambda(b)$ and
$\lambda(c)=(2^{p-i)}$. The corresponding summand is
$\bigoplus_{i=0}^p(S_{(2^i)}V_1\otimes
S_{(2^{p-i})}V_1)^{SL\,V_1}\otimes(S_{(2^i)}V_2^*\otimes
S_{(2^{p-i})}V_2)^{SL\,V_2^*}$ and it is spanned by the
coefficients of $\varphi^{i}\psi^{i}$ in
$$
det\left(\begin{array}{cc} \varphi V(a) & V(c)\\
V(\sigma(c)) & \psi V(b)\end{array}\right),
$$
semi-invariants of weight $(2,-2)$. In particular for $i=0$ we
have $(det\,V(c))^2$ and for $i=p$ we have $det\,V(a)\cdot
det\,V(b)$.\\
\textbf{O)} The ring of orthogonal semi-invariants is
$$
OSI(\widetilde{A}^{2,0,2}_{2,0},ph)=\bigoplus_{{\lambda(a),\lambda(b)\in
EC\Lambda\atop \lambda(c)\in\Lambda}}(S_{\lambda(a)}V_1\otimes
S_{\lambda(c)}V_1)^{SL\,V_1}\otimes (S_{\lambda(b)}V_2^*\otimes
S_{\lambda(c)}V_2^*)^{SL\,V_2}.
$$
By proposition \ref{i2} we have
\begin{equation}\label{4}
\left\{\begin{array}{c}
\lambda(a)_j+\lambda(c)_{p+j-1}=k_1\\
\lambda(b)_j+\lambda(c)_{p+j-1}=k_2
\end{array}\right.
\end{equation}
for some $k_1,k_2\in\mathbb{Z}_{\geq 0}$ and for every $0\leq
j\leq
p$.\\
We consider the summands in which $k_1=0,1$ and $k_2=0,1$ because
the other ones are generated by the monomials of these. Let $p$ be
even. If $k_1=1$ and $k_2=0$ we have $\lambda(b)=0=\lambda(c)$ and
so the summand is $(S_{(1^p)}V_1)^{SL\,V_1}$ which is generated by
a semi-invariant of weight $(1,0)$, i.e. $pf\,V(a)$. If $k_1=0$
and $k_2=1$ as before we obtain the generator of ring of
semi-invariant $pf\,V(b)$ of weight $(0,-1)$. If $k_1=1=k_2$, the
solutions of (\ref{4}) are $\lambda(a)=(1^{2i})=\lambda(b)$ and
$\lambda(c)=(1^{p-2i})$ with $0\leq i\leq\frac{p}{2}$. So the
summand is $\bigoplus_{i=0}^{\frac{p}{2}}(S_{(1^{2i})}V_1\otimes
S_{(1^{p-2i})}V_1)^{SL\,V_1}\otimes(S_{(1^{2i})}V_2^*\otimes
S_{(1^{p-2i})}V_2^*)^{SL\,V_2}$ which is generated by the
coefficients of $\varphi^{i}\psi^{i}$ in
$$
pf\left(\begin{array}{cc} \varphi V(a) & V(c)\\
V(\sigma(c)) & \psi V(b)\end{array}\right),
$$
semi-invariants of weight $(1,-1)$. In particular for $i=0$ we
have $det\,V(c)=det\,V(\sigma(c))$ and for $i=p$ we have
$pf\,V(a)\cdot pf\,V(b)$. Let $p$ be odd. In this case the summand
$(S_{(1^p)}V_1)^{SL\,V_1}$ (respectively
$(S_{(1^p)}V_2)^{SL\,V_2}$) doesn't exist since $\lambda(a)$
(respectively $\lambda(b)$) must have even columns. If
$k_1=1=k_2$, the solutions of \ref{4} are
$\lambda(a)=(1^{2i})=\lambda(b)$ and $\lambda(c)=(1^{p-2i})$ with
$0\leq i\leq\frac{p-1}{2}$. So the summand is
$\bigoplus_{i=0}^{\frac{p-1}{2}}(S_{(1^{2i})}V_1\otimes
S_{(1^{p-2i})}V_1)^{SL\,V_1}\otimes(S_{(1^{2i})}V_2^*\otimes
S_{(1^{p-2i})}V_2^*)^{SL\,V_2}$ which is generated by the
coefficients of $\varphi^{i}\psi^{i}$ in
$$
pf\left(\begin{array}{cc} \varphi V(a) & V(c)\\
V(\sigma(c)) & \psi V(b)\end{array}\right),
$$
semi-invariants of weight $(1,-1)$. In particular for $i=0$ we get
$det\,V(c)=det\,V(\sigma(c))$. $\Box$
\subsection{$\widetilde{A}^{0,2}_{k,l}$ for dimension vector
$ph$}
\begin{teorema}
Let $(Q,\sigma)$ be a symmetric quiver of type $(0,2,k,l)$ with
orientation
$$
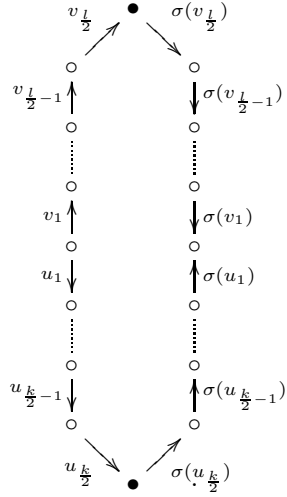
\xymatrix@-1pc{ &\bullet\ar[dr]^{\sigma(v_{\frac{l}{2}})}&\\
\circ\ar[ur]^{v_{\frac{l}{2}}}&&\circ\ar[d]^{\sigma(v_{\frac{l}{2}-1})}\\
\circ\ar[u]^{v_{\frac{l}{2}-1}}\ar@{.}[d]&&\circ\\
\circ&&\circ\ar[d]^{\sigma(v_{1})}\ar@{.}[u]\\
\circ\ar[d]_{u_1}\ar[u]^{v_1}&&\circ\\
\circ&&\circ\ar[u]_{\sigma(u_{1})}\ar@{.}[d]\\
\circ\ar[d]_{u_{\frac{k}{2}-1}}\ar@{.}[u]&&\circ\\
 \circ\ar[dr]_{u_{\frac{k}{2}}}&&\circ\ar[u]_{\sigma(u_{\frac{k}{2}-1})}\\
 &\bullet\ar[ur]_{\sigma(u_{\frac{k}{2}})}&.}
$$
Then\\
O) $OSI(Q,ph)$ is generated by the following indecomposable
semi-invariants:
\begin{itemize}
\item[a)] $det\,V(u_j)$ with $j\in\{1,\ldots,\frac{k}{2}\}$;
\item[b)] $det\,V(v_j)$ with $j\in\{1,\ldots,\frac{l}{2}\}$;
\item[c)] the coefficients $c_i$ of $\varphi^{p-i}\psi^i$, $0\leq i\leq p$, in $det(\psi V(\sigma(\bar{a})\bar{a})+\varphi
V(\sigma(\bar{b})\bar{b}))$, where $\bar{a}=\sigma(v_{1})\cdots
\sigma(v_{\frac{l}{2}}) v_{\frac{l}{2}}\cdots v_1$ and
$\bar{b}=\sigma(u_{1})\cdots \sigma(u_{\frac{k}{2}})
u_{\frac{k}{2}}\cdots u_1$.
\end{itemize}
Sp) $SpSI(Q,ph)$ is generated by the following indecomposable
semi-invariants:\\
if $p$ is even,
\begin{itemize}
\item[a)] $det\,V(u_j)$ with $j\in\{1,\ldots,\frac{k}{2}\}$;
\item[b)] $det\,V(v_j)$ with $j\in\{1,\ldots,\frac{l}{2}\}$;
\item[c)] the coefficients $c_i$ of $\varphi^{\frac{p}{2}-i}\psi^{i}$, $0\leq i\leq \frac{p}{2}$, in $pf(\psi V(\sigma(\bar{a})\bar{a})+\varphi
V(\sigma(\bar{b})\bar{b}))$, where $\bar{a}=\sigma(v_{1})\cdots
\sigma(v_{\frac{l}{2}}) v_{\frac{l}{2}}\cdots v_1$ and
$\bar{b}=\sigma(u_{1})\cdots \sigma(u_{\frac{k}{2}})
u_{\frac{k}{2}}\cdots u_1$;
\end{itemize}
if $p$ is odd, $SpSI(Q,ph)=\mathbb{K}$.
\end{teorema}
\textit{Proof.} We proceed by induction on
$\frac{k}{2}+\frac{l}{2}$. The smallest case is
$\widetilde{A}^{0,2}_{2,2}$
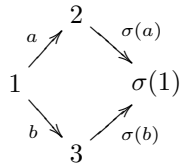
$$
\xymatrix@-1pc{&2\ar[dr]^{\sigma(a)}&\\
1\ar[dr]_b\ar[ur]^a&&\sigma(1)\\
&3\ar[ur]_{\sigma(b)}&}
$$
and so it's enough to study the semi-invariants of
$\widetilde{A}^{0,2}_{2,2}$.\\
The induction step follows by lemma \ref{cl}, so it's enough to prove the theorem for $\widetilde{A}^{0,2}_{2,2}$.\\
\textbf{O)} The ring of orthogonal semi-invariants is
$$
\bigoplus_{\lambda(a),\lambda(b)\in\Lambda}(S_{\lambda(a)}V_1\otimes
S_{\lambda(b)}V_1)^{SL\,V_1}\otimes
(S_{\lambda(a)}V_2)^{SO\,V_2}\otimes(S_{\lambda(b)}V_3)^{SO\,V_3}.
$$
By proposition \ref{i2} we have
\begin{equation}\label{5}
\lambda(a)_j+\lambda(b)_{p-j+1}=k_1
\end{equation}
for every $0\leq j\leq p$ and for some $k_1\in\mathbb{Z}_{\geq
0}$. By proposition \ref{i3} we have $\lambda(a)=2\mu+(l^p)$ and
$\lambda(b)=2\nu+(m^p)$ for some $\mu,\nu\in\Lambda$ and for some
$l,m\in\mathbb{Z}_{\geq 0}$. We consider the summands in which
$k_1=1,2$ because the other ones are generated by products of
powers of the generators of this summands. If $k_1=1$ the only
solutions of (\ref{5}) are $\lambda(a)=(1^p)$, $\lambda(b)=0$ and
$\lambda(a)=0$, $\lambda(b)=(1^p)$. Respectively, the summand
$(S_{(1^p)}V_1)^{SL\,V_1}\otimes(S_{(1^p)}V_2)^{SO\,V_2}$ is
generated by a semi-invariant of weight $(1,0,0)$, i.e
$det\,V(a)=det\,V(\sigma(a))$, and the summand
$(S_{(1^p)}V_1)^{SL\,V_1}\otimes(S_{(1^p)}V_3)^{SO\,V_3}$ is
generated by a semi-invariant of weight $(1,0,0)$, i.e
$det\,V(b)=det\,V(\sigma(b))$. If $k_1=2$, the solutions of
(\ref{5}) are $\lambda(a)=(2^i)$, $\lambda(b)=(2^{p-i})$ with
$0\leq i\leq p$. So the summand is
$$
\bigoplus_{i=0}^p(S_{(2^i)}V_1\otimes
S_{(2^{p-i})}V_1)^{SL\,V_1}\otimes
(S_{(2^i)}V_2)^{SO\,V_2}\otimes(S_{(2^{p-i})}V_3)^{SO\,V_3}
$$
which is generated by the coefficients of $\varphi^{p-i}\psi^i$ in
$det(\psi V(\sigma(a)a)+\varphi V(\sigma(b)b))$, semi-invariants
of weight $(2,0,0)$. In particular for $i=0$ we have
$det\,V(\sigma(b)b)$ and for $i=p$ we have $det\,V(\sigma(a)a)$.\\
\textbf{Sp)} The ring of symplectic semi-invariants is
$$
\bigoplus_{\lambda(a),\lambda(b)\in\Lambda}(S_{\lambda(a)}V_1\otimes
S_{\lambda(b)}V_1)^{SL\,V_1}\otimes
(S_{\lambda(a)}V_2)^{Sp\,V_2}\otimes(S_{\lambda(b)}V_3)^{Sp\,V_3}.
$$
By proposition \ref{i2} we have
\begin{equation}\label{6}
\lambda(a)_j+\lambda(b)_{p-j+1}=k_1
\end{equation}
for every $0\leq j\leq p$ and for some $k_1\in\mathbb{Z}_{\geq
0}$. By proposition \ref{i3} $\lambda(a)$ and $\lambda(b)$ have to
be
in $EC\Lambda$.\\
Let $p$ be even. We consider the summands in which $k_1=1$ because
the other ones are generated by products of powers of the
generators of this summands. The solutions of (\ref{6}) are
 $\lambda(a)=(1^{2i})$, $\lambda(b)=(1^{p-2i})$ with $0\leq
i\leq \frac{p}{2}$. So the summand is
$$
\bigoplus_{i=0}^{\frac{p}{2}}(S_{(1^{2i})}V_1\otimes
S_{(1^{p-2i})}V_1)^{SL\,V_1}\otimes
(S_{(1^{2i})}V_2)^{SO\,V_2}\otimes(S_{(1^{p-2i})}V_3)^{SO\,V_3}
$$
which is generated by the coefficients of
$\varphi^{\frac{p}{2}-i}\psi^i$ in $pf(\psi V(\sigma(a)a)+\varphi
V(\sigma(b)b))$, semi-invariants of weight $(1,0,0)$. In
particular for $i=0$ we have
$pf\,V(\sigma(b)b)=\sqrt{det\,V(\sigma(b)b)}=\sqrt{det\,V(\sigma(b))\cdot
det\,V(b)}=\sqrt{(det\,V(b))^2}=det\,V(b)$ and for
$i=\frac{p}{2}$ we have $pf\,V(\sigma(a)a)=det\,V(a)$.\\
If $p$ is odd there not exist any non-trivial symplectic
representations because a symplectic space of dimension odd
doesn't exist. So we have $SpSI(Q,ph)=\mathbb{K}$. $\Box$
\subsection{$\widetilde{A}^{1,1}_{k,l}$ for dimension vector $ph$}
\begin{teorema}
Let $(Q,\sigma)$ be a symmetric quiver of type $(1,1,k,l)$ with
orientation
$$
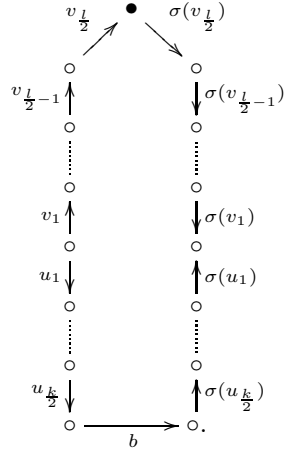
\xymatrix@-1pc{ &\bullet\ar[dr]^{\sigma(v_{\frac{l}{2}})}&\\
\circ\ar[ur]^{v_{\frac{l}{2}}}&&\circ\ar[d]^{\sigma(v_{\frac{l}{2}-1})}\\
\circ\ar[u]^{v_{\frac{l}{2}-1}}\ar@{.}[d]&&\circ\\
\circ&&\circ\ar[d]^{\sigma(v_{1})}\ar@{.}[u]\\
\circ\ar[d]_{u_1}\ar[u]^{v_1}&&\circ\\
\circ&&\circ\ar[u]_{\sigma(u_{1})}\ar@{.}[d]\\
\circ\ar[d]_{u_{\frac{k}{2}}}\ar@{.}[u]&&\circ\\
 \circ\ar[rr]_{b}&&\circ\ar[u]_{\sigma(u_{\frac{k}{2}})}.}
$$
Then\\
O) $OSI(Q,ph)$ is generated by the following indecomposable
semi-invariants:\\
if $p$ is even,
\begin{itemize}
\item[a)] $det\,V(u_j)$ with $j\in\{1,\ldots,\frac{k}{2}\}$;
\item[b)] $det\,V(v_j)$ with $j\in\{1,\ldots,\frac{l}{2}\}$;
\item[c)] $pf\,V(b)$
\item[d)] the coefficients $c_i$ of $\varphi^{p-2i}\psi^{2i}$, $0\leq i\leq \frac{p}{2}$, in $det(\psi V(\sigma(\bar{a})\bar{a})+\varphi
V(\bar{b}))$, where $\bar{a}=v_{\sigma(1)}\cdots
\sigma(v_{\frac{l}{2}}) v_{\frac{l}{2}}\cdots v_1$ and
$\bar{b}=\sigma(u_{1})\cdots \sigma(u_{\frac{k}{2}})b
u_{\frac{k}{2}}\cdots u_1$;
\end{itemize}
if $p$ is odd,
\begin{itemize}
\item[a)] $det\,V(u_j)$ with $j\in\{1,\ldots,\frac{k}{2}\}$;
\item[b)] $det\,V(v_j)$ with $j\in\{1,\ldots,\frac{l}{2}\}$;
\item[c)] the coefficients $c_i$ of $\varphi^{p-2i}\psi^{2i}$, $0\leq i\leq \frac{p-1}{2}$, in $det(\psi V(\sigma(\bar{a})\bar{a})+\varphi
V(\bar{b}))$, where $\bar{a}=v_{\sigma(1)}\cdots
\sigma(v_{\frac{l}{2}}) v_{\frac{l}{2}}\cdots v_1$ and
$\bar{b}=\sigma(u_{1})\cdots \sigma(u_{\frac{k}{2}})b
u_{\frac{k}{2}}\cdots u_1$.
\end{itemize}
Sp) $SpSI(Q,ph)$ is generated by the following indecomposable
semi-invariants:\\
if $p$ is even,
\begin{itemize}
\item[a)] $det\,V(u_j)$ with $j\in\{1,\ldots,\frac{k}{2}\}$;
\item[b)] $det\,V(v_j)$ with $j\in\{1,\ldots,\frac{l}{2}\}$;
\item[c)] $det\,V(b)$
\item[d)] the coefficients $c_i$ of $\varphi^{p-2i}\psi^{2i}$, $0\leq i\leq \frac{p}{2}$, in $det(\psi V(\sigma(\bar{a})\bar{a})+\varphi
V(\bar{b}))$, where $\bar{a}=v_{\sigma(1)}\cdots
\sigma(v_{\frac{l}{2}}) v_{\frac{l}{2}}\cdots v_1$ and
$\bar{b}=\sigma(u_{1})\cdots \sigma(u_{\frac{k}{2}})b
u_{\frac{k}{2}}\cdots u_1$;
\end{itemize}
if $p$ is odd, $SpSI(Q,ph)=\mathbb{K}$.
\end{teorema}
\textit{Proof.} We proceed by induction on
$\frac{k}{2}+\frac{l}{2}$. The smallest case is
$\widetilde{A}^{1,1}_{0,2}$
$$
\xymatrix@-1pc{&2\ar[dr]^{\sigma(a)}&\\
1\ar[rr]_b\ar[ur]^a&&\sigma(1)}
$$
and so it's enough to study the semi-invariants of
$\widetilde{A}^{1,1}_{0,2}$.\\
The induction step follows by lemma \ref{cls} and by lemma \ref{cl}, so it's enough to prove the theorem for $\widetilde{A}^{1,1}_{0,2}$.\\
\textbf{O)} The ring of orthogonal semi-invariants is
$$
\bigoplus_{{\lambda(a)\in\Lambda \atop\lambda(b)\in EC
\Lambda}}(S_{\lambda(a)}V_1\otimes
S_{\lambda(b)}V_1)^{SL\,V_1}\otimes (S_{\lambda(a)}V_2)^{SO\,V_2}.
$$
By proposition \ref{i2} we have
\begin{equation}\label{7}
\lambda(a)_j+\lambda(b)_{p-j+1}=k_1
\end{equation}
for every $0\leq j\leq p$ and for some $k_1\in\mathbb{Z}_{\geq
0}$. By proposition \ref{i3} we have $\lambda(a)=2\mu+(l^p)$ for
some $\mu\in\Lambda$ and for some $l\in\mathbb{Z}_{\geq 0}$. We
consider the summands in which $k_1=1,2$ because the other ones
are generated by products of powers of the generators of this
summands. Let $p$ be even. If $k_1=1$ the only solutions of
(\ref{7}) are $\lambda(a)=(1^p)$, $\lambda(b)=0$ and
$\lambda(a)=0$, $\lambda(b)=(1^p)$. Respectively, the summand
$(S_{(1^p)}V_1)^{SL\,V_1}\otimes(S_{(1^p)}V_2)^{SO\,V_2}$ is
generated by a semi-invariant of weight $(1,0)$, i.e
$det\,V(a)=det\,V(\sigma(a))$, and the summand
$(S_{(1^p)}V_1)^{SL\,V_1}$ is generated by a semi-invariant of
weight $(1,0)$, i.e $pf\,V(b)$. If $k_1=2$, the solutions of
(\ref{7}) are $\lambda(a)=(2^{2i})$, $\lambda(b)=(2^{p-2i})$ with
$0\leq i\leq \frac{p}{2}$. So the summand is
$$
\bigoplus_{i=0}^{\frac{p}{2}}(S_{(2^{2i})}V_1\otimes
S_{(2^{p-2i})}V_1)^{SL\,V_1}\otimes (S_{(2^{2i})}V_2)^{SO\,V_2}
$$
which is generated by the coefficients of
$\varphi^{p-2i}\psi^{2i}$ in $det(\psi V(\sigma(a)a)+\varphi
V(b))$, semi-invariants of weight $(2,0)$. In particular for $i=0$
we have
$det\,V(b)$ and for $i=\frac{p}{2}$ we have $det\,V(\sigma(a)a)$.\\
Let $p$ be odd. If $k_1=1$ the only solutions of (\ref{7}) are
$\lambda(a)=(1^p)$, $\lambda(b)=0$. The summand
$(S_{(1^p)}V_1)^{SL\,V_1}\otimes(S_{(1^p)}V_2)^{SO\,V_2}$ is
generated by a semi-invariant of weight $(1,0)$, i.e
$det\,V(a)=det\,V(\sigma(a))$. If $k_1=2$, the solutions of
(\ref{7}) are $\lambda(b)=(2^{2i})$, $\lambda(a)=(2^{p-2i})$ with
$0\leq i\leq \frac{p-1}{2}$. So the summand is
$$
\bigoplus_{i=0}^{\frac{p-1}{2}}(S_{(2^{p-2i})}V_1\otimes
S_{(2^{2i})}V_1)^{SL\,V_1}\otimes (S_{(2^{p-2i})}V_2)^{SO\,V_2}
$$
which is generated by the coefficients of
$\varphi^{2i}\psi^{p-2i}$ in $det(\psi V(\sigma(a)a)+\varphi
V(b))$, semi-invariants of weight $(2,0)$. In particular for $i=\frac{p-1}{2}$ we have $det\,V(\sigma(a)a)$.\\
\textbf{Sp)} The ring of symplectic semi-invariants is
$$
\bigoplus_{{\lambda(a)\in\Lambda \atop\lambda(b)\in ER
\Lambda}}(S_{\lambda(a)}V_1\otimes
S_{\lambda(b)}V_1)^{SL\,V_1}\otimes (S_{\lambda(a)}V_2)^{Sp\,V_2}.
$$
By proposition \ref{i2} we have
\begin{equation}\label{8}
\lambda(a)_j+\lambda(b)_{p-j+1}=k_1
\end{equation}
for every $0\leq j\leq p$ and for some $k_1\in\mathbb{Z}_{\geq
0}$. By proposition \ref{i3} we have $\lambda(a)\in EC\Lambda$. We
consider the summands in which $k_1=1,2$ because the other ones
are generated by products of powers of the generators of this
summands. Let $p$ be even. If $k_1=1$ the only solutions of
(\ref{8}) are $\lambda(a)=(1^p)$, $\lambda(b)=0$. The summand
$(S_{(1^p)}V_1)^{SL\,V_1}\otimes(S_{(1^p)}V_2)^{Sp\,V_2}$ is
generated by a semi-invariant of weight $(1,0)$, i.e
$det\,V(a)=det\,V(\sigma(a))=pf\,V(\sigma(a)a)$. If $k_1=2$, the
solutions of (\ref{8}) are $\lambda(a)=(2^{2i})$,
$\lambda(b)=(2^{p-2i})$ with $0\leq i\leq \frac{p}{2}$. So the
summand is
$$
\bigoplus_{i=0}^{\frac{p}{2}}(S_{(2^{2i})}V_1\otimes
S_{(2^{p-2i})}V_1)^{SL\,V_1}\otimes (S_{(2^{2i})}V_2)^{Sp\,V_2}
$$
which is generated by the coefficients of
$\varphi^{p-2i}\psi^{2i}$ in $det(\psi V(\sigma(a)a)+\varphi
V(b))$, semi-invariants of weight $(2,0)$. In particular for $i=0$
we have
$det\,V(b)$ and for $i=\frac{p}{2}$ we have $det\,V(\sigma(a)a)$.\\
If $p$ is odd there not exist any non-trivial symplectic
representations because a symplectic space of dimension odd
doesn't exist. So we have $SpSI(Q,ph)=\mathbb{K}$. $\Box$
\subsection{$\widetilde{A}^{0,0}_{k,k}$ for dimension vector
$ph$}
\begin{teorema}
Let $(Q,\sigma)$ be a symmetric quiver of type $(0,0,k,k)$ with
orientation
$$
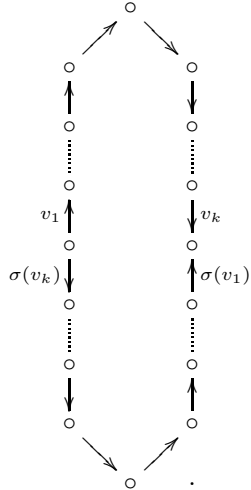
\xymatrix@-1pc{ &\circ\ar[dr]&\\
\circ\ar[ur]&&\circ\ar[d]\\
\circ\ar[u]\ar@{.}[d]&&\circ\\
\circ&&\circ\ar[d]^{v_{k}}\ar@{.}[u]\\
\circ\ar[d]_{\sigma(v_{k})}\ar[u]^{v_1}&&\circ\\
\circ&&\circ\ar[u]_{\sigma(v_{1})}\ar@{.}[d]\\
\circ\ar[d]\ar@{.}[u]&&\circ\\
 \circ\ar[dr]&&\circ\ar[u]\\
 &\circ\ar[ur]&.}
$$
Then\\
$OSI(Q,ph)=SpSI(Q,ph)$ is generated by the following
indecomposable semi-invariants:
\begin{itemize}
\item[a)] $det\,V(v_j)$ with $j\in\{1,\ldots,k\}$;
\item[b)] $pf(V(\bar{a})+
V(\sigma(\bar{a})))$
\item[c)] the coefficients $c_i$ of $\varphi^{p-i}\psi^i$, $0\leq i\leq p$, in $det(\psi V(\bar{a})+\varphi
V(\sigma(\bar{a})))$, where $\bar{a}=v_{k}\cdots v_{1}$.
\end{itemize}
\end{teorema}
\textit{Proof.} We proceed by induction on
$\frac{k}{2}+\frac{h}{2}$. The smallest case is
$\widetilde{A}^{0,0}_{2,2}$
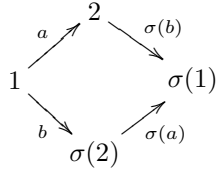
$$
\xymatrix@-1pc{&2\ar[dr]^{\sigma(b)}&\\
1\ar[dr]_b\ar[ur]^a&&\sigma(1)\\
&\sigma(2)\ar[ur]_{\sigma(a)}&}
$$
and so it's enough to study the semi-invariants of
$\widetilde{A}^{0,0}_{2,2}$.\\
The induction step follows by lemma \ref{cl}, so it's enough to prove the theorem for $\widetilde{A}^{0,0}_{2,2}$.\\
In this case we have $ORep(Q,ph)=SpRep(Q,ph)$ and so
$OSI(Q,ph)=SpSI(Q,ph)$. The ring of semi-invariants is
$$
\bigoplus_{\lambda(a),\lambda(b)\in\Lambda}(S_{\lambda(a)}V_1\otimes
S_{\lambda(b)}V_1)^{SL\,V_1}\otimes (S_{\lambda(a)}V_2^*\otimes
S_{\lambda(b)}V_2 )^{SL\,V_2}.
$$
By proposition \ref{i2} we have
\begin{equation}\label{9}
\left\{\begin{array}{c} \lambda(a)_j+\lambda(b)_{p-j+1}=k_1\\
\lambda(a)_j=\lambda(b)_{j}+k_2\end{array}\right.
\end{equation}
for every $0\leq j\leq p$ and for some $k_1,k_2\in\mathbb{Z}_{\geq
0}$. We consider the summands in which $k_1=1,2$ and $k_2=0,1$
because the other ones are generated by products of powers of the
generators of this summands. Let $p$ even. If $k_1=1$ and $k_2=0$
the only solution of (\ref{9}) are $\lambda(a)=(1^{\frac{p}{2}})$,
$\lambda(b)=(1^{\frac{p}{2}})$. The summand
$(S_{(1^{\frac{p}{2}})}V_1\otimes
S_{(1^{\frac{p}{2}})}V_1)^{SL\,V_1}\otimes
(S_{(1^{\frac{p}{2}})}V_2^*\otimes S_{(1^{\frac{p}{2}})}V_2
)^{SL\,V_2}$ is generated by a semi-invariant of weight $(1,0)$,
i.e. $pf(V(\sigma(b)a)+ V(\sigma(a)b))$. If $k_1=2$ and $k_2=0$,
the solutions of (\ref{9}) are $\lambda(a)=(2^i,1^{p-2i})$,
$\lambda(b)=(2^i,1^{p-2i})$ with $0\leq i\leq \frac{p}{2}$. So the
summand is
$$
\bigoplus_{i=0}^{\frac{p}{2}}(S_{(2^i,1^{p-2i})}V_1\otimes
S_{(2^i,1^{p-2i})}V_1)^{SL\,V_1}\otimes
(S_{(2^i,1^{p-2i})}V_2^*\otimes S_{(2^i,1^{p-2i})}V_2 )^{SL\,V_2}
$$
which is generated by the coefficients of $\varphi^{p-i}\psi^{i}$
with $0\leq i\leq \frac{p}{2}$ in $det(\psi V(\sigma(b)a)+\varphi
V(\sigma(a)b))$, semi-invariants of weight $(2,0)$. In particular
for $i=0$ we have $det\,V(\sigma(b)a)=det\,V(\sigma(a)b)$. Let $p$
be odd. If $k_1=1$ and $k_2=0$ we don't have any solutions of
(\ref{9}). If $k_1=2$ and $k_2=0$, the solutions of (\ref{9}) are
$\lambda(a)=(2^i,1^{p-2i})$, $\lambda(b)=(2^i,1^{p-2i})$ with
$0\leq i\leq \frac{p-1}{2}$. So the summand is
$$
\bigoplus_{i=0}^{\frac{p-1}{2}}(S_{(2^i,1^{p-2i})}V_1\otimes
S_{(2^i,1^{p-2i})}V_1)^{SL\,V_1}\otimes
(S_{(2^i,1^{p-2i})}V_2^*\otimes S_{(2^i,1^{p-2i})}V_2 )^{SL\,V_2}
$$
which is generated by the coefficients of $\varphi^{p-i}\psi^{i}$
with $0\leq i\leq \frac{p-1}{2}$ in $det(\psi
V(\sigma(b)a)+\varphi V(\sigma(a)b))$, semi-invariants of weight
$(2,0)$. In particular for $i=0$ we have
$det\,V(\sigma(b)a)=det\,V(\sigma(a)b)$.\\
If $k_2=1$, in both cases $p$ even or odd, $k_1$ can't be 0
otherwise we have $\lambda(b)_j+\lambda(b)_{p-j+1}=-1$ but this is
impossible. So $k_1=1$ and the only solutions of (\ref{9}) are
$\lambda(a)=(1^p)$, $\lambda(b)=0$ and $\lambda(a)=0$,
$\lambda(b)=(1^p)$; respectively we have the summand
$(S_{(1^p)}V_1)^{SL\,V_1}\otimes(S_{(1^p)}V_2^*)^{SL\,V_2}$
generated by the semi-invariant $det\,V(a)$ of weight $(1,-1)$ and
the summand
$(S_{(1^p)}V_1)^{SL\,V_1}\otimes(S_{(1^p)}V_2)^{SL\,V_2}$
generated by the semi-invariant $det\,V(b)$ of weight $(1,-1)$.
$\Box$
\subsection{$\widetilde{D}^{1,0}_n$ for dimension vector $ph$}
\begin{teorema}
Let $(Q,\sigma)$ be a symmetric quiver of type
$\widetilde{D}^{1,0}_n$ with orientation
$$
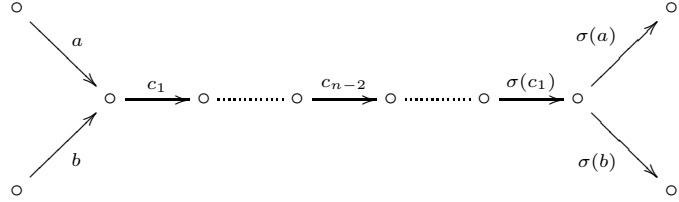
\xymatrix{\circ\ar[dr]^a&&&&&&&\circ\\
&\circ\ar[r]^{c_1}&\circ\ar@{.}[r]&\circ\ar[r]^{c_{n-2}}&\circ\ar@{.}[r]&\circ\ar[r]^{\sigma(c_1)}&\circ\ar[ur]^{\sigma(a)}\ar[dr]_{\sigma(b)}&\\
\circ\ar[ur]_b&&&&&&&\circ}
$$
and let $\bar{c}=\sigma(c_1)\cdots c_{n-2}\cdots c_1$.
Then\\
Sp) $SpSI(Q,ph)$ is generated by the following indecomposable
semi-invariants:\\
\begin{itemize}
\item[a)] $det\,V(c_j)$ with $j\in\{1,\ldots,n-2\}$
\item[b)] $det\left(V(a),V(b)\right)=det\left(\begin{array}{c}V(\sigma(a))\\
V(\sigma(b))\end{array}\right)$
\item[c)] $det\,V(\sigma(a)\bar{c}a)$
\item[d)] $det\,V(\sigma(b)\bar{c}b)$
\item[e)] $det\,V(\sigma(b)\bar{c}a)=det\,V(\sigma(a)\bar{c}b)$
\item[f)] the coefficients $c_i$ of $\varphi^{i}\psi^i$, $0\leq i\leq p$, in
$$det\left(\begin{array}{cc} \varphi
V(\sigma(a)\bar{c}a)&V(\sigma(b)\bar{c}a)\\
V(\sigma(a)\bar{c}b)&\psi
V(\sigma(b)\bar{c}b)\end{array}\right).$$
\end{itemize}
O) $OSI(Q,ph)$ is generated by the following indecomposable
semi-invariants:\\
if $p$ is even,
\begin{itemize}
\item[a)] $det\,V(c_j)$ with $j\in\{1,\ldots,n-2\}$;
\item[b)] $det\left(V(a),V(b)\right)=det\left(\begin{array}{c}V(\sigma(a))\\
V(\sigma(b))\end{array}\right)$
\item[c)] $pf\,V(\sigma(a)\bar{c}a)$
\item[d)] $pf\,V(\sigma(b)\bar{c}b)$
\item[e)] the coefficients $c_i$ of $\varphi^{i}\psi^i$, $0\leq i\leq \frac{p}{2}$, in
$$pf\left(\begin{array}{cc} \varphi
V(\sigma(a)\bar{c}a)&V(\sigma(b)\bar{c}a)\\
V(\sigma(a)\bar{c}b)&V(\sigma(b)\bar{c}b)\end{array}\right);$$
\end{itemize}
if $p$ is odd,
\begin{itemize}
\item[a)] $det\,V(c_j)$ with $j\in\{1,\ldots,n-2\}$
\item[b)] $det\left(V(a),V(b)\right)=det\left(\begin{array}{c}V(\sigma(a))\\
V(\sigma(b))\end{array}\right)$
\item[c)] the coefficients $c_i$ of $\varphi^{i}\psi^i$, $0\leq i\leq \frac{p-1}{2}$, in
$$pf\left(\begin{array}{cc} \varphi
V(\sigma(a)\bar{c}a)&V(\sigma(b)\bar{c}a)\\
V(\sigma(a)\bar{c}b)&\psi
V(\sigma(b)\bar{c}b)\end{array}\right).$$
\end{itemize}
\end{teorema}
\textit{Proof.} We proceed by induction on $n$. The smallest case
is $(\widetilde{D}^{1,0}_{3})^{eq}$
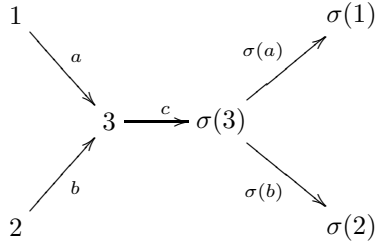
$$
\xymatrix{1\ar[dr]^a&&&\sigma(1)\\
&3\ar[r]^{c}&\sigma(3)\ar[ur]^{\sigma(a)}\ar[dr]_{\sigma(b)}&\\
2\ar[ur]_b&&&\sigma(2)}
$$
The induction step follows by lemma \ref{cls}, so it's enough to prove the theorem for $(\widetilde{D}^{1,0}_{3})^{eq}$.\\
Let $V$ be a representation of $(\widetilde{D}^{1,0}_{3})^{eq}$ of
dimension $ph$ for some $p\in\mathbb{Z}_{\geq 0}$, in this case $h=(1,1,2)$.\\
\textbf{Sp)} The ring of symplectic semi-invariants is
$$
SpSI(\widetilde{D}^{1,0}_{3},ph)=\bigoplus_{{\lambda(a),\lambda(b)\in
\Lambda\atop\lambda(c)\in
ER\Lambda}}(S_{\lambda(a)}V_1)^{SL\,V_1}\otimes
S_{\lambda(b)}V_2)^{SL\,V_2}\otimes
$$
$$
(S_{\lambda(a)}V_3^*\otimes S_{\lambda(b)}V_3^*\otimes
S_{\lambda(c)}V_3)^{SL\,V_3}.
$$
By proposition \ref{i1} we have $\lambda(a)=(k_1^p)$,
$\lambda(b)=(k_2^p)$, for some $k_1,k_2\in\mathbb{Z}_{\geq 0}$,
and by proposition \ref{clrr} we have
\begin{equation}\label{10}
S_{(k_1^p)}V_3^*\otimes S_{(k_2^p)}V_3^*=\bigoplus_{i=0}^p
S_{\nu_i}V_3^*
\end{equation}
where
$$
v_i=(k_1+\lambda_1,\ldots,k_1+\lambda_{p-i},\underbrace{k_1,\ldots,k_1}_i,\underbrace{k_2,\ldots,k_2}_i,k_2-\lambda_{p-i},\ldots,k_2-\lambda_1)
$$
with $0\leq\lambda_{p-i}\leq\ldots\leq\lambda_1\leq k_2$ and for
every $0\leq i\leq p$. Moreover we have
$$
(S_{\nu_i}V_3^*\otimes S_{\lambda(c)}V_3)^{SL\,V_3}\neq
0\Leftrightarrow \lambda(c)=v_i+(k_3^{2p})
$$
for some $k_3\in\mathbb{Z}_{\geq 0}$.\\
We consider the summands in which $k_1=0,1,2$ and $k_2=0,1,2$
because the other ones are generated by products of powers of the
generators of these summands.\\
If $\lambda(c)=0$, then $\lambda(a)=(k_1^p)\neq 0
\neq\lambda(b)=(k_2^p)$ because otherwise if for example
$\lambda(a)=0$ we have $(S_{(k_2^p)}V_3^*)^{SL\,V_3}=0$. We
consider the summand in which $\lambda(c)=0$ and $k_1=1=k_2$, the
only $\nu_i$ such that $(S_{\nu_i}V_3^*)^{SL\,V_3}\neq 0$ is
$\nu_p=(1^{2p})$. So $(S_{(1^p)}V_1)^{SL\,V_1}\otimes
(S_{(1^p)}V_2)^{SL\,V_2}\otimes (S_{(1^{2p})}V_3)^{SL\,V_3}$ is
generated by a semi-invariant of weight $(1,1,-1)$, i.e. $det\left(V(a),V(b)\right)=det\left(\begin{array}{c}V(\sigma(a))\\
V(\sigma(b))\end{array}\right)$. Now we suppose $\lambda(c)\neq
0$. We can't consider $k_1=1$, $k_2=0$ and $k_1=0$, $k_2=1$
because otherwise we haven't $\lambda(c)$ with even rows. If
$k_1=2$, $k_2=0$ and $k_3=0$ the summand
$(S_{(2^p)}V_1)^{SL\,V_1}\otimes(S_{(2^p)}V_3^*\otimes
S_{(2^p)}V_3)^{SL\,V_3}$ is generated by a semi-invariant of
weight $(2,0,0)$, i.e. $det\,V(\sigma(a)\bar{c}a)$. If $k_1=0$,
$k_2=2$ and $k_3=0$ the summand
$(S_{(2^p)}V_2)^{SL\,V_2}\otimes(S_{(2^p)}V_3^*\otimes
S_{(2^p)}V_3)^{SL\,V_3}$ is generated by a semi-invariant of
weight $(0,2,0)$, i.e. $det\,V(\sigma(b)\bar{c}b)$. If
$k_1=0=k_2$, then $k_3$ has to be even. So, considering $k_3=2$,
$(S_{(2^{2p})}V_3)^{SL\,V_3}$ is generated by a semi-invariant of
weight $(0,0,2)$, i.e. $det\,V(c)$. If $k_1=k_2=1$, by (\ref{10}),
$\lambda(c)=(2^p)$. So
$(S_{(1^p)}V_1)^{SL\,V_1}\otimes(S_{(1^p)}V_2)^{SL\,V_2}\otimes(S_{(2^p)}V_3^*\otimes
S_{(2^p)}V_3)^{SL\,V_3}$ is generated by a semi-invariant of
weight $(1,1,0)$, i.e.
$det\,V(\sigma(b)\bar{c}a)=det\,V(\sigma(a)\bar{c}b)$. Finally if
$k_1=k_2=2$, considering $k_3=0$, the summand is
$$
(S_{(2^p)}V_1)^{SL\,V_1}\otimes(S_{(2^p)}V_2)^{SL\,V_2}\otimes(\bigoplus_{i=0}^p
S_{(4^{p-2i},2^{4i})}V_3^*\otimes
S_{(4^{p-2i},2^{4i})}V_3)^{SL\,V_3}
$$
which is generated by the coefficients of $\varphi^i\psi^i$ in
$$det\left(\begin{array}{cc} \varphi
V(\sigma(a)\bar{c}a)&V(\sigma(b)\bar{c}a)\\
V(\sigma(a)\bar{c}b)&\psi
V(\sigma(b)\bar{c}b)\end{array}\right),$$ semi-invariants of
weight $(2,2,0)$. In particular for $i=0$ we have\\
$(det\,V(\sigma(b)\bar{c}a))^2$ and for $i=p$ we have
$det\,V(\sigma(a)\bar{c}a)\cdot det\,V(\sigma(b)\bar{c}b)$.\\
\textbf{O)} The ring of orthogonal semi-invariants is
$$
SpSI(\widetilde{D}^{1,0}_{3},ph)=\bigoplus_{{\lambda(a),\lambda(b)\in
\Lambda\atop\lambda(c)\in
EC\Lambda}}(S_{\lambda(a)}V_1)^{SL\,V_1}\otimes
S_{\lambda(b)}V_2)^{SL\,V_2}\otimes
$$
$$
(S_{\lambda(a)}V_3^*\otimes S_{\lambda(b)}V_3^*\otimes
S_{\lambda(c)}V_3)^{SL\,V_3}.
$$
By proposition \ref{i1} we have $\lambda(a)=(k_1^p)$,
$\lambda(b)=(k_2^p)$, for some $k_1,k_2\in\mathbb{Z}_{\geq 0}$,
and by proposition \ref{clrr} we have
\begin{equation}\label{11}
S_{(k_1^p)}V_3^*\otimes S_{(k_2^p)}V_3^*=\bigoplus_{i=0}^p
S_{\nu_i}V_3^*
\end{equation}
where
$$
v_i=(k_1+\lambda_1,\ldots,k_1+\lambda_{p-i},\underbrace{k_1,\ldots,k_1}_i,\underbrace{k_2,\ldots,k_2}_i,k_2-\lambda_{p-i},\ldots,k_2-\lambda_1)
$$
with $0\leq\lambda_{p-i}\leq\ldots\leq\lambda_1\leq k_2$ and for
every $0\leq i\leq p$. Moreover we have
$$
(S_{\nu_i}V_3^*\otimes S_{\lambda(c)}V_3)^{SL\,V_3}\neq
0\Leftrightarrow \lambda(c)=v_i+(k_3^{2p})
$$
for some $k_3\in\mathbb{Z}_{\geq 0}$. Since $\lambda(c)\in EC\Lambda$, also $\nu_i\in EC\Lambda$ for every $i$.\\
We consider the summands in which $k_1=0,1$ and $k_2=0,1$ because
the other ones are generated by products of powers of the
generators of these summands.\\
As before if $\lambda(c)=0$, the only $\nu_i$ such that
$(S_{\nu_i}V_3^*)^{SL\,V_3}\neq 0$ is $\nu_p=(1^{2p})$. So
$(S_{(1^p)}V_1)^{SL\,V_1}\otimes (S_{(1^p)}V_2)^{SL\,V_2}\otimes
(S_{(1^{2p})}V_3)^{SL\,V_3}$ is
generated by a semi-invariant of weight $(1,1,-1)$, i.e. $det\left(V(a),V(b)\right)=det\left(\begin{array}{c}V(\sigma(a))\\
V(\sigma(b))\end{array}\right)$. Now we suppose $\lambda(c)\neq
0$.\\
 Let $p$ be even. If $k_1=1$, $k_2=0$ and $k_3=0$ the summand
$(S_{(1^p)}V_1)^{SL\,V_1}\otimes(S_{(1^p)}V_3^*\otimes
S_{(1^p)}V_3)^{SL\,V_3}$ is generated by a semi-invariant of
weight $(1,0,0)$, i.e. $pf\,V(\sigma(a)ca)$. If $k_1=0$, $k_2=1$
and $k_3=0$ the summand
$(S_{(1^p)}V_2)^{SL\,V_2}\otimes(S_{(1^p)}V_3^*\otimes
S_{(1^p)}V_3)^{SL\,V_3}$ is generated by a semi-invariant of
weight $(0,1,0)$, i.e. $pf\,V(\sigma(b)cb)$. If $k_1=0=k_2$, then
$k_3$ has to be not zero. So, considering $k_3=1$,
$(S_{(1^{2p})}V_3)^{SL\,V_3}$ is generated by a semi-invariant of
weight $(0,0,1)$, i.e. $pf\,V(c)$. Finally if $k_1=k_2=1$,
considering $k_3=0$, the summand is
$$
(S_{(1^p)}V_1)^{SL\,V_1}\otimes(S_{(1^p)}V_2)^{SL\,V_2}\otimes(\bigoplus_{i=0}^{\frac{p}{2}}
S_{(2^{p-2i},1^{4i})}V_3^*\otimes
S_{(2^{p-2i},1^{4i})}V_3)^{SL\,V_3}
$$
which is generated by the coefficients of $\varphi^i\psi^i$ in
$$pf\left(\begin{array}{cc} \varphi
V(\sigma(a)ca)&V(\sigma(b)ca)\\
V(\sigma(a)cb)&\psi V(\sigma(b)cb)\end{array}\right),$$
semi-invariants of weight $(1,1,0)$. In particular for $i=0$ we
have\\ $det\,V(\sigma(b)ca)=det\,V(\sigma(a)cb)$ and for
$i=\frac{p}{2}$ we have
$pf\,V(\sigma(a)ca)\cdot pf\,V(\sigma(b)cb)$.\\
Let $p$ be odd. In this case we can't consider $k_1=1$, $k_2=0$,
$k_3=0$ and $k_1=0$, $k_2=1$, $k_3=0$  because otherwise we have
$\lambda(c)=(1^p)$ with $p$ odd but $\lambda(c)$ has to be in
$EC\Lambda$. As before, if $k_1=0=k_2$, then $k_3$ has to be not
zero. So, considering $k_3=1$, $(S_{(1^{2p})}V_3)^{SL\,V_3}$ is
generated by a semi-invariant of weight $(0,0,1)$, i.e.
$pf\,V(c)$. Finally if $k_1=k_2=1$, considering $k_3=0$, the
summand is
$$
(S_{(1^p)}V_1)^{SL\,V_1}\otimes(S_{(1^p)}V_2)^{SL\,V_2}\otimes(\bigoplus_{i=0}^{\frac{p-1}{2}}
S_{(2^{p-(2i+1)},1^{4i+2})}V_3^*\otimes
S_{(2^{p-(2i+1)},1^{4i+2})}V_3)^{SL\,V_3}
$$
which is generated by the coefficients of $\varphi^i\psi^i$ in
$$pf\left(\begin{array}{cc} \varphi
V(\sigma(a)ca)&V(\sigma(b)ca)\\
V(\sigma(a)cb)&\psi V(\sigma(b)cb)\end{array}\right),$$
semi-invariants of weight $(1,1,0)$. In particular for $i=0$ we
have $det\,V(\sigma(b)ca)=det\,V(\sigma(a)cb)$. $\Box$
\subsection{$\widetilde{D}^{0,1}_n$ for dimension vector $ph$}
\begin{teorema}
Let $(Q,\sigma)$ be a symmetric quiver of type
$\widetilde{D}^{0,1}_n$ with orientation
$$
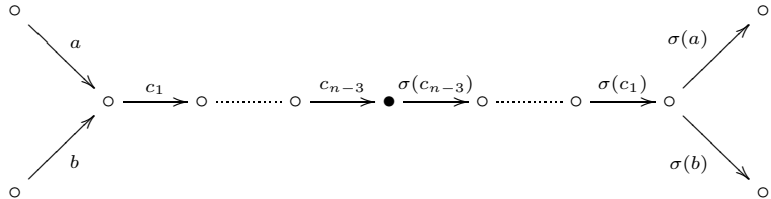
\xymatrix{\circ\ar[dr]^a&&&&&&&&\circ\\
&\circ\ar[r]^{c_1}&\circ\ar@{.}[r]&\circ\ar[r]^{c_{n-3}}&\bullet\ar[r]^{\sigma(c_{n-3})}&\circ\ar@{.}[r]&\circ\ar[r]^{\sigma(c_1)}&\circ\ar[ur]^{\sigma(a)}\ar[dr]_{\sigma(b)}&\\
\circ\ar[ur]_b&&&&&&&&\circ}
$$
and let $\bar{c}=\sigma(c_1)\cdots \sigma(c_{n-3})c_{n-3}\cdots
c_1$.
Then\\
O) $OSI(Q,ph)$ is generated by the following indecomposable
semi-invariants:\\
\begin{itemize}
\item[a)] $det\,V(c_j)$ with $j\in\{1,\ldots,n-3\}$
\item[b)] $det\left(V(a),V(b)\right)=det\left(\begin{array}{c}V(\sigma(a))\\
V(\sigma(b))\end{array}\right)$
\item[c)] $det\,V(\sigma(a)\bar{c}a)$
\item[d)] $det\,V(\sigma(b)\bar{c}b)$
\item[e)] $det\,V(\sigma(b)\bar{c}a)=det\,V(\sigma(a)\bar{c}b)$
\item[f)] the coefficients $c_i$ of $\varphi^{i}\psi^i$, $0\leq i\leq p$, in
$$det\left(\begin{array}{cc} \varphi
V(\sigma(a)\bar{c}a)&V(\sigma(b)\bar{c}a)\\
V(\sigma(a)\bar{c}b)&\psi
V(\sigma(b)\bar{c}b)\end{array}\right).$$
\end{itemize}
Sp) $SpSI(Q,ph)$ is generated by the following indecomposable
semi-invariants:\\
if $p$ is even,
\begin{itemize}
\item[a)] $det\,V(c_j)$ with $j\in\{1,\ldots,n-2\}$;
\item[b)] $det\left(V(a),V(b)\right)=det\left(\begin{array}{c}V(\sigma(a))\\
V(\sigma(b))\end{array}\right)$
\item[c)] $pf\,V(\sigma(a)\bar{c}a)$
\item[d)] $pf\,V(\sigma(b)\bar{c}b)$
\item[e)] $det\,V(\sigma(b)\bar{c}a)=det\,V(\sigma(a)\bar{c}b)$
\item[f)] the coefficients $c_i$ of $\varphi^{i}\psi^i$, $0\leq i\leq \frac{p}{2}$, in
$$pf\left(\begin{array}{cc} \varphi
V(\sigma(a)\bar{c}a)&V(\sigma(b)\bar{c}a)\\
V(\sigma(a)\bar{c}b)&\psi
V(\sigma(b)\bar{c}b)\end{array}\right);$$
\end{itemize}
if $p$ is odd,
\begin{itemize}
\item[a)] $det\,V(c_j)$ with $j\in\{1,\ldots,n-2\}$
\item[b)] $det\left(V(a),V(b)\right)=det\left(\begin{array}{c}V(\sigma(a))\\
V(\sigma(b))\end{array}\right)$
\item[c)] $det\,V(\sigma(b)\bar{c}a)=det\,V(\sigma(a)\bar{c}b)$
\item[d)] the coefficients $c_i$ of $\varphi^{i}\psi^i$, $0\leq i\leq \frac{p-1}{2}$, in
$$pf\left(\begin{array}{cc} \varphi
V(\sigma(a)\bar{c}a)&V(\sigma(b)\bar{c}a)\\
V(\sigma(a)\bar{c}b)&\psi
V(\sigma(b)\bar{c}b)\end{array}\right).$$
\end{itemize}
\end{teorema}
\textit{Proof.} We proceed by induction on $n$. The smallest case
is ($\widetilde{D}^{0,1}_{3})^{eq}$
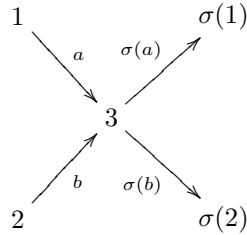
$$
\xymatrix{1\ar[dr]^a&&\sigma(1)\\
&3\ar[ur]^{\sigma(a)}\ar[dr]_{\sigma(b)}&\\
2\ar[ur]_b&&\sigma(2)}
$$
The induction step follows by lemma \ref{cl}, so it's enough to prove the theorem for $(\widetilde{D}^{0,1}_{3})^{eq}$.\\
Let $V$ be a representation of $(\widetilde{D}^{0,1}_{3})^{eq}$ of
dimension $ph$ for some $p\in\mathbb{Z}_{\geq 0}$, in this case $h=(1,1,2)$.\\
\textbf{O)} The ring of orthogonal semi-invariants is
$$
OSI(\widetilde{D}^{0,1}_{3},ph)=\bigoplus_{\lambda(a),\lambda(b)\in
\Lambda}(S_{\lambda(a)}V_1)^{SL\,V_1}\otimes
(S_{\lambda(b)}V_2)^{SL\,V_2}\otimes (S_{\lambda(a)}V_3^*\otimes
S_{\lambda(b)}V_3^*)^{SO\,V_3}.
$$
By proposition \ref{i1} we have $\lambda(a)=(k_1^p)$,
$\lambda(b)=(k_2^p)$, for some $k_1,k_2\in\mathbb{Z}_{\geq 0}$,
and by proposition \ref{clrr} we have
\begin{equation}\label{12}
S_{(k_1^p)}V_3^*\otimes S_{(k_2^p)}V_3^*=\bigoplus_{i=0}^p
S_{\nu_i}V_3^*
\end{equation}
where
$$
v_i=(k_1+\lambda_1,\ldots,k_1+\lambda_{p-i},\underbrace{k_1,\ldots,k_1}_i,\underbrace{k_2,\ldots,k_2}_i,k_2-\lambda_{p-i},\ldots,k_2-\lambda_1)
$$
with $0\leq\lambda_{p-i}\leq\ldots\leq\lambda_1\leq k_2$ and for
every $0\leq i\leq p$. Moreover we have
\begin{equation}\label{13}
(S_{\nu_i}V_3^*)^{SO\,V_3}\neq 0\Leftrightarrow
v_i=2\mu_i+(k_3^{2p})
\end{equation}
for some $k_3\in\mathbb{Z}_{\geq 0}$ and for some $\mu_i\in\Lambda$.\\
We consider the summands in which $k_1=0,1,2$ and $k_2=0,1,2$
because the other ones are generated by products of powers of the
generators of these summands.\\
We can't consider $k_1=1$, $k_2=0$ and $k_1=0$, $k_2=1$ because
otherwise we haven't $v_i=2\mu_i+(k_3^{2p})$. If $k_1=2$, $k_2=0$
the summand
$(S_{(2^p)}V_1)^{SL\,V_1}\otimes(S_{(2^p)}V_3^*)^{SO\,V_3}$ is
generated by a semi-invariant of weight $(2,0,0)$,\\ i.e.
$det\,V(\sigma(a)a)$. If $k_1=0$, $k_2=2$ the summand
$(S_{(2^p)}V_2)^{SL\,V_2}\otimes(S_{(2^p)}V_3^*)^{SO\,V_3}$ is
generated by a semi-invariant of weight $(0,2,0)$, i.e.
$det\,V(\sigma(b)b)$. If $k_1=k_2=1$, by (\ref{12}) and by
(\ref{13}), $\nu_i=(2^p)$ or $\nu_i=(1^{2p})$. So we have
$(S_{(1^p)}V_3^*\otimes
S_{(1^p)}V_3^*)^{SO\,V_3}=(S_{(2^p)}V_3^*\oplus
S_{(1^{2p})}V_3^*)^{SO\,V_3}$. Now
$$(S_{(1^p)}V_1)^{SL\,V_1}\otimes(S_{(1^p)}V_2)^{SL\,V_2}\otimes(S_{(1^{2p})}V_3^*)^{SO\,V_3}$$
is generated by a semi-invariant of weight $(1,1,0)$, i.e. $det\left(V(a),V(b)\right)=det\left(\begin{array}{c}V(\sigma(a))\\
V(\sigma(b))\end{array}\right)$ and
$(S_{(1^p)}V_1)^{SL\,V_1}\otimes(S_{(1^p)}V_2)^{SL\,V_2}\otimes(S_{(2^{p})}V_3^*)^{SO\,V_3}$
is generated by a semi-invariant of weight $(1,1,0)$, i.e.
$det\,V(\sigma(b)a)=det\,V(\sigma(a)b)$. Finally if $k_1=k_2=2$
the summand is
$$
(S_{(2^p)}V_1)^{SL\,V_1}\otimes(S_{(2^p)}V_2)^{SL\,V_2}\otimes(\bigoplus_{i=0}^p
S_{(4^{p-2i},2^{4i})}V_3^*)^{SO\,V_3}
$$
which is generated by the coefficients of $\varphi^i\psi^i$ in
$$det\left(\begin{array}{cc} \varphi
V(\sigma(a)a)&V(\sigma(b)a)\\
V(\sigma(a)b)&\psi V(\sigma(b)b)\end{array}\right),$$
semi-invariants of weight $(2,2,0)$. In particular for $i=0$ we
have $(det\,V(\sigma(b)a))^2$ and for $i=p$ we have
$det\,V(\sigma(a)a)\cdot det\,V(\sigma(b)b)$.\\
\textbf{Sp)} The ring of symplectic semi-invariants is
$$
SpSI(\widetilde{D}^{0,1}_{3},ph)=\bigoplus_{\lambda(a),\lambda(b)\in
\Lambda}(S_{\lambda(a)}V_1)^{SL\,V_1}\otimes
(S_{\lambda(b)}V_2)^{SL\,V_2}\otimes (S_{\lambda(a)}V_3^*\otimes
S_{\lambda(b)}V_3^*)^{Sp\,V_3}.
$$
By proposition \ref{i1} we have $\lambda(a)=(k_1^p)$,
$\lambda(b)=(k_2^p)$, for some $k_1,k_2\in\mathbb{Z}_{\geq 0}$,
and by proposition \ref{clrr} we have
\begin{equation}\label{14}
S_{(k_1^p)}V_3^*\otimes S_{(k_2^p)}V_3^*=\bigoplus_{i=0}^p
S_{\nu_i}V_3^*
\end{equation}
where
$$
v_i=(k_1+\lambda_1,\ldots,k_1+\lambda_{p-i},\underbrace{k_1,\ldots,k_1}_i,\underbrace{k_2,\ldots,k_2}_i,k_2-\lambda_{p-i},\ldots,k_2-\lambda_1)
$$
with $0\leq\lambda_{p-i}\leq\ldots\leq\lambda_1\leq k_2$ and for
every $0\leq i\leq p$. Moreover we have
\begin{equation}\label{15}
(S_{\nu_i}V_3^*)^{Sp\,V_3}\neq 0\Leftrightarrow v_i\in EC\Lambda.
\end{equation}
We consider the summands in which $k_1=0,1$ and $k_2=0,1$ because
the other ones are generated by products of powers of the
generators of these summands.\\
Let $p$ be even. If $k_1=1$, $k_2=0$ the summand
$(S_{(1^p)}V_1)^{SL\,V_1}\otimes(S_{(1^p)}V_3^*)^{Sp\,V_3}$ is
generated by a semi-invariant of weight $(1,0,0)$, i.e.
$pf\,V(\sigma(a)a)$. If $k_1=0$, $k_2=1$ the summand
$(S_{(1^p)}V_2)^{SL\,V_2}\otimes(S_{(1^p)}V_3^*)^{Sp\,V_3}$ is
generated by a semi-invariant of weight $(0,1,0)$, i.e.
$pf\,V(\sigma(b)b)$. Finally if $k_1=k_2=1$, the summand is
$$
(S_{(1^p)}V_1)^{SL\,V_1}\otimes(S_{(1^p)}V_2)^{SL\,V_2}\otimes(\bigoplus_{i=0}^{\frac{p}{2}}
S_{(2^{p-2i},1^{4i})}V_3^*)^{Sp\,V_3}
$$
which is generated by the coefficients of $\varphi^i\psi^i$ in
$$pf\left(\begin{array}{cc} \varphi
V(\sigma(a)a)&V(\sigma(b)a)\\
V(\sigma(a)b)&\psi V(\sigma(b)b)\end{array}\right),$$
semi-invariants of weight $(1,1,0)$. In particular for $i=0$ we
have\\ $det\,V(\sigma(b)a)=det\,V(\sigma(a)b)$ and for
$i=\frac{p}{2}$ we have
$pf\,V(\sigma(a)a)\cdot pf\,V(\sigma(b)b)$.\\
Let $p$ be odd. In this case we can't consider $k_1=1$, $k_2=0$
and $k_1=0$, $k_2=1$  because otherwise, by (\ref{15}),
$(S_{(1^p)}V_3)^{Sp\,V_3}=0$. Finally if $k_1=k_2=1$, the summand
is
$$
(S_{(1^p)}V_1)^{SL\,V_1}\otimes(S_{(1^p)}V_2)^{SL\,V_2}\otimes(\bigoplus_{i=0}^{\frac{p-1}{2}}
S_{(2^{p-(2i+1)},1^{4i+2})}V_3^*)^{SL\,V_3}
$$
which is generated by the coefficients of $\varphi^i\psi^i$ in
$$pf\left(\begin{array}{cc} \varphi
V(\sigma(a)a)&V(\sigma(b)a)\\
V(\sigma(a)b)&\psi V(\sigma(b)b)\end{array}\right),$$
semi-invariants of weight $(1,1,0)$. In particular for $i=0$ we
have $det\,V(\sigma(b)a)=det\,V(\sigma(a)b)$. $\Box$

\subsection{End of the proof of conjecture \ref{mt1} and \ref{mt2}
for dimension vector $ph$}\label{fineph} First of all we note
that, by definition of $c^W$ and $pf^W$, when we have it, are not
zero if $0=\langle \underline{dim}\,W,ph\rangle=p\langle
\underline{dim}\,W,h\rangle=-p\langle
h,\underline{dim}\,W\rangle$, so we have to consider only regular
representations $W$. Moreover it is enough to consider only simple
regular representations $W$, because the other regular
representations are extensions of simple regular ones and so, by
lemma \ref{cVcV'}, we obtain the $c^W$ and $pf^W$ with non-simple
regular $W$ as products of those with  simple regular $W$. Now we
check only for $\widetilde{A}^{2,0,1}_{k,l}$ and
$\widetilde{D}^{1,0}_{n}$ that the generators found for
$SpSI(Q,ph)$ and $OSI(Q,d)$ are of type $c^{W}$, for some simple
regular $W$, and $pf^W$, for some simple regular $W$ satisfying
property \textit{(Op)} in symplectic case and \textit{(Spp)} in
orthogonal case (see lemma \ref{ss1}). For the
other types of quivers it is similar (see also [D, section 4.1]).\\
We use notation of section \ref{Qtt}. For
$\widetilde{A}^{2,0,1}_{k,l}$, by definition of $c^W$ and $pf^W$,
\begin{itemize}
\item[\textbf{Sp)}] if $V$ is a symplectic representation, we have $c^{E_0}(V)=det(V(v_{\frac{l}{2}}))=det(V(v_1))=c^{E_1}(V)$,
$c^{E_i}(V)=det(V(v_i))=det(V(v_{\sigma(i)}))=c^{E_{\sigma(i)}}(V)$
for every $i\in\{2,\ldots,l\}\setminus\{\frac{l}{2}+1\}$,
$c^{E_{\frac{l}{2}+1}}(V)=det(V(a))$,
$c^{E'_0}(V)=det(V(u_{\frac{k}{2}}))=det(V(u_1))=c^{E'_1}(V)$,
$c^{E'_i}(V)=det(V(u_i))=det(V(u_{\sigma(i)}))=c^{E'_{\sigma(i)}}(V)$
for every $i\in\{2,\ldots,k\}\setminus\{\frac{k}{2}+1\}$,
$c^{E'_{\frac{k}{2}+1}}(V)=det(V(b))$ and
$c_{V_{(\varphi,\psi)}}(V)=det(\psi V(a)+\varphi V(b)$;
\item[\textbf{O)}] if $V$ is an orthogonal representation, the only differences with the symplectic case are,
when $p$ is even, we have $pf^{E_{\frac{l}{2}+1}}(V)=pf(V(a))$,
$pf^{E'_{\frac{k}{2}+1}}(V)=pf(V(b))$ and
$pf^{V_{(\varphi,\psi)}}(V)=pf(\psi V(a)+\varphi V(b)$, in fact
$E_{\frac{l}{2}+1}$, $E'_{\frac{k}{2}+1}$ and $V_{(\varphi,\psi)}$
satisfy property \textit{(Spp)}.
\end{itemize}
For $\widetilde{D}^{1,0}_n$, by definition of $c^W$ and $pf^W$,
\begin{itemize}
\item[\textbf{Sp)}] if $V$ is a symplectic representation, we have
$c^{E_0}(V)=det\left(\begin{array}{c}V(\sigma(a))\\V(\sigma(b))\end{array}\right)=det(V(a),V(b))=c^{E_1}(V)$,
$c^{E_i}(V)=det(V(c_{i-1}))=det(V(c_{\sigma(i-1)}))=c^{E_{\sigma(i)}}(V)$
for every $i\in\{2,\ldots,2n-3\}$,
$c^{E'_0}(V)=det(V(\sigma(b)\bar{c}a))=det(V(\sigma(a)\bar{c}b))=c^{E'_1}(V)$,
$c^{E''_0}(V)=det(V(\sigma(a)\bar{c}a))$,\\
$c^{E''_1}(V)=det(V(\sigma(b)\bar{c}b))$ and
$$c^{V_{(\varphi,\psi)}}(V)=det\left(\begin{array}{cc}\varphi
V(\sigma(a)\bar{c}a)&V(\sigma(b)\bar{c}a)\\V(\sigma(a)\bar{c}b)&\psi
V(\sigma(b)\bar{c}b)\end{array}\right);$$
\item[\textbf{O)}] if $V$ is an orthogonal representation, the
only differences with the symplectic case are that we have
$$pf^{V_{(\varphi,\psi)}}(V)=pf\left(\begin{array}{cc}\varphi
V(\sigma(a)\bar{c}a)&V(\sigma(b)\bar{c}a)\\V(\sigma(a)\bar{c}b)&\psi
V(\sigma(b)\bar{c}b)\end{array}\right),$$ since
$V_{(\varphi,\psi)}$ satisfies property \textit{(Spp)} and
$c^{E'_0}(V)=det(V(\sigma(b)\bar{c}a))=det(V(\sigma(a)\bar{c}b))=c^{E'_1}(V)$
is the coefficient of $\varphi^0\psi^0=1$ in
$pf^{V_{(\varphi,\psi)}}(V)$; moreover, if $p$ is even, we have
$pf^{E_{n-1}}(V)=pf(V(c_{n-2}))$,
$pf^{E''_0}(V)=pf(V(\sigma(a)\bar{c}a))$ and
$pf^{E''_1}(V)=pf(V(\sigma(b)\bar{c}b))$, because $E_{n-1}$,
$E''_{0}$ and $E''_1$ satisfy property \textit{(Spp)}.
\end{itemize}

\section{Semi-invariants of symmetric
quivers of tame type for any regular dimension vector} In this
section we prove theorems \ref{mt1} and \ref{mt2} for symmetric
quiver of tame type and any regular symmetric dimension vector
$d$.\\
We will use the same notation of section 3.1. For the type
$\widetilde{A}$ we call $a_0=tv_1=tu_1$, $x_i=hv_i$ for every
$i\in\{1,\dots,\frac{l}{2}\}$ and $y_i=hv_i$ for every
$i\in\{1,\dots,\frac{k}{2}\}$. For the type $\widetilde{D}$ we
call $t_1=ta$, $t_2=tb$ and $z_i=tc_i$ for every $i$ such that
$c_i\in (Q_1^+\sqcup Q_1^{\sigma})\setminus\{a,b\}$.\\
First we consider the canonical decomposition of $d$ for the
symmetric quivers.\\
Let $(Q,\sigma)$ be a symmetric quiver of tame type and let
$\Delta=\{e_i|\,i\in I=\{0,\ldots,u\}\}$, $\Delta'=\{e'_i|\,i\in
I'= \{0,\ldots,v\}\}$ and $\Delta''=\{e''_i|\,i\in
I''=\{0,\ldots,w\}\}$ be the three $\tau^+$-orbits of
nonhomogeneous simple regular representations of the underlying
quiver $Q$ (see proposition \ref{orbitetau}).\\
We shall call $I_{\delta}=\{i\in I|\,e_i=\delta e_i\}$
(respectively $I'_{\delta}$ and $I''_{\delta}$).
\begin{lemma}\label{dI}
Let $[x]:=max\{z\in\mathbb{N}|\,z\leq x\}$ is the floor of
$x\in\mathbb{R}$.\\
(1) For $\widetilde{A}_{k,l}^{2,0,1}$, we have:
\begin{itemize}
\item[(1.1)] decomposition $I=I_+\sqcup I_{\delta}\sqcup I_-$
where $I_+=\{2,\ldots,\frac{l}{2}+1\}$, $I_{\delta}=\{1\}$ and
$I_-=I\setminus (I_+\sqcup I_{\delta})$;
\item[(1.2)] decomposition $I'=I'_+\sqcup
I'_{\delta}\sqcup I'_-$ where $I'_+=\{2,\ldots,\frac{k}{2}+1\}$,
$I'_{\delta}=\{1\}$ and $I'_-=I'\setminus (I'_+\sqcup
I'_{\delta})$;
\item[(1.3)] $I''=\emptyset$.
\end{itemize}
(2) For $\widetilde{A}_{k,l}^{2,0,2}$, we have:
\begin{itemize}
\item[(2.1)] decomposition $I=I_+\sqcup I_{\delta}\sqcup I_-$ where
$I_+=\{2,\ldots,[\frac{l+1}{2}]+2\}$, $I_{\delta}=\emptyset$ and
$I_-=I\setminus I_+$;
\item[(2.2)] decomposition $I'=I'_+\sqcup
I'_{\delta}\sqcup I'_-$ where
$I'_+=\{2,\ldots,[\frac{k-1}{2}]+1\}$,
$I'_{\delta}=\{1,[\frac{k-1}{2}]+2\}$ and $I'_-=I'\setminus
(I'_+\sqcup I'_{\delta})$;
\item[(2.3)] $I''=\emptyset$.
\end{itemize}
(3) For $\widetilde{A}_{k,l}^{0,2}$, we have:
\begin{itemize}
\item[(3.1)] decomposition $I=I_+\sqcup I_{\delta}\sqcup I_-$
where  $I_+=\{2,\ldots,[\frac{l-1}{2}]+1\}$,
$I_{\delta}=\{1,[\frac{l-1}{2}]+2\}$ and $I_-=I\setminus
(I_+\sqcup I_{\delta})$;
\item[(3.2)] decomposition $I'=I'_+\sqcup
I'_{\delta}\sqcup I'_-$ where
$I'_+=\{2,\ldots,[\frac{k-1}{2}]+1\}$,
$I'_{\delta}=\{1,[\frac{k-1}{2}]+2\}$ and $I'_-=I'\setminus
(I'_+\sqcup I'_{\delta})$;
\item[(3.3)] $I''=\emptyset$.
\end{itemize}
(4) For $\widetilde{A}_{k,l}^{1,1}$, we have:
\begin{itemize}
\item[(4.1)] decomposition $I=I_+\sqcup I_{\delta}\sqcup I_-$
where  $I_+=\{2,\ldots,[\frac{l-1}{2}]+1\}$,
$I_{\delta}=\{1,[\frac{l-1}{2}]+2\}$ and $I_-=I\setminus
(I_+\sqcup I_{\delta})$;
\item[(4.2)] decomposition $I'=I'_+\sqcup
I'_{\delta}\sqcup I'_-$ where $I'_+=\{2,\ldots,\frac{k}{2}+1\}$,
$I'_{\delta}=\{1\}$ and $I'_-=I'\setminus (I'_+\sqcup
I'_{\delta})$;
\item[(4.3)] $I''=\emptyset$.
\end{itemize}
(5) For $\widetilde{A}_{k,k}^{0,0}$, we have:
\begin{itemize}
\item[(5.1)] $\Delta=\delta \Delta'$ and so $I=I'$;
\item[(5.2)] $I''=\emptyset$.
\end{itemize}
(6) For $(\widetilde{D}_{n}^{1,0})^{eq}$, we have:
\begin{itemize}
\item[(6.1)] decomposition $I=I_+\sqcup I_{\delta}\sqcup I_-$
where  $I_+=\{2,\ldots,[\frac{2n-4}{2}]+1\}$, $I_{\delta}=\{1\}$
and $I_-=I\setminus (I_+\sqcup I_{\delta})$;
\item[(6.2)] $I'=I'_{\delta}=\{0,1\}$ and $I'_-=I'_+=\emptyset$;
\item[(6.3)] decomposition $I''=I''_+\sqcup I''_-$ where $I''_+=\{0\}$, $I''_{\delta}=\emptyset$ and $I''_-=I''\setminus I''_+$.
\end{itemize}
(7) For $(\widetilde{D}_{n}^{0,1})^{eq}$, we have:
\begin{itemize}
\item[(7.1)] decomposition $I=I_+\sqcup I_{\delta}\sqcup I_-$
where  $I_+=\{2,\ldots,[\frac{2n-5}{2}]+1\}$,
$I_{\delta}=\{1,[\frac{2n-5}{2}]+2\}$ and $I_-=I\setminus
(I_+\sqcup I_{\delta})$;
\item[(7.2)] $I'=I'_{\delta}=\{0,1\}$ and $I'_-=I'_+=\emptyset$;
\item[(7.3)] decomposition $I''=I''_+\sqcup I''_-$ where $I''_+=\{0\}$, $I''_{\delta}=\emptyset$ and $I''_-=I''\setminus I''_+$.
\end{itemize}
\end{lemma}
\textit{Proof.} We prove \textit{(1), (2), (3), (4), (6)} and
\textit{(7)}.  By [DR, section 6, page 40] and by [DR, section 6,
pages 40 and 46] we note type by type that we have
$|I_{\delta}|=0,1,2$ (respectively $|I'_{\delta}|=0,1,2$ and
$|I''_{\delta}|=0$). Now
\begin{itemize}
\item[i)] if $|I_{\delta}|=0$ we have $e_3=\delta e_0$, $e_2=\delta e_1$
and $e_i=\delta e_{u-i+4}$ for every
$i\in\{4,\ldots,[\frac{u}{2}]+2\}$,
\item[ii)] if $|I_{\delta}|=1$ we have $e_2=\delta e_0$, $e_1=\delta e_1$
and $e_i=\delta e_{u-i+3}$ for every
$i\in\{3,\ldots,[\frac{u}{2}]+1\}$,
\item[iii)] if $|I_{\delta}|=2$ we have $e_2=\delta e_0$, $e_1=\delta
e_1$,  $e_i=\delta e_{u-i+3}$ for every
$i\in\{3,\ldots,[\frac{u}{2}]+1\}$ and $e_{[\frac{u}{2}]+2}=\delta
e_{[\frac{u}{2}]+2}$
\end{itemize}
We define $I_+\subseteq I$ such that
\begin{itemize}
\item[i)] $I_+=\{2,\ldots,[\frac{u}{2}]+2\}\Leftrightarrow
|I_{\delta}|=0$,
\item[ii)] $I_+=\{2,\ldots,[\frac{u}{2}]+1\}\Leftrightarrow
|I_{\delta}|=1$,
\item[iii)] $I_+=\{2,\ldots,[\frac{u}{2}]+1\}\Leftrightarrow
|I_{\delta}|=2$.
\end{itemize}
So respectively decompositions of $I$ of the statement follow. One
proceeds similarly for $I'$ and $I''$.\\
\textit{(5)} follows by the symmetry and considering [DR,
section 6]. $\Box$\\
\\
We note that in part \textit{(5)} of previous lemma we can
consider
$I_{\delta}=I_-=I'_{\delta}=I'_-=I''_{\delta}=I''_-=I''_+=\emptyset$
and so $I_+=I=I'=I'_+$.
\begin{proposizione}\label{dc}
Let $(Q,\sigma)$ be a symmetric quiver of tame type and let $I_+$,
$I_{\delta}$, $I'_+$, $I'_{\delta}$, $I''_+$ and $I''_{\delta}$ be
as above. Any regular symmetric dimension vector can be written
uniquely in the following form:
\begin{equation}\label{fdc}
d=ph+\sum_{i\in I_+}p_i(e_i+\delta e_i)+\sum_{i\in I_{\delta}}p_i
e_i+\sum_{i\in I'_+}p'_i(e'_i+\delta e'_i)+\sum_{i\in
I'_{\delta}}p'_i e'_i+\sum_{i\in I''_+}p''_i(e''_i+\delta e''_i)
\end{equation}
for some non-negative $p,p_i,p'_i,p''_i$ with at least one
coefficient in each family $\{p_i|\;i\in I_+\sqcup I_{\delta}\}$,
$\{p'_i|\;i\in I_+'\sqcup I_{\delta}'\}$, $\{p''_i|\;i\in I''_+\}$
being zero. In particular, in the symplectic case,
\begin{itemize}
\item[i)] if $Q$ has one $\sigma$-fixed vertex and one
$\sigma$-fixed arrow (i.e. $Q=\widetilde{A}_{k,l}^{1,1}$), then
$p_{[\frac{l-1}{2}]+2}$ and $p'_1$ have to be even,
\item[ii)] if $Q$ has one or two $\sigma$-fixed vertices and it
has not any $\sigma$-fixed arrows (i.e.
$Q=\widetilde{A}_{k,l}^{0,2}$ or $\widetilde{D}_{n}^{0,1}$), then
both $p_i$'s and $p'_j$'s, with $i\in I_{\delta}$ and $j\in
I'_{\delta}$, have to be even.
\end{itemize}
\end{proposizione}
\textit{Proof.} It follows by lemma \ref{dI} and by decomposition
of any regular dimension vector of the underlying quiver of
$(Q,\sigma)$. In particular, since symplectic spaces with odd
dimension don't exist, it implies \textit{i)} and
\textit{ii)}. $\Box$\\
\\
Graphically we can represent $\Delta$ (similarly $\Delta'$ and
$\Delta''$) as the polygons
$$
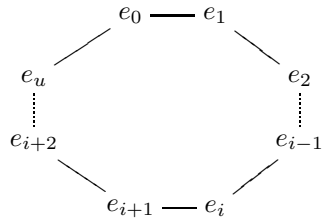
\xymatrix@-1pc{
&e_0\ar@{-}[r]\ar@{-}[dl]&e_1\ar@{-}[dr]&\\
e_u\ar@{.}[d]&&&e_2\ar@{.}[d]\\
e_{i+2}\ar@{-}[dr]&&&\ar@{-}[dl]e_{i-1}\\
&e_{i+1}\ar@{-}[r]&e_i&}
$$
if $Q=\widetilde{A}_{k,k}^{0,0}$ and
$$
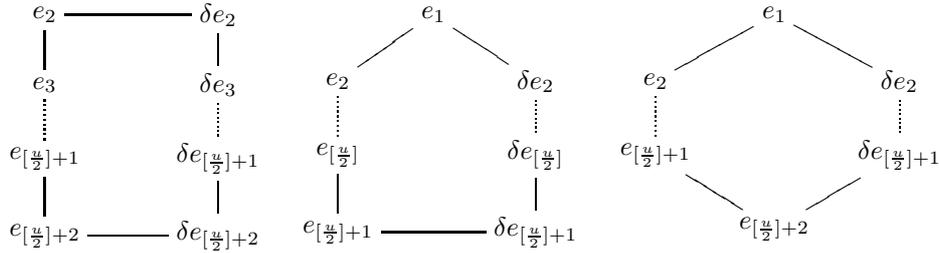
\begin{array}{ccc}
\xymatrix@-1pc{
 e_2\ar@{-}[rr]\ar@{-}[d]&&\delta e_2\ar@{-}[d]\\
 e_3\ar@{.}[d]&&\delta e_3\ar@{.}[d]\\
 e_{[\frac{u}{2}]+1}\ar@{-}[d]&&\delta e_{[\frac{u}{2}]+1}\ar@{-}[d]\\
  e_{[\frac{u}{2}]+2}\ar@{-}[rr]&&\delta e_{[\frac{u}{2}]+2}}&
\xymatrix@-1pc{
&e_1\ar@{-}[dr]\ar@{-}[dl]&\\
 e_2\ar@{.}[d]&&\delta e_2\ar@{.}[d]\\
 e_{[\frac{u}{2}]}\ar@{-}[d]&&\delta e_{[\frac{u}{2}]}\ar@{-}[d]\\
  e_{[\frac{u}{2}]+1}\ar@{-}[rr]&&\delta e_{[\frac{u}{2}]+1}}&
\xymatrix@-1pc{
&e_1\ar@{-}[dr]\ar@{-}[dl]&\\
 e_2\ar@{.}[d]&&\delta e_2\ar@{.}[d]\\
e_{[\frac{u}{2}]+1}\ar@{-}[dr]&&\delta e_{[\frac{u}{2}]+1}\ar@{-}[dl]\\
&e_{[\frac{u}{2}]+2}&}
 \end{array}
$$
with a reflection respect to a central vertical line, in the other
cases.
\begin{definizione}
We define an involution $\sigma_I$ on the set of indices $I$ such
that $e_{\sigma_I(i)}=\delta e_i$ for every $i\in I$. Hence
$\sigma_I(I)=I'$ for $\widetilde{A}_{k,k}^{0,0}$ and
$\sigma_{I}I_+=I_-$, $\sigma_I I_{\delta}=I_{\delta}$ for the
other cases. Similarly we define an involution $\sigma_{I'}$ and
an involution $\sigma_{I''}$ respectively on $I'$ and on $I''$.
\end{definizione}
\begin{lemma}\label{regsport}
(1) For $\widetilde{A}^{2,0,1}_{k,l}$, no one indecomposable
regular representation is orthogonal. The following indecomposable
regular representations are symplectic
\begin{itemize}
\item[(1.1)] $E_{i,\sigma_I(i)}$ such that $\sum_{k=i}^{\sigma_I(i)} e_k\neq h$ and $E_{i,\sigma_I(i)}$
 of dimension $h$ containing
$E_{\frac{l}{2}+1}$.
\item[(1.2)] $E'_{i,\sigma_I'(i)}$ such that $\sum_{k=i}^{\sigma_{I'}(i)} e'_k\neq h$ and $E'_{i,\sigma_I'(i)}$
 of dimension $h$ containing
$E'_{\frac{k}{2}+1}$.
\end{itemize}
(2) For $\widetilde{A}^{2,0,2}_{k,l}$, no one indecomposable
regular representation is orthogonal. The following indecomposable
regular representations are symplectic
\begin{itemize}
\item[(2.1)] $E_{i,\sigma_I(i)}$ such that $\sum_{k=i}^{\sigma_I(i)} e_k\neq h$, $E_{i,\sigma_I(i)}$
 of dimension $h$ containing
$E_{0}$ and $E_{i,\sigma_I(i)}$
 of dimension $h$ containing
$E_{[\frac{l+1}{2}]+1}$.
\item[(2.2)] $E'_{i,\sigma_{I'}(i)}$ such that $\sum_{k=i}^{\sigma_{I'}(i)} e'_k\neq h$.
\end{itemize}
(3) For $\widetilde{A}^{0,2}_{k,l}$, no one indecomposable regular
representations is symplectic. The following indecomposable
regular representations are orthogonal
\begin{itemize}
\item[(3.1)] $E_{i,\sigma_{I}(i)}$ such that $\sum_{k=i}^{\sigma_I(i)} e_k\neq h$.
\item[(3.2)] $E'_{i,\sigma_{I'}(i)}$ such that $\sum_{k=i}^{\sigma_{I'}(i)} e'_k\neq h$.
\end{itemize}
(4) For $\widetilde{A}^{0,2}_{k,l}$, the following indecomposable
regular representations are orthogonal
\begin{itemize}
\item[(4.1.1)] $E_{i,\sigma_{I}(i)}$, with $i\leq\sigma_I(i)$, such that $\sum_{k=i}^{\sigma_I(i)} e_k\neq h$.
\item[(4.1.2)] $E'_{i,\sigma_{I'}(i)}$, with $i\geq\sigma_{I'}(i)$, such that $\sum_{k=i}^{\sigma_{I'}(i)} e'_k\neq h$.
\end{itemize}
The following indecomposable regular representations are
symplectic
\begin{itemize}
\item[(4.2.1)] $E_{i,\sigma_{I}(i)}$, with $i\geq\sigma_I(i)$, such that $\sum_{k=i}^{\sigma_I(i)} e_k\neq h$.
\item[(4.2.2)] $E'_{i,\sigma_{I'}(i)}$, with $i\leq\sigma_{I'}(i)$, such that $\sum_{k=i}^{\sigma_{I'}(i)} e'_k\neq
h$ and $E'_{i,\sigma_{I'}(i)}$, with $i\leq\sigma_{I'}(i)$, of
dimension $h$ containing $E'_{\frac{k}{2}+1}$.
\end{itemize}
(5) For $\widetilde{A}^{0,0}_{k,k}$, no one indecomposable regular
representation is symplectic
or orthogonal.\\
\\
(6) For $(\widetilde{D}^{1,0}_{n})^{eq}$, no one indecomposable
regular representation is othogonal. The following indecomposable
regular representations are symplectic
\begin{itemize}
\item[(6.1)] $E_{i,\sigma_I(i)}$ such that $\sum_{k=i}^{\sigma_I(i)} e_k\neq h$ and $E_{i,\sigma_I(i)}$
 of dimension $h$ containing
$E_{n-1}$.
\item[(6.2)] $E'_0$ and $E'_1$.
\item[(6.3)] $E''_{0,1}$ and $E''_{1,0}$.
\end{itemize}
(7) For $(\widetilde{D}^{0,1}_{n})^{eq}$, no one indecomposable
regular representation is symplectic. The orthogonal
indecomposable regular representations are
\begin{itemize}
\item[(7.1)] $E_{i,\sigma_I(i)}$ such that $\sum_{k=i}^{\sigma_I(i)} e_k\neq h$.
\item[(7.2)] $E'_0$ and $E'_1$.
\item[(7.3)] $E''_{0,1}$ and $E''_{1,0}$.
\end{itemize}
\end{lemma}
\textit{Proof.} We check only part (1.1), similarly one proves the
other parts. Let $Q=\widetilde{A}^{2,0,1}_{k,l}$. The only
$E_{i,j}$ such that $\delta
\underline{dim}E_{i,j}=\underline{dim}E_{i,j}$ are
$E_{i,\sigma_I(i)}$. We have three cases.
\begin{itemize}
\item[(i)] If $\sum_{k=i}^{\sigma_I(i)} e_k\neq h$ and $i<
\sigma_I(i)$ then we have for $j\in Q_0$
$$
E_{i,\sigma_I(i)}(j)=\left\{\begin{array}{ll}\mathbb{K} &
j=x_s,\sigma(x_s)\;\textrm{with}\;i-1\leq s\leq \frac{l}{2}\\0 &
\textrm{otherwise}\end{array}\right.
$$
and for $c\in Q_1$
$$
E_{i,\sigma_I(i)}(c)=\left\{\begin{array}{ll}Id &
c=v_s,\sigma(v_s),a\;\textrm{with}\;i\leq s\leq \frac{l}{2}\\0 &
\textrm{otherwise}.\end{array}\right.
$$
So we note that we can define on such $E_{i,\sigma_I(i)}$ a
symplectic structure.
\item[(ii)] If $\sum_{k=i}^{\sigma_I(i)} e_k\neq h$ and $i\geq
\sigma_I(i)$ then we have for $j\in Q_0$
$$
E_{i,\sigma_I(i)}(j)=\left\{\begin{array}{ll}0 &
j=x_s,\sigma(x_s)\;\textrm{with}\;i\leq s\leq
\frac{l}{2}\\\mathbb{K} & \textrm{otherwise}\end{array}\right.
$$
and for $c\in Q_1$
$$
E_{i,\sigma_I(i)}(c)=\left\{\begin{array}{ll}0 &
c=v_s,\sigma(v_s),a\;\textrm{with}\;i\leq s\leq \frac{l}{2}\\Id &
\textrm{otherwise}.\end{array}\right.
$$
So we note that we can define on such $E_{i,\sigma_I(i)}$ a
symplectic structure.
\item[(iii)] If $\sum_{k=i}^{\sigma_I(i)} e_k= h$ and
$E_{i,\sigma_I(i)}$ contains $E_{\frac{l}{2}+1}$, then, by AR
quiver of $Q$, we note the following almost split sequence
$$
0\longrightarrow
E_{\frac{l}{2}+1,\sigma_I(\frac{l}{2})}\longrightarrow
E_{i,\sigma_I(i)}\oplus
E_{\frac{l}{2},\sigma_I(\frac{l}{2})}\longrightarrow
E_{\frac{l}{2},\sigma(\frac{l}{2})-1}\longrightarrow 0.
$$
So we have for every $j\in Q_0$, $E_{i,\sigma_I(i)}(j)=\mathbb{K}$
and for $c\in Q_1$
$$
E_{i,\sigma_I(i)}(c)=\left\{\begin{array}{ll}0 & c=a\\Id &
\textrm{otherwise}.\end{array}\right.
$$
Finally, we note that we can define on such $E_{i,\sigma_I(i)}$ a
symplectic structure. $\Box$
\end{itemize}
In the remainder of the section, we shall call
\begin{equation}\label{d'}
d'=\sum_{i\in I_+}p_i(e_i+\delta e_i)+\sum_{i\in
I_{\delta}}p_i e_i+\sum_{i\in I'_+}p'_i(e'_i+\delta
e'_i)+\sum_{i\in I'_{\delta}}p'_i e'_i+\sum_{i\in
I''_+}p''_i(e''_i+\delta e''_i).
\end{equation}
\begin{proposizione}\label{oa}
If $d$ is regular with decomposition (\ref{fdc}) such that $d= d'$
or $d$ is not regular then $SpRep(Q, d)$ (respectively
$ORep(Q,d)$) has an open $Sp(Q, d)$-orbit (respectively
$O(Q,d)$-orbit).
\end{proposizione}
\textit{Proof.} If $d=d'$, we have no indecomposable of dimension
vector $ph$ and so there are finitely many orbits.
If $d$ is not regular, it follows from [R2, theorem 3.2]. $\Box$\\
\\
In the next $d$ shall be a regular symmetric dimension vector with
decomposition (\ref{fdc}) with $p\geq 1$ and $p\neq 0$. Now we
shall describe the generators of $SpSI(Q,d)$ and $OSI(Q,d)$. To do
this the following theorem, which we prove later, is useful.
\begin{teorema}\label{dtso} Let $(Q,\sigma)$ be a symmetric quiver of
tame type and the decomposition (\ref{fdc}) of a regular symmetric
dimension vector with $p\geq 1$ and $d'\neq 0$. There exist
isomorphisms of algebras
\begin{equation}SpSI(Q,d)\stackrel{\Phi_d}{\rightarrow}\bigoplus_{\chi\in
char(Sp(Q,d))}SpSI(Q,ph)_{\chi}\otimes
SpSI(Q,d')_{\chi'}\end{equation} and
\begin{equation}OSI(Q,d)\stackrel{\Psi_d}{\rightarrow}\bigoplus_{\chi\in
char(O(Q,d))}OSI(Q,ph)_{\chi}\otimes
OSI(Q,d')_{\chi'},\end{equation}
where $\chi'=\chi|_{d'}$, i.e.
the restriction of the weight $\chi$ to the support of $d'$.
\end{teorema}
By proposition \ref{oa} $Sp(Q,d')$ (respectively $O(Q,d')$) acting
on $SpRep(Q,d')$ (respectively on $ORep(Q,d')$) has an open orbit
so , by lemma \ref{sk}, dimension of $SpSI(Q,d')_{\chi'}$
(respectively dimension of $OSI(Q,d')_{\chi'}$) is 0 or 1. This
allows us to identify one non-zero element of $SpSI(Q,d)_{\chi}$
(respectively of $OSI(Q,d)_{\chi}$) with the element of
$SpSI(Q,ph)_{\chi}$ (respectively of $OSI(Q,ph)_{\chi}$) to which
it restricts.\\
We proceed now to describe the generators of the algebra $SpSI(Q,
d)$ (respectively $OSI(Q,d)$). If the corresponding $I,\,I',\,I''$
are not empty, we label the vertices $e_i,\, e'_i,\, e''_i$ of the
polygons $\Delta,\, \Delta',\, \Delta''$ with the coefficients
$p_i,\, p'_i,\, p''_i$. We recall that
\begin{itemize}
\item[a)] we have to label with $p_i$ (respectively with $p'_i$ and $p''_i$)
both vertices $e_i$ and $\delta e_i$, i.e $p_i=p_{\sigma_I(i)}$
(respectively $p'_i=p'_{\sigma_I''(i)}$ and
$p''_i=p''_{\sigma_I''(i)}$), if $e_i\neq\delta e_i$.
\end{itemize}
and in the symplectic case, by \textit{i)} and \textit{ii)} of
proposition \ref{dc}
\begin{itemize}
\item[b)] for $\widetilde{A}_{k,l}^{1,1}$, $p_{[\frac{u}{2}]+2}$ and $p'_1$
have to be even,
\item[c)] for $\widetilde{A}_{k,l}^{0,2}$ and
$\widetilde{D}_{n}^{0,1}$,  $p_i\in I_{\delta}$ and $p'_i\in
I'_{\delta}$ have to be even.
\end{itemize}
We shall call these labelled polygons respectively $\Delta(d),\,
\Delta'(d),\, \Delta''(d)$.
\begin{definizione}
We shall say that the labelled arc
$\xymatrix@-1pc{p_i\ar@{-}[r]&\ar@{.}[r]&\ar@{-}[r]&p_j}$ (in
clockwise orientation) of the labelled polygon $\Delta(d)$ is
admissible if $p_i=p_j$ and $p_i<p_k$ for every its interior
labels $p_k$. We denote such a labelled arc
$\xymatrix@-1pc{p_i\ar@{-}[r]&\ar@{.}[r]&\ar@{-}[r]&p_j}$ by
$[i,j]$, and we define $p_i=p_j$ the index $ind[i,j]$ of $[i,j]$.
Similarly we define admissible arcs and their indexes for the
labelled polygons $\Delta'(d)$ and $\Delta''(d)$.
\end{definizione}
We denote by $\mathcal{A}(d)$, $\mathcal{A}'(d)$,
$\mathcal{A}''(d)$ the sets of all admissible labelled arcs in the
polygons $\Delta(d)$, $\Delta'(d)$, $\Delta''(d)$ respectively. In
particular we note that if $d=ph$, then the polygons $\Delta(d)$,
$\Delta'(d)$, $\Delta''(d)$ are labelled by zeros and so
$\mathcal{A}(d)$, $\mathcal{A}'(d)$, $\mathcal{A}''(d)$ consist of
all edges of respective polygons. With these notations we have the
following
\begin{proposizione}\label{fij}
For each arc $[i,j]$ from $\mathcal{A}(d)$ (respectively
$\mathcal{A}'(d)$ and $\mathcal{A}''(d)$) there exists in
$SpSI(Q,d)$ and in $OSI(Q,d)$ a non zero semi-invariant
\begin{itemize}
\item[(i)] of type $c^{E_{i,j-1}}$ (respectively $c^{E'_{i,j-1}}$
and $c^{E''_{i,j-1}}$) or of type $c^{V_{(\varphi,\psi)}}$, with
$(\varphi,\psi)\in\{(1,0),(0,1),(1,1)\}$;
\item[(ii)] of type $pf^{E_{i,j-1}}$ (respectively $pf^{E'_{i,j-1}}$
and $pf^{E''_{i,j-1}}$) or of type $pf^{V_{(\varphi,\psi)}}$, with
$(\varphi,\psi)\in\{(1,0),(0,1),(1,1)\}$, if $E_{i,j-1}$,
$E'_{i,j-1}$, $E''_{i,j-1}$ and $V_{(\varphi,\psi)}$ satisfy
property \textit{(Op)} in the symplectic case and property
\textit{(Spp)} in the orthogonal case.
\end{itemize}
\end{proposizione}
Let $c_0,\ldots,c_t$, with $t=\frac{p-1}{2}$, $\frac{p}{2}$ and
$p$, defined case by case in section 3.1. The generators of
algebras $SpSI(Q,d)$ and $OSI(Q,d)$ are described by the following
theorem
\begin{teorema}\label{stq}
Let $(Q,d)$ a symmetric quiver of tame type and $d=ph+d'$ the
decomposition of a regular symmetric dimension vector $d$ with
$p\geq 1$. Then $SpSI(Q,d)$ (respectively $OSI(Q,d)$) is generated
by
\begin{itemize}
\item[(i)] $c_0,\ldots,c_t$;
\item[(ii)] $c^{E_{i,j-1}}$, $c^{E'_{r,s-1}}$, $c^{E''_{t,m-1}}$ and $c^{V_{(\varphi,\psi)}}$ with $[i,j]\in\mathcal{A}(d)$,
$[r,s]\in\mathcal{A}'(d)$, $[t,m]\in\mathcal{A}''(d)$ and
$(\varphi,\psi)\in\{(1,0),(0,1),(1,1)\}$;
\item[(iii)] $pf^{E_{i,j-1}}$, $pf^{E'_{r,s-1}}$, $pf^{E''_{t,m-1}}$ and $pf^{V_{(\varphi,\psi)}}$ with $[i,j]\in\mathcal{A}(d)$,
$[r,s]\in\mathcal{A}'(d)$, $[t,m]\in\mathcal{A}''(d)$ and
$(\varphi,\psi)\in\{(1,0),(0,1),(1,1)\}$, if $E_{i,j-1}$,
$E'_{i,j-1}$, $E''_{i,j-1}$ and $V_{(\varphi,\psi)}$ satisfy
property \textit{(Op)} (respectively property \textit{(Spp)}).
\end{itemize}
\end{teorema}
First we note that $\langle h,d\rangle=0$ and further we have the
following
\begin{lemma}
For every regular dimension vector $d$
$$
\langle\underline{dim} E_{i,j-1},d\rangle=0\Leftrightarrow
p_i=p_j.
$$
\end{lemma}
\textit{Proof.} See [D, section 4.3]. $\Box$\\
\\
So theorem \ref{stq} is equivalent to conjectures \ref{mt1} and
\ref{mt2}.

\subsection{Proof of theorem \ref{stq} and \ref{dtso}}
In this section we prove the theorem \ref{stq} and theorem
\ref{dtso}. For theorem \ref{stq}, by proposition \ref{oA},
proposition \ref{Q=Q'} and lemma \ref{kac}, we can reduce the
proof to the orientation of $\widetilde{A}$ as in proposition
\ref{oA} and to the equiorientation for $\widetilde{D}$. In the
proof we use the notion of generic decomposition of the symmetric
dimension vector $d$ (see [K1], [K2], [KR]).
\begin{definizione}
A decomposition $\alpha=\beta_1\oplus\cdots\oplus\beta_q$ of a
dimension vector $\alpha$ is called generic if there is a Zariski
open subset $\mathcal{U}$ of $Rep(Q,\alpha)$ such that each
$U\in\mathcal{U}$ decomposes in $U=\bigoplus_{i=1}^q U_i$ with
$U_i$ indecomposable representation of dimension $\beta_i$, for
every $i\in\{1,\ldots,q\}$.
\end{definizione}
\begin{definizione}
\begin{itemize}
\item[(1)] A decomposition $\alpha=\beta_1\oplus \cdots\oplus\beta_q$ of a symmetric
dimension vector $\alpha$ is called symplectic generic if there is
a Zariski open subset $\mathcal{U}$ of $SpRep(Q,\alpha)$ such that
each $U\in\mathcal{U}$ decomposes in $U=\bigoplus_{i=1}^q U_i$
with $U_i$ indecomposable symplectic representation of dimension
$\beta_i$, for every $i\in\{1,\ldots,q\}$.
\item[(2)] A decomposition $\alpha=\beta_1\oplus\cdots\oplus\beta_q$ of a symmetric
dimension vector $\alpha$ is called orthogonal generic if there is
a Zariski open subset $\mathcal{U}$ of $ORep(Q,\alpha)$ such that
each $U\in\mathcal{U}$ decomposes in $U=\bigoplus_{i=1}^q U_i$
with $U_i$ indecomposable orthogonal representation of dimension
$\beta_i$, for every $i\in\{1,\ldots,q\}$.
\end{itemize}
\end{definizione}
For tame quivers the generic decomposition of any regular
dimension vector is given by
results of [R2, section 3].\\
We describe this decomposition explicitly for a symmetric regular
dimension vector $d$ with decomposition (\ref{fdc}).\\
In the remainder of this section we set
\begin{equation}\label{d'1}
\bar{d}=\sum_{i\in I_+}p_i(e_i+\delta e_i)+\sum_{i\in
I_{\delta}}p_i e_i,
\end{equation}
\begin{equation}
\bar{d}'=\sum_{i\in I'_+}p'_i(e'_i+\delta e'_i)+\sum_{i\in
I'_{\delta}}p'_i e'_i
\end{equation}
\begin{equation}
\bar{d}''=\sum_{i\in I''_+}p''_i(e''_i+\delta e''_i).
\end{equation}
\begin{oss}
\begin{itemize}\label{p_i=0}
\item[(i)] We remember that at least one coefficient in each family $\{p_i|\;i\in
I_+\sqcup I_{\delta}\}$, $\{p'_i|\;i\in I_+'\sqcup I_{\delta}'\}$,
$\{p''_i|\;i\in I''_+\}$ is zero.
\item[(ii)] We can assume $p_i=0$ for $i\in I_{\delta}$ or $p_{i}=0$, for $i\in I_+$, and so
$p_{\sigma_I(i)}=0$.
\end{itemize}
\end{oss}
\begin{definizione}
We divide the polygon $\Delta(\bar{d})$ in two parts:
\begin{itemize}
\item[(i)] the up part $\Delta_{up}(\bar{d})$ is the part of $\Delta(\bar{d})$ from $p_{i-1}$ to
$p_{\sigma_I(i-1)}$;
\item[(ii)] the down part $\Delta_{down}(\bar{d})$ is the part of $\Delta(\bar{d})$ from $p_{i+1}$ to
$p_{\sigma_I(i+1)}$.
\end{itemize}
Similarly for $\Delta'$ and $\Delta''$.
\end{definizione}
\begin{oss}
We note that if $p_i=0$ with $i\in I_{\delta}$, then we have only
the part $\Delta_{up}$ or the part $\Delta_{down}$.
\end{oss}
We consider $\Delta$, similarly one proceeds for $\Delta'$ and
$\Delta''$.
\begin{definizione}
We shall call symmetric arc, an arc invariant under $\sigma_I$,
i.e. an arc of type $[i,\sigma_I(i)]$.
\end{definizione}
\begin{oss}
By the division of $\Delta$ in $\Delta_{up}$ and $\Delta_{down}$,
we note that all symmetric arcs pass through the same
$\sigma_I$-fixed vertex of $\Delta$ or through the same
$\sigma_I$-fixed edge of $\Delta$.
\end{oss}
\begin{lemma}
Let $(Q,\sigma)$ be a symmetric quiver of tame type.
\begin{itemize}
\item[(i)] If $n=\sigma_I(n)$ then either there exists unique $x\in
Q_0^{\sigma}$ such that $e_n(x)\neq 0$ or there exists unique
$a\in Q_1^{\sigma}$ such that $e_n(ta)\neq 0$.
\item[(ii)] If $\xymatrix@-1pc{n\ar@{-}[r]&\sigma_I(n)}$ is
a $\sigma_I$-fixed edge in $\Delta$, then there exists unique
$a\in Q_1^{\sigma}$ such that $e_n(ta)\neq  0$.
\end{itemize}
\end{lemma}
\textit{Proof.} One proceeds type by type. We consider
$Q=\widetilde{A}^{2,0,1}_{k,l}$ since
for the other types one proves similarly.\\
(i) By lemma \ref{dI}, the only $\sigma_I$-fixed vertex of
$\Delta$ is $1$ and $b$ is the unique arrow in $Q_1^{\sigma}$ such that $e_1(tb)\neq 0$.\\
(ii) The only $\sigma_I$-fixed edge of $\Delta$ is
$\xymatrix@-1pc{{\frac{l}{2}+1}\ar@{-}[r]&{\sigma_I(\frac{l}{2}+1)}}$
and $a$ is the unique arrow in $Q_1^{\sigma}$ such that
$e_{\frac{l}{2}+1}(ta)\neq 0$. $\Box$
\begin{definizione}
\begin{itemize}
\item[(i)] If $n=\sigma_I(n)$, we call $x(n)$ the unique $x\in
Q_0^{\sigma}$ such that $e_n(x)\neq 0$.
\item[(ii)] If $n=\sigma_I(n)$ or $\xymatrix@-1pc{n\ar@{-}[r]&\sigma_I(n)}$ is
a $\sigma_I$-fixed edge in $\Delta$, we call $a(n)$ the unique
$a\in Q_1^{\sigma}$ such that $e_n(ta)\neq  0$.
\end{itemize}
\end{definizione}
\begin{definizione}
For every arc $[i,j]$ in $\Delta$, we define
$$
e_{[i,j]}=\sum_{k\in [i,j]}e_k.
$$
\end{definizione}
\begin{definizione}
\begin{itemize}
\item[(i)] $\mathcal{A}_+(\bar{d}):=\{[i,j]\in\mathcal{A}(\bar{d})|\;[i,j]\subset I_+\}$
\item[(ii)] $\mathcal{A}^k_+(\bar{d}):=\{[i,j]\in\mathcal{A}(\bar{d})|\;[i,j]\subset I_+,ind[i,j]=k\}$.
\item[(iii)] $\mathcal{A}^k_{\sigma_I}(\bar{d})=\{[i,j]=\sigma_I[i,j]\in\mathcal{A}(\bar{d})|\;ind[i,j]=k\}$.
\end{itemize}
\end{definizione}
\begin{oss}
$[i,j]\subset I_+$ if and only if
$[\sigma_I(j),\sigma_I(i)]\subset I_{-}$ and
$ind[i,j]=ind[\sigma_I(j),\sigma_I(i)]$.
\end{oss}
First we consider all the admissible arcs in
$\mathcal{A}^r_{\sigma_I}(\bar{d})\cup \mathcal{A}^r_{+}(\bar{d})$
such that $r=max\{p_k\}$. So we get
$$
=\sum_{i\in I_+}\tilde{p}_i(\tilde{e}_i+\delta
\tilde{e}_i)+\sum_{i\in I_{\delta}}\tilde{p}_i\tilde{e}_i=
$$
\begin{equation}\label{barp}
\sum_{i\in I_+}p_i(e_i+\delta e_i)+\sum_{i\in
I_{\delta}}p_ie_i-\left(\bigoplus_{[i,j]\in
\mathcal{A}^r_{+}(\bar{d})}(e_{[i,j]}+\delta
e_{[i,j]})+\bigoplus_{[i,\sigma_I(i)]\in
\mathcal{A}^r_{\sigma_I}(\bar{d})}e_{[i,\sigma_I(i)]}\right),
\end{equation}
where $max(\bar{p}_i)=r-1$. Then we repeat the procedure for
(\ref{barp}) and so on we have
$$
\sum_{i\in I_+}p_i(e_i+\delta e_i)+\sum_{i\in I_{\delta}}p_ie_i=
$$
\begin{equation}\label{vj,k}
\bigoplus_{k=1}^{r}\left(\bigoplus_{[i,j]\in
\mathcal{A}^k_{+}(\bar{d})}(e_{[i,j]}+\delta
e_{[i,j]})+\bigoplus_{[i,\sigma_I(i)]\in
\mathcal{A}^k_{\sigma_I}(\bar{d})}e_{[i,\sigma_I(i)]}\right).
\end{equation}
\begin{oss}
\begin{itemize}
\item[(i)] If $[i,j]$ and $[i',j']$ are two admissible arcs in $\mathcal{A}(\bar{d})$ such that $[i,j]\supseteq [i',j']$, then
$ind[i,j]\leq ind[i',j']$.
\item[(ii)] If there not exists $[i,j]\in\mathcal{A}^h_{\sigma_I}(\bar{d})\cup
\mathcal{A}^h_{+}(\bar{d})$ such that $[i,j]\supseteq [i',j']$ for
some $[i',j']\in\mathcal{A}^k_{\sigma_I}(\bar{d})\cup
\mathcal{A}^k_{+}(\bar{d})$, then the symmetric dimension vector
corresponding to $[i',j']$ appears $k$-times in the decomposition
(\ref{vj,k}), with $1\leq h< k$.
\end{itemize}
\end{oss}
\begin{definizione}
Let $[i_1,j_1],\ldots,[i_k,j_k]$ be the admissible arcs such that
$[i_1,j_1]\supseteq\cdots\supseteq[i_k,j_k]$, with $k\geq 1$. We
define $q_{[i_h,j_h]}=ind[i_h,j_h]-ind[i_{h-1},j_{h-1}]$ for every
$1\leq h\leq k$, where $ind[i_0,i_0]=0$.
\end{definizione}
We note that for every $[i,j]\in
\mathcal{A}^k_{\sigma_I}(\bar{d})\cup \mathcal{A}^k_{+}(\bar{d})$,
$q_{[i,j]}$ is the multiplicity of the symmetric dimension vector
corresponding to $[i,j]$ in the
decomposition (\ref{vj,k}).\\
Finally we have
$$
\sum_{i\in I_+}p_i(e_i+\delta e_i)+\sum_{i\in I_{\delta}}p_ie_i=
$$
\begin{equation}\label{qi,j}
\bigoplus_{[i,j]\in\mathcal{A}_+(\bar{d}) }(e_{[i,j]}+\delta
e_{[i,j]})^{\oplus q_{[i,j]}} +
\bigoplus_{[i,\sigma_I(i)]\in\mathcal{A}(\bar{d})
}(e_{[i,\sigma_I(i)]})^{\oplus q_{[i,\sigma_I(i)]}}.
\end{equation}
\begin{esempi}\label{esempioe}
If $\Delta$ is of the form
\begin{equation}\label{esempio}
\xymatrix{&e_1=\delta e_1\ar@{-}[dl]\ar@{-}[dr]&\\
 e_2\ar@{-}[d]&&\delta e_2=e_{\sigma_I(2)}\ar@{-}[d]\\
 e_3\ar@{-}[dr]&&\delta e_3=e_{\sigma_I(3)}\ar@{-}[dl]\\
&e_4=\delta e_4&}
\end{equation}
and $p_1=2$, $p_2=3$, $p_3=0$ and $p_4=2$, then
$[2,\sigma_I(2)]=\{2,1,\sigma_I(2)\}\subset I_+\sqcup
I_{\delta}\sqcup I_-$ with
$q_{[2,\sigma_I(2)]}=ind[2,\sigma_I(2)]=2$, $[2,2]=\{2\}\in I_+$
with $q_{[2,2]}=ind[2,2]-ind[2,\sigma_I(2)]=1$ and $[4,4]=\{4\}\in
I_{\delta}$ with $q_{[4,4]}=ind[4,4]=2$. So we have
$$
\sum_{i\in I_+}p_i(e_i+\delta e_i)+\sum_{i\in I_{\delta}}p_ie_i=
((e_2+\delta e_2)+e_1)^{\oplus 2}\oplus(e_2+\delta e_2)\oplus
(e_4)^{\oplus 2}.
$$
\end{esempi}
Similarly we proceed with the decomposition of $\bar{d}'$ and
$\bar{d}''$. So we have the following
\begin{proposizione}
Let $(Q,\sigma)$ be a symmetric quiver of tame type and let $d$ be
a symmetric dimension vector of a representation of the underlying
quiver $Q$ with decomposition (\ref{fdc}). Then
$$
d=\bigoplus_{i=1}^p h+\bigoplus_{[i,j]\in\mathcal{A}_+(\bar{d})
}(e_{[i,j]}+\delta e_{[i,j]})^{\oplus q_{[i,j]}} +
\bigoplus_{[i,\sigma_I(i)]\in\mathcal{A}(\bar{d})}(e_{[i,\sigma_I(i)]})^{\oplus
q_{[i,\sigma_I(i)]}}+
$$
$$
\bigoplus_{[i,j]\in\mathcal{A}'_+(\bar{d}')}(e'_{[i,j]}+\delta
e'_{[i,j]})^{\oplus q'_{[i,j]}}
+\bigoplus_{[i,\sigma_{I'}(i)]\in\mathcal{A}'(\bar{d}')
}(e'_{[i,\sigma_{I'}(i)]})^{\oplus q'_{[i,\sigma_{I'}(i)]}}+
$$
\begin{equation}\label{dgvs}
\bigoplus_{[i,j]\in\mathcal{A}''_+(\bar{d}'')}(e''_{[i,j]}+\delta
e''_{[i,j]})^{\oplus q''_{[i,j]}}
+\bigoplus_{[i,\sigma_{I''}(i)]\in\mathcal{A}''(\bar{d}'')
}(e''_{[i,\sigma_{I''}(i)]})^{\oplus q''_{[i,\sigma_{I''}(i)]}}
\end{equation}
is the generic decomposition of $d$.
\end{proposizione}
We restrict to dimension vectors of regular symplectic
representations and of regular orthogonal representations. We
modify generic decomposition (\ref{dgvs}) of $d=(d_i)_{i\in Q_0}$
to get symplectic generic decomposition of $d$ or orthogonal
generic decomposition of $d$.\\
Let $[i,j]$ be an arc in $\Delta_{up}$ and let $[h,k]$ be an arc
in $\Delta_{down}$. If $E_{[i,j]}$ is the regular indecomposable
symplectic (respectively orthogonal) representation of
$(Q,\sigma)$ corresponding to $[i,j]$ and $E_{[h,k]}$ is the
regular indecomposable symplectic (respectively orthogonal)
representation of $(Q,\sigma)$ corresponding to $[h,k]$, then
$$
Hom_Q(E_{[i,j]},E_{[h,k]})=0=Hom_Q(E_{[h,k]},E_{[i,j]})
$$
and
$$
Ext_Q^1(E_{[i,j]},E_{[h,k]})=0=Ext_Q^1(E_{[h,k]},E_{[i,j]}).
$$
So we deal separately with $\Delta_{up}$ and $\Delta_{down}$. We
consider $I=I^{up}\sqcup I^{down}$, $I_+=I_+^{up}\sqcup
I_+^{down}$ and $I_{\delta}=I_{\delta}^{up}\sqcup
I_{\delta}^{down}$. We have the decomposition
$\bar{d}=\bar{d}_{up}+\bar{d}_{down}$, where
\begin{equation}\label{d'1up}
\bar{d}_{up}=\sum_{i\in I^{up}_+}p_i(e_i+\delta e_i)+\sum_{i\in
I^{up}_{\delta}}p_i e_i
\end{equation}
and
\begin{equation}
\bar{d}_{down}=\sum_{i\in I^{down}_+}p_i(e_i+\delta
e_i)+\sum_{i\in I^{down}_{\delta}}p_i e_i.
\end{equation}
By what has be said, the symplectic (respectively orthogonal)
generic decomposition of $\bar{d}$ is direct sum of the symplectic
(respectively orthogonal) generic decomposition of $\bar{d}_{up}$
and the symplectic (respectively orthogonal) generic decomposition
of $\bar{d}_{down}$.
\begin{oss}
\begin{itemize}
\item[(i)] In the symplectic case, since $\bar{d}_{x}$ has to
be even for every $x\in Q_0^{\sigma}$, we have to modify the
symmetric dimension vectors corresponding to the arcs passing
through the $\sigma_I$-fixed vertex $n$ such that there exists
$x=x(n)\in Q_0^{\sigma}$ such that $e_n(x)\neq 0$.
\item[(ii)] In the orthogonal case, we have to modify the symmetric
dimension vectors corresponding to the arcs passing through the
$\sigma_I$-fixed vertex $n$ such that $\bar{d}_{ta(n)}$ is even
and those corresponding to the arcs passing through the
$\sigma_I$-fixed edge $\xymatrix@-1pc{n\ar@{-}[r]&\sigma_I(n)}$
such that $\bar{d}_{ta(n)}$ is even.
\item[(iii)] We have to modify also $ph+e_{[i,\sigma_I(i)]}$, with
$p$ odd, if $[i,\sigma_I(i)]$ is like in part (i) (respectively
part (ii)), since $h+e_{[i,\sigma_I(i)]}$ is the dimension vector
of regular indecomposable symplectic (respectively orthogonal)
representation.
\end{itemize}
\end{oss}
\begin{definizione}
\begin{itemize}
\item[(i)] $\mathcal{A}^{up}(\bar{d})=\{[i,j]\in
\mathcal{A}(\bar{d})|\;[i,j]\subset I^{up}\}$.
\item[(ii)] $\mathcal{A}_+^{up}(\bar{d})=\{[i,j]\in
\mathcal{A}(\bar{d})|\;[i,j]\subset I_+^{up}\}$.
\item[(iii)] $\mathcal{A}^{down}(\bar{d})=\{[i,j]\in
\mathcal{A}(\bar{d})|\;[i,j]\subset I^{down}\}$.
\item[(iv)] $\mathcal{A}_+^{down}(\bar{d})=\{[i,j]\in
\mathcal{A}(\bar{d})|\;[i,j]\subset I_+^{down}\}$.
\end{itemize}
\end{definizione}
Let $\bar{d}=\bar{d}_{up}+\bar{d}_{down}$ be a regular symplectic
dimension vector. We consider $\Delta_{up}$. $\Delta_{up}$
contains either a $\sigma_I$-fixed vertex $n_{up}$ or a
$\sigma_I$-fixed edge
$\xymatrix@-1pc{n_{up}\ar@{-}[r]&\sigma_I(n_{up})}$. Starting from
generic decomposition (\ref{dgvs}) of $\bar{d}_{up}$ we modify it
as follows.
\begin{itemize}
\item[(1)] We keep the summands $(e_{[i,j]}+\delta e_{[i,j]})^{\oplus q_{[i,j]}}$ corresponding to the
arc $[i,j]\subset I_+^{up}$.
\item[(2)] If $n_{up}$ is
such that there exists $a=a(n_{up})\in Q_1^{\sigma}$, then we keep
the summands $(e_{[i,\sigma_I(i)]})^{\oplus q_{[i,\sigma_I(i)]}}$
corresponding to the symmetric arcs $[i,\sigma_I(i)]$ of
$\Delta_{up}$.
\item[(3)] If $n_{up}$ is
such that there exists $x=x(n_{up})\in Q_0^{\sigma}$, we have the
symmetric dimension vectors
$$e_{[i_1,\sigma_I(i_1)]}
,\ldots,e_{[i_{2s},\sigma_I(i_{2s})]}$$ corresponding to the arcs
$[i_1,\sigma_I(i_1)],\ldots, [i_{2s},\sigma_I(i_{2s})]$ such that
$[i_1,\sigma_I(i_1)]\supseteq\cdots\supseteq
[i_{2s},\sigma_I(i_{2s})]$. Then we divide them into pairs
$$([i_{2k},\sigma_I(i_{2k})],[i_{2k-1},\sigma_I(i_{2k-1})]),$$ with
$1\leq k\leq s$. For each pair we consider
$[i_{2k},\sigma_I(i_{2k-1})]\cup [i_{2k-1},\sigma_I(i_{2k})]$ and
we substitute $e_{[i_{2k},\sigma_I(i_{2k})]}\oplus e_{
[i_{2k-1},\sigma_I(i_{2k-1})]}$ for
$$e_{[i_{2k},\sigma_I(i_{2k-1})]}+ e_{
[i_{2k-1},\sigma_I(i_{2k})]}.$$
\end{itemize}
So, by equation \ref{qi,j}, in the symplectic case we get
\begin{itemize}
\item[(i)] if $n_{up}$ is
such that there exists $a=a(n_{up})\in Q_1^{\sigma}$,
\begin{equation}\label{supqi,j1}
\bar{d}_{up}=\bigoplus_{[i,j]\in\mathcal{A}_+^{up}(\bar{d})
}(e_{[i,j]}+\delta e_{[i,j]})^{\oplus q_{[i,j]}} +
\bigoplus_{[i,\sigma_I(i)]\in\mathcal{A}^{up}(\bar{d})
}(e_{[i,\sigma_I(i)]})^{\oplus q_{[i,\sigma_I(i)]}};
\end{equation}
\item[(ii)] If $n_{up}$ is
such that there exists $x=x(n_{up})\in Q_0^{\sigma}$,
\begin{equation}\label{supqi,j2}
\bar{d}_{up}=\bigoplus_{[i,j]\in\mathcal{A}_+^{up}(\bar{d})
}(e_{[i,j]}+\delta e_{[i,j]})^{\oplus q_{[i,j]}}
+\bigoplus_{k=1}^{s}(e_{[i_{2k},\sigma_I(i_{2k-1})]}+ e_{
[i_{2k-1},\sigma_I(i_{2k})]}).
\end{equation}
\end{itemize}
Similarly one proceeds for $\Delta_{down}$.\\
Finally we have to modify like in (3) the dimension vector
$ph+e_{[i,\sigma_I(i)]}$ if $p$ is odd and $[i,\sigma_I(i)]$
passes through $n_{up}$ such that there exists $x=x(n_{up})\in
Q_0^{\sigma}$.
\begin{esempi}
Let $(Q,\sigma)$ be the symmetric quiver
$\widetilde{A}^{1,1}_{0,6}$. We recall that
$x_{\frac{l}{2}}=\sigma(x_{\frac{l}{2}})$. $\Delta$ has the form
(\ref{esempio}).\\
As in example \ref{esempioe}, let $p_1=2$, $p_2=3$, $p_3=0$ and
$p_4=2$. The $\sigma_I$-fixed vertex $4$ is such that
$e_4(x_{\frac{l}{2}})\neq 0$. The only symmetric arc passing
through 4 is $[4,4]$. Thus we substitute $(e_4)^{\oplus 2}$ for
$2e_4$. So, in the symplectic case we get
$$
\sum_{i\in I_+}p_i(e_i+\delta e_i)+\sum_{i\in I_{\delta}}p_ie_i=
((e_2+\delta e_2)+e_1)^{\oplus 2}\oplus(e_2+\delta e_2)\oplus 2
e_4.
$$
\end{esempi}
Similarly we proceed with the decomposition of
$\bar{d}'$ and $\bar{d}''$.\\
Let $\bar{d}=\bar{d}_{up}+\bar{d}_{down}$ be a regular orthogonal
dimension vector. We consider $\Delta_{up}$. Starting from generic
decomposition (\ref{dgvs}) of $\bar{d}_{up}$ we modify it as
follows.
\begin{itemize}
\item[(1)] We keep the summands $(e_{[i,j]}+\delta e_{[i,j]})^{\oplus q_{[i,j]}}$ corresponding to the
arc $[i,j]\subset I_+^{up}$.
\item[(2)] If $n_{up}$ is
such that there exists $a=a(n_{up})\in Q_1^{\sigma}$ such that
$\bar{d}_{ta}$ is odd or $n_{up}$ is such that there exist
$x=x(n_{up})\in Q_0^{\sigma}$, then we keep the summands
$(e_{[i,\sigma_I(i)]})^{\oplus q_{[i,\sigma_I(i)]}}$ corresponding
to the symmetric arcs $[i,\sigma_I(i)]$ of $\Delta_{up}$.
\item[(3)] If $n_{up}$ is
such that there exists $a=a(n_{up})\in Q_1^{\sigma}$ such that
$\bar{d}_{ta}$ is even, we have the symmetric dimension vectors
$$e_{[i_1,\sigma_I(i_1)]}
,\ldots,e_{[i_{2s},\sigma_I(i_{2s})]}$$ corresponding to the arcs
$[i_1,\sigma_I(i_1)],\ldots, [i_{2s},\sigma_I(i_{2s})]$ such that
$[i_1,\sigma_I(i_1)]\supseteq\cdots\supseteq
[i_{2s},\sigma_I(i_{2s})]$. Then we divide them into pairs
$$([i_{2k},\sigma_I(i_{2k})],[i_{2k-1},\sigma_I(i_{2k-1})]),$$ with
$1\leq k\leq s$. For each pair we consider
$[i_{2k},\sigma_I(i_{2k-1})]\cup [i_{2k-1},\sigma_I(i_{2k})]$ and
we substitute $e_{[i_{2k},\sigma_I(i_{2k})]}\oplus e_{
[i_{2k-1},\sigma_I(i_{2k-1})]}$ for
$$e_{[i_{2k},\sigma_I(i_{2k-1})]}+ e_{
[i_{2k-1},\sigma_I(i_{2k})]}.$$
\end{itemize}
So, by equation \ref{qi,j}, in the orthogonal case we get
\begin{itemize}
\item[(i)] if $n_{up}$ is
such that there exists $a=a(n_{up})\in Q_1^{\sigma}$ such that
$\bar{d}_{ta}$ is odd or $n_{up}$ is such that there exist
$x=x(n_{up})\in Q_0^{\sigma}$,
\begin{equation}\label{supqi,j1}
\bar{d}_{up}=\bigoplus_{[i,j]\in\mathcal{A}_+^{up}(d')
}(e_{[i,j]}+\delta e_{[i,j]})^{\oplus q_{[i,j]}} +
\bigoplus_{[i,\sigma_I(i)]\in\mathcal{A}^{up}(d')
}(e_{[i,\sigma_I(i)]})^{\oplus q_{[i,\sigma_I(i)]}};
\end{equation}
\item[(ii)] if $n_{up}$ is
such that there exists $a=a(n_{up})\in Q_1^{\sigma}$ such that
$\bar{d}_{ta}$ is even,
\begin{equation}\label{supqi,j2}
\bar{d}_{up}=\bigoplus_{[i,j]\in\mathcal{A}_+^{up}(d')
}(e_{[i,j]}+\delta e_{[i,j]})^{\oplus q_{[i,j]}}
+\bigoplus_{k=1}^{s}(e_{[i_{2k},\sigma_I(i_{2k-1})]}+ e_{
[i_{2k-1},\sigma_I(i_{2k})]}).
\end{equation}
\end{itemize}
Similarly one proceeds for $\Delta_{down}$.\\
Finally we have to modify like in (3) the dimension vector
$ph+e_{[i,\sigma_I(i)]}$ if $p$ is odd and $[i,\sigma_I(i)]$
passes through $n_{up}$ such that there exists $a=a(n_{up})\in
Q_1^{\sigma}$ such that $\bar{d}_{ta}$ is even.
\begin{esempi}
Let $(Q,\sigma)$ be the symmetric quiver
$\widetilde{A}^{1,1}_{0,6}$. We recall that $b=\sigma(b)$.
$\Delta$ has the form
(\ref{esempio}).\\
As in example \ref{esempioe}, let $p_1=2$, $p_2=3$, $p_3=0$ and
$p_4=2$. The $\sigma_I$-fixed vertex 1 is such that $e_1(tb)\neq
0$ and $\bar{d}_{tb}$ is 2. The only symmetric arc passing through
1 is $[2,\sigma_I(2)]$. Thus we substitute $((e_2+\delta
e_2))+e_1)^{\oplus 2}$ for $2((e_2+\delta e_2))+e_1)$. So, in the
orthogonal case we get
$$
\sum_{i\in I_+}p_i(e_i+\delta e_i)+\sum_{i\in I_{\delta}}p_ie_i=
(e_4)^{\oplus 2}\oplus(e_2+\delta e_2)\oplus 2((e_2+\delta e_2)+
e_1).
$$
\end{esempi}
Similarly we proceed with the decomposition of $\bar{d}'$ and
$\bar{d}''$.\\
In general we have
\begin{proposizione}\label{dg}
Let $(Q,\sigma)$ be a symmetric quiver of tame type.
\begin{itemize}
\item[(1)] If $d$ is a regular symplectic dimension vector with decomposition
(\ref{fdc}). Then
\begin{equation}\label{fdg}
d=\bigoplus_{i=1}^p
h\oplus\bar{d}_{up}\oplus\bar{d}_{down}\oplus\bar{d}'_{up}\oplus\bar{d}'_{down}\oplus\bar{d}''_{up}\oplus\bar{d}''_{down}
\end{equation}
is the symplectic generic decomposition of $d$.
\item[(2)] If $d$ is a regular orthogonal dimension vector with decomposition
(\ref{fdc}). Then (\ref{fdc}). Then
\begin{equation}\label{fdgo}
d=\bigoplus_{i=1}^p
h\oplus\bar{d}_{up}\oplus\bar{d}_{down}\oplus\bar{d}'_{up}\oplus\bar{d}'_{down}\oplus\bar{d}''_{up}\oplus\bar{d}''_{down}
\end{equation}
is the orthogonal generic decomposition of $d$.
\end{itemize}
\end{proposizione}
For the proof, we need two propositions. We state and prove these
propositions only for regular indecomposable symplectic
(respectively orthogonal) representations related to polygon
$\Delta$, because for those related to polygon $\Delta'$ and to
polygon $\Delta''$ the statement and the proof are similar.
\begin{proposizione}\label{extreg1}
Let $(Q,\sigma)$ be a symmetric quiver of tame tape. Let $V_1\neq
V_2$ be two regular indecomposable symplectic (respectively
orthogonal) representations of $(Q,\sigma)$ with symmetric
dimension vector corresponding respectively to the arc $[i,j]$ and
the arc $[h,k]$ of $\Delta$ ($\Delta'$ or $\Delta''$). Moreover we
suppose that $[i,j]$ and $[h,k]$ don't satisfy the following
properties
\begin{itemize}
\item[(i)] $[i,j]\cap [h,k]\neq \emptyset$ and $[i,j]$ doesn't
contain $[h,k]$;
\item[(ii)] $[i,j]\cap [h,k]\neq \emptyset$ and $[h,k]$ doesn't
contain $[i,j]$;
\item[(iii)] $[i,j]$ and $[h,k]$ are linked by one edge of
$\Delta$ (respectively $\Delta'$ or $\Delta''$).
\end{itemize}
Then $Ext^1_Q(V_1,V_2)=0$.
\end{proposizione}
\textit{Proof.} We restrict to decomposition $\bar{d}_j=\sum_{i\in
I_+} p_i^j(e_i+\delta e_i)+\sum_{i\in I_{\delta}}p_i^j e_i$, for
$j=1,2$. We have nine cases:
\begin{itemize}
\item[(1)] $V_1=E_{i,\sigma_I(i)}$, $V_2=E_{j,\sigma_I(j)}$ and $V_1=E_{\sigma_I(j),j}$, $V_2=E_{\sigma_I(i),i}$  with $i,j\in I_+\sqcup
I_{\delta}$.
\item[(2)] $V_1=E_{i,\sigma_I(i)}$, $V_2= E_{\sigma_I(j),j}$ and $V_1= E_{\sigma_I(j),j}$, $V_2= E_{i,\sigma_I(i)}$ with $i,j\in I_+\sqcup
I_{\delta}$ such that $j> i+1$.
\item[(3)] $V_1= E_{i,j}\oplus E_{\sigma_I(j),\sigma_I(i)}$, $V_2= E_{k,\sigma_I(k)}$ and $V_1=E_{k,\sigma_I(k)}$, $V_2= E_{i,j}\oplus
E_{\sigma_I(j),\sigma_I(i)}$ with $i,j,k\in I_+\sqcup I_{\delta}$
such that either $j> k+1$ or $k\geq i$.
\item[(4)] $V_1= E_{i,j}\oplus
E_{\sigma_I(j),\sigma_I(i)}$, $V_2= E_{\sigma_I(k),k}$ or
 $V_1=E_{\sigma_I(k),k}$, $V_2= E_{i,j}\oplus
E_{\sigma_I(j),\sigma_I(i)}$ with $i,j,k\in I_+\sqcup I_{\delta}$
such that either $j\geq k$ or $k>i+1$.
\item[(5)] $V_1= E_{i,j}\oplus
E_{\sigma_I(j),\sigma_I(i)}$ and $V_2= E_{h,k}\oplus
E_{\sigma_I(k),\sigma_I(h)}$ with $i,j,k,h\in I_+\sqcup
I_{\delta}$ such that either $k\leq j$ and $i\leq h$ or $k\geq j$
and $i\geq h$.
\item[(6)] $V_1= E_{i,\sigma_I(j)}\oplus
E_{j,\sigma_I(i)}$, $V_2= E_{h,\sigma_I(k)}\oplus
E_{k,\sigma_I(h)}$ and $V_1= E_{\sigma_I(j),i}\oplus
E_{\sigma_I(i),j}$ and $V_2= E_{\sigma_I(k),h}\oplus
E_{\sigma_I(h),k}$) with $i,j,k\in I_+\sqcup I_{\delta}$.
\item[(7)] $V_1= E_{i,\sigma_I(j)}\oplus
E_{j,\sigma_I(i)}$, $V_2= E_{\sigma_I(k),k}$ (resp.
$V_1=E_{\sigma_I(k),k}$, $V_2= E_{i,\sigma_I(j)}\oplus
E_{j,\sigma_I(i)}$) with $i,j,k\in I_+\sqcup I_{\delta}$ such that
$k> i+1$ and $i>j$ and $V_1= E_{\sigma_I(j),i}\oplus
E_{\sigma_I(i),j}$, $V_2=E_{k,\sigma_I(k)}$ (resp.
$V_1=E_{k,\sigma_I(k)}$, $V_2= E_{\sigma_I(j),i}\oplus
E_{\sigma_I(i),j}$) with $i,j,k\in I_+\sqcup I_{\delta}$ such that
$i> k+1$ and $i<j$.
\item[(8)] $V_1= E_{i,\sigma_I(j)}\oplus
E_{j,\sigma_I(i)}$, $V_2= E_{h,k}\oplus
E_{\sigma_I(k),\sigma_I(h)}$ (resp. $V_1= E_{h,k}\oplus
E_{\sigma_I(k),\sigma_I(h)}$, $V_2= E_{i,\sigma_I(j)}\oplus
E_{j,\sigma_I(i)}$) with $i,j,k\in I_+\sqcup I_{\delta}$ such that
$i>j$ and either $k> i+1$ or $i\geq h$ and $V_1=
E_{\sigma_I(j),i}\oplus E_{\sigma_I(i),j}$, $V_2= E_{h,k}\oplus
E_{\sigma_I(k),\sigma_I(h)}$ (resp. $V_1= E_{h,k}\oplus
E_{\sigma_I(k),\sigma_I(h)}$, $V_2= E_{\sigma_I(j),i}\oplus
E_{\sigma_I(i),j}$) with $i,j,k\in I_+\sqcup I_{\delta}$ such that
$i<j$ and either $k\geq i$ or $i> h+1$.
\item[(9)] $V_1= E_{i,\sigma_I(j)}\oplus
E_{j,\sigma_I(i)}$ and $V_2= E_{\sigma_I(k),h}\oplus
E_{\sigma_I(h),k}$ (resp. $V_1= E_{\sigma_I(k),h}\oplus
E_{\sigma_I(h),k}$ and $V_2= E_{i,\sigma_I(j)}\oplus
E_{j,\sigma_I(i)}$) with $i,j,k\in I_+\sqcup I_{\delta}$ such that
$h> i+1$, $i>j$ and $h<k$.
\end{itemize}
We consider (1). By [D, lemma 4.1],
$$Hom_Q(E_{i,\sigma_I(i)},E_{j,\sigma_I(j)})=0=Hom_Q(E_{\sigma_I(j),j},E_{\sigma_I(i),i})$$
and by lemma \ref{eEe},
$$\langle\underline{dim}(E_{i,\sigma_I(i)}),\underline{dim}(E_{j,\sigma_I(j)})\rangle=0=
\langle\underline{dim}(E_{\sigma_I(j),j}),\underline{dim}(E_{\sigma_I(i),i})\rangle.$$
So we get
$$Ext^1_Q(E_{i,\sigma_I(i)},E_{j,\sigma_I(j)})=0=Ext^1_Q(E_{\sigma_I(j),j},E_{\sigma_I(i),i}).$$
Similarly for (2), by [D, lemma 4.1] and by lemma \ref{eEe}, we
get
$Ext^1_Q(V_1,V_2)=0.$\\
We consider (3). We suppose $j> k+1$.  By [D, lemma 4.1], we have
$$
Hom_Q(E_{i,j}, E_{k,\sigma_I(k)})=0=
Hom_Q(E_{\sigma_I(j),\sigma_I(i)}, E_{k,\sigma_I(k)})
$$
and so
$$
Hom_Q(E_{i,j}\oplus E_{\sigma_I(j),\sigma_I(i)},
E_{k,\sigma_I(k)})
$$
$$
=Hom_Q(E_{i,j}, E_{k,\sigma_I(k)})\oplus
Hom_Q(E_{\sigma_I(j),\sigma_I(i)}, E_{k,\sigma_I(k)})=0.
$$
Moreover, by lemma \ref{eEe}
$$
\langle\underline{dim}(E_{i,j}),\underline{dim}(E_{k,\sigma_I(k)})\rangle=0=\langle
\underline{dim}(E_{\sigma_I(j),\sigma_I(i)}),\underline{dim}(E_{k,\sigma_I(k)})\rangle
$$
and hence
$$
\langle\underline{dim}(E_{i,j}\oplus
E_{\sigma_I(j),\sigma_I(i)}),\underline{dim}(E_{k,\sigma_I(k)})\rangle=
$$
$$
\langle\underline{dim}(E_{i,j}),\underline{dim}(E_{k,\sigma_I(k)})\rangle+\langle
\underline{dim}(E_{\sigma_I(j),\sigma_I(i)}),\underline{dim}(E_{k,\sigma_I(k)})\rangle=0.
$$
So we have
$$
Ext_Q^1(E_{i,j}\oplus E_{\sigma_I(j),\sigma_I(i)},
E_{k,\sigma_I(k)})=0.
$$
Similarly to (3), one proceeds for the other cases.  $\Box$

\begin{proposizione}\label{autoext}
Let $(Q,\sigma)$ be a symmetric quiver of tame tape. Let $V$ be a
regular indecomposable symplectic (respectively orthogonal)
representation of $(Q,\sigma)$ such that $\underline{dim}(V)=h$ or
$\bar{d}$. Moreover we suppose $V\neq E_{i,j}\oplus
E_{\sigma_I(j),\sigma_I(i)}$ with $i,j\in I_+$ such that
$e_i(ta)\neq 0$ or $e_j(ta)\neq 0$ for $a\in Q_1^{\sigma}$. Then,
for every non-trivial short exact sequence
$$
0\rightarrow V\rightarrow W\rightarrow V\rightarrow 0,
$$
$W$ is not symplectic (respectively it is not orthogonal).
\end{proposizione}
\textit{Proof.} We give a proof for
$(Q=\widetilde{A}_{k,l}^{2,0,1},\sigma)$ for the symplectic case,
one proves similarly the other cases.\\
(i) Let $\underline{dim}(V)=h$. By lemma \ref{regsport}, the
regular indecomposable symplectic representation of dimension $h$
is $E_{i,\sigma_I(i)}$ containing $E_{\frac{l}{2}+1}$, i.e. the
representation $V$ defined by $V(x)=\mathbb{K}$ for every $x\in
Q_0$ and
$$
V(c)=\left\{\begin{array}{ll} 0& \textrm{if}\;c=a\\Id
&\textrm{otherwise},\end{array}\right.
$$
for $c\in Q_1$.\\
By [D, lemma 4.1], $Hom_Q(V,V)=\mathbb{K}$ and since $\langle h,
h\rangle=0$, then $Ext^1_Q(V,V)=\mathbb{K}$. One non-trivial
auto-extension $W$ of $V$ is defined by $W(x)=\mathbb{K}^2$ for
every $x\in Q_0$, and
$$
W(c)=\left\{\begin{array}{cl} \left(\begin{array}{cc}0& 1\\0&
0\end{array}\right)&
\textrm{if}\;c=a\\&\\\left(\begin{array}{cc}1& 0\\0&
1\end{array}\right) &\textrm{otherwise},\end{array}\right.
$$
for $c\in Q_1$. Finally we note that $W$ is not symplectic,
because $W(a)$ is not symmetric. Since $Ext_Q^1(V,V)=\mathbb{K}$,
the non-trivial auto-extensions of $V$ is not symplectic.\\
(ii) Let $\underline{dim}(V)=\bar{d}$. The only regular
indecomposable symplectic representations which we have to
consider are $E_{i,\sigma_I(j)}\oplus E_{j,\sigma_I(i)}$ and
$E_{\sigma_I(j),i}\oplus E_{\sigma_I(i),j}$ with $i,j\in I_+\sqcup
I_{\delta}$.\\
Let $V=E_{i,\sigma_I(j)}\oplus E_{j,\sigma_I(i)}$, with $j<i$.
$$
V(x)=V(\sigma(x))\left\{\begin{array}{ll}\mathbb{K}&\textrm{if}\;
x\in\{x_r|\;e_m(x_r)\neq 0,\;m\in\{j+1,\ldots,i\}\}\\
0&\textrm{if}\;x\in\{x_r|\;e_m(x_r)=0,m\in I_+\}\\
\mathbb{K}^2&\textrm{otherwise}\end{array}\right.
$$
for $x\in Q_0$ and
$$
V(c)=-V(\sigma(c))^t=\left\{\begin{array}{ll}1 &
\textrm{if}\;c\in\{v_r|\;e_m(tv_r)\neq 0,\;m\in\{j+1,\ldots,i\}\}\\
(1,1) & \textrm{if}\;c=v_r\;\textrm{s.t.}\;e_j(tv_r)\neq 0\\
0 & \textrm{if}\;c\in\{v_r|\;e_m(tv_r)=0,m\in I_+\}\cup\{a\}\\
Id_{2\times 2} & \textrm{otherwise}\end{array}\right.
$$
for $c\in Q_1^+$ and $V(b)=Id_{2\times 2}$.\\
By [D, lemma 4.1],
$$
dim_{\mathbb{K}}(Hom_Q(E_{i,\sigma_I(j)}\oplus
E_{j,\sigma_I(i)},E_{i,\sigma_I(j)}\oplus E_{j,\sigma_I(i)}))=3
$$
and by lemma \ref{eEe},
$$
\langle\underline{dim}(E_{i,\sigma_I(j)}\oplus
E_{j,\sigma_I(i)}),\underline{dim}(E_{i,\sigma_I(j)}\oplus
E_{j,\sigma_I(i)})\rangle=2.
$$
So we have
$$
Ext_Q^1(E_{i,\sigma_I(j)}\oplus
E_{j,\sigma_I(i)},E_{i,\sigma_I(j)}\oplus
E_{j,\sigma_I(i)})=\mathbb{K}.
$$
Let
$$
\begin{array}{cc}A=\left(\begin{array}{cccc}1&1&0&0\\0&0&1&1\end{array}\right)\quad\textrm{and}\quad
B=\left(\begin{array}{cccc}1&0&0&1\\0&1&0&0\\0&0&1&0\\0&0&0&1\end{array}\right).\end{array}
$$
One non-trivial auto-extension $W$ of $V$ is defined by
$$
W(x)=W(\sigma(x))\left\{\begin{array}{ll}\mathbb{K}^2&\textrm{if}\;
x\in\{x_r|\;e_m(x_r)\neq 0,\;m\in\{j+1,\ldots,i\}\}\\
0&\textrm{if}\;x\in\{x_r|\;e_m(x_r)=0,m\in I_+\}\\
\mathbb{K}^4&\textrm{otherwise}\end{array}\right.
$$
for $x\in Q_0$ and
$$
W(c)=-W(\sigma(c))^t=\left\{\begin{array}{ll}Id_{2\times 2} &
\textrm{if}\;c\in\{v_r|\;e_m(tv_r)\neq 0,\;m\in\{j+1,\ldots,i\}\}\\
A & \textrm{if}\;c=v_r\;\textrm{s.t.}\;e_j(tv_r)\neq 0\\
0 & \textrm{if}\;c\in\{v_r|\;e_m(tv_r)=0,m\in I_+\}\cup\{a\}\\
Id_{4\times 4} & \textrm{otherwise},\end{array}\right.
$$
for $c\in Q_1^+$ and $W(b)=B$. Finally we note that $W$ is not
symplectic because $W(b)$ is not symmetric. Since
$Ext^1_Q(V,V)=\mathbb{K}$, this concludes the
proof for $(\widetilde{A}_{k,l}^{2,0,1},\sigma)$. $\Box$\\
\\
\textit{Proof of \ref{dg}.} \textit{(1)} Let $d$ be a symplectic
regular dimension vector with decomposition (\ref{fdg}). First we
note that the symmetric dimension vectors appearing in
decomposition (\ref{dg}) are not dimension vectors of the regular
indecomposable symplectic representations which are exceptions of
proposition \ref{extreg1} and \ref{autoext}. Let $\mathcal{O}(d)$
be the open orbit of the regular symplectic representations of
dimension $d$. By [Bo1] and [Z], we obtain each
representation $V$ in $\mathcal{O}(d)$ as follows.\\
There are representations $M_i$, $U_i$, $V_i$ and short exact
sequences
$$
0\rightarrow U_i\rightarrow M_i\rightarrow V_i\rightarrow 0
$$
such that $M_{i+1}=U_i\oplus V_i$ and $V=U_{n+1}\oplus V_{n+1}$,
with $1\leq i\leq n$ for some $n\in\mathbb{N}$.\\
By propositions \ref{extreg1} and \ref{autoext}, we have
\begin{itemize}
\item[(i)] If $U_i\neq V_i$, then $Ext^1_Q(V_i,U_i)=0$.
\item[(ii)] If $U_i=V_i$, then either $Ext^1_Q(U_i,U_i)=0$ or
no one non-trivial auto-extension of $U_i$ is symplectic. So, if
$Ext^1_Q(U_i,U_i)\neq 0$ then $U_i$ doesn't appear in
decomposition of a symplectic representation.
\end{itemize}
Hence $V$ decomposes in regular indecomposable symplectic
representations of dimension $\beta_i$, where $\beta_i$ are
regular symmetric dimension
vectors appearing in decomposition (\ref{fdg}) of $d$.\\
\textit{(2)} One proves similarly to \textit{(1)}. $\Box$\\
\\
Let $d$ be a regular symmetric vector with a decomposition
(\ref{fdg}) or (\ref{fdgo}). We note that if $d=d_1+d_2$ with
$d_1$ and $d_2$ summands of this generic decomposition, we have
canonical embeddings
\begin{equation}\label{Spe}
SpSI(Q,d)\stackrel{\Phi_d}{\rightarrow}\bigoplus_{\chi\in
char(Sp(Q,d))}SpSI(Q,d_1)_{\chi|_{d_1}}\otimes
SpSI(Q,d_2)_{\chi_{d_2}}
\end{equation}
and
\begin{equation}\label{Oe}
OSI(Q,d)\stackrel{\Psi_d}{\rightarrow}\bigoplus_{\chi\in
char(O(Q,d))}OSI(Q,d_1)_{\chi|_{d_1}}\otimes
OSI(Q,d_2)_{\chi_{d_2}},
\end{equation}
induced by the restriction homomorphism. We prove theorem
\ref{stq} by induction on the number of the summands
$e_{[i,j]}+\delta e_{[i,j]}$, $e_{[i,\sigma_I(i)]}$,
$e_{[i_{2k},\sigma_I(i_{2k-1})]}+e_{[i_{2k-1},\sigma_I(i_{2k})]}$
and respective summands corresponding to the admissible arcs in
$\mathcal{A}'(d)$ and in $\mathcal{A}''(d)$. If this number is 0,
then $d=ph$ and it was already proved. We suppose that the generic
decomposition of $d$ contains one of those summands and, without
loss of generality, we can assume that this summand is one of
those corresponding to the arcs in $\mathcal{A}(d)$. In particular
we suppose that this summand is $e_{[s,\sigma_I(s)]}$ (one
proceeds similarly for the other types), with $s\in I_+\sqcup
I_{\delta}$, and we can assume $ind[s,\sigma_I(s)]=r=max\{p_k\}$.
We call $d_2=e_{[s,\sigma_I(s)]}$ and so
$d_1=d-e_{[s,\sigma_I(s)]}$. Now we compare the generators of the
algebras $SpSI(Q,d)$ and $SpSI(Q,d_1)$ (respectively  $OSI(Q,d)$
and $OSI(Q,d_1)$). By induction the generators of $SpSI(Q,d_1)$
(respectively of $OSI(Q,d_1)$) are described by theorem \ref{stq}.
Since $\Delta'(d)=\Delta'(d_1)$ and $\Delta''(d)=\Delta''(d_1)$,
the generators $c_0,\ldots,c_t$ (with $t=\frac{p}{2}$,
$\frac{p-1}{2}$ or $p$), those corresponding to the arcs from
$\mathcal{A}'(d)$ and those corresponding to the arcs from
$\mathcal{A}''(d)$ occur. So it's enough to study the behavior of
the semi-invariants corresponding to the arcs from
$\mathcal{A}(d)$. We describe the link between the admissible arcs
of the polygons $\Delta(d)$ and $\Delta(d_1)$. We have
$$
d_1=ph+\sum_{i\in I_+\setminus(I_+\cap
[s,\sigma_I(s)])}p_i(e_i+\delta e_i)+\sum_{i\in
I_{\delta}\setminus(I_{\delta}\cap [s,\sigma_I(s)])}p_i e_i+
$$
$$
\sum_{i\in I_+\cap [s,\sigma_I(s)]}p_i(e_i+\delta e_i)+\sum_{i\in
I_{\delta}\cap [s,\sigma_I(s)]}p_i e_i+
$$
$$
\sum_{i\in I'_+}p'_i(e'_i+\delta e'_i)+\sum_{i\in I'_{\delta}}p'_i
e'_i+\sum_{i\in I''_+}p''_i(e''_i+\delta e''_i).
$$
We have two cases
\begin{itemize}
\item[(1)] $p_{s-1}=p_{\sigma_I(s)+1}<r-1$ with $s-1\in I_+$,
\item[(2)] $p_{s-1}=p_{\sigma_I(s)+1}=r-1$ with $s-1\in I_+$.
\end{itemize}
in the case (1) the only difference between the structure of
$\mathcal{A}(d)$ and $\mathcal{A}(d_1)$ is that the admissible
arcs $[s,s+1],[s+1,s+2],\ldots,[\sigma_I(s)-1,\sigma(s)]$ are of
index $r$ in $\mathcal{A}(d)$ and of index $r-1$ in
$\mathcal{A}(d_1)$. In the case (2) we have the admissible arc
$[s-1,\sigma_I(s)+1]$ of index $r-1$. The admissible arcs
$[s,s+1],[s+1,s+2],\ldots,[\sigma_I(s)-1,\sigma_I(s)]$ are of
index $s$ in $\mathcal{A}(d)$ and the admissible arcs
$[s-1,s],[s,s+1],\ldots,[\sigma_I(s)-1,\sigma_I(s)],[\sigma_I(s),\sigma_I(s)+1]$
are of index $r-1$ in $\mathcal{A}(d_1)$.\\
Now we prove that the embeddings $\Phi_d$ and $\Psi_d$ are
isomorphisms and this will be done in two steps. The first step is
to show case by case that the semi-invariants corresponding to the
admissible arcs $[i,j]$ are non zero $c^V$ for some $V\in Rep(Q)$
and, if $V$ satisfy property \textit{(Spp)} or \textit{(Op)}, they
are non zero $pf^V$. The second step is to give an explicit
description of the generators of the algebras on the right hand
side of $\Phi_d$ and $\Psi_d$. This is based on the knowledge,
given by inductive hypothesis, of the algebra $SpSI(Q,d_1)$
(respectively $OSI(Q,d_1)$). We can describe explicitly the
generators of the algebra $SpSI(Q,d_2)$ (respectively
$OSI(Q,d_2)$) and we can note that they are determinants or
pfaffians, knowing that the group $Sp(Q,d_2)$ (respectively
$O(Q,d_2)$) has an open orbit in $SpRep(Q,d_2)$ (respectively
$ORep(Q,d_2)$) and hence that $SpSI(Q,d_2)$ (respectively
$OSI(Q,d_2)$) is a polynomial ring (lemma \ref{sk}). At this point
we know the generators of the algebras on the right hand side of
$\Phi_d$ and $\Psi_d$. Now, using the fact that these are
determinants or pfaffians, we prove that they actually are in
$SpSI(Q,d)$ (respectively in $OSI(Q,d)$) and that the
embeddings $\Phi_d$ and $\Psi_d$ are isomorphisms.\\
We will consider case by case the semi-invariants corresponding to
each admissible arc $[i,j]$. To simplify the notation we shall
call $a$ both the arrow $a\in Q_1$ and the linear map $V(a)$
defined on $a$, where $V$ is a representation of $Q$.

\subsubsection{3.2.1.1\quad$\widetilde{A}_{k,l}^{2,0,1}$}\addcontentsline{toc}{subsection}{3.2.1.1\quad$\widetilde{A}_{k,l}^{2,0,1}$} We have at most two
$\tau^+$-orbits $\Delta$ and $\Delta'$ of the dimension vectors of
nonhomogeneous simple regular representation. We assume $n\geq 2$
and we consider the $\tau$-orbit $\{e_1=\delta
e_1,e_2,\ldots,e_{[\frac{l}{2}]+1},\delta
e_{[\frac{l}{2}]+1},\ldots,\delta e_2\}$. Let
$[i,j]\in\mathcal{A}(d)$. If we consider the arc $[1,1]$ of index
0, i.e. $p_1=0,p_2\ne 0,\ldots,p_{[\frac{l}{2}]+1}\ne 0$, we have
the minimal projective resolution of $V_{(0,1)}$
$$
0\longrightarrow
P_{\sigma(a_0)}\stackrel{d^{V_{(0,1)}}_{min}}{\longrightarrow}
P_{a_0}\longrightarrow V_{(0,1)}\longrightarrow 0
$$
where
$d^{V_{(0,1)}}_{min}=\sigma(v_1)\cdots\sigma(v_{\frac{l}{2}})av_{\frac{l}{2}}\cdots
v_1$ and so
$$
c^{V_{(0,1)}}=det(Hom_Q(d^{V_{(0,1)}}_{min},\cdot))=det(\sigma(v_1)\cdots\sigma(v_{\frac{l}{2}})av_{\frac{l}{2}}\cdots
v_1)
$$
in the symplectic case and
$pf^{V_{(0,1)}}=pf(\sigma(v_1)\cdots\sigma(v_{\frac{l}{2}})av_{\frac{l}{2}}\cdots
v_1)$ in the orthogonal case, since in this case $a$ is
skew-symmetric and $\sigma(v_i)=-(v_i)^t$. If we consider the arc
$[\sigma_I(2),2]=[0,2]$ of index 0, i.e.
$p_{\sigma_I(2)}=0=p_2,p_1\ne 0$, we have the minimal projective
resolution of $V_{(1,0)}$
$$
0\longrightarrow
P_{\sigma(a_0)}\stackrel{d^{V_{(1,0)}}_{min}}{\longrightarrow}
P_{a_0}\longrightarrow V_{(1,0)}\longrightarrow 0
$$
where
$d^{V_{(1,0)}}_{min}=\sigma(u_1)\cdots\sigma(u_{\frac{k}{2}})bu_{\frac{k}{2}}\cdots
u_1$ and so
$$
c^{V_{(1,0)}}=det(Hom_Q(d^{V_{(1,0)}}_{min},\cdot))=
det(\sigma(u_1)\cdots\sigma(u_{\frac{k}{2}})bu_{\frac{k}{2}}\cdots
u_1)
$$
in the symplectic case and
$pf^{V_{(1,0)}}=pf(\sigma(u_1)\cdots\sigma(u_{\frac{k}{2}})bu_{\frac{k}{2}}\cdots
u_1)$ in the orthogonal case, since in this case $b$ is
skew-symmetric and $\sigma(u_i)=-(u_i)^t$. We note that for $l=2$
we have only the admissible arcs $[1,1]$ an $[\sigma_I(2),2]$. We
assume now that $l\geq 4$ ($l$ is even) and $[i,j]$ is not an
admissible arc considered above. If $1\leq i<j\leq\frac{l}{2}+1$,
then we identify $[i,j]$ with the path $v_{j-1}\cdots v_{i}$ in
$Q$ and we have the minimal projective resolution of $E_{i,j-1}$
$$
0\longrightarrow
P_{x_{j-1}}\stackrel{d^{E_{i,j-1}}_{min}}{\longrightarrow}
P_{x_{i-1}}\longrightarrow E_{i,j-1}\longrightarrow 0
$$
where $d^{E_{i,j-1}}_{min}=v_{j-1}\cdots v_{i}$ and so
$$
c^{E_{i,j-1}}=det(Hom_Q(d^{E_{i,j-1}}_{min},\cdot))=
det(v_{j-1}\cdots v_{i}).
$$
We note that
$$
c^{\tau^-\nabla
E_{i,j-1}}=c^{E_{\sigma_I(j),\sigma_I(i)-1}}=det(\sigma(v_{i})\cdots
\sigma(v_{j-1}))=det(v_{j-1}\cdots v_{i})=c^{E_{i,j-1}}.
$$
If $j=\sigma_I(i)$ then in the symplectic case we get
$c^{E_{i,\sigma_I(i)-1}}=det(\sigma(v_{i})\cdots a\cdots v_{i})$
and in the orthogonal case, we get
$pf^{E_{i,\sigma_I(i)-1}}=pf(\sigma(v_{i})\cdots a\cdots
v_{i})$,\quad since $\sigma(v_{i})\cdots a\cdots v_{i}$ is
skew-symmetric. Now we consider the arcs $[i,j]$ which have $e_1$
as internal vertex. For these arcs, $2\leq j< i-1<l$ and $[i,j]$
can be identify with the path in $Q$ consisting of the path
$v_l\cdots
 v_{i-1}=\sigma(v_1)\cdots v_{i-1}$, then coming back by
$\sigma(u_1)\cdots b\cdots u_1$ and at last passing for
$v_{j-1}\cdots v_1$. We have the minimal projective resolution of
$E_{i,j-1}$
$$
0\longrightarrow P_{\sigma(a_0)}\oplus
P_{x_{j-1}}\stackrel{d^{E_{i,j-1}}_{min}}{\longrightarrow}
P_{a_0}\oplus P_{x_{i-2}}\longrightarrow E_{i,j-1}\longrightarrow
0
$$
where
$d^{E_{i,j-1}}_{min}=\left(\begin{array}{cc}\sigma(u_1)\cdots
b\cdots u_1
&v_{j-1}\cdots v_1\\
 \sigma(v_1)\cdots v_{i-1}& 0 \end{array}\right)$ and so
$$
c^{E_{i,j-1}}=det(Hom_Q(d^{E_{i,j-1}}_{min},\cdot))=det\left(\begin{array}{cc}\sigma(u_1)\cdots
b\cdots u_1
&\sigma(v_1)\cdots v_{i-1}\\
v_{j-1}\cdots v_1 & 0 \end{array}\right).
$$
In particular we note that if $i=\sigma_I(j)$, in the orthogonal
case, we get
$$pf^{E_{\sigma_I(j),j-1}}=pf\left(\begin{array}{cc}\sigma(u_1)\cdots
b\cdots u_1
&\sigma(v_1)\cdots \sigma(v_{j-1})\\
v_{j-1}\cdots v_1 & 0
\end{array}\right),$$
since $b$ is skew-symmetric and $\sigma(v_i)=-(v_i)^t$. Finally we note that $V_{(0,1)}$, $V_{(1,0)}$,
${E_{i,\sigma_I(i)-1}}$ and ${E_{\sigma_I(j),j-1}}$ satisfy
property \textit{(Spp)}. Similarly we define the semi-invariants
for the admissible arcs $[i,j]$ in $\mathcal{A'}(d)$, exchanging
the upper paths of $\widetilde{A}_{k,l}^{2,0,1}$ with the lower
ones.

\subsubsection{3.2.1.2\quad$\widetilde{A}_{k,l}^{2,0,2}$}\addcontentsline{toc}{subsection}{3.2.1.2\quad$\widetilde{A}_{k,l}^{2,0,2}$} We have at most two
$\tau^+$-orbits $\Delta$ and $\Delta'$ of the dimension vectors of
nonhomogeneous simple regular representation. We assume $n\geq 2$
and we consider the $\tau$-orbit
$$\{e_2,\ldots,e_{[\frac{l}{2}]+2},\delta
e_{[\frac{l}{2}]+2},\ldots,\delta e_2=e_1\}.$$ Let
$[i,j]\in\mathcal{A}(d).$ If we consider the arc
$[\sigma_I(2),2]=[1,2]$ of index 0, i.e. $p_2=0,p_3\ne 0,
\ldots,p_{[\frac{l}{2}]+2}\ne 0$, we have the minimal projective
resolution of $V_{(1,0)}$
$$
0\longrightarrow
P_{\sigma(a_0)}\stackrel{d^{V_{(1,0)}}_{min}}{\longrightarrow}
P_{a_0}\longrightarrow V_{(1,0)}\longrightarrow 0
$$
where
$d^{V_{(1,0)}}_{min}=\sigma(v_1)\cdots\sigma(v_{\frac{l}{2}})av_{\frac{l}{2}}\cdots
v_1$ and so
$$
c^{V_{(1,0)}}=det(Hom_Q(d^{V_{(1,0)}}_{min},\cdot))=
det(\sigma(v_1)\cdots\sigma(v_{\frac{l}{2}})av_{\frac{l}{2}}\cdots
v_1)
$$
in the symplectic case and
$pf^{V_{(1,0)}}=pf(\sigma(v_1)\cdots\sigma(v_{\frac{l}{2}})av_{\frac{l}{2}}\cdots
v_1)$ in the orthogonal case, since in this case $a$ is
skew-symmetric and $\sigma(v_i)=-(v_i)^t$. If we consider the arc
$[\sigma_I(3),3]=[0,3]$ of index 0, i.e. $p_3=0,p_2\ne 0$, we have
the minimal projective resolution of $V_{(0,1)}$
$$
0\longrightarrow P_{y_{\frac{k}{2}}}\oplus
P_{\sigma(a_0)}\stackrel{d^{V_{(0,1)}}_{min}}{\longrightarrow}
P_{\sigma(y_{\frac{k}{2}})}\oplus P_{a_0}\longrightarrow
V_{(0,1)}\longrightarrow 0
$$
where $d^{V_{(0,1)}}_{min}=\left(\begin{array}{cc}
b & \sigma(u_1)\cdots\sigma(u_{\frac{k}{2}})\\
u_{\frac{k}{2}}\cdots u_1 & 0\end{array}\right)$ and so
$$
c^{V_{(0,1)}}=det(Hom_Q(d^{V_{(0,1)}}_{min},\cdot))=det\left(\begin{array}{cc}
b &
 u_{\frac{k}{2}}\cdots u_1
\\ \sigma(u_1)\cdots\sigma(u_{\frac{k}{2}})& 0\end{array}\right)
$$
in the symplectic case and
$$pf^{V_{(0,1)}}=pf\left(\begin{array}{cc} b &
 u_{\frac{k}{2}}\cdots u_1
\\ \sigma(u_1)\cdots\sigma(u_{\frac{k}{2}})& 0\end{array}\right)$$ in the orthogonal case, since $b$ is
skew-symmetric and $\sigma(u_i)=-(u_i)^t$. We note that for $l=2$
we have only the admissible arcs $[\sigma_I(2),2]$ an
$[\sigma_I(3),3]$. We assume now that $l\geq 4$ and $[i,j]$ is not
an admissible arc considered above. If $2\leq i<j\in I\leq
\frac{l}{2}+2$, then we identify $[i,j]$ with the path
$v_{j-2}\cdots v_{i-1}$ in $Q$ and we have the minimal projective
resolution of $E_{i,j-1}$
$$
0\longrightarrow
P_{x_{j-2}}\stackrel{d^{E_{i,j-1}}_{min}}{\longrightarrow}
P_{x_{i-2}}\longrightarrow E_{i,j-1}\longrightarrow 0
$$
where $d^{E_{i,j-1}}_{min}=v_{j-2}\cdots v_{i-1}$ and so
$$
c^{E_{i,j-1}}=det(Hom_Q(d^{E_{i,j-1}}_{min},\cdot))=
det(v_{j-2}\cdots v_{i-1}).
$$
We note that
$$
c^{\tau^-\nabla
E_{i,j-1}}=c^{E_{\sigma_I(j),\sigma_I(i)-1}}=det(\sigma(v_{i-1})\cdots
\sigma(v_{j-2}))=det(v_{j-2}\cdots v_{i-1})=c^{E_{i,j-1}}.
$$
Moreover, if $j=\sigma_I(i)$ then, only in the orthogonal case, we
get\quad $pf^{E_{i,\sigma_I(i)-1}}=pf(\sigma(v_{i-1})\cdots
a\cdots v_{i-1})$ since $\sigma(v_{i-1})\cdots a\cdots v_{i-1}$ is
skew-symmetric. Now we consider the arcs $[i,j]$ which have
$\delta e_2=e_1$ and $e_2$ as internal vertex. For these arcs,
$3\leq j< i-1< l+1$ and we have the minimal projective resolution
of $E_{i,j-1}$
$$
0\longrightarrow P_{y_{\frac{k}{2}}}\oplus P_{\sigma(a_0)}\oplus
P_{x_{j-2}} \stackrel{d^{E_{i,j-1}}_{min}}{\longrightarrow}
P_{\sigma(y_{\frac{k}{2}})}\oplus P_{a_0}\oplus
P_{x_{i-3}}\longrightarrow E_{i,j-1}\longrightarrow 0
$$
where $d^{E_{i,j-1}}_{min}=\left(\begin{array}{ccc} b
&\sigma(u_1)\cdots\sigma(u_{\frac{k}{2}})
 &0\\
u_{\frac{k}{2}}\cdots u_1 & 0 & v_{j-2}\cdots v_1\\
0 & \sigma(v_1)\cdots v_{i-2}&0
 \end{array}\right)$ and so
$$
c^{E_{i,j-1}}=det\left(\begin{array}{ccc} b &
u_{\frac{k}{2}}\cdots
u_1 &0\\
\sigma(u_1)\cdots\sigma(u_{\frac{k}{2}}) & 0 & \sigma(v_1)\cdots
v_{i-2}\\
0 & v_{j-2}\cdots v_1&0
 \end{array}\right).
$$
In particular we note that if $i=\sigma_I(j)$, in the orthogonal
case, we get
$$pf^{E_{\sigma_I(j),j-1}}=pf\left(\begin{array}{ccc} b &
u_{\frac{k}{2}}\cdots
u_1 &0\\
\sigma(u_1)\cdots\sigma(u_{\frac{k}{2}}) & 0 & \sigma(v_1)\cdots
\sigma(v_{j-2})\\
0 & v_{j-2}\cdots v_1&0
 \end{array}\right),$$
 since $b$ is skew-symmetric, $\sigma(v_i)=-(v_i)^t$ and $\sigma(u_i)=-(u_i)^t$. Finally we note
 that $V_{(0,1)}$, $V_{(1,0)}$, ${E_{i,\sigma_I(i)-1}}$ and ${E_{\sigma_I(j),j-1}}$
satisfy property \textit{(Spp)}. Similarly we define the
semi-invariants for the admissible arcs $[i,j]$ in
$\mathcal{A}'(d)$, exchanging the upper paths of
$\widetilde{A}_{k,l}^{2,0,2}$ with the lower ones.

\subsubsection{3.2.1.3\quad$\widetilde{A}_{k,l}^{0,2}$}\addcontentsline{toc}{subsection}{3.2.1.3\quad$\widetilde{A}_{k,l}^{0,2}$} We have at most two
$\tau^+$-orbits $\Delta$ and $\Delta'$ of the dimension vectors of
nonhomogeneous simple regular representation. We assume $n\geq 2$
and we consider the $\tau$-orbit $$\{e_1=\delta
e_1,e_2,\ldots,e_{[\frac{l-1}{2}]+1},e_{[\frac{l-1}{2}]+2}=\delta
e_{[\frac{l-1}{2}]+2} ,\delta e_{[\frac{l-1}{2}]+1},\ldots,\delta
e_2\}.$$ Let $[i,j]\in\mathcal{A}(d)$. If we consider the arc
$[1,1]$ of index 0, i.e. $p_1=0,p_2\ne
0,\ldots,p_{[\frac{l-1}{2}]+2}\ne 0$, we have the minimal
projective resolution of $V_{(0,1)}$
$$
0\longrightarrow
P_{\sigma(a_0)}\stackrel{d^{V_{(0,1)}}_{min}}{\longrightarrow}
P_{a_0}\longrightarrow V_{(0,1)}\longrightarrow 0
$$
where
$d^{V_{(0,1)}}_{min}=\sigma(v_1)\cdots\sigma(v_{\frac{l}{2}})v_{\frac{l}{2}}\cdots
v_1$ and so
$$
c^{V_{(0,1)}}=det(Hom_Q(d^{V_{(0,1)}}_{min},\cdot))=det(\sigma(v_1)\cdots\sigma(v_{\frac{l}{2}})v_{\frac{l}{2}}\cdots
v_1)
$$
in the orthogonal case and
$pf^{V_{(0,1)}}=pf(\sigma(v_1)\cdots\sigma(v_{\frac{l}{2}})v_{\frac{l}{2}}\cdots
v_1)$ in the symplectic case, since by definition of symplectic
representation\\
$\sigma(v_1)\cdots\sigma(v_{\frac{l}{2}})v_{\frac{l}{2}}\cdots
v_1$ is skew-symmetric. If we consider the arc
$[\sigma_I(2),2]=[0,2]$ of index 0, i.e.
$p_{\sigma_I(2)}=0=p_2,p_1\ne 0$, we have the minimal projective
resolution of $V_{(1,0)}$
$$
0\longrightarrow
P_{\sigma(a_0)}\stackrel{d^{V_{(1,0)}}_{min}}{\longrightarrow}
P_{a_0}\longrightarrow V_{(1,0)}\longrightarrow 0
$$
where
$d^{V_{(1,0)}}_{min}=\sigma(u_1)\cdots\sigma(u_{\frac{k}{2}})u_{\frac{k}{2}}\cdots
u_1$ and so
$$
c^{V_{(1,0)}}=det(Hom_Q(d^{V_{(1,0)}}_{min},\cdot))=
det(\sigma(u_1)\cdots\sigma(u_{\frac{k}{2}})u_{\frac{k}{2}}\cdots
u_1)
$$
in the orthogonal case and
$pf^{V_{(1,0)}}=pf(\sigma(u_1)\cdots\sigma(u_{\frac{k}{2}})u_{\frac{k}{2}}\cdots
u_1)$ in the symplectic case, since by definition of symplectic
representation\\
$\sigma(u_1)\cdots\sigma(u_{\frac{k}{2}})u_{\frac{k}{2}}\cdots
u_1$ is skew-symmetric. We note that for $l=2$ we have only the
admissible arcs $[1,1]$ an $[\sigma_I(2),2]$. We assume now that
$l\geq 4$ ($l$ is even) and $[i,j]$ is not an admissible arc
considered above. If $1\leq i<j\leq l+1$, then we identify $[i,j]$
with the path $v_{j-1}\cdots v_{i}$ in $Q$ and we have the minimal
projective resolution of $E_{i,j-1}$
$$
0\longrightarrow
P_{x_{j-1}}\stackrel{d^{E_{i,j-1}}_{min}}{\longrightarrow}
P_{x_{i-1}}\longrightarrow E_{i,j-1}\longrightarrow 0
$$
where $d^{E_{i,j-1}}_{min}=v_{j-1}\cdots v_{i}$ and so
$$
c^{E_{i,j-1}}=det(Hom_Q(d^{E_{i,j-1}}_{min},\cdot))=
det(v_{j-1}\cdots v_{i}).
$$
We note that $$c^{\tau^-\nabla
E_{i,j-1}}=c^{E_{\sigma_I(j),\sigma_I(i)-1}}=det(\sigma(v_{i})\cdots
\sigma(v_{j-1}))=det(v_{j-1}\cdots v_{i})=c^{E_{i,j-1}}.$$
Moreover, if $j=\sigma_I(i)$ then, only in the symplectic case, we
get $pf^{E_{i,\sigma_I(i)-1}}=pf(\sigma(v_{i})\cdots v_{i})$,
since $\sigma(v_{i})\cdots v_{i}$ is skew-symmetric. Now we
consider the arcs $[i,j]$ which have $e_1$ as internal vertex. For
these arcs, $2\leq j< i-1< l$ and we have the minimal projective
resolution of $E_{i,j-1}$
$$
0\longrightarrow P_{\sigma(a_0)}\oplus
P_{x_{j-1}}\stackrel{d^{E_{i,j-1}}_{min}}{\longrightarrow}
P_{a_0}\oplus P_{x_{i-1}}\longrightarrow E_{i,j-1}\longrightarrow
0
$$
where
$d^{E_{i,j-1}}_{min}=\left(\begin{array}{cc}\sigma(u_1)\cdots  u_1
&v_{j-1}\cdots v_1\\
 \sigma(v_1)\cdots v_{i}& 0 \end{array}\right)$ and so
$$
c^{E_{i,j-1}}=det(Hom_Q(d^{E_{i,j-1}}_{min},\cdot))=det\left(\begin{array}{cc}\sigma(u_1)\cdots
u_1
&\sigma(v_1)\cdots v_{i}\\
v_{j-1}\cdots v_1 & 0 \end{array}\right).
$$
In particular we note that if $i=\sigma_I(j)$, in the symplectic
case, we get
$$pf^{E_{\sigma_I(j),j-1}}=pf\left(\begin{array}{cc}\sigma(u_1)\cdots
 u_1
&\sigma(v_1)\cdots \sigma(v_{j-1})\\
v_{j-1}\cdots v_1 & 0
\end{array}\right),$$
since $\sigma(u_1)\cdots
u_1$ and $\sigma(v_i)=-(v_i)^t.$ Finally we note that $V_{(0,1)}$,
$V_{(1,0)}$, ${E_{i,\sigma_I(i)-1}}$ and ${E_{\sigma_I(j),j-1}}$
satisfy \textit{(Op)}. Similarly we define the semi-invariants for
the admissible arcs $[i,j]$ in $\mathcal{A}'(d)$, exchanging the
upper paths of $\widetilde{A}_{k,l}^{0,2}$ with the lower ones.

\subsubsection{3.2.1.4\quad$\widetilde{A}_{k,l}^{1,1}$}\addcontentsline{toc}{subsection}{3.2.1.4\quad$\widetilde{A}_{k,l}^{1,1}$}\label{A11} We have at
most two $\tau^+$-orbits $\Delta$ and $\Delta'$ of the dimension
vectors of nonhomogeneous simple regular representation. We assume
$n\geq 2$ and we consider the $\tau$-orbit $$\{e_1=\delta
e_1,e_2,\ldots,e_{[\frac{l-1}{2}]+1},e_{[\frac{l-1}{2}]+2}=\delta
e_{[\frac{l-1}{2}]+2} ,\delta e_{[\frac{l-1}{2}]+1},\ldots,\delta
e_2\}.$$ Let $[i,j]\in\mathcal{A}(d)$. If we consider the arc
$[1,1]$ of index 0, i.e. $p_1=0,p_2\ne
0,\ldots,p_{[\frac{l-1}{2}]+2}\ne 0$, we have the minimal
projective resolution of $V_{(0,1)}$
$$
0\longrightarrow
P_{\sigma(a_0)}\stackrel{d^{V_{(0,1)}}_{min}}{\longrightarrow}
P_{a_0}\longrightarrow V_{(0,1)}\longrightarrow 0
$$
where
$d^{V_{(0,1)}}_{min}=\sigma(v_1)\cdots\sigma(v_{\frac{l}{2}})v_{\frac{l}{2}}\cdots
v_1$ and so
$$
c^{V_{(0,1)}}=det(Hom_Q(d^{V_{(0,1)}}_{min},\cdot))=det(\sigma(v_1)\cdots\sigma(v_{\frac{l}{2}})v_{\frac{l}{2}}\cdots
v_1)
$$
in the orthogonal case and
$pf^{V_{(0,1)}}=pf(\sigma(v_1)\cdots\sigma(v_{\frac{l}{2}})v_{\frac{l}{2}}\cdots
v_1)$ in the symplectic case, since by definition of symplectic
representation\\
$\sigma(v_1)\cdots\sigma(v_{\frac{l}{2}})v_{\frac{l}{2}}\cdots
v_1$ is skew-symmetric . If we consider the arc
$[\sigma_I(2),2]=[0,2]$ of index 0, i.e.
$p_{\sigma_I(2)}=0=p_2,p_1\ne 0$, then we have the minimal
projective resolution of $V_{(1,0)}$
$$
0\longrightarrow
P_{\sigma(a_0)}\stackrel{d^{V_{(1,0)}}_{min}}{\longrightarrow}
P_{a_0}\longrightarrow V_{(1,0)}\longrightarrow 0
$$
where
$d^{V_{(1,0)}}_{min}=\sigma(u_1)\cdots\sigma(u_{\frac{k}{2}})bu_{\frac{k}{2}}\cdots
u_1$ and so
$$
c^{V_{(1,0)}}=det(Hom_Q(d^{V_{(1,0)}}_{min},\cdot))=
det(\sigma(u_1)\cdots\sigma(u_{\frac{k}{2}})bu_{\frac{k}{2}}\cdots
u_1)
$$
in the symplectic case and
$pf^{V_{(1,0)}}=pf(\sigma(u_1)\cdots\sigma(u_{\frac{k}{2}})bu_{\frac{k}{2}}\cdots
u_1)$ in the orthogonal case, since $b$ is skew-symmetric and
$\sigma(u_i)=-(u_i)^t$. We note that for $l=2$ we have only the
admissible arcs $[1,1]$ an $[\sigma_I(2),2]$. We assume now that
$l\geq 4$ ($l$ is even) and $[i,j]$ is not an admissible arc
considered above. If $1\leq i<j\leq l+1$, then we identify $[i,j]$
with the path $v_{j-1}\cdots v_{i}$ in $Q$ and we have the minimal
projective resolution of $E_{i,j-1}$
$$
0\longrightarrow
P_{x_{j-1}}\stackrel{d^{E_{i,j-1}}_{min}}{\longrightarrow}
P_{x_{i-1}}\longrightarrow E_{i,j-1}\longrightarrow 0
$$
where $d^{E_{i,j-1}}_{min}=v_{j-1}\cdots v_{i}$ and so
$$
c^{E_{i,j-1}}=det(Hom_Q(d^{E_{i,j-1}}_{min},\cdot))=
det(v_{j-1}\cdots v_{i}).
$$
We note that
$$c^{\tau^-\nabla
E_{i,j-1}}=c^{E_{\sigma_I(j),\sigma_I(i)-1}}=det(\sigma(v_{i})\cdots
\sigma(v_{j-1}))=det(v_{j-1}\cdots v_{i})=c^{E_{i,j-1}}.$$
Moreover, if $j=\sigma_I(i)$ then, only in the symplectic case, we
get $pf(\sigma(v_{i})\cdots v_{i})=pf^{E_{i,\sigma_I(i)-1}}$ since
$\sigma(v_{i})\cdots v_{i}$ is skew-symmetric. Now we consider the
arcs $[i,j]$ which have $e_1$ as internal vertex. For these arcs,
$2\leq j< i-1< l$ and we have the minimal projective resolution of
$E_{i,j-1}$
$$
0\longrightarrow P_{\sigma(a_0)}\oplus
P_{x_{j-1}}\stackrel{d^{E_{i,j-1}}_{min}}{\longrightarrow}
P_{a_0}\oplus P_{x_{i-1}}\longrightarrow E_{i,j-1}\longrightarrow
0
$$
where
$d^{E_{i,j-1}}_{min}=\left(\begin{array}{cc}\sigma(u_1)\cdots
b\cdots u_1
&\sigma(v_1)\cdots v_{i}\\
v_{j-1}\cdots v_1 & 0 \end{array}\right)$ and so
$$
c^{E_{i,j-1}}=det(Hom_Q(d^{E_{i,j-1}}_{min},\cdot))=det\left(\begin{array}{cc}\sigma(u_1)\cdots
b\cdots u_1
&\sigma(v_1)\cdots v_{i}\\
v_{j-1}\cdots v_1 & 0 \end{array}\right).
$$
In particular we note that if $i=\sigma_I(j)$, in the orthogonal
case, we get
$$pf^{E_{\sigma_I(j),j-1}}=pf\left(\begin{array}{cc}\sigma(u_1)\cdots
b\cdots u_1
&\sigma(v_1)\cdots \sigma(v_{j-1})\\
v_{j-1}\cdots v_1 & 0
\end{array}\right),$$ since $b$ is skew-symmetric, $\sigma(v_i)=-(v_i)^t$ and
 $\sigma(u_i)=-(u_i)^t$. Finally we note that
$V_{(0,1)}$, ${E_{i,\sigma_I(i)-1}}$ satisfy \textit{(Op)} and
$V_{(1,0)}$, ${E_{\sigma_I(j),j-1}}$ satisfy property
\textit{(Spp)}. Similarly we define the semi-invariants for the
admissible arcs $[i,j]$ in $\mathcal{A}'(d)$, exchanging the upper
paths of $\widetilde{A}_{k,l}^{1,1}$ with the lower ones and
tracing out the procedure done for $\widetilde{A}_{k,l}^{2,0,1}$.

\subsubsection{3.2.1.5\quad$\widetilde{A}_{k,k}^{0,0}$}\addcontentsline{toc}{subsection}{3.2.1.5\quad$\widetilde{A}_{k,k}^{0,0}$} We have at most two
$\tau^+$-orbits $\Delta$ and $\Delta'$ of the dimension vectors of
nonhomogeneous simple regular representation but in this case
$\Delta=\delta \Delta'$ so it's enough to study the
semi-invariants associated to the arcs in $\mathcal{A}(d)$,
because these are equal to those ones associated to the arcs in
$\mathcal{A}'(d)$. We assume $k\geq 2$ and we consider the
$\tau$-orbit $\{e_0, e_1,e_2,\ldots,e_{k-1}\}$. Let
$[i,j]\in\mathcal{A}(d)$. If we consider the arc $[1,1]$ of index
0, i.e. $p_1=0,p_2\ne 0,\ldots,p_{k-1}\ne 0$, we have the minimal
projective resolution of $V_{(0,1)}$
$$
0\longrightarrow
P_{\sigma(a_0)}\stackrel{d^{V_{(0,1)}}_{min}}{\longrightarrow}
P_{a_0}\longrightarrow V_{(0,1)}\longrightarrow 0
$$
where $d^{V_{(0,1)}}_{min}=v_k\cdots v_1$ and so
$$
c^{V_{(0,1)}}=det(Hom_Q(d^{V_{(0,1)}}_{min},\cdot))=det(v_k\cdots
v_1).
$$
If we consider the arc $[0,2]$ of index 0, i.e.
$p_{0}=0=p_2,p_1\ne 0$, then we have the minimal projective
resolution of $V_{(1,0)}$
$$
0\longrightarrow
P_{\sigma(a_0)}\stackrel{d^{V_{(1,0)}}_{min}}{\longrightarrow}
P_{a_0}\longrightarrow V_{(1,0)}\longrightarrow 0
$$
where $d^{V_{(0,1)}}_{min}=u_k\cdots u_1$ and so
$$
c^{V_{(1,0)}}=det(Hom_Q(d^{V_{(1,0)}}_{min},\cdot))=det(u_k\cdots
u_1).
$$
We note that for $k=2$ we have only the admissible arcs $[1,1]$ an
$[0,2]$. We assume now that $k\geq 3$ and $[i,j]$ is not an
admissible arc considered above. If $1\leq i<j\leq k$, then we
identify $[i,j]$ with the path $v_{j-1}\cdots v_{i}$ in $Q$ and we
have the minimal projective resolution of $E_{i,j-1}$
$$
0\longrightarrow
P_{x_{j-1}}\stackrel{d^{E_{i,j-1}}_{min}}{\longrightarrow}
P_{x_{i-1}}\longrightarrow E_{i,j-1}\longrightarrow 0
$$
where $d^{E_{i,j-1}}_{min}=v_{j-1}\cdots v_{i}$ and so
$$
c^{E_{i,j-1}}=det(Hom_Q(d^{E_{i,j-1}}_{min},\cdot))=
det(v_{j-1}\cdots v_{i}).
$$
Now we consider the arcs $[i,j]$ which have $e_1$ as internal
vertex. For these arcs, $2\leq j< i-1< k-1$ and we have the
minimal projective resolution of $E_{i,j-1}$
$$
0\longrightarrow P_{\sigma(a_0)}\oplus
P_{x_{j-1}}\stackrel{d^{E_{i,j-1}}_{min}}{\longrightarrow}
P_{a_0}\oplus P_{x_{i-1}}\longrightarrow E_{i,j-1}\longrightarrow
0
$$
where $d^{E_{i,j-1}}_{min}=\left(\begin{array}{cc}u_k\cdots u_1
&v_{j-1}\cdots v_1 \\
v_k\cdots v_{i}& 0 \end{array}\right)$ and so
$$
c^{E_{i,j-1}}=det(Hom_Q(d^{E_{i,j-1}}_{min},\cdot))=det\left(\begin{array}{cc}u_k\cdots
u_1
&v_k\cdots v_{i}\\
v_{j-1}\cdots v_1 & 0 \end{array}\right).
$$

\subsubsection{3.2.1.6\quad$\widetilde{D}_{n}^{1,0}$}\addcontentsline{toc}{subsection}{3.2.1.6\quad$\widetilde{D}_{n}^{1,0}$} In this case there are
three $\tau$-orbit $\Delta=\{e_1=\delta
e_1,e_2,\ldots,e_{n-1},\delta e_{n-1},\ldots,\delta e_2\}$,
$\Delta'=\{e'_0=\delta e'_0,e'_1=\delta e'_1\}$ and
$\Delta''=\{e''_0=\delta e''_1\}$. The only admissible arcs in
$\Delta'(d)$ and $\Delta''(d)$ are $[0,0]$ and $[1,1]$, recalling
that $e'_0+e'_1=h=e''_0+e''_1$. For such arcs in $\Delta'$ we have
the minimal projective resolution of $E'_{0,1}$
$$
0\longrightarrow P_{\sigma(t_1)}\oplus P_{\sigma(t_2)}
\stackrel{d^{E'_{0,1}}_{min}}{\longrightarrow}
 P_{t_1}\oplus P_{t_2}\longrightarrow E'_{0,1}\longrightarrow 0
$$
where $d^{E'_{0,1}}_{min}=\left(\begin{array}{cc}\sigma(a)\bar{c}a&0\\
\sigma(a)\bar{c}b&\sigma(b)\bar{c}b\end{array}\right)$, similarly
for $E'_{1,0}$ and so
$$
c^{E'_{1,0}}=c^{E'_{0,1}}=det(Hom_Q(d^{E'_{0,1}}_{min},\cdot))=
det\left(\begin{array}{cc}\sigma(a)\bar{c}a&\sigma(a)\bar{c}b\\
0&\sigma(b)\bar{c}b\end{array}\right).
$$
We note that the matrices $\sigma(a)\bar{c}b$, $\sigma(a)\bar{c}a$
and $\sigma(b)\bar{c}b)$ have different size for $[0,0]$ and for
$[1,1]$. Whereas in $\Delta''$ we have we have the minimal
projective resolution of $c^{E''_{0,1}}=c^{E''_{1,0}}$
$$
0\longrightarrow P_{\sigma(t_1)}\oplus P_{\sigma(t_2)}
\stackrel{d^{E''_{0,1}}_{min}}{\longrightarrow}
 P_{t_1}\oplus P_{t_2}\longrightarrow E''_{0,1}\longrightarrow 0
$$
where $d^{E''_{0,1}}_{min}=\left(\begin{array}{cc}0&\sigma(b)\bar{c}a\\
\sigma(a)\bar{c}b&\sigma(b)\bar{c}b\end{array}\right)$ and so
$$
c^{E''_{1,0}}=c^{E''_{0,1}}=det(Hom_Q(d^{E''_{0,1}}_{min},\cdot))=
det\left(\begin{array}{cc}0&\sigma(a)\bar{c}b\\
\sigma(b)\bar{c}a&\sigma(b)\bar{c}b\end{array}\right)
$$
in the symplectic case and
$$
pf^{E''_{0,1}}=pf^{E''_{1,0}}=pf\left(\begin{array}{cc}0&\sigma(a)\bar{c}b\\
\sigma(b)\bar{c}a&\sigma(b)\bar{c}b\end{array}\right)
$$
in the orthogonal case, since $\bar{c}$ is skew-symmetric,
$\sigma(b)=-b^t$ and $\sigma(a)=-a^t$. We assume $n\geq 3$ and we
take $[i,j]\in\mathcal{A}(d)$. If we consider the arc $[1,1]$, we
have the minimal projective resolution $V_{(1,1)}$
$$
0\longrightarrow P_{\sigma(t_1)}\oplus P_{\sigma(t_2)}
\stackrel{d^{V_{(1,1)}}_{min}}{\longrightarrow}
 P_{t_1}\oplus P_{t_2}\longrightarrow V_{(1,1)}\longrightarrow 0
$$
where $d^{V_{(1,1)}}_{min}=\left(\begin{array}{cc}\sigma(a)\bar{c}a&\sigma(b)\bar{c}a\\
\sigma(a)\bar{c}b&\sigma(b)\bar{c}b\end{array}\right)$ and so
$$
c^{V_{(1,1)}}=det(Hom_Q(d^{V_{(1,1)}}_{min},\cdot))=det\left(\begin{array}{cc}\sigma(a)\bar{c}a&\sigma(a)\bar{c}b\\
\sigma(b)\bar{c}a&\sigma(b)\bar{c}b\end{array}\right)
$$
in the symplectic case and
$$
pf^{V_{(1,1)}}=pf\left(\begin{array}{cc}\sigma(a)\bar{c}a&\sigma(a)\bar{c}b\\
\sigma(b)\bar{c}a&\sigma(b)\bar{c}b\end{array}\right)
$$
in the orthogonal case. If $[i,j]$ doesn't contain $e_1$ as an
internal vertex, then we have $1\leq i<j\leq 2n$ and we have the
minimal projective resolution of $E_{i,j-1}$
$$
0\longrightarrow
P_{z_{j-1}}\stackrel{d^{E_{i,j-1}}_{min}}{\longrightarrow}
P_{z_{i-1}}\longrightarrow E_{i,j-1}\longrightarrow 0
$$
where $d^{E_{i,j-1}}_{min}=c_{j-2}\cdots c_{i-1}$ and so
$$
c^{E_{i,j-1}}=det(Hom_Q(d^{E_{i,j-1}}_{min},\cdot))=
det(c_{j-2}\cdots c_{i-1}),
$$
where $c_0=(a,b)$ and $c_{2n-1}=\sigma(c_0)$. In particular in the
orthogonal case if $j=\sigma_I(i)$ then
$pf^{E_{i,\sigma_I(i)-1}}=pf(\sigma(c_{i-1})\cdots c_{i-1})$,
since in this case $\sigma(c_{i})=-(c_{i})^t$ and $c_{n-2}$ is
skew-symmetric. If $[i,j]$ contains $e_1$ as an internal vertex,
i.e. $2\leq j<i\leq 2n-1$ and we have the minimal projective
resolution of $E_{i,j-1}$
$$
0\longrightarrow  P_{z_{j-1}}\oplus P_{\sigma(t_1)}\oplus
P_{\sigma(t_2)}\stackrel{d^{E_{i,j-1}}_{min}}{\longrightarrow}
 P_{z_{i-1}}\oplus P_{t_1}\oplus P_{t_2}
\longrightarrow E_{i,j-1}\longrightarrow 0
$$
where
$d^{E_{i,j-1}}_{min}=\left(\begin{array}{ccc} 0&\sigma(a)c_{2n-3,i-1}&\sigma(b)c_{2n-3,i-1}\\
c_{j-2,1}a&\sigma(a)c_{2n-3,1}a&\sigma(b)c_{2n-3,1}a\\
c_{j-2,1}b&\sigma(a)c_{2n-3,1}b&\sigma(b)c_{2n-3,1}b\end{array}\right)$
and so
$$
c^{E_{i,j-1}}=\left(\begin{array}{ccc} 0&c_{j-2,1}a&c_{j-2,1}b\\
\sigma(a)c_{2n-3,i-1}&\sigma(a)c_{2n-3,1}a&\sigma(a)c_{2n-3,1}b\\
\sigma(b)c_{2n-3,i-1}&\sigma(b)c_{2n-3,1}a&\sigma(b)c_{2n-3,1}b\end{array}\right)
$$
where $c_{k,l}=c_k\cdots c_l$ and $c_{0,1}=id$. If $\sigma_I(i)=j$
then, only in the orthogonal case, we have
$$
pf^{E_{\sigma_I(j),j-1}}=pf\left(\begin{array}{ccc} 0&c_{j-2,1}a&c_{j-2,1}b\\
\sigma(a)\sigma(c_{j-2,1})&\sigma(a)c_{2n-3,1}a&\sigma(a)c_{2n-3,1}b\\
\sigma(b)\sigma(c_{j-2,1})&\sigma(b)c_{2n-3,1}a&\sigma(b)c_{2n-3,1}b\end{array}\right),
$$
since $\sigma(c_{j-2,1})=-(c_{j-2,1})^t$, $\sigma(a)=-a^t$, $\sigma(b)=-b^t$ and $c_{2n-3,1}$ is skew-symmetric.\\
Finally we note that $E''_{1,0}$, $V_{(1,1)}$,
${E_{i,\sigma_I(i)-1}}$ and ${E_{\sigma_I(j),j-1}}$ satisfy
property \textit{(Spp)}.

\subsubsection{3.2.1.7\quad$\widetilde{D}_{n}^{0,1}$}\addcontentsline{toc}{subsection}{3.2.1.7\quad$\widetilde{D}_{n}^{0,1}$} There are again three
$\tau$-orbit $\Delta=\{e_1=\delta e_1,e_2,\ldots,e_{n-1}=\delta
e_{n-1},\ldots,\delta e_2\}$, $\Delta'=\{e'_0=\delta
e'_0,e'_1=\delta e'_1\}$ and $\Delta''=\{e''_0=\delta e''_1\}$.
The only admissible arcs in $\Delta'(d)$ and $\Delta''(d)$ are
$[0,0]$ and $[1,1]$, recalling that $e'_0+e'_1=h=e''_0+e''_1$. For
such arcs in $\Delta'$ we have the minimal projective resolution
of $E'_{0,1}$
$$
0\longrightarrow P_{\sigma(t_1)}\oplus P_{\sigma(t_2)}
\stackrel{d^{E'_{0,1}}_{min}}{\longrightarrow}
 P_{t_1}\oplus P_{t_2}\longrightarrow E'_{0,1}\longrightarrow 0
$$
where $d^{E'_{0,1}}_{min}=\left(\begin{array}{cc}\sigma(a)\bar{c}a&0\\
\sigma(a)\bar{c}b&\sigma(b)\bar{c}b\end{array}\right)$, similarly
for $E'_{1,0}$ and so
$$
c^{E'_{1,0}}=c^{E'_{0,1}}=det(Hom_Q(d^{E'_{0,1}}_{min},\cdot))=
det\left(\begin{array}{cc}\sigma(a)\bar{c}a&\sigma(a)\bar{c}b\\
0&\sigma(b)\bar{c}b\end{array}\right).
$$
We note that the matrices $\sigma(a)\bar{c}b$, $\sigma(a)\bar{c}a$
and $\sigma(b)\bar{c}b)$ have different size for $[0,0]$ and for
$[1,1]$. Whereas in $\Delta''$ we have the minimal projective
resolution of $c^{E''_{0,1}}=c^{E''_{1,0}}$
$$
0\longrightarrow P_{\sigma(t_1)}\oplus P_{\sigma(t_2)}
\stackrel{d^{E''_{0,1}}_{min}}{\longrightarrow}
 P_{t_1}\oplus P_{t_2}\longrightarrow E''_{0,1}\longrightarrow 0
$$
where $d^{E''_{0,1}}_{min}=\left(\begin{array}{cc}0&\sigma(b)\bar{c}a\\
\sigma(a)\bar{c}b&\sigma(b)\bar{c}b\end{array}\right)$ and so
$$
c^{E''_{1,0}}=c^{E''_{0,1}}=det(Hom_Q(d^{E''_{0,1}}_{min},\cdot))=
det\left(\begin{array}{cc}0&\sigma(a)\bar{c}b\\
\sigma(b)\bar{c}a&\sigma(b)\bar{c}b\end{array}\right)
$$
in the orthogonal case and
$$
pf^{E''_{0,1}}=pf^{E''_{1,0}}=pf\left(\begin{array}{cc}0&\sigma(a)\bar{c}b\\
\sigma(b)\bar{c}a&\sigma(b)\bar{c}b\end{array}\right)
$$
in the symplectic case, since $\bar{c}$ is skew-symmetric,
$\sigma(b)=-b^t$ and $\sigma(a)=-a^t$. We assume $n\geq 3$ and we
take $[i,j]\in\mathcal{A}(d)$. If we consider the arc $[1,1]$, we
have the minimal projective resolution $V_{(1,1)}$
$$
0\longrightarrow P_{\sigma(t_1)}\oplus P_{\sigma(t_2)}
\stackrel{d^{V_{(1,1)}}_{min}}{\longrightarrow}
 P_{t_1}\oplus P_{t_2}\longrightarrow V_{(1,1)}\longrightarrow 0
$$
where $d^{V_{(1,1)}}_{min}=\left(\begin{array}{cc}\sigma(a)\bar{c}a&\sigma(b)\bar{c}a\\
\sigma(a)\bar{c}b&\sigma(b)\bar{c}b\end{array}\right)$ and so
$$
c^{V_{(1,1)}}=det(Hom_Q(d^{V_{(1,1)}}_{min},\cdot))=det\left(\begin{array}{cc}\sigma(a)\bar{c}a&\sigma(a)\bar{c}b\\
\sigma(b)\bar{c}a&\sigma(b)\bar{c}b\end{array}\right)
$$
in the orthogonal case and
$$
pf^{V_{(1,1)}}=pf\left(\begin{array}{cc}\sigma(a)\bar{c}a&\sigma(a)\bar{c}b\\
\sigma(b)\bar{c}a&\sigma(b)\bar{c}b\end{array}\right)
$$
in the symplectic case. If $[i,j]$ doesn't contain $e_1$ as an
internal vertex, then we have $1\leq i<j\leq 2n-3$ and we have the
minimal projective resolution of $E_{i,j-1}$
$$
0\longrightarrow
P_{z_{j-1}}\stackrel{d^{E_{i,j-1}}_{min}}{\longrightarrow}
P_{z_{i-1}}\longrightarrow E_{i,j-1}\longrightarrow 0
$$
where $d^{E_{i,j-1}}_{min}=c_{j-2}\cdots c_{i-1}$ and so
$$
c^{E_{i,j-1}}=det(Hom_Q(d^{E_{i,j-1}}_{min},\cdot))=
det(c_{j-2}\cdots c_{i-1}),
$$
where $c_0=(a,b)$ and $c_{2n-4}=\sigma(c_0)$. In particular in the
symplectic case if $j=\sigma_I(i)$ then
$pf^{E_{i,\sigma_I(i)-1}}=pf(\sigma(c_{i-1})\cdots c_{i-1})$. If
$[i,j]$ contains $e_1$ as an internal vertex, i.e. $2\leq j<i\leq
2n-4$ and we have the minimal projective resolution of $E_{i,j-1}$
$$
0\longrightarrow  P_{z_{j-1}}\oplus P_{\sigma(t_1)}\oplus
P_{\sigma(t_2)}\stackrel{d^{E_{i,j-1}}_{min}}{\longrightarrow}
 P_{z_{i-1}}\oplus P_{t_1}\oplus P_{t_2}
\longrightarrow E_{i,j-1}\longrightarrow 0
$$
where
$d^{E_{i,j-1}}_{min}=\left(\begin{array}{ccc} 0&\sigma(a)c_{2n-6,i-1}&\sigma(b)c_{2n-6,i-1}\\
c_{j-2,1}a&\sigma(a)c_{2n-6,1}a&\sigma(b)c_{2n-6,1}a\\
c_{j-2,1}b&\sigma(a)c_{2n-6,1}b&\sigma(b)c_{2n-6,1}b\end{array}\right)$
and so
$$
c^{E_{i,j-1}}=\left(\begin{array}{ccc} 0&c_{j-2,1}a&c_{j-2,1}b\\
\sigma(a)c_{2n-6,i-1}&\sigma(a)c_{2n-6,1}a&\sigma(a)c_{2n-6,1}b\\
\sigma(b)c_{2n-6,i-1}&\sigma(b)c_{2n-6,1}a&\sigma(b)c_{2n-6,1}b\end{array}\right).
$$
If $\sigma_I(i)=j$ then, only in the symplectic case, we have
$$
pf^{E_{\sigma_I(j),j-1}}=pf\left(\begin{array}{ccc} 0&c_{j-2,1}a&c_{j-2,1}b\\
\sigma(a)\sigma(c_{j-2,1})&\sigma(a)c_{2n-6,1}a&\sigma(a)c_{2n-6,1}b\\
\sigma(b)\sigma(c_{j-2,1})&\sigma(b)c_{2n-6,1}a&\sigma(b)c_{2n-6,1}b\end{array}\right),
$$
since $\sigma(c_{j-2,1})=-(c_{j-2,1})^t$, $\sigma(a)=-a^t$, $\sigma(b)=-b^t$ and $c_{2n-6,1}$ is skew-symmetric.\\
Finally we note that $E''_{1,0}$, $V_{(1,1)}$,
${E_{i,\sigma_I(i)-1}}$ and ${E_{\sigma_I(j),j-1}}$ satisfy
property \textit{(Op)}.

\subsubsection{3.2.1.8\quad End of proof of theorem \ref{stq}, theorem
\ref{dtso} and proposition \ref{fij}} We prove the second step of
proof of theorem \ref{stq}. By the analysis case by case we note
that if $[i,j]$ is admissible then the semi-invariants associated
to $[i,j]$ define a nonzero element of $SpSI(Q,d)$ (respectively
of $OSI(Q,d)$).\\
For a symmetric dimension vector $d$ we denote
\begin{equation}\label{SpGamma}
Sp\Gamma(Q,d)=\{\chi\in\mathbb{Z}^{Q_0}\cup\frac{1}{2}\mathbb{Z}^{Q_0}|\,SpSI(Q,d)_{\chi}\neq
0\}
\end{equation}
and
\begin{equation}\label{OGamma}
O\Gamma(Q,d)=\{\chi\in\mathbb{Z}^{Q_0}\cup\frac{1}{2}\mathbb{Z}^{Q_0}|\,OSI(Q,d)_{\chi}\neq
0\}
\end{equation}
the semigroup of weights of symplectic (respectively orthogonal)
semi-invariants. We note that (\ref{SpGamma}) and (\ref{OGamma})
involve also $\frac{1}{2}\mathbb{Z}^{Q_0}$ because in $SpSI(Q,d)$
and in $OSI(Q,d)$ also pfaffians can appear. To simplify the
notation, we shall call $\chi_{[i,j]}$, $\chi'_{[i,j]}$ and
$\chi''_{[i,j]}$ be respectively the weights of the
semi-invariants associated to admissible arcs $[i,j]$ respectively
from $\mathcal{A}(d)$, $\mathcal{A}'(d)$ and $\mathcal{A}''(d)$.
In the next the following proposition will be useful. We will
state it only for $\Delta$, because for $\Delta'$ and $\Delta''$
the statements are similar. Let $d$ be a regular symmetric
dimension vector with canonical decomposition $d=ph+d'$ with
$p\geq 1$.
\begin{proposizione}\label{pe}
Let $(Q,\sigma)$ be a symmetric quiver of tame type. Let $d_2$ be
of type $e_{[s,\sigma_I(s)]}$, $e_{[s,t]}+\delta e_{[s,t]}$ or
$e_{[i_{2k},\sigma_I(i_{2k-1})]}+e_{[i_{2k-1},\sigma_I(i_{2k})]}$.\\
(i) If $d_2=e_{[s,\sigma_I(s)]}$, then
\begin{itemize}
\item[(a)] For every arc $[i,j]$ of $\Delta'$ and $\Delta''$ we
have $\chi'_{[i,j]}|_{supp(d_2)},\chi''_{[i,j]}|_{supp(d_2)}\in
Sp\Gamma(Q,d_2)$ (respectively in $O\Gamma(Q,d_2)$).
\item[(b)] For every arc $[i,j]$ of $\Delta$ that doesn't intersect
 $[s,\sigma_I(s)]$ or contains $[s-1,\sigma_I(s)+1]$ we
have $\chi_{[i,j]}|_{supp(d_2)}\in Sp\Gamma(Q,d_2)$ (respectively
in $O\Gamma(Q,d_2)$).
\item[(c)] Let $\rho_1,\ldots,\rho_r$ be the weights of generators of the polynomial
algebra $SpSI(Q,d_2)$ (respectively $OSI(Q,d_2)$). Then $r\geq
n'-s$, where $n'\in I_+\sqcup I_{\delta}$ is either a
$\sigma_I$-fixed vertex or the extremal vertex of a
$\sigma_I$-fixed edge, and $\rho_1,\ldots,\rho_r$ can be reordered
such that
$\rho_1=\chi_{[s,s+1]},\ldots,\rho_{n'-s}=\chi_{[n'-1,n']}$ and
for every $m>n'-s$ we have $\langle\rho_m,e_n\rangle=0$ for
$n=s,\ldots,n'$.
\end{itemize}
(ii) Let $d_2=e_{[s,t]}+\delta e_{[s,t]}$, then
\begin{itemize}
\item[(a)] For every arc $[i,j]$ of $\Delta'$ and $\Delta''$ we
have $\chi'_{[i,j]}|_{supp(d_2)},\chi''_{[i,j]}|_{supp(d_2)}\in
Sp\Gamma(Q,d_2)$ (respectively in $O\Gamma(Q,d_2)$).
\item[(b)] For every symmetric arc $[i,j]$ of $\Delta$ that doesn't intersect
 $[s,t]\cup[\sigma_I(t),\sigma_I(s)]$ or contains $[s-1,\sigma_I(s-1)]$ or $[\sigma_I(t+1),t+1]$, we have
$\chi_{[i,j]}|_{supp(d_2)}\in Sp\Gamma(Q,d_2)$ (respectively in
$O\Gamma(Q,d_2)$).
\item[(c)] For every arc $[i,j]\subset I_+$ (respectively
$[i,j]\subset I_-$) that doesn't intersect $[s,t]$ (respectively
$[\sigma_I(t),\sigma_I(s)]$ or contains $[s-1,t+1]$ we have
$\chi_{[i,j]}|_{supp(d_2)}\in Sp\Gamma(Q,d_2)$ (respectively in
$O\Gamma(Q,d_2)$)..
\item[(d)] Let $\rho_1,\ldots,\rho_r$ be the weights of generators of the polynomial
algebra $SpSI(Q,d_2)$ (respectively $OSI(Q,d_2)$). Then $r\geq
t-s$ and $\rho_1,\ldots,\rho_r$ can be reordered such that
$\rho_1=\chi_{[s,s+1]},\ldots,\rho_{t-s}=\chi_{[t-1,t]}$ and for
every $m>t-s$ we have $\langle\rho_m,e_n\rangle=0$ for
$n=s,\ldots,t$.
\end{itemize}
(iii) Let
$d_2=e_{[i_{2k},i_{\sigma_I(i_{2k-1})}]}+e_{[i_{2k-1},i_{\sigma_I(i_{2k})}]}$,
then
\begin{itemize}
\item[(a)] For every arc $[i,j]$ of $\Delta'$ and $\Delta''$ we
have $\chi'_{[i,j]}|_{supp(d_2)},\chi''_{[i,j]}|_{supp(d_2)}\in
Sp\Gamma(Q,d_2)$ (respectively in $O\Gamma(Q,d_2)$).
\item[(b)] For every arc $[i,j]$ of $\Delta$ that doesn't intersect
 $[i_{2k-1},\sigma_I(i_{2k-1})]$ or contains $[i_{2k-1}-1,\sigma_I(i_{2k-1})+1]$ we
have $\chi_{[i,j]}|_{supp(d_2)}\in Sp\Gamma(Q,d_2)$ (respectively
in $O\Gamma(Q,d_2)$).
\item[(c)] Let $\rho_1,\ldots,\rho_r$ be the weights of generators of the polynomial
algebra $SpSI(Q,d_2)$ (respectively $OSI(Q,d_2)$). Then $r\geq
n'-s$, where $n'\in I_+\sqcup I_{\delta}$ is either a
$\sigma_I$-fixed vertex or the extremal vertex of a
$\sigma_I$-fixed edge, and $\rho_1,\ldots,\rho_r$ can be reordered
such that
$\rho_1=\chi_{[s,s+1]},\ldots,\rho_{n'-s}=\chi_{[n'-1,n']}$ and
for every $m>n'-s$ we have $\langle\rho_m,e_n\rangle=0$ for
$n=s,\ldots,n'$.
\end{itemize}
\end{proposizione}
\textit{Proof.} It proceeds type by type analysis, considering the
description of the weights of symplectic and orthogonal
semi-invariants done above. We recall that
$\gamma\chi_{[i,j]}=\chi_{[\sigma_I(j),\sigma_I(i)]}$ and we
observe that if $x$ is a $\sigma$-fixed vertex and $\chi$ is a
weight, then $\chi(x)=0$. We prove only the symplectic case for
$Q=\widetilde{A}^{1,1}_{k,l}$ and for $d_2=e_{[s,\sigma_I(s)]}$,
because the procedure to prove all other cases is similar. We
order the vertices of $\widetilde{A}^{1,1}_{k,l}$ such that the
only source is 1 (so the only sink is $\sigma(1)$), $hv_{i-1}=i$
for every $i\in\{2,\ldots,\frac{l}{2}+1\}$, $hu_i=\frac{l}{2}+i+1$
for every $i\in\{1,\ldots,\frac{k}{2}\}$ and then the respective
conjugates by $\sigma$ of these. We shall call
$w_{(t^1)_{i_1},\ldots,(t^f)_{i_f}}$, where
$t^1,\ldots,t^f\in\mathbb{Z}\cup\frac{1}{2}\mathbb{Z}$ and
$\{i_1,\ldots,i_f\}$ is an ordered subset of
$\{1,\ldots,\frac{l}{2}+\frac{k}{2}+1,\sigma(\frac{l}{2}+\frac{k}{2}+1),\ldots,\sigma(1)\}$,
the vector such that
$$
w_{(t^1)_{i_1},\ldots,(t^f)_{i_f}}(y)=\left\{\begin{array}{ll}(t^j)_{i_j}
& y=i_j,\forall j=1,\ldots,f\\
0&\textrm{otherwise}.\end{array}\right.
$$
Moreover we can associate in bijective way the vertex
$i\in\{2,\ldots,\frac{l}{2}\}\subset(\widetilde{A}^{1,1}_{k,l})_0^+$
to $i\in I_+$, the vertex $\frac{l}{2}+i+1$ of
$\widetilde{A}^{1,1}_{k,l}$ to $i+1\in I'_+$ and
the vertex $\frac{l}{2}$ to $[\frac{l-1}{2}]+2\in I_{\delta}$.\\
\textit{(a)} By section 3.2.1.4 we have
$$
\chi'_{[i,j]}=w_{(1)_{\frac{l}{2}+i+1},(-1)_{\frac{l}{2}+j+1}}\quad\textrm{for}\quad
1\leq i<j\leq \frac{k}{2}+1,
$$
if $[i,j]$ has not $e_1$ as internal vertex;
$$
\chi'_{[i,j]}=w_{(1)_1,(-1)_{\frac{l}{2}+j+1},(1)_{\frac{l}{2}+i+1},(-1)_{\sigma(1)}}\quad\textrm{for}\quad
j<i-1
$$
if $[i,j]$ has $e_1$ as internal vertex and in particular if
$j=\sigma_I(i)$ we have
$$
\chi'_{[i,j]}=w_{(\frac{1}{2})_1,(-\frac{1}{2})_{\frac{l}{2}+i+1},(\frac{1}{2})_{\sigma(\frac{l}{2}+i+1)},(-\frac{1}{2})_{\sigma(1)}}.
$$
Now if $\langle\chi'_{[i,j]},e_{[s,\sigma_I(s)]}\rangle\neq 0$
then $\chi'_{[i,j]}\not\in SpSI(Q,d_2)$, but we note that
$\langle\chi'_{[i,j]},e_{[s,\sigma_I(s)]}\rangle=0$ for every $i$
and $j$, so we have \textit{(a)}.\\
\textit{(b)} By section 3.2.1.4 we have
$$
\chi_{[i,j]}=w_{(1)_{i},(-1)_{j}}\quad\textrm{for}\quad 1\leq
i<j\leq
\frac{l}{2}\quad\textrm{and}\quad\chi_{[\frac{l}{2}+1,\sigma(\frac{l}{2}+1)]}=w_{(\frac{1}{2})_{\frac{l}{2}+1},(-\frac{1}{2})_{\sigma(\frac{l}{2}+1)}}
$$
if $[i,j]$ has not $e_1$ as internal vertex;
$$
\chi_{[i,j]}=w_{(1)_1,(-1)_{j},(1)_{i},(-1)_{\sigma(1)}}\quad\textrm{for}\quad
j<i-1.
$$
if $[i,j]$ has $e_1$ as internal vertex.\\
Now we note that
$\langle\chi_{[i,j]},e_{[s,\sigma_I(s)]}\rangle\neq 0$ if
$[i,j]\cap[s,t]\neq\emptyset$ and
$[i,j]\nsupseteq[s-1,\sigma(s-1)=\sigma(s)+1]$, so we have
\textit{(b)}.\\
\textit{(c)} First we note that we can choose symmetric arcs of
each length from a fixed vertex of $\Delta$, because the result of
theorem \ref{stq} is invariant respect to the Coxeter
transformation $\tau^+$. We note that $[s,\sigma_I(s)]$ has $e_1$
as internal vector. The generators of $SpSI(Q,d_2)$ associated to
$\Delta(d_2)$ are $c^{E_i}=det(v_i)$ of weight
$\chi_{[i,i+1]}=w_{(1)_{i},(-1)_{i+1}}$ for every
$i\in\{1,\ldots,s-1\}$ and
$c^{E_{s,\sigma_I(s)-1}}=det\left(\begin{array}{cc}\sigma(u_1)\cdots
b\cdots u_1
&\sigma(v_1)\cdots \sigma(v_{s})\\
v_{s-1}\cdots v_1 & 0 \end{array}\right)$ of weight
$\chi_{[s,\sigma_I(s)]}=w_{(1)_{1},(-1)_{s},(1)_{\sigma(s)},(-1)_{\sigma(1)}}$.
So we call $\rho_i=\chi_{[i,i+1]}$  for every
$i\in\{1,\ldots,s-1\}$ and $\rho_{n'-s}=\chi_{[s,\sigma_I(s)]}$,
where in this case $n'=\frac{[l-1]}{2}+2$. The other generators
are associated to $\Delta'(d_2)$ and so, as done in the part
\textit{(a)} of this proposition, their weight $\rho_m$, for
$m\in\{n'-s+1,\ldots,r\}$, are such that
$\langle\rho_m,e_n\rangle=0$ for
$n\in\{s,\ldots,n'\}$. $\Box$ \\
\\
We assume now that $d=d_1+d_2$ where $d_1=ph+d'_1$ with $p\geq 1$
and $d_2=e_{[s,\sigma_I(s)]}$, $e_{[s,t]}+\delta e_{s,t]}$ or
$e_{[i_{2k},i_{\sigma_I(i_{2k-1})}]}+e_{[i_{2k-1},i_{\sigma_I(i_{2k})}]}$..
So we take the corresponding arc in a chosen position (for which
we proved proposition \ref{pe}).
\begin{proposizione}\label{d1d2}
Let $d,d_1,d_2$ be as above. We suppose that the semigroup
$Sp\Gamma(Q,d_1)$ (respectively $O\Gamma(Q,d_1)$) is generated by
the weights $\chi_{[i,j]}$, $\chi'_{[i,j]}$, $\chi''_{[i,j]}$ for
admissible arcs $[i,j]$ of the labelled polygons $\Delta(d_1)$,
$\Delta'(d_1)$, $\Delta''(d_1)$. Then $Sp\Gamma(Q,d_1)\cap
Sp\Gamma(Q,d_2)$ (respectively $O\Gamma(Q,d_1)\cap
O\Gamma(Q,d_2)$) is generated by the weights $\chi_{[i,j]}$,
$\chi'_{[i,j]}$, $\chi''_{[i,j]}$ for admissible arcs $[i,j]$ of
the labelled polygons $\Delta(d)$, $\Delta'(d)$, $\Delta''(d)$.
\end{proposizione}
\textit{Proof.} We prove it only for the othogonal case and for $d_2=e_{[s,\sigma_I(s)]}$, because the symplectic case is similar.\\
We are two cases.\\
(1) Assume $p_{s-1}=p_{\sigma_I(s)+1}<r-1$. The admissible arcs of
$\Delta(d_1)$, $\Delta'(d_1)$, $\Delta''(d_1)$ and $\Delta(d)$,
$\Delta'(d)$, $\Delta''(d)$ are the same. By proposition \ref{pe}
$O\Gamma(Q,d_2)$ contains
$\chi_{[s,s+1]},\ldots,\chi_{[\sigma_I(s)-1,\sigma_I(s)]}$ and all
the other weights corresponding to the admissible arcs of
$\Delta(d)$, $\Delta'(d)$ and $\Delta''(d)$.\\
(2) Assume $p_{s-1}=p_{\sigma_I(s)+1}=r-1$. We prove that
$O\Gamma(Q,d_1)\cap O\Gamma(Q, d_2)$ is generated by
$\chi'_{[i,j]}$ for every admissible arc $[i,j]$ of
$\Delta'(d_1)=\Delta'(d)$, $\chi''_{[i,j]}$ for every admissible
arc $[i,j]$ of $\Delta''(d_1)=\Delta''(d)$ and $\chi_{[i,j]}$ for
every admissible arc $[i,j]$ of $\Delta(d_1)$ of index smaller
than $r-1$ or not intersecting $[s,\sigma_I(s)]$, i.e.
$\chi_{[s,s+1]},\ldots,\chi_{[\sigma_I(s)-1,\sigma_I(s)]}$ and
$\chi_{[s-1,\sigma_I(s)+1]}=\chi_{[s-1,s]}+\cdots+\chi_{[\sigma_I(s),\sigma_I(s)+1]}$.
Let
$$
\chi=\sum_{[i,j]\in\mathcal{A}(d_1)}n_{i,j}\chi_{[i,j]}+\sum_{[i,j]\in\mathcal{A}'(d_1)}n'_{i,j}\chi'_{[i,j]}
+\sum_{[i,j]\in\mathcal{A}''(d_1)}n''_{i,j}\chi''_{[i,j]},
$$
with $n_{i,j},n'_{i,j},n''_{i,j}\geq 0$, be an element of
$O\Gamma(Q,d_1)$. We assume that $\chi$ is also in
$O\Gamma(Q,d_2)$. By proposition \ref{pe}, we note that all the
generators of $O\Gamma(Q,d_1)$ except of $\chi_{[s-1,s]}$ and
$\chi_{[\sigma_I(s),\sigma_I(s)+1]}$ are also in $O\Gamma(Q,d_2)$.
Hence, if $\chi$ contains neither $\chi_{[s-1,s]}$ nor
$\chi_{[\sigma_I(s),\sigma_I(s)+1]}$, then $\chi$ is a linear
combination of desired generators. So we have to prove that if
$\chi$ contains $\chi_{[s-1,s]}$ (resp.
$\chi_{[\sigma_I(s),\sigma_I(s)+1]}$) with positive coefficient,
then it contains
$\chi_{[s,s+1]},\ldots,\chi_{[\sigma_I(s),\sigma_I(s)+1]}$ (resp.
$\chi_{[s-1,s]},\ldots,\chi_{[\sigma_I(s)-1,\sigma_I(s)]}$). Thus
we can subtract $\chi_{[s-1,\sigma_I(s)+1]}$ from
$\chi$.\\
We assume that $\chi$ contains $\chi_{[s-1,s]}$ with positive
coefficient (the proof is similar for
$\chi_{[\sigma_I(s),\sigma_I(s)+1]}$). We note that
$\langle\chi_{[s-1,s]},e_s\rangle=-1$ and, by proposition
\ref{pe}, the other generators of $O\Gamma(Q,d_1)$, except
$\chi_{[s,s+1]}$, have zero product scalar with $e_s$. Moreover,
$\chi\in O\Gamma(Q,d_2)$ and so, by proposition \ref{pe},
$\langle\chi,e_s\rangle\geq 0$. Hence $\chi$ contains
$\chi_{[s,s+1]}$ with positive coefficient. By proposition
\ref{pe}, it follows that $\langle\chi,e_s+e_{s+1}\rangle\geq 0$.
But $\langle\chi_{[s-1,s]}+\chi_{[s,s+1]},e_s+e_{s+1}\rangle=-1$
and $\chi_{[s+1,s+2]}$ is the only generator of $O\Gamma(Q,d_1)$
with positive scalar product with $e_s+e_{s+1}$. Continuing in
this way, we check that $\chi$ contains
$\chi_{[s-1,s]},\chi_{[s,s+1]},\ldots,\chi_{[\sigma_I(s)-1,\sigma_I(s)]},\chi_{[\sigma_I(s),\sigma_I(s)+1]}$
with positive coefficients. So we can subtract
$\chi_{[s-1.\sigma_I(s)+1]}$ from $\chi$ and continue. In this way
we complete the proof. $\Box$\\
\\
Now we can finish the proof of theorem \ref{stq}. Since  theorem
\ref{stq} is equivalent to conjectures \ref{mt1} and \ref{mt2} for
tame type and regular dimension vectors, then, in this way, we
finish also the proof of conjectures \ref{mt1} and \ref{mt2}.\\
Again we consider the embeddings
\begin{equation}\label{Spe1}
SpSI(Q,d)\stackrel{\Phi_d}{\rightarrow}\bigoplus_{\chi\in
char(Sp(Q,d))}SpSI(Q,d_1)_{\chi|_{d_1}}\otimes
SpSI(Q,d_2)_{\chi|_{d_2}}
\end{equation}
and
\begin{equation}\label{Oe1}
OSI(Q,d)\stackrel{\Psi_d}{\rightarrow}\bigoplus_{\chi\in
char(O(Q,d))}OSI(Q,d_1)_{\chi|_{d_1}}\otimes
OSI(Q,d_2)_{\chi|_{d_2}}
\end{equation}
where $Q$, $d$, $d_1$ and $d_2$ are as above. The semigroup of
weights of the right hand side of $\Phi_d$ and $\Psi_d$ are
respectively $Sp\Gamma(Q,d_1)\cap Sp\Gamma(Q,d_2)$ and
$O\Gamma(Q,d_1)\cap O\Gamma(Q,d_2)$. These are generated by
$\chi_{[i,j]}$, $\chi'_{[i,j]}$, $\chi''_{[i,j]}$ for admissible
arcs $[i,j]$ of the labelled polygons $\Delta(d)$, $\Delta'(d)$,
$\Delta''(d)$, by proposition \ref{d1d2}. So the algebras on the
right hand side of $\Phi_d$ and $\Psi_d$ are generated by the
semi-invariants of weights $\chi_{[i,j]}$, $\chi'_{[i,j]}$,
$\chi''_{[i,j]}$ and by the semi-invariants of weights $\langle
h,\cdot\rangle$ (or $\frac{1}{2}\langle h,\cdot\rangle$).\\
Finally, we note that the embeddings $\Phi_d$ and $\Psi_d$ are
isomorphisms because they are also isomorphisms in the weight
$\langle h,\cdot\rangle$ (or $\frac{1}{2}\langle h,\cdot\rangle$)
and so we completed the proof of theorem \ref{stq}. Moreover, in
that way, we also proved proposition \ref{fij}, expliciting the
semi-invariants of type $c^V$ for every admissible arc $[i,j]$,
and theorem \ref{dtso}, by isomorphisms $\Phi_d$ and $\Psi_d$
considering $d_1=ph$ and $d_2=d'$.

\appendix
\chapter{Representations of $GL$ and invariant theory}
\section{Highest weight theory for $GL$ and Schur modules}\label{A1}
    We recall the basics of representation theory of general
    linear group.\\
    We fix an algebraically closed field $\mathbb{K}$.
\begin{definizione}
Let $G$ be an algebraic group. $(V,\rho)$ is a rational
representation if $V$ is a vector space of dimension $m$,
$\rho:G\times V\longrightarrow V$ such that $\rho(g,v)=g\cdot v$
is a rational action, i.e.
\begin{itemize}
\item[a)] $g\cdot(h\cdot v)=(gh)\cdot v$ for every $g,h\in G$ and $v\in V$,
\item[b)] $e\cdot v=v$ for every $v\in V$ where $e$ is identity in $G$,
\item[c)] $\rho$ is a morphism of varieties.
\end{itemize}
\end{definizione}
\begin{definizione}
 $G$ is linearly reductive if and only if every rational linear representation of $G$ is semisimple.
\end{definizione}
Let $G$ be a linearly reductive group and let $\rho:G\rightarrow
GL(V)$ be a finite dimensional rational representation of $G$. Let
$H$ be a \textit{maximal torus of $G$}, i.e. a maximal subgroup of
$G$ isomorphic to $(\mathbb{K}^*)^h$ for some $h\in\mathbb{N}$,
restricting $\rho$ to $H$ we obtain a rational representation of
$H$. So we can decompose $V$ into the direct sum of eigenspaces
$$
V=\bigoplus_{\chi\in char(H)}V_{\chi}
$$
where $char(H)=\{\textrm{homomorphisms of algebraic groups}\;
\chi:H\rightarrow\mathbb{K}^*\}$ is the set of characters of $H$
and $V_{\chi}=\{v\in V|\rho(t)(v)=\chi(t)v,\;\forall t\in H\}$.
The elements $\chi\in char(H)$ such that $V_{\chi}\ne 0$ are
called weights of $\rho$, $V_{\chi}$ is called \textit{weight
space of weight $\chi$} and $dim\,V_{\chi}$ is called
\textit{multiplicity of the weight $\chi$}. The set of weights
$char(G)$ forms a free abelian group $\mathscr{X}=char(G)$. Let
$\Phi=\Phi(G,H)$ be the set of roots of $G$ relative to $H$.
$\Phi$ is an abstract root system in a real vector space $E$. Let
$\Delta$ be a base of $\Phi$. So $\mathscr{X}$ has a dual base by
the inner product on $E$ defined by Cartan matrix of $\Phi$ (see
[Hu, Appendix]. A weight is called \textit{dominant weight} if it
is a linear combination of elements of a such base of
$\mathscr{X}$ with integer non-negative coefficients.
\begin{teorema}\label{B}
Let $B$ be a \textit{Borel subgroup of $G$}, i.e. a closed,
connected and solvable subgroup of $G$ which is maximal for these
properties, containing $H$.
\begin{itemize}
\item[(a)] For every irreducible rational representation $V$ of
$G$ there exists a unique $B$-stable 1-dimensional subspace which
is a weight space $V_{\mu}$, for some dominant weight $\mu$ of
multiplicity 1 ($\mu$ is called the \textit{highest weight of $V$}
and any generator of $V_{\mu}$ is called highest weight vector ).
\item[(b)] For every dominant weight $\mu\in char(H)$ there exists an
irreducible rational representation $V$ of $G$ with highest weight
$\mu$ (called the \textit{highest weight representation of $G$})
which is unique up to isomorphism, i.e. if $V'$ is another
irreducible rational representation of $G$ with highest weight
$\mu'$ then $V$ is isomorphic to $V'$ if and only if $\mu$ equals
$\mu'$.
\end{itemize}
\end{teorema}
\textit{Proof.} See [Hu, theorem 31.3]. $\Box$\\
\\
The groups $GL(n)$ and $SL(n)$ are linearly reductive (see [GW,
theorem 2.4.5]. Hence for $GL(n)=GL(E)$, where $E=\mathbb{K}^n$
with $\mathbb{K}$ an algebraically closed field of characteristic
0,
it's enough to classify irreducible rational representations.\\
If $V$ is a vector space of dimension $m$, a rational
representation $\rho:GL(E)\rightarrow GL(V)$ is called polynomial
if and only if the entries $\rho_{ij}(g)$ of $\rho$ (for $1\leq
i,j\leq m$) are polynomials in $\{g_{ij}\}_{1\leq i,j\leq n}$,
where $g=(g_{ij})_{1\leq i,j\leq n} \in GL(E)$. A polynomial
representation $\rho:GL(E)\rightarrow GL(V)$ is homogeneous of
degree $d$ if and only if the entries $\rho_{ij}(g)$ of $\rho$
(for $1\leq i,j\leq m$) are homogeneous of degree $d$ in
$\{g_{ij}\}_{1\leq i,j\leq n}$.
\begin{proposizione}
\begin{itemize}
\item[a)] Every rational representation $V$ of $GL(E)$ is of the form
$V=V'\otimes(\bigwedge^n E)^{\otimes t}$ for some $t$, where $V'$
is a polynomial representation and $\bigwedge^n E$ is the $n$-th
exterior power of $E$.
\item[b)] Every polynomial representation of $GL(E)$ is a direct
sum of homogeneous representations.
\end{itemize}
\end{proposizione}
\textit{Proof.} See [FH, sec. 15.5]. $\Box$\\
\\
Hence it's enough to classify irreducible
homogeneous representations of degree $d$.\\
Let $\lambda$ be a partition of $d$, i.e.
$\lambda=(\lambda_1,\ldots,\lambda_k)$ with
$\lambda=\lambda_1\geq\ldots\geq  \lambda_k\geq 0$ and
$\lambda_1+\ldots+\lambda_k=d$.   We identify partitions
$(\lambda_1,\ldots,\lambda_k,0)$ with  $
(\lambda_1,\ldots,\lambda_k)$. We shall denote $d=|\lambda|$ and
we shall call the \textit{height of} $\lambda$, denoted by
$ht(\lambda)$, the number $k$ of nonzero components of $\lambda$.
Graphically we represent $\lambda$ as a set of boxes with
$\lambda_i$ boxes in the i-th row (called \textit{Young diagram
of} $\lambda$), so $|\lambda|$ and $ht(\lambda)$ are,
respectively, the number of boxes and the number of rows of the
diagram of $\lambda$. For example, if $\lambda=(4,3,1)$, then the
Young diagram of $\lambda$ is:
$$
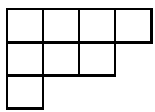
\begin{Young}
& & &\cr & & \cr \cr
\end{Young}
$$

For a partition $\lambda$ we denote  its conjugate (or transpose)
partition $\lambda'=(\lambda_1',\ldots,\lambda_t')$, where
$\lambda_j'$ is the number of boxes in the $j$-th column of the
Young diagram of $\lambda$. For example, if $\lambda=(4,3,1)$ then
$\lambda'=(3,2,2,1)$ and the Young diagram of $\lambda'$ is:
$$
  \begin{Young} &&\cr & \cr & \cr \cr
\end{Young}.
$$
 Let $T$ be a \textit{tableau} of shape
$\lambda$, i.e. a filling of the Young diagram of $\lambda$ with
numbers $1,\ldots,d$. We define the \textit{Young idempotent}
$e_T$ to be an element of the group ring $\mathbb{K}[S_d]$. In the
symmetric group $S_d$ we define the subgroups $R_T$ and $C_T$ to
be the sets of permutations in $S_d$ preserving respectively the
rows and the columns of $T$. We define
$$
e_T= \sum_{\sigma\in R_T,\tau\in C_T}sgn(\tau)\sigma\tau.
$$
Finally we define the \textit{Schur module}
$$
S_{\lambda}V:=e_T V^{\otimes d},
$$
where $V$ is a finite dimensional vector space, $dimV=n$. If $T$
and $T'$ are two tableaux of the same partition $\lambda$, then
$e_T V^{\otimes d}$ and $e_T' V^{\otimes d}$ are isomorphic as
$GL(V)$-modules [W, lemma 2.2.13]; thus $S_{\lambda}V=e_T
V^{\otimes d}$ depends on the partition $\lambda$ and not on the
tableau $T$. The representations $S_{\lambda}V$ give all
irreducible representations of $GL(V)$ homogeneous of degree $d$
[P, chap. 9
sec. 8.1].\\
For  the Schur modules sometimes we shall use the notation $S_{\lambda}V$ and sometimes the notation $S_{(\lambda_1,\ldots,\lambda_k)}V$, it depends if we want to consider or not the components of $\lambda$.\\
Now we give two examples of Schur modules. If $V$ is finite
dimensional vector space we shall call $S_n(V)$ the $n$-th
symmetric power of $V$, so the symmetric algebra of $V$ is
$S(V)=\bigoplus_{n\geq 0}S_n(V)$, and $\bigwedge ^n(V)$ the $n$-th
exterior power of $V$, so the exterior algebra of $V$ is
$\bigwedge(V)=\bigoplus_{n\geq 0}\bigwedge^n(V)$.
\begin{esempi}
Let $V$ be an $n$-dimensional vector space
\begin{itemize}
\item[(a)] If $\lambda=(d,\overbrace{0,\ldots,0}^{n-1})=(d,0^{n-1})$ then $S_{(d,0^{n-1})}V$ is just the $d$-th symmetric power $S_d(V)$.
\item[(b)] If $\lambda=(\overbrace{1,\ldots,1}^d,\overbrace{0,\ldots,0}^{n-d})=(1^d,0^{n-d})$ then  $S_{(1^d,0^{n-d})}V$ is just the $d$-th exterior power $\bigwedge^d(V)$; in particular if $d=dim\,V$, $S_{(1^{dim\,V})}V=\bigwedge^{dim\,V}(V):=D$ is called a \textit{determinant representation of $G$}.
\item[c)] If $k>n$ and $\lambda_k>0$, we have $S_{(\lambda_1,\ldots,\lambda_k)}V=0$.
\end{itemize}
\end{esempi}
Introducing the convention
$\bigwedge^n(V^*)=S_{(\underbrace{-1,\ldots,-1}_n)}V$ and
$S_{(\lambda_1,\ldots,\lambda_n)}V\otimes\bigwedge^n(V^*)=S_{(\lambda_{1}-1,\ldots,\lambda_{n}-1)}V$,
we see that there is a bijective correspondence between rational
irreducible representations of $GL(n)$ and vectors
$(\lambda_1,\ldots,\lambda_n)\in\mathbb{Z}^n$ such that
$\lambda_1\geq\cdots\geq\lambda_n$.\\
We give an alternative description of Schur modules equivalent to
that already given [W, lemma 2.2.13]. Let $V$ be an
$n$-dimensional vector space. Let
$$
m:\bigwedge^r V\otimes \bigwedge^s V\rightarrow\bigwedge^{r+s} V,
$$
such that
$$
m(u_1\wedge\ldots\wedge u_r\otimes v_1\wedge\ldots\wedge
v_s)=u_1\wedge\ldots\wedge u_r\wedge v_1\wedge\ldots\wedge v_s,
$$
be the multiplication in the exterior algebra $\bigwedge V$ and
let
$$
\Delta:\bigwedge^{r+s} V\rightarrow\bigwedge^r V\otimes
\bigwedge^s V,
$$
such that
$$
\Delta(u_1\wedge\ldots\wedge u_{r+s})=\sum_{\sigma\in
S_{r+s}^{r,s}}(-1)^{sgn(\sigma)}u_{\sigma(1)}\wedge\ldots\wedge
u_{\sigma(r)}\otimes u_{\sigma(r+1)}\wedge\ldots\wedge
u_{\sigma(r+s)}
$$
where $S_{r+s}^{r,s}=\{\sigma\in S_{r+s}|\sigma(1)<\cdots
<\sigma(r);\sigma(r+1)<\cdots <\sigma(r+s)\}$, be the
comultiplication in the exterior algebra $\bigwedge V$. We
consider $\lambda=(\lambda_1,\ldots,\lambda_k)$ a partition of
$d$. We can define the Schur module as
$$
S_{\lambda}V:=\bigwedge^{\lambda_1}V\otimes\cdots\otimes\bigwedge^{\lambda_k}V/
R(\lambda,V),
$$
where
$$
R(\lambda,V)=\sum_{1\leq a\leq
k-1}\bigwedge^{\lambda_1}V\otimes\cdots\otimes\bigwedge^{\lambda_{a-1}}V\otimes
R_{a,a+1}(V)\otimes\bigwedge^{\lambda_{a+2}}V\otimes\cdots\otimes\bigwedge^{\lambda_k}V
$$
where $R_{a,a+1}(V)$ is the submodule spanned by the images of the
following maps $\theta(\lambda,a,u,v;V)$ with $u+v<\lambda_{a+1}$:
$$
\begin{array}{c}
\bigwedge^{u}V\otimes\bigwedge^{\lambda_{a}-u+\lambda_{a+1}-v}V\otimes\bigwedge^{v}V\\

\downarrow \scriptstyle{1\otimes\Delta\otimes 1}\\

 \bigwedge^{u}V\otimes\bigwedge^{\lambda_{a}-u}\otimes\bigwedge^{\lambda_{a+1}-v}V\otimes\bigwedge^{v}V\\

\;\;\;\downarrow \scriptstyle{m_{12}\otimes m_{34}}\\

 \bigwedge^{\lambda_{a}}V\otimes\bigwedge^{\lambda_{a+1}}V.
 \end{array}
$$
Let us choose an ordered basis $\{e_1,\ldots,e_n\}$ of V. If $T$
is a tableau of shape $\lambda$ with entries in $\{1,\ldots,n\}$,
we associate to $T$ the element in $S_{\lambda}V$
$$
e_{T(1,1)}\wedge\ldots\wedge
e_{T(1,\lambda_1)}\otimes\ldots\otimes
e_{T(k,1)}\wedge\ldots\wedge e_{T(k,\lambda_k)}+R(\lambda,V),
$$
where $T(i,j)$ is the entry of $T$ in the $i$-th row and $j$-th
column of the Young diagram of $\lambda$.\\
We recall some properties and some known results about Schur modules.\\
A filling of the Young diagram of a partition $\lambda$ with the
numbers $1,\ldots,n$ weakly increasing along each row and strictly
increasing along each column  is called \textit{column standard
tableau corresponding to the basis $\{e_1,\ldots,e_n\}$}.
\begin{teorema}
Let $\{e_1,\ldots,e_n\}$ be a basis of $V$. The column standard
tableaux corresponding to this basis form a basis of
$S_{\lambda}V$
\end{teorema}
\textit{Proof.} See [W, prop. 2.1.4]. $\Box$\\
\\
If $V$ is an $n$-dimensional vector space, a Borel subgroup of
$GL(V)=GL(n)$ is the subgroup of all upper triangular matrices,
the maximal torus $H$ of $GL(n)$ is the subgroup of diagonal
matrices and the sequences $(\lambda_1,\ldots,\lambda_n)$, with
$\lambda_i\in\mathbb{Z}$ and $\lambda_1\geq\ldots\geq\lambda_n$,
are the \textit{dominant integral weights for} $GL(n)$; we shall
write $x=diag(x_1,\ldots,x_n)$ in $H$ for the diagonal matrix with
these entries. The decomposition of $V$ into direct sum of weight
spaces is
$$
\bigoplus_{a=(a_1,\ldots,a_n)\in\mathbb{Z}^n}V_a=\{v\in V|x\cdot
v=\prod_{i=1}^n x_i^{a_i}v\;\forall x\in H\},
$$
see [B, chap. 3 sec. 8].
\begin{teorema}
Let $V$ be an $n$-dimensional vector space.
\begin{itemize}
\item[1)] If $\lambda$ is a partition with at most $n$ components
then the representation $S_{\lambda}V$ of $GL(n)$ is an
irreducible representation of highest weight
$\lambda=(\lambda_1,\ldots,\lambda_n)$.
\item[2)] For any $\mu=(\mu_1,\ldots,\mu_n)$ with
$\mu_1\geq\cdots\geq\mu_n$ integers, there is a unique irreducible
representation of $GL(n)$ with highest weight $\mu$, which can be
realized as $S_{\lambda}V\otimes D^{\otimes k}$, for any
$k\in\mathbb{Z}$ and where $\lambda_i=\mu_i-k\geq 0$ for every
$i\in\{1,\ldots,n\}$.
\end{itemize}
\end{teorema}
\textit{Proof.} See [F, sec. 8.2 theorem 2]. $\Box$\\
\\
By theorem \ref{B} and by the previous one, every irreducible
rational representation is a Schur module tensored with a power of
a determinant representation.
\begin{teorema}[Properties of Schur modules]\label{pS}
Let $V$ be vector space of dimension $n$ and $\lambda$ be the
highest weight for $GL(n)$.
\begin{itemize}
\item[(i)] $S_{\lambda}V=0\Leftrightarrow ht(\lambda)> n$.
\item[(ii)] $dim\;S_{\lambda}V=1\Leftrightarrow\lambda=(\overbrace{k,\ldots,k}^n)=(k^n)$ for some $k\in\mathbb{Z}$.
\item[(iii)] $\big(S_{(\lambda_1,\ldots,\lambda_n)}V\big)^*\cong S_{(\lambda_1,\ldots,\lambda_n)}V^*\cong S_{(-\lambda_n,\ldots,-\lambda_1)}V$.
\end{itemize}
\end{teorema}
\textit{Proof.} See [FH, theorem 6.3].
\begin{teorema}[Cauchy formulas]\label{c}
Let $V$ and $W$ be two finite dimensional vector spaces.
\begin{itemize}
\item[a)]  As a representation of $GL(V)\times GL(W)$, $S_d(V\otimes W)$ decomposes as
$$
S_d(V\otimes W)=\bigoplus_{|\lambda|=d}S_{\lambda}V\otimes
S_{\lambda}W;
$$
\item[b)] As a representation of $GL(V)\times GL(W)$, $S_d(V\otimes W)$ decomposes as
$$
\bigwedge^d(V\otimes W)=\bigoplus_{|\lambda|=d}S_{\lambda}V\otimes
S_{\lambda'}W;
$$
\item[c)]  As a representation of $GL(V)$, $S_d(S_2(V))$ decomposes as
$$
S_d(S_2(V))=\bigoplus_{|\lambda|=d}S_{2\lambda}V,
$$
where $2\lambda=(2\lambda_1,\ldots,2\lambda_k)$ if
$\lambda=(\lambda_1,\ldots,\lambda_k)$;
\item[d)] As a representation of $GL(V)$ the ring  $S_d(\bigwedge^2(V))$ decomposes as
$$
S_d(\bigwedge^2(V))=\bigoplus_{|\lambda|=d}S_{2\lambda'}V.
$$
\end{itemize}
\end{teorema}
\textit{Proof.} See [P] chap. 9 sec. 6.3 and sec. 8.4 , chap. 11
sec. 4.5.\\
\\
Finally we consider the tensor product of Schur modules
\begin{lemma}
$$
S_{\lambda}V\otimes
S_{\mu}V=\bigoplus_{\nu}c^{\nu}_{\lambda\mu}S_{\nu}V,
$$
where $c^{\nu}_{\lambda\mu}$'s are called Littlewood-Richardson
coefficients.
\end{lemma}
There is a combinatorial formula to calculate
$c^{\nu}_{\lambda\mu}$.\\
Let
$$
D_{\lambda}=\{(i,j)|\,1\leq i\leq k\,1\leq j\leq\lambda_i\}
$$
be the Young diagram of $\lambda$ and let
$f:D_{\nu/\lambda}\rightarrow \{1,\ldots,n\}$ be a column standard
tableau. We denote $CST(\nu/\lambda,\{1,\ldots,n\})$ the set of
column standard tableaux of shape $\nu/\lambda$ with values in
$\{1,\ldots,n\}$. We define $cont(f)$, the \textit{content of}
$f$, to be the sequence $\{|f^{-1}(1)|,\ldots,|f^{-1}(n)|\}$. We
define $w(f)$ to be the word we get from $f$ when we read it by
rows, starting with the first row, from right to left in each row.
A word $w=(w_1,\ldots,w_m)$ on the alphabet $\{1,\ldots,n\}$ is a
\textit{lattice permutation} if for each $1\leq u\leq m$ and for
each $1\leq i\leq n-1$ we have
$$
|\{1\leq j\leq u|\,w_j=i\}|\geq|\{1\leq j\leq u|\,w_j=i+1\}|.
$$
Finally we define the set
$$
LR^{\nu}_{\lambda,\mu}=\{f\in
CST(\nu/\lambda,\{1,\ldots,n\})|\,cont(f)=(\mu_1,\ldots,\mu_n),\,w(f)\;\textrm{is
a lattice permutation}\}.
$$

\begin{teorema}[Littlewood-Richardson rule]\label{lrr}
Let $\lambda,\mu,\nu$ be partitions, then
$$
c^{\nu}_{\lambda\mu}=|LR^{\nu}_{\lambda,\mu}|.
$$
\end{teorema}
\textit{Proof.} See [P, chap. 12 sec. 5.3]. $\Box$
\begin{cor}\label{clrr}
If $\lambda=(l^s)$ and $\mu=(m^t)$, then $S_{\lambda} V\otimes
S_{\mu}V$ is multiplicity free, i.e. $S_{\lambda} V\otimes
S_{\mu}V=\bigoplus_{\nu}S_{\nu}V$. Moreover if $s\geq t$ then
$\nu=(\nu_1,\ldots,\nu_{s+t})$ with $\nu_i=l+c_i$ for $1\leq i\leq
t$, $\nu_i=l$ for $t< i\leq s$ and $\nu_{s+i}=m-c_{t-i+1}$ for
$1\leq i\leq t$, where $m\geq c_1\geq\ldots\geq c_t\geq 0$ and
$l+c_t\geq m$.
\end{cor}
\textit{Proof.} We note that we can suppose in the statement
$s\geq t$, since the tensor product is commutative. The proof is a
consequence of Littlewood-Richardson rule. $\Box$

\section{Invariant theory}
In this section we recall definitions and fundamental results of invariant theory.\\
If $G$ is a group which acts on a finite dimensional vector space
$V$, we shall call $V^G=\{v\in V|g\cdot v=v\;\forall g\in G\}$ the
space of invariants of $V$ and we have a general lemma
\begin{lemma}\label{VxW}
Let $G$ be a group which acts on two finite dimensional vector
space $V$ and $W$. If $G$ acts trivially on $V$ , then $(V\otimes
W)^G=V\otimes W^G$.
\end{lemma}

If $G$ is an algebraic group and $V$ is a rational representation
of $G$, then $G$ acts on the coordinate ring of $V$
$\mathbb{K}[V]$ as follows: if $f\in\mathbb{K}[V]$ and $g\in G$,
$$(g\cdot f)(v)=f(g^{-1}\cdot v).$$ The ring of $G$-invariants in
$\mathbb{K}[V]$ is $$\mathbb{K}[V]^G=\{f\in\mathbb{K}[V]|g\cdot
f=f\;\forall g\in G\}.$$

\begin{teorema}[Hilbert]
If $G$ is linearly reductive and acts rationally on a finite
dimensional vector space $V$ then $\mathbb{K}[V]^G$ is finitely
generated.
\end{teorema}
\textit{Proof.} See [P, chap. 14 sec. 1.1].\\
\\
Now we formulate the first fundamental theorem for the linear
group.
\begin{teorema}[FFT for $GL$]\label{fft}
Let $V$ be a finite dimensional vector space. We take the space
$(V^*)^p\times
V^q=\big\{(\alpha_1,\ldots,\alpha_p,v_1,\ldots,v_q)|\alpha_j\in
V^*,v_i\in V\;\forall j\in\{1,\ldots,p\}\;and\;\forall
i\in\{1,\ldots,q\}\big\}$ as a representation of $GL(V)$. On this
space we consider the $pq$ polynomial functions $u_{ij}(
\alpha_1,\ldots,\alpha_p,v_1,\ldots,v_q)=\alpha_j(v_i)$ which are
$GL(V)$-invariant. Then
 $$
 \mathbb{K}[(V^*)^p\times V^q]^{GL(V)}=\mathbb{K}[u_{ij}]_{{1\leq i\leq q  \atop 1\leq j\leq p}}
 $$
 \end{teorema}
 \textit{Proof.} See [P, chap. 9 sec.1.4].\\
 \\
 Now we give the definition of  semi-invariant and of character of an algebraic group.
 \begin{definizione}
 Let $G$ be an algebraic group and let $V$ be a rational representation of $G$.
 \begin{itemize}
 \item[(i)] $\chi:G\rightarrow \mathbb{K}^*$ is a character of $G$ if it is a homomorphism of algebraic groups;
 \item[(ii)] $f\in\mathbb{K}[V]$ is a semi-invariant of weight $\chi$ of the action of $G$ on $V$ if for every $g\in G$, $g\cdot f=\chi(g) f$ where $\chi$ is a character of $G$.
 \end{itemize}
 \end{definizione}
 If $char(G)$ is the set of characters of $G$, then the ring of semi-invariants of the action of $G$ on $V$ is
 $$
 SI(G,V)=\bigoplus_{\chi\in char(G)}SI(G,V)_{\chi}
 $$
 where $SI(G,V)_{\chi}=\{f\in\mathbb{K}[V]|\forall g\in G, g\cdot f=\chi(g)f\}$ is called weight space.
 In general we have the following lemma proved in [SK].
\begin{lemma}[Sato-Kimura]\label{sk}
Let $G$ be a linear algebraic group acting rationally on the
vector space $V$. If there is a Zariski open $G$-orbit in $V$ then
the ring $SI(G, V )$ spanned by the semi-invariants is a
polynomial ring:
$$
SI(G, V )=k[f_1,\ldots,f_s]
$$
 for some collection of algebraically independent and
irreducible semi-invariants $f_1,\ldots, fs$. Moreover if $f_i\in
SI(G, V )_{\chi_i}$ then the $\chi_i$ are linearly independent
over $\mathbb{Z}$ in the space of characters of $G$.
\end{lemma}
\begin{cor}
Under the assumptions of the lemma \ref{sk}, the set of characters
$\chi$ such that $SI(G, V )_{\chi}\neq 0$ forms a free abelian
semigroup, isomorphic to $\mathbb{N}^s$. In particular, if $f$ is
any semi-invariant of weight $\chi$, then $f = uf^{a_1}_1\cdots
f^{a_s}_s$ , where $u$ is a unit in $\mathbb{K}$ and the $a_i\geq
0$ are the unique integers such that $\chi =\sum_{i=1}^s
a_i\chi_i$ in the space of characters of $G$. Thus $SI(G,V)$ is a
polynomial ring.
\end{cor}
If $G=GL(n)$, there exists an isomorphism $\mathbb{Z}\cong
char(GL(n))$ which sends an element $a$ of $\mathbb{Z}$ in
$(det)^a$ (where $det$ associates to $g\in GL(n)$ its
determinant). So we have
 $$
 SI(G,V)=\mathbb{K}[V]^{SL(V)}.
 $$
 Finally other two results on Schur modules and invariant theory.
 \begin{proposizione}\label{i1}
 Let $V$ be a finite dimensional vector space of dimension $n$.
 $$
 (S_{\lambda}V)^{SL(V)}\neq 0\Longleftrightarrow\lambda=(k^n)
 $$
 for some $k$ and in this case $S_{\lambda}V$, and so also $(S_{\lambda}V)^{SL(V)}$, have dimension one.
 \end{proposizione}
 \begin{proposizione}\label{i2}
  Let $V$ be a finite dimensional vector space of dimension $n$ and let $\lambda$ and $\mu$ be two dominant integral weights. Then\\
  $$
  S_{\lambda}V\otimes S_{\mu}V\; \textrm{contains a
  semi-invariant}
  $$
  $$
  \Longleftrightarrow
 $$

 $$
 \begin{array}{ccc}
 \lambda_1-\lambda_2&=&\mu_{n-1}-\mu_n\\
  \lambda_2-\lambda_3&=&\mu_{n-2}-\mu_{n-1}\\
  & \vdots & \\
  \lambda_{n-1}-\lambda_n&=&\mu_1-\mu_2
  \end{array}
  $$
  and in this case the semi-invariant is unique (up to a non zero scalar) and has weight
  $\lambda_1+\mu_n=\lambda_2+\mu_{n-1}=\cdots=\lambda_n+\mu_1$.
  \end{proposizione}
  \textit{Proof.} It is a corollary of (5.6) in [M, I.5]. $\Box$\\
  \\
  Let $Sp(2n)=\{A\in GL(2n)|AJA=J\}$ be
  the simplectic group, let $O(n)=\{A\in GL(n)|A^tA=I\}$
  be the orthogonal group and $SO(n)=\{A\in O(n)|det\,A=1\}$ be the special orthogonal group,
   where I is the identity matrix and
  $J=\left(\begin{array}{cc} 0 & I\\
  -I & 0 \end{array}\right)$.
  \begin{proposizione}\label{i3}
  Let $V$ be an orthogonal space of dimension $n$ and let $W$ be a
  symplectic space of dimension $2n$.
  \begin{itemize}
  \item[(a)] $ dim\,(S_{\lambda}V)^{O(V)}=\left\{\begin{array}{ll}1 &
  \textrm{if}\quad \lambda=2\mu\\
  0 & \textrm{otherwise}
  \end{array}\right.$,
\item[(b)] $ dim\,(S_{\lambda}V)^{SO(V)}=\left\{\begin{array}{ll}1 &
  \textrm{if}\quad \lambda=2\mu+(k^n)\\
  0 & \textrm{otherwise}\end{array}\right.$,
  \item[(c)] $ dim\,(S_{\lambda}W)^{Sp(W)}=\left\{\begin{array}{ll}1 &
  \textrm{if}\quad \lambda=2\mu'\\
  0 & \textrm{otherwise}\end{array}\right.$
  \end{itemize}
  for some partition $\mu$ and for some $k\in\mathbb{Z}_{\geq 0}$.
  \end{proposizione}
  \textit{Proof.} See [P] chap. 11 cor. 5.2.1 and 5.2.2. $\Box$\\
  \\
  We end this section recalling definition and properties of the
  \textit{Pfaffian} of a skew-symmetric matrix.\\
  Let $A=(a_{ij})_{1\leq i,j\leq 2n}$ be a skew-symmetric $2n\times 2n$
  matrix. Given $2n$ vectors $x_1,\ldots,x_{2n}$ in $\mathbb{K}^{2n}$, with $\mathbb{K}$ an algebraically closed field with characteristic 0,  we define
  $$
  F_A(x_1,\ldots,x_{2n})=\frac{1}{n!2^n}\sum_{s\in
  S_{2n}}sgn(s)\prod_{i=1}^n(x_{s(2i-1)},x_{s(2i)}),
  $$
  where $S_{2n}$ is the symmetric group on $2n$ elements, $sgn(s)$
  is the sign of permutation $s$ and $(\cdot,\cdot)$ is the
  skew-symmetric bilinear form associated to $A$. So $F_A$ is a
  skew-symmetric multilinear function of $x_1,\ldots,x_{2n}$.
  Since, up to a scalar, the only one skew-symmetric
  multilinear function of $2n$ vectors in $\mathbb{K}^{2n}$ is the
  determinant, there is a complex number $Pf(A)$, called
  \textit{Pfaffian of} $A$, such that
  $$
  F_A(x_1,\ldots,x_{2n})=Pf(A)det[x_1,\ldots,x_{2n}]
  $$
  where $[x_1,\ldots,x_{2n}]$ is the matrix which has the vector
  $x_i$ for $i$-th column. In particular one proves that
  $$
  Pf(A)=\frac{1}{n!2^n}\sum_{s\in
  S_{2n}\setminus B_n}sgn(s)\prod_{i=1}^n a_{s(2i-1)s(2i)}
  $$
  where $B_n$ is a subgroup of $S_{2n}$ isomorphic to the
  semidirect product $S_n\ltimes(\mathbb{Z}_2)^n$. We can write
  the Pfaffian of $A$ avoiding to sum on all possible
  permutations,
  $$
 Pf(A)=\sum_{{i_1<j_1,\ldots,i_n<j_n \atop i_1<\ldots<i_n}}sgn(s)a_{1_1j_1}\cdots
  a_{i_nj_n}
  $$
  where $s$ is the permutation $\left[\begin{array}{ccccc} 1 & 2
  &\ldots & 2n-1 & 2n \\
  i_1 & j_1 & \ldots & i_n & j_n \end{array}\right]$.
  \begin{proposizione}
  Let $A$ be a skew-symmetric $2n\times 2n$ matrix.
  \begin{itemize}
  \item[(i)] For every invertible $2n\times 2n$ matrix $B$,
  $$
  Pf(BAB^t)=det(B)Pf(A);
  $$
  \item[(ii)] $det(A)=Pf(A)^2$.
  \end{itemize}
  \end{proposizione}
  \textit{Proof.} See [P, chap. 5 sec. 3.6]. $\Box$

  \chapter{Quiver representations and semi-invariants}

\section{Auslander-Reiten theory}\label{B1}

A quiver $Q$ is a pair $(Q_0,Q_1)$ where $Q_0$ is the set of
vertices and $Q_1$ is the set of arrows. Let
$$
a:ta\longrightarrow ha,\quad ta,ha\in Q_0
$$
be an arrow in $Q_1$. We shall call $ta$ the tail of the arrow $a$
and $ha$ the head of the arrow $a$. A path $p$ in $Q$ is a
sequence of arrows $p=a_1\cdots a_n$ such that $ha_i=ta_{i+1}$,
$(1\leq i\leq n-1)$. For every $x\in Q_0$ we also have a trivial
path $e_x$ such that $he_x=te_x=x$.
We say that $Q$ \textit{has no oriented cycles} if there are no paths  $p=a_1\cdots a_n$ such that $ta_1=ha_n$.\\
We fix an algebraically closed field $\mathbb{K}$. A
representation $V$ of $Q$ is a family of finite dimensional vector
spaces
 $\{V(x)|x\in Q _0\}$ and of linear maps
$\{V(a):V(ta)\rightarrow V(ha)\}_{a\in Q_1}$. The dimension vector
of $V$ is a function
$\underline{dim}(V):Q_0\rightarrow\mathbb{Z}_{\geq 0}$ defined by
$\underline{dim}(V)(x):=dim V(x)$.\\
A morphism $f:V\rightarrow W$ of two representations is a family
of linear maps $\{f(x):V(x)\rightarrow
W(x)|\,f(ha)V(a)=W(a)f(ta)\forall a\in Q_1\}_{x\in Q_0}$. We
denote the space of morphisms from $V$ to $W$ by $Hom_Q(V,W)$ and
the space of extensions of $V$ by $W$ by $Ext_Q^1(V,W)$.
\begin{definizione}\label{fe}
The non symmetric bilinear form on the space $\mathbb{Z}^{Q_0}$ of
dimension vectors given by
$$
\langle\alpha,\beta\rangle=\sum_{x\in
Q_0}\alpha(x)\beta(x)-\sum_{a\in Q_1}\alpha(ta)\beta(ha)
$$
is the Euler form of $Q$, where $\alpha,\beta\in\mathbb{Z}^{Q_0}$.
\end{definizione}
If $\underline{dim}\,V=\alpha$ and $\underline{dim}\,W=\beta$, we
have
$$
\langle\alpha,\beta\rangle=dim\,Hom_Q(V,W)-dim\,Ext^1_Q(V,W)
$$
We shall call $Rep(Q,\alpha)$ the variety of representations of
$Q$ of dimension vector $\alpha$.
\begin{definizione}
Let $Q$ be a quiver and let $\alpha$ be a dimension vector. A
general representation of $Q$ is a representation from some
nonempty Zariski open set in $Rep(Q,\alpha)$.
\end{definizione}
We recall the definitions of simple, projective and injective
representation of a quiver $Q=(Q_0,Q_1)$. For each vertex $x$, a
simple representation $S_x$ is the representation for which
$S_x(x)=\mathbb{K}$, $S_x(y)=0$ for every $y\in Q_0\setminus\{x\}$
and $S_x(a)$ is the zero map for every $a\in Q_1$. For every $x\in
Q_0$ we define an indecomposable projective representation $P_x$
as follows:
$$
P_x(y)=[x,y]\;\textrm{and}\;P_x(a):=a\circ:[x,ta]\rightarrow
[x,ha]
$$
with $x,y\in Q_0$ and $a\in Q_1$, where $[x,y]$ is a vector space
over $\mathbb{K}$ with a basis labelled by all paths from x to y
in $Q$ and $a\circ$ is the map which sends the path $p$ to the
path $a\circ p$. Every indecomposable projective representation of
$Q$ is isomorphic to $P_x$ for some $x\in Q_0$ and moreover we
have $Hom_Q(P_x,V)\cong V(x)$ for every representation $V$ of $Q$,
see [ARS, sec III.1]. Similarly every indecomposable injective
representation of $Q$ is isomorphic to $I_x$, where $I_x$ is
defined as follows:
$$
I_x(y)=[y,x]^*\;\textrm{and}\;I_x(a):=(\circ
a)^*:[ta,x]^*\rightarrow [ha,x]^*
$$
with $x,y\in Q_0$ and $a\in Q_1$, where $[y,x]^*$ is the dual
space of $[y,x]$ and $\circ a:[ha,x]\rightarrow [ta,x]$ is the map
which sends $p$ to $p\circ a$. In this case we have
$Hom_Q(V,I_x)\cong V(x)^*$ for every representation $V$ of $Q$,
where $V(x)^*$ is the dual space of $V(x)$.\\
Now we recall some definitions and results of Auslander-Reiten
 Theory, for deepening see [ARS] and [ASS].\\
 We define \textit{the path algebra} $\mathbb{K}Q$ of a quiver
 $Q$, the $\mathbb{K}$-algebra which has the paths of
 $Q$ as basis. The multiplication in $\mathbb{K}Q$ is defined by
 $$
 p\cdot q=\left\{\begin{array}{cc}
 pq & \textrm{if}\;tp=hq\\
 0 & \textrm{otherwise}.
 \end{array}\right.
 $$
 \begin{proposizione}
 \begin{itemize}

\item[1)] $\mathbb{K}Q$ is a finite-dimensional
 $\mathbb{K}$-algebra if and only if $Q$ has no oriented cycles.
 \item[2)] The categories $Rep(Q)$ of representations of $Q$ and
 $\mathbb{K}Q-mod$ of left $\mathbb{K}Q$-modules are equivalent.
 \end{itemize}
\end{proposizione}
\textit{Proof.} See [ARS, sec. 3.1 prop. 1.1 and prop. 1.3] and
[ASS, sec. II.1 lemma 1.4(c) and sec. III.1 cor. 1.7].
$\Box$\\
\\
Let $A$ be a finite-dimensional $\mathbb{K}$-algebra, a morphism
$f:V\rightarrow W$ in the category of left $A$-modules $A-mod$ is
called a \textit{retraction} if there exists $g:W\rightarrow V$
such that $fg=id_W$ and it is called a \textit{section} if there
exists $g:W\rightarrow V$ such that $gf=id_V$.
\begin{definizione}
Let $f:V\rightarrow W$ be a morphism in $A-mod$.
\begin{itemize}
\item[(a)] $f$ is called
minimal right almost split if
\begin{itemize}
\item[(i)] every endomorphism $h:V\rightarrow V$ such that $fh=f$,
is an  isomorphism (right minimal morphism),
\item[(ii)] $f$ is not a retraction,
\item[(iii)] for every $g:V'\rightarrow W$ which is not a
retraction there exists $g':V'\rightarrow V$ such that $fg'=g$.
\end{itemize}
\item[(b)] $f$ is called irreducible if it is neither a section
nor a retraction and if $f=ts$, for some $s:V\rightarrow X$ and
$t:X\rightarrow W$, then $s$ is a section or $t$ is a retraction.
\end{itemize}
\end{definizione}
Now we are able to define the Auslander-Reiten quiver and the
almost split sequences.
\begin{definizione}\label{arQ}
Let $Q$ be a quiver and $\mathbb{K}Q$ be the path algebra of $Q$.
The quiver $AR(Q)=(AR(Q)_0,AR(Q)_1)$, where the set of vertices
$AR(Q)_0$ is the set of indecomposables of $\mathbb{K}Q$ and the
set of arrows $AR(Q)_1$ is the set of the irreducible morphisms
not zero between indecomposables, is called Auslander-Reiten
quiver of $Q$.
\end{definizione}
\begin{teorema}
If $W$ is an indecomposable non-projective $A$-module
(respectively $V$ is an indecomposable non-injective $A$-module)
then there exists an exact sequence $0\rightarrow
V\stackrel{f}{\rightarrow}Z\stackrel{g}{\rightarrow}W\rightarrow
0$ such that $f$ and $g$ are both irreducible, called almost split
sequence.
\end{teorema}
\textit{Proof.} See [ARS, sec. 5.1 theorem 1.15]. $\Box$\\
\\
If $V$ is an $A$-module, a right minimal morphism $p:P\rightarrow
V$, with $P$ projective, is called a \textit{projective cover of}
$V$. One can prove that every $A$-module $V$ has a minimal
projective presentation
$P_1\stackrel{p_1}{\rightarrow}P_0\stackrel{p_0}{\rightarrow}V\rightarrow
0$, i.e. an exact sequence where $p_0$ is a projective cover of
$V$ and $p_1$ is a projective cover of $Ker\,p_0$ ([ARS, sec. 1.4
theorem 4.2] and [ASS, sec. I.5 theor. 5.8]).\\
Let $V\in  A-mod$, we assume that $V$ has no projective summands
and let
$P_1\stackrel{p_1}{\rightarrow}P_0\stackrel{p_0}{\rightarrow}V\rightarrow
0$ be a minimal presentation of $V$. Applying the functor
$Hom_A(\cdot,A)$ on it, we obtain a minimal presentation
$$
Hom(P_0,A)\stackrel{Hom(p_1,A)}{\longrightarrow}Hom(P_1,A)\longrightarrow
Coker(Hom(p_1,A))\longrightarrow 0.
$$
We define $coKer(Hom(p_1,A)):=Tr(V)$, the \textit{transpose of}
$V$. Thus the transpose is a contravariant functor
$Tr:A-mod\rightarrow mod-A$ ($mod-A$ is the category of right
$A$-modules) which equals zero on projective modules. We can
define also $Tr:mod-A\rightarrow A-mod$ considering
$$
mod-A\cong A^{op}-mod\stackrel{Tr}{\longrightarrow}mod-A^{op}\cong
A-mod.
$$
\begin{proposizione}
If $A=\mathbb{K}Q$ and $V$ is a representation of $Q$ without
projective direct summands, then $Tr(V)=Ext^1_{A}(V,A)$.
\end{proposizione}
\textit{Proof.} See [ARS, sec. 4.1 corollary 1.14]. $\Box$
\begin{definizione}
The functor
$$
\tau^+:=\nabla\circ
Tr:A-mod\stackrel{Tr}{\longrightarrow}mod-A\cong
A^{op}-mod\stackrel{\nabla}{\longrightarrow}A-mod,
$$
where $\nabla$ is the duality functor sending the representation
$V$ to $V^*$, is called \textit{Auslander-Reiten translation
(AR-translation)}. Similarly we can define the functor
$\tau^-:=Tr\circ\nabla$.
\end{definizione}
We note that, by definition, $\nabla\tau^-=\tau^+\nabla$ and
$\nabla\tau^+=\tau^-\nabla$.\\
The following theorem records an important property of the
AR-translation.
\begin{teorema}[Auslander-Reiten duality]\label{ard}
Let $A=\mathbb{K}Q$ and let $V$ and $W$ be two $A$-modules.
\begin{itemize}
\item[(a)] If $V$ has no projective summands, then there exist
isomorphisms of vector spaces
$$
Hom_Q(W,\tau^+V)\cong
Ext^1_Q(V,W)^*\;\textit{and}\;Ext^1_Q(W,\tau^+V)\cong
Hom_Q(V,W)^*.
$$
\item[(b)] If $V$ has no injective summands, then there exist
isomorphisms of vector spaces
$$
Hom_Q(\tau^-V,W)\cong
Ext^1_Q(W,V)^*\;\textit{and}\;Ext^1_Q(\tau^-V,W)\cong
Hom_Q(W,V)^*.
$$
\end{itemize}
\end{teorema}
\textit{Proof} See [ASS, sec. IV.2 cor. 2.14]. $\Box$
\begin{cor}
Let $A=\mathbb{K}Q$ and let $V$ and $W$ be two $A$-modules.
\begin{itemize}
\item[(a)] If $V$ and $W$ have no projective summands, then there exist
isomorphisms of vector spaces
$$
Hom_Q(\tau^+V,\tau^+W)\cong Hom_Q(V,W)
$$
and
$$
Ext^1_Q(\tau^+V,\tau^+W)\cong Ext^1_Q(V,W).
$$
\item[(b)] If $V$ and  has no injective summands, then there exist
isomorphisms of vector spaces
$$
Hom_Q(\tau^-V,\tau^-W)\cong Hom_Q(V,W)
$$
and
$$
Ext^1_Q(\tau^-V,\tau^-W)\cong Ext^1_Q(V,W).
$$
\end{itemize}
\end{cor}
\textit{Proof.} It is an immediate consequence of theorem 1.9.
$\Box$\\
\\
By AR-duality, if we consider $\tau^+$ and $\tau^-$ as linear
transformations on the space of dimension vectors, i.e. if $V$ is
a representation of a quiver with dimension $\alpha$ then
$\tau^\pm\alpha:=\underline{dim}\,\tau^\pm V$, we have, for every
$\alpha$ and $\beta$ dimension vectors, then
\begin{itemize}
\item[(i)]
$\langle\alpha,\beta\rangle=-\langle\tau^-\beta,\alpha\rangle$
\item[(ii)] $\langle\alpha,\beta\rangle=-\langle\beta,\tau^+\alpha\rangle$
\item[(iii)]
$\langle\alpha,\beta\rangle=\langle\tau^\pm\alpha,\tau^\pm\beta\rangle$.
\end{itemize}
At last another result about the existence of the almost split
sequences.
\begin{teorema}[Auslander-Reiten 1975]\label{ar}
\begin{itemize}
\item[1)] For every finitely generated indecomposable non-projective
module $V$ there is an almost split sequence $0\rightarrow\tau^+
V\rightarrow X\rightarrow V\rightarrow 0$ in $A-mod$ with finitely
generated modules.
\item[2)] For every finitely generated indecomposable non-injective
module $V$ there is an almost split sequence $0\rightarrow
V\rightarrow Z\rightarrow \tau^-V\rightarrow 0$ in $A-mod$ with
finitely generated modules.
\end{itemize}
\end{teorema}
\textit{Proof.} It is a direct consequence of the theorem 1.8, see
also [ASS, sec. IV.3 theor. 3.1]. $\Box$

\section{Quivers of tame type}\label{Qtt}
\begin{definizione}
A quiver $Q$ is called of tame type if the underlying graph of $Q$
is of type $\widetilde{A},\widetilde{D}$ or $\widetilde{E}$.
\end{definizione}
For all of the next results we refer to [DR].
\begin{proposizione}\label{h}
Let $Q$ be a quiver of tame type, then the quadratic form
$q_Q:\mathbb{Z}^{Q_0}\rightarrow\mathbb{Z}$ defined by
$$
q_Q(\alpha):=\sum_{x\in Q_0}\alpha(x)^2-\sum_{a\in
Q_1}\alpha(ta)\alpha(ha)
$$
is positive semi-definite and there exists a unique vector
$h\in\mathbb{N}^{Q_0}$ such that $\mathbb{Z}h$ is the radical of
$q_Q$ or, equivalently, such that $\tau^+h=h$ and $|h|:=\sum_{x\in
Q_0}h(x)$ is minimum in $\mathbb{N}^{Q_0}$. For quivers of type
$\widetilde{A}$ and $\widetilde{D}$ the vector $h$ has the
following form
\begin{equation}
\begin{array}{cccccccccccc}
&&&&1&\cdots&1&\\
\widetilde{A}:&&&1&&&&1\\
&&&&1&\cdots&1&
\end{array}
\end{equation}
\\
\begin{equation}
\begin{array}{cccccccccccc}
&&&1&&&&1\\
\widetilde{D}:&&&&2&\cdots&2&\\
&&&1&&&&1
\end{array}
\end{equation}
\end{proposizione}
\begin{definizione}
Let $V$ be an indecomposable representation of $Q$.
\begin{itemize}
\item[(i)] $V$ is preprojective if and only if $(\tau^+)^iV=0$ for
$i>>0$.
\item[(ii)] $V$ is preinjective if and only if $(\tau^-)^iV=0$ for
$i>>0$.
\item[(iii)] $V$ is regular if and only if $(\tau^+)^iV\neq 0$
for every $i\in\mathbb{Z}$.
\end{itemize}
\end{definizione}
\begin{definizione}
Let $V$ be a representation of $Q$. The linear map
$$
\partial:\mathbb{N}^{Q_0}\longrightarrow\mathbb{Z}
$$
defined by $\partial(\underline{dim}\,V):=\langle
h,\underline{dim}\,V\rangle$ is called defect of $V$.
\end{definizione}
\begin{lemma}\label{difetto}
Let $V$ an indecomposable representation of $Q$. $V$ is
preprojective, preinjective or regular if and only if the defect
of $V$ is respectively negative, positive or zero.
\end{lemma}
The regular representations of $Q$ form an Abelian category
$Reg_{\mathbb{K}}(Q)$. Moreover $Reg_{\mathbb{K}}(Q)$ is serial,
i.e. every indecomposable regular representation has only one
regular composition series and so it is only determined by its
regular socle and by its regular length.
\begin{definizione}
A simple regular module $E$ is called homogeneous if and only if
$\underline{dim}\,E=h$.
\end{definizione}
\begin{proposizione}\label{orbitetau}
Let $Q$ be a quiver of tame type. Then there exist at most three
$\tau^+$-orbits $\Delta=\{e_i|\,i\in I=\{0,\ldots,u\}\}$,
$\Delta'=\{e'_i|\,i\in I'=\{0,\ldots,v\}\}$,
$\Delta''=\{e''_i|\,i\in I''=\{0,\ldots,w\}\}$, of dimension
vectors of non-homogeneous simple regular representations of $Q$
($I$, $I'$, $I''$ could be empty). We can assume that
$\tau^+(e_i)=e_{i+1}$ for $i\in I$ ($e_{u+1}=e_0$),
$\tau^+(e'_i)=e'_{i+1}$ for $i\in I'$ ($e'_{v+1}=e'_0$) and
$\tau^+(e''_i)=e''_{i+1}$ for $i\in I''$ ($e''_{w+1}=e_0$).
\end{proposizione}
We denote the set of all regular representations of $Q$ with
$\mathcal{D}_r$. Every vector $d\in\mathcal{D}_r$ can be
decomposed as
\begin{equation}\label{dcsolita}
d=ph+\sum_{i\in I}p_ie_i+\sum_{i\in I'}p'_ie'_i+\sum_{i\in
I''}p''_ie''_i
\end{equation}
for some $p,p_i,p'_i,p''_i\in\mathbb{N}$ such that at least one of
coefficients in each family $\{p_i|\,i\in I\}$, $\{p'_i|\,i\in
I'\}$, $\{p''_i|\,i\in I''\}$ is zero. The decomposition
(\ref{dcsolita}) is called \textit{canonical decomposition of
$d$}. It is unique because the only linear relations between $h$,
$e_i$, $e'_i$ and $e''_i$ are
$$
h=\sum_{i\in I}e_i=\sum_{i\in I'}e'_i=\sum_{i\in I''}e''_i.
$$
We observe that the category $Reg_{\mathbb{K}}(Q)$ can be
decomposed as direct sum of categories $\mathcal{R}_t$, with
$t=(\varphi,\psi)\in\mathbb{P}_1(\mathbb{K})$. In all categories
$\mathcal{R}_t$, but at most three of these, there is only one
simple object $V_t$ which is necessarily homogeneous.
\begin{definizione}\label{indregij}
(1) We call $E_i$, $E_i'$ and $E_i''$ the simple non-homogeneous
regular representations respectively of dimension $e_i$, $e_i'$
and $e_i''$.\\
(2) We call $V_{(\varphi,\psi)}$, where $(\varphi,\psi
)\in\mathbb{P}_1(\mathbb{K})$, the indecomposable regular
representation of dimension $h$.\\
(3) We define $E_{i,j}$ to be the indecomposable regular
representations with socle $E_i$ and dimension $\sum_{k=i}^j e_k$,
where $e_k$ are vertices of the arc with clockwise orientation
$\xymatrix@-1pc{e_i\ar@{-}[r]&\ar@{.}[r]&\ar@{-}[r]&e_j}$ in
$\Delta$, without repetitions of $e_k$. We denote $E_i:=E_{i,i}$
and similarly we define $E'_{i,j}$ and $E''_{i,j}$.
\end{definizione}
\begin{lemma}\label{eEe}
$$
\langle e_i,e_j\rangle=\left\{\begin{array}{ll}1&\textrm{if}\;i=j\\
-1&\textrm{if}\;i=j-1\\
0&\textrm{otherwise}.\end{array}\right.
$$
\end{lemma}
\textit{Proof.} By Schur's lemma, we have
$$
dim_{\mathbb{K}}(Hom_Q(E_i,E_j))=\left\{\begin{array}{ll}1&\textrm{if}\;i=j\\
0&\textrm{otherwise}.\end{array}\right.
$$
By [DR, lemma 3.3], we have $dim_{\mathbb{K}}(Ext^1_Q(E_i,E_j))=0$
for every $i\neq j-1$. So by the relation
$$
\langle
e_i,e_j\rangle=dim_{\mathbb{K}}(Hom_Q(E_i,E_j))-dim_{\mathbb{K}}(Ext^1_Q(E_i,E_j)),
$$
we obtain the thesis. $\Box$

\section{Reflection functors and Coxeter functors}\label{B2}
\begin{definizione}
Let $Q$ be a quiver.
\begin{itemize}
\item[a)] The vertex $x\in Q_0$ is a sink if there are no arrows $a\in Q_1$ such that $ta=x$.
\item[b)] The vertex $x\in Q_0$ is a source if there are no arrows $a\in Q_1$ such that $ha=x$.
\end{itemize}
\end{definizione}

Let $Q$ be a quiver and let $x\in Q_0$ be a sink (respectively a
source). We define the quiver $c_x(Q)$   in which the direction of
the arrows connecting to $x$ are reversed.
\begin{definizione}
Let $\{a_1,\ldots,a_k\}$ be the set of arrows in $Q$ whose head
(respectively tail) equals $x$. We put
$$
c_x(Q)_0=Q_0
$$
$$
c_x(Q)_1=\{c_x(a);a\in Q_1\}
$$
where $tc_x(a_i)=ha_i$,   $hc_x(a_i)=ta_i$ for every
$i\in\{1,\ldots,k\}$ and $tc_x(b)=tb$, $hc_x(b)=hb$ for every
$b\in Q_1\setminus\{a_1,\ldots,a_k\}$.
\end{definizione}
Now we define the functors $C_x^+$ and $C_x^-$ from $Rep(Q)$ to
$Rep(c_x(Q))$.
\begin{definizione}

 Let $Q$ be a quiver and $x\in Q_0$ be a sink. Let $\{a_1,\ldots,a_k\}$ be the set of arrows in $Q$ whose head equals $x$.  Let $V\in Rep(Q)$. We define the representation $C_x^+(V):=W\in Rep(c_x(Q))$ as follows.
$$
W(y)=\left\{\begin{array}{ll}
V(y) & if \quad x\ne y\\
Ker\big(\bigoplus_{i=1}^k V(ta_i)\stackrel{h}{ \longrightarrow}
V(x)\big) & otherwise,
\end{array}\right.
$$
where $h(v_1,\ldots,v_k)=V(a_1)(v_1)+\cdots+V(a_k)(v_k)$ with
$(v_1,\ldots,v_k)\in\bigoplus_{i=1}^k V(ta_i)$.
$$
W(c_x(a))=\left\{\begin{array}{ll}
V(a) & if\quad ha\ne x\\
W(x)\hookrightarrow \bigoplus_{i=1}^k
V(ta_i)\stackrel{p_j}{\longrightarrow} V(ta_j) & if\quad a=a_j
\end{array}\right.
$$
where $p_j$ denotes the projection on the $j$-th factor.
\end{definizione}
\begin{definizione}
 Let $Q$ be a quiver and $x\in Q_0$ be a source. Let $\{b_1,\ldots,b_l\}$ be the set of arrows in $Q$ whose tail equals $x$.  Let $V\in Rep(Q)$. We define the representation $C_x^-(V):=W\in Rep(c_x(Q))$ as follows.
$$
W(y)=\left\{\begin{array}{ll}
V(y) & if \quad x\ne y\\
Coker\big( V(x)\stackrel{\tilde{h}}{
\longrightarrow}\bigoplus_{i=1}^l V(hb_i)\big) & otherwise,
\end{array}\right.
$$
where $\tilde{h}(v)=(V(b_1)(v),\ldots,V(b_l)(v))$ with $v\in
V(x)$.
$$
W(c_x(a))=\left\{\begin{array}{ll}
V(a) & if\quad ta\ne x\\
V(hb_j)\stackrel{i_j}{\longrightarrow} \bigoplus_{i=1}^l
V(hb_i)\twoheadrightarrow W(x) & if\quad a=b_j
\end{array}\right.
$$
where $i_j$ denotes the immersion of the $j$-th factor.
\end{definizione}
Let $f=(f_y)_{y\in Q_0}:V\rightarrow W$ be a morphism in
$Rep(Q)$.\\
If $x$ is a sink and $\{a_1,\ldots,a_k\}$ is the set of arrows
whose head equals $x$, we define $C^+_xf=((C^+_xf)_y)_{y\in
Q_0}:C^+_x V\rightarrow C^+_x W$ a morphism in $Rep(c_xQ)$ as
follows. For every $y\neq x$, we have $f_y=(C^+_xf)_y$, whereas
$(C^+_xf)_x$ is the unique $\mathbb{K}$-linear map which makes the
diagram
$$
\begin{array}{ccccccc}
0 & \longrightarrow & (C^+_xV)_x & \longrightarrow &
\bigoplus_{i=1}^kV_{ta_i} & \stackrel{h}{\longrightarrow} & V_x\\
& & \downarrow\scriptstyle{(C^+_xf)_x} & & \downarrow\scriptstyle{\bigoplus_{i=1}^kf_{ta_i}} & & \downarrow\scriptstyle{f_x}\\
0 & \longrightarrow & (C^+_xW)_x & \longrightarrow &
\bigoplus_{i=1}^kW_{ta_i} & \stackrel{h'}{\longrightarrow} & W_x
\end{array}
$$
commutative.\\
If $x$ is a source and $\{b_1,\ldots,b_l\}$ is the set of arrows
whose tail equals $x$, we define $C^-_xf=((C^-_xf)_y)_{y\in
Q_0}:C^-_x V\rightarrow C^-_x W$ a morphism in $Rep(c_xQ)$ as
follows. For every $y\neq x$, we have $f_y=(C^-_xf)_y$, whereas
$(C^-_xf)_x$ is the unique $\mathbb{K}$-linear map which makes the
diagram
$$
\begin{array}{ccccccc}
 V_x & \stackrel{\tilde{h}}{\longrightarrow}
&
\bigoplus_{i=1}^lV_{tb_i} & \longrightarrow & (C^-_xV)_x & \longrightarrow & 0\\
 \downarrow\scriptstyle{f_x} & & \downarrow\scriptstyle{\bigoplus_{i=1}^lf_{tb_i}} & & \downarrow\scriptstyle{(C^-_xf)_x} & &\\
 W_x & \stackrel{\tilde{h'}}\longrightarrow &
\bigoplus_{i=1}^lW_{tb_i} & \longrightarrow & (C^-_xW)_x &
\longrightarrow & 0
\end{array}
$$
commutative.\\
In particular, by definition, we have $Hom(V,W)=0$ if and only
if\\
$Hom(C^+_xV,C^+_xW)=0$, with $x$ a sink and $Hom(V,W)=0$ if and
only if $Hom(C^-_xV,C^-_xW)=0$, with $x$ a source.\\
$C^+_x$, for every $x$ sink, and $C^-_x$, for every $x$ source,
are called \textit{reflection functors}.\\
We state the main result about reflection functors.
\begin{teorema}[Bernstein-Gelfand-Ponomarev]\label{BGP}
\begin{itemize}
\item[1)] Let $x\in Q_0$ be a sink. Let $V\in Rep(Q)$ be an indecomposable representation of dimension $\alpha$. Then we have two possibilities
\begin{itemize}
\item[a)] $V=S_x$ and then $C_x^+(V)=0$,
\item[b)] $C_x^+(V)$ is indecomposable and $C_x^-C_x^+(V)\cong V$ and the dimension of $C_x^+(V)$ equals $c_x(\alpha)$ where
$$
c_x(\alpha)(y)=\left\{\begin{array}{ll}
\alpha(y) & if\quad y\ne x\\
\sum_{i=1}^k \alpha(ta_i)-\alpha(x) & otherwise.
\end{array}\right.
$$
\end{itemize}
\item[2)] Let $x\in Q_0$ be a source. Let $V\in Rep(Q)$ be an indecomposable representation of dimension $\alpha$. Then we have two possibilities
\begin{itemize}
\item[a)] $V=S_x$ and then $C_x^-(V)=0$,
\item[b)] $C_x^-(V)$ is indecomposable and $C_x^+C_x^-(V)\cong V$ and the dimension of $C_x^-(V)$ equals $c_x(\alpha)$ where
$$
c_x(\alpha)(y)=\left\{\begin{array}{ll}
\alpha(y) & if\quad y\ne x\\
\sum_{i=1}^l \alpha(hb_i)-\alpha(x) & otherwise.
\end{array}\right.
$$
\end{itemize}
\item[3)] Let $V_1,V_2\in Rep(Q)$
$$
C_x^{\pm}(V_1\oplus V_2)=C_x^{\pm}(V_1)\oplus C_x^{\pm} (V_2).
$$
\end{itemize}

\end{teorema}
\textit{Proof.} See [BGP, theorem 1.1].\\
 \begin{definizione}
 A sequence $x_1,\ldots,x_m$ of vertices of $Q$ is an admissible sequence of sinks (respectively of sources) if $x_{i+1}$ is a sink (respectively a source) in $c_{x_i}\cdots c_{x_1}(Q)$ for $i=0,1,\ldots,m-1$.
 \end{definizione}
  \begin{cor}
  Let $Q$ be a quiver and let $x_1,\ldots,x_m$ be an admissible sequence of sinks.
  \begin{itemize}
  \item[1)] For every $i=1,\ldots,m$, $C_{x_1}^-\cdots C_{x_{i-1}}^-(S_{x_i})$ is either 0 or indecomposable (here $S_{x_i}\in Rep(c_{x_{i-1}}\cdots c_{x_1}(Q))$).
  \item[2)] Let $V\in Rep(Q)$ be an indecomposable. We assume $C_{x_k}\cdots C_{x_1}(V)=0$ for some $k$. Then there exists $i\in\{0,\ldots,k-1\}$ such that $V\cong C_{x_1}^-\cdots C_{x_{i-1}}^-(S_{x_i})$.
  \end{itemize}
  \end{cor}
  \textit{Proof.} Follows by induction from theorem 1.7.\\
\begin{definizione}
Let $Q$ be a quiver with $n$ vertices without oriented cycles. We
choose the numbering $(x_1,\ldots,x_n)$ of vertices such that
$ta>ha$ for every $a\in Q_1$. We define
$$
C^+:=C^+_{x_n}\cdots C^+_{x_1}\quad\textrm{and}\quad
C^-:=C^-_{x_1}\cdots C^-_{x_n}.
$$
The functors $C^+,C^-:Rep(Q)\rightarrow Rep(Q)$ are called Coxeter
functors. \end{definizione} These functors don't depend on the
choice of numbering of vertices because of the following
interpretation of the Coxeter functors in terms of the
Auslander-Reiten functors.
\begin{lemma}\label{C=tau}
Let $\mathbb{K}Q$ be the path algebra of a quiver $Q$ without
oriented cycles and $(x_1,\ldots,x_n)$ be an admissible numbering
of vertices.
\begin{itemize}
\item[(i)] If $V$ is an indecomposable nonprojective
$\mathbb{K}Q$-module, then there are isomorphisms
$C^+V\cong\tau^+V$ and $C^-C^+V\cong V$.
\item[(ii)] If $W$ is an indecomposable noninjective
$\mathbb{K}Q$-module, then there are isomorphisms
$C^-W\cong\tau^-W$ and $C^+C^-W\cong W$.
\end{itemize}
\end{lemma}
\textit{Proof.} See [ASS, chap. VII lemma 5.8]. $\Box$

\section{Semi-invariants of quivers without oriented
cycles}\label{B3} For a dimension vector $\alpha$ we have
$$
Rep(Q,\alpha):=\bigoplus_{a\in
Q_1}Hom(\mathbb{K}^{\alpha(ta)},\mathbb{K}^{\alpha(ha)}),
$$
the space of $\alpha$-dimensional representations of $Q$. Moreover
we define the group
$$
GL(Q,\alpha):=\prod_{x\in Q_0}GL(\mathbb{K},\alpha(x))
$$
and its subgroup
$$
SL(Q,\alpha):=\prod_{x\in Q_0}SL(\mathbb{K},\alpha(x)).
$$
 These groups act on $Rep(Q,\alpha)$ as follows: if $V\in Rep(Q,\alpha)$ and $g=(g_x)_{x\in Q_0}\in  GL(Q,\alpha)$, then $g\cdot V=\{g_{ha}V(a)g_{ta}^{-1}\}_{a\in Q_1}$. Finally we denote the ring of semi-invariants by
 $$
 SI(Q,\alpha):=\mathbb{K}[Rep(Q,\alpha)]^{SL(Q,\alpha)}=\{f\in Rep(Q,\alpha)|\forall g\in SL(Q,\alpha) g\cdot f=f\},
 $$
 where the action of $GL(Q,\alpha)$ on $\mathbb{K}[Rep(Q,\alpha)]$, the coordinate ring of polynomial functions on $Rep(Q,\alpha)$, is induced  by the action of $GL(Q,\alpha)$ on $Rep(Q,\alpha)$ by the rule
 $$
 (g\cdot f)(V):=f(g^{-1}\cdot V),
 $$
 with $g\in GL(Q,\alpha)$, $f\in \mathbb{K}[Rep(Q,\alpha)]$ and $V\in
 Rep(Q,\alpha)$.
 \begin{definizione}
 If $f$ is a semi-invariant of a quiver $Q$, we call $Z(f)$ the
 vanishing set of $f$.
 \end{definizione}
 \begin{lemma}\label{zeridif}
 Let $f$ and $f'$ be two semi-invariants of a quiver $Q$ such that $Z(f)=Z(f')$ is
 irreducible. Then $f=\lambda\cdot f'$ for some non zero
 $\lambda\in\mathbb{K}$.
 \end{lemma}
 \textit{Proof.} Since $Z(f)$ is irreducible, also $f$ is an
 irreducible polynomial. From $Z(f)=Z(f')$ it follows that $f'|f$
 and so $f=\lambda\cdot f'$ for some non zero
 $\lambda\in\mathbb{K}$. $\Box$.
 \begin{oss}\label{genirr}
 Let $\alpha$ be a dimension vector. Any set $S$ of generators of
 $SI(Q,\alpha)$ contains a subset of irreducible generators.
 Indeed if $f\in S$ is a reducible polynomial, then it can be expressed
 as a product of irreducible elements from $S$.
 \end{oss}
Now we define the semi-invariants which appear in the principal
  theorem.
\begin{lemma}\label{dVW}
The spaces $Hom_Q(V,W)$ and $Ext_Q^1(V,W)$ are respectively the
kernel and the cokernel of the following linear map
$$
d^V_W:\bigoplus_{x\in
Q_0}Hom(V(x),W(x))\longrightarrow\bigoplus_{a\in
Q_1}Hom(V(ta),W(ha))
$$
where $d^V_W$ is given by
$$
\{f(x)|x\in Q_0\}\longmapsto\{f(ha)V(a)-W(a)f(ta)|a\in Q_1\}.
$$
\end{lemma}
\textit{Proof.} See [R].\\
\\
If a representation $V$ has dimension vector $\alpha$, then
$d^V_W$ can be seen as the $\mathbb{K}$-linear map which sends
$\bigoplus_{x\in Q_0}W(x)^{\alpha(x)}$ to $\bigoplus_{a\in
Q_1}W(ha)^{\alpha(ta)}$.\\
For every representation $V$ of a quiver $Q$ without oriented
cycles of dimension $\alpha$, we can construct a projective
resolution, called \textit{Ringel resolution of} $V$:
\begin{equation}\label{Rr}
0\longrightarrow\bigoplus_{a\in Q_1}V(ta)\otimes
P_{ha}\stackrel{d^V}{\longrightarrow}\bigoplus_{x\in
Q_0}V(x)\otimes
P_{x}\stackrel{p_V}{\longrightarrow}V\longrightarrow 0
\end{equation}
where $P_x$ is the indecomposable projective associated to vertex
$x$ for every $x\in Q_0$ (see section B.1 of appendix), $d^V$
restricted to $V(ta)\otimes P_{ha}$ sends $v\otimes e_{ha}$ to
$V(a)(v)\otimes e_{ha}-v\otimes a$ and $p_V$ restricted to
$V(x)\otimes P_{x}$ sends $v$ to $v\otimes e_{x}$, see [R].
Moreover, applying the functor $Hom_Q(\cdot,W)$ to Ringel
resolution of $V$, we have $Hom_Q(d^V,W)=d^V_W$ for every
representation $W$ of $Q$.\\
Any character $\tau$ of $GL(Q,\alpha)$ has the form
 $$
 \tau:\{g_x\in GL(\alpha(x))|x\in Q_0\}\mapsto\prod_{x\in Q_0}(detg_x)^{\chi(e_x)}
$$
with $e_x$ a dimension vector, defined by $e_x(x)=1$ and
$e_x(y)=0$ if $x=y$, and $\chi(e_x)\in\mathbb{Z}\;\forall x\in
Q_0$. A vector $\chi\in\mathbb{Z}^{|Q_0|}$ is called
\textit{weight}.\\
The ring $SI(Q,\alpha)$ decomposes in graded components as
  $$
  SI(Q,\alpha)=\bigoplus_{\tau\in char(GL(Q,\alpha)}SI(Q,\alpha)_{\tau}
  $$
 where $SI(Q,\alpha)_{\tau} =\big\{f\in\mathbb{K}[Rep(Q,\alpha)]|g\cdot f=\tau(g)f\;\forall g\in GL(Q,\alpha)\big\}$.
 \begin{oss}
 \begin{itemize}
 \item[(1)] Each vector $\chi\in\mathbb{Z}^{|Q_0|}$ determines a unique character $\tau_{\chi}$.
\item[(2)] A character $\tau$ for some semi-invariant
might not uniquely determine the weight of the semi-invariant,
e.g. if $\alpha(x) = 0$, then $g_x$ is a $0 \times 0$ matrix, in
which case $det(g_x) = 1$, therefore for any
$\chi(x)\in\mathbb{Z}$, $det(g_x)^{\chi(x)} = det(g_x) = 1$.
\end{itemize}
 \end{oss}
If $\alpha$ and $\beta$ are dimension vectors such that $\langle
\alpha ,\beta \rangle=0$,
 $V\in Rep(Q,\alpha)$ and $W\in Rep(Q,\beta)$, then
 the matrix of $d^V_W$ is a square
matrix.
\begin{definizione}
We define the semi-invariant $c(V,W) := det\,d^V_W$ of the action
of $GL(Q,\alpha)\times GL(Q,\beta)$ on $Rep(Q,\alpha)\times
Rep(Q,\beta)$ (see [S]). For a fixed $V$ the restriction of $c$ to
$\{ V\}\times Rep(Q,\beta)$ defines a semi-invariant
$c^V=c(V,\cdot)$ in $SI(Q, \beta)$ of weight $\langle\alpha ,
\cdot\rangle$ [S, lemma 1.4]. Similarly, for a fixed $W$ the
restriction of $c$ to $Rep(Q,\alpha)\times\{ W\}$ defines a
semi-invariant $c_W=c(\cdot,W)$ in $SI(Q, \alpha)$ of weight
$-\langle\cdot, \beta\rangle$ [S, lemma 1.4]. These
semi-invariants are called Schofield semi-invariants.
\end{definizione}
These semi-invariants have the following properties.
\begin{lemma}\label{cVcV'}
Suppose that $V'$, $V$, $V''$ and $W'$, $W$, $W''$ are
representations of $Q$, that
$\langle\underline{dim}(V),\underline{dim}(W)\rangle=0$ and that
there are exact sequences
$$
0\rightarrow V'\rightarrow V\rightarrow V''\rightarrow
0,\qquad0\rightarrow W'\rightarrow W\rightarrow W''\rightarrow 0
$$
then
\begin{itemize}
\item[(i)] If $\langle\underline{dim}(V'),\underline{dim}(W)\rangle<0$, then $c^V(W)=0$
\item[(ii)] If $\langle\underline{dim}(V'),\underline{dim}(W)\rangle=0$, then $c^{V'}(W)=c^{V''}(W)
c^V(W)$
\item[(iii)] If $\langle\underline{dim}(V),\underline{dim}(W')\rangle>0$, then $c^V(W)=0$
\item[(iv)] If $\langle\underline{dim}(V),\underline{dim}(W')\rangle=0$, then $c^V(W)=c^V(W')c^V(W'')$
\end{itemize}
and similarly for $c_W$.
\end{lemma}
\textit{Proof.} See [DW1, lemma 1]. $\Box$
\begin{oss}
A consequence of lemma \ref{dVW} in [S] is that any projective
resolution of $V$ (respectively injective coresolution of $W$) can
be used to calculate $c^V$ (respectively $c_W$).So if $P$ is a
projective module and $I$ is an injective module then $c^P=0$ and
$c_I=0$.
\end{oss}
Now we formulate the result of Derksen and Weyman about the set of
generators of the ring of semi-invariants $SI(Q,\alpha)$, defined
in section 1.1, where $Q$ is a quiver without oriented cycles and
$\alpha$ is a dimension vector. So we assume throughout this
section that there are no oriented cycles in $Q$.
\begin{teorema}[Derksen-Weyman]\label{dw}
Let $Q$ be a quiver without oriented cycles and let $\beta$ be a
dimension vector. The ring $SI(Q,\beta)$ is spanned by
semi-invariants of the form $c^V$ of weight $\langle
\underline{dim}(V),\cdot\rangle$, for which $\langle
\underline{dim}(V),\beta\rangle=0$. It is also spanned by
semi-invariants of the form $c_W$ of weight $-\langle\cdot,
\underline{dim}(W)\rangle$,  for which $\langle \beta,
\underline{dim}(W)\rangle=0$.
\end{teorema}
\textit{Proof.} See [DW1, theorem 1]. $\Box$
\begin{oss}\label{cV=0}
If $\langle \underline{dim}(V), \underline{dim}(W)\rangle=0$ then
we have $c(V,W)=c^V(W)=c_W(V)=0$ if and only if $Hom_Q (V, W)\ne
0$ which is equivalent to $Ext_Q^1(V, W)\ne 0$ by lemma \ref{dVW}.
\end{oss}
\begin{oss}\label{cV+V'}
\begin{itemize}
\item[i)] If $V,V'\in Rep(Q)$ and $V\cong V'$ then $c^V$ and $c^{V'}$ are equal up to a scalar.
\item[ii)] If $V=V'\oplus V''$ is decomposable then, by lemma \ref{cVcV'}, we have $c^V=0$ in $SI(Q,\beta)$
if $\langle\underline{dim}(V'),\beta\rangle\ne 0$, and
$c^V=c^{V'}c^{V''}$ in $SI(Q,\beta)$ if
$\langle\underline{dim}(V'),\beta\rangle=0$.
\end{itemize}
So the algebra $SI(Q,\beta)$ is generated by all $c^V$ where $V$
is indecomposable and $\langle\underline{dim}\,V,\beta\rangle=0$.
\end{oss}
Moreover in [DW1] Derksen and Weyman show the following
\begin{cor}[Reciprocity]
Let $\alpha$ and $\beta$ be the dimension vectors satisfying
$\langle\alpha,\beta\rangle=0$. Then
$$
dim SI(Q,\beta)_{\langle \alpha,\cdot\rangle}= dim
SI(Q,\alpha)_{-\langle \cdot,\beta\rangle}.
$$
\end{cor}

\section{$c^V$, reflection functors and duality functor}\label{dVWrf}
The following results show the relation between $c^V$ and $C^+_x$
(respectively $C^-_x$).
\begin{lemma}\label{cVW=c+xVW}
Let $V$ be an indecomposable representation of $Q$ of dimension
$\alpha$ such that $Z(c^V)$ is irreducible and let $x$ be a sink
of $Q$. Then
$$
c^V=\lambda\cdot (c^{C^+_xV}\circ C^+_x)
$$
on $Rep(Q,\beta)$ such that $\langle\alpha,\beta\rangle=0$ and for
some non zero $\lambda\in\mathbb{K}$.
\end{lemma}
\textit{Proof.} First we note that, by remark \ref{genirr} and by
theorem \ref{dw}, it's not restrictive to suppose $Z(c^V)$ is
irreducible. By remark \ref{cV=0}, the vanishing set of $c^V$ is
the hypersurface
$$
Z(c^V)=\{W\in Rep(Q,\beta)|\, Hom_Q(V,W)\neq 0\}
$$
and the vanishing set of $c^{C^+_xV}$ is the hypersurface
$$
Z(c^{C^+_xV})=\{C^+_xW\in Rep(c_x(Q),c_x(\beta))|\,
Hom_Q(C^+_xV,C^+_xW)\neq 0\}.
$$
By definition of reflection functor, for every $W\in
Rep(Q,\beta)$,
$$
Hom_Q(V,W)\neq 0\Leftrightarrow Hom_Q(C^+_xV,C^+_xW)\neq 0.
$$
Hence $Z(c^V)=Z(c^{C^+_xV})$.\\
So, by lemma \ref{zeridif}, we conclude that there exist non zero
$\lambda\in\mathbb{K}$ such
that $c^V=\lambda\cdot (c^{C^+_xV}\circ C^+_x)$. $\Box$\\
\\
Similarly one proves the following
\begin{lemma}\label{cVW=c-xVW}
Let $V$ be an indecomposable representation of $Q$ of dimension
$\alpha$ such that $Z(c^V)$ is irreducible and let $x$ be a source
of $Q$. Then
$$
c^V=\lambda\cdot (c^{C^-_xV}\circ C^-_x)
$$
on $Rep(Q,\beta)$ such that $\langle\alpha,\beta\rangle=0$ and for
some non zero $\lambda\in\mathbb{K}$.
\end{lemma}
Next we study the relation between $c^V$ and duality functor
$\nabla$.
\begin{lemma}\label{cV=cVnabla1}
  Let $(Q,\sigma)$ be a symmetric quiver. For every representation
  $V$ of the underlying quiver $Q$ such that $Z(c^V)$ is irreducible, we have
  \begin{eqnarray}
  c^{V}=\lambda\circ (c^{\tau^-\nabla V}\circ\nabla)
  \end{eqnarray}
  for some non zero $\lambda\in\mathbb{K}$.
  \end{lemma}
  \textit{Proof.} First we note that, by remark \ref{genirr} and by
theorem \ref{dw}, it's not restrictive to suppose $Z(c^V)$ is
irreducible. Let $\beta$ be a dimension vector such that
  $\langle\underline{dim}\,V,\beta\rangle=0$. By equation (\ref{alphabetataunabla}) we note
  that, for every $W\in Rep(Q,\beta)$,
\begin{equation}\label{homtau}
Hom_Q(V,W)= 0\Leftrightarrow Hom_Q(\nabla W,\nabla V)=
0\Leftrightarrow Hom_Q(\tau^-\nabla V,\nabla W)=0.
\end{equation}
Thus, by remark \ref{cV=0}, the vanishing set of $c^V$ is the
hypersurface
$$
Z(c^V)=\{W\in Rep(Q,\beta)|\, Hom_Q(V,W)\neq 0\}
$$
and the vanishing set of $c^{\tau^-\nabla V}$ is the hypersurface
$$
Z(c^{\tau^-\nabla V})=\{\nabla W\in Rep(Q,\delta\beta)|\,
Hom_Q(\nabla W,\nabla V)\neq 0\}.
$$
Finally, by equation (\ref{homtau}), $Z(c^V)=Z(c^{\tau^-\nabla V})$.\\
So, by lemma \ref{zeridif}, we conclude that there exist non zero
$\lambda\in\mathbb{K}$ such that $c^V=\lambda\cdot
(c^{\tau^-\nabla V}\circ \nabla)$. $\Box$
\section{$c^{V}$'s, weights and partitions}
\begin{lemma}\label{SIcxQ}
Let $Q$ be a quiver, let $x$ be a sink and let $\alpha$ be a
vector dimension.
\begin{itemize}
\item[(i)] If $V$ is indecomposable not projective such that
$C_x^+V$ is not projective and
$0=\langle\underline{dim}V,\alpha\rangle(=\langle
c_x\underline{dim}V,c_x\alpha\rangle)$, then $c^V\in SI(Q,\alpha)$
and $c^{C_x^+V}\in SI(c_xQ,c_x\alpha)$.
\item[(ii)] If $V=S_x$ and
$\langle\underline{dim}S_x,c_x\alpha\rangle=0$, then we have
$c^V\in SI(c_xQ,c_x\alpha)$, where $S_x$ is considered as
representation of $c_xQ$, but $c^V$ is zero for $Q$.
\item[(iii)] If $V=C^-S_x$ and
$\langle\underline{dim}C^-S_x,\alpha\rangle=0$, then we have
$c^V\in SI(Q,\alpha)$ but $c^{C^+_xV}$ is zero for $c_xQ$.
\end{itemize}
\end{lemma}
\textit{Proof.} First of all we observe that if $x$ is a sink and
$V\neq S_x$ is projective then $C^+_xV$ is projective since $C^+$
doesn't depend on any admissible numbering of vertices. Moreover
$\langle\underline{dim}S_x,c_x\alpha\rangle=0$ and
$\langle\underline{dim}C^-S_x,\alpha\rangle=0$
 are not both zero. By theorem \ref{ard} and
since $x$ is a sink,
$0=\langle\underline{dim}C^-S_x,\alpha\rangle=-\langle\alpha,\underline{dim}S_x\rangle=-\alpha_x+\sum_{a\in
Q_1 :ha=x}\alpha_{ta}$ and
$0=\langle\underline{dim}S_x,c_x\alpha\rangle=(c_x\alpha)_x-\sum_{a\in
c_x(Q)_1}(c_x\alpha)_{ha}=\sum_{a\in Q_1
:ha=x}\alpha_{ta}-\alpha_x-\sum_{a\in Q_1
:ha=x}\alpha_{ta}=-\alpha_x$ and so $\sum_{a\in Q_1
:ha=x}\alpha_{ta}=0$ which is an absurd unless $\alpha_{ta}=0$ for
every $a$ such that $ha=x$ but in such case $c^{S_x}=0$ for $c_xQ$
and $c^{C^-S_x}=0$ for $Q$.\\
\textit{Proof of (i)}. Since
$\langle\underline{dim}V,\alpha\rangle=0$, by theorem \ref{dw},
the $c^V$'s are generators of $SI(Q,\alpha)$ and $c^{C^+_xV}$'s
are generators of $SI(c_xQ,c_x\alpha)$.
Moreover we note that the number of generators of $SI(Q,\alpha)$ is equal to the number of generators of $SI(c_xQ,c_x\alpha)$.\\
\textit{Proof of (ii)}. We can study $S_x$ since if $V\neq S_x$ is
projective, by remark above, we have $c^V=0$ and also $c^{C^+_x
V}=0$. $S_x$ is projective in $Q$ and so $c^{S_x}$ is zero in
$SI(Q,\alpha)$ but $S_x$, considered as a representation of
$c_xQ$, is injective. So, if
$\langle\underline{dim}S_x,c_x\alpha\rangle=0$ then $c^{S_x}\in
SI(c_xQ,c_x\alpha)$.\\
\textit{Proof of (iii)}. $C^-S_x$ is not projective otherwise
$S_x=C^+(C^-S_x)=0$ which is an absurd. Thus if
$\langle\underline{dim}C^-S_x,\alpha\rangle=0$ then $c^{C^-S_x}\in
SI(Q,\alpha)$. Moreover $C^+C^+_xC^-S_x=C^+_xC^+C^-S_x=C^+_xS_x=0$
hence $C^+_xC^-S_x$ is projective and so $c^{C^+_xC^-S_x}=0$ in
$SI(c_xQ,c_x\alpha)$. $\Box$\\
\\
We recall that if $Q$ is Dynkin, then $SI(Q,\alpha)$ has a finite
number of generators by remark \ref{cV+V'}.
\begin{cor}\label{cSIcxQ}
Let $Q$ be a Dynkin quiver and let $x$ be a sink. We call
$N(Q,\alpha)$ the number of generators of $SI(Q,\alpha)$ and
$N(c_xQ,c_x\alpha)$ the number of generators of
$SI(c_xQ,c_x\alpha)$. We have three possibilities.
\begin{itemize}
\item[(a)] $N(Q,\alpha)=N(c_xQ,c_x\alpha)$\; if
$\langle\underline{dim}S_x,c_x\alpha\rangle\neq 0$ and
$\langle\underline{dim}C^-S_x,\alpha\rangle\neq 0$;
\item[(b)] $N(Q,\alpha)+1=N(c_xQ,c_x\alpha)$\; if
$\langle\underline{dim}S_x,c_x\alpha\rangle=0$;
\item[(c)] $N(c_xQ,c_x\alpha)+1=N(Q,\alpha)$\; if
$\langle\underline{dim}C^-S_x,\alpha\rangle=0$.
\end{itemize}
\end{cor}
\textit{Proof.} (a) follows directly from (i) of the previous
lemma. (b): the generators of $SI(c_xQ,c_x\alpha)$ are those of
$SI(Q,\alpha)$ and $c^{S_x}$. (c): the generators of
$SI(Q,\alpha)$ are those of $SI(c_xQ,c_x\alpha)$ and $c^{C^-S_x}$.
$\Box$\\
\\
Now we study weights of a quiver $A_n$ and associated partitions.
We denote vertices of $A_n$ with $\{1,\ldots, n\}$ in increasing
way from left to right and we call $a_i$ the arrow which has $i$
on the left and $i+1$ on the right. Let $V_{i,j}$ be the
indecomposable of $A_n$ with dimension vector
$$
(v_{i,j})_h=\left\{\begin{array}{cc} 1 &
\textrm{if}\quad i\leq h \leq j\\
0& \textrm{otherwise}.\end{array}\right.
$$
Let $E=(E_{i,j})_{1\leq i,j \leq n}$ be the Euler matrix of a
quiver $Q$ , i.e the matrix associated to the Euler form
$\langle\cdot,\cdot\rangle$. In general we have
$$
E_{i,j}=\left\{\begin{array}{ll} 1 & \textrm{if}\;i=j\\
\sharp\{a\in Q_1|ta=i,ha=j\} & \textrm{otherwise}.
\end{array}\right.
$$
If $Q=A_n$
$$
E_{i,j}=\left\{\begin{array}{ll} 1 & \textrm{if}\;i=j\\
-1 & \textrm{if}\; i\rightarrow j\\
0 & \textrm{otherwise}.
\end{array}\right.
$$
Let $\langle v_{i,j},\cdot\rangle=v_{i,j}E=\chi=(\chi_l)_{1\leq
l\leq n}$ be the weight of
$c^{V_{i,j}}$.\\
 We consider the following notation for $A_n$, let $s,p\geq 1$ be respectively the number
 of sources and the number of sinks in $A_n$ (there are at least one
source and one sink, which occurs in the equioriented case).
$$
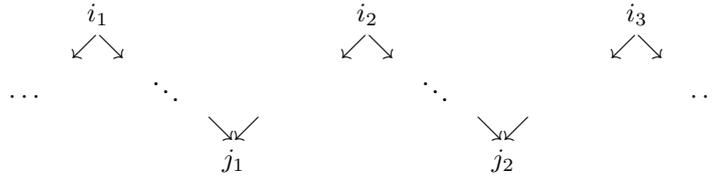
\begin{array}{cccccccccccccccc}
& i_1 & & & & & i_2 & & & & & i_3 &\\
 & \swarrow  \searrow & & & & &  \swarrow  \searrow& & & & & \swarrow  \searrow
 &\\
 \cdots &  & \ddots & & &  & & \ddots & & & & & \cdots\\
  & & &\searrow   \swarrow & & & & &\searrow   \swarrow& & & &\\
  & & & j_1 & & & & & j_2 & & & &
 \end{array}
 $$
where $i_k$ and $j_h$ in $\{1,\ldots,n\}$ with $1\leq k\leq s$ and
$1\leq h\leq p$ are respectively sources and sinks of $Q$. By the
previous picture we note that in $A_n$ sinks and sources
alternate.\\
Let $K=\{k\in\{1,\ldots,s\}|i\leq i_k\leq j\}$ and
$H=\{h\in\{1,\ldots,p\}|i\leq j_h\leq j\}$
\begin{lemma}
The weight of $c^{V_{i,j}}$ is
$\chi=(\chi_l)_{l\in\{1,\ldots,n\}}$ such that
$$
\chi_l=\left\{\begin{array}{ll} 1 & l=i_k\;\textrm{with}\; k\in
K\; \textrm{or}\; l=i\;\textrm{and}\quad ta_i=i\;\textrm{or}\;
l=j\;\textrm{and}\; ta_{j-1}=j\\
-1 & l=j_h\;\textrm{with}\; h\in H\; \textrm{or}\;
l=i-1\;\textrm{and}\; ha_{i-1}=i-1\;\textrm{or}\;
l=j+1\;\textrm{and}\; ha_{j}=j+1\\
0 & \textrm{otherwise}.\end{array}\right.
$$
\end{lemma}
\textit{Proof}. Since $v_{i,j}E=\chi=(\chi_l)_{1\leq l\leq n}$ is
the weight of $c^{V_{ij}}$ then
$\chi_l=E_{i,l}+E_{i+i,l}+\cdots+E_{j,l}$ for every
$l\in\{1,\ldots,n\}$. So
$$
\chi_l=\left\{\begin{array}{ll} E_{l-1,l}+E_{l,l}+E_{l+1,l} &
l\in\{i+1,\ldots,j-1\}\\
E_{l+1,l} & l=i-1\\
E_{l-1,l} & l=j+1\\
E_{l,l}+E_{l+1,l} & l=i\\
E_{l-1,l}+E_{l,l} & l=j\\
0 & \textrm{otherwise}.
\end{array}\right.
$$
Hence $\chi_{l}=0$ for every $l\in
\{1,\ldots,i-2\}\cup\{j+2,\ldots,n\}$,
$$
\chi_{i-1}=\left\{\begin{array}{ll} -1 & i-1\leftarrow i\\
0 & \textrm{otherwise},
\end{array}\right.
$$
$$
\chi_{j+1}=\left\{\begin{array}{ll} -1 & j\rightarrow j+1\\
0 & \textrm{otherwise},
\end{array}\right.
$$
$$
\chi_{i}=\left\{\begin{array}{ll} 1 & i\rightarrow i+1\\
0 & \textrm{otherwise},
\end{array}\right.
$$
$$
\chi_{j}=\left\{\begin{array}{ll} 1 & j-1\leftarrow j\\
0 & \textrm{otherwise}
\end{array}\right.
$$
and for every $l\in\{i+1,\ldots,j-1\}$
$$
\chi_{l}=\left\{\begin{array}{ll} 1 & l-1\leftarrow l\rightarrow l+1\\
-1 & l-1\rightarrow l\leftarrow l+1\\
0 & \textrm{otherwise}.\quad\quad\quad\quad\Box
\end{array}\right.
$$
\begin{cor}
Let $Q=A_n$ and let $w$ be the weight of $c^{V_{i,j}}$.
\begin{itemize}
\item[(i)] Let $\chi_l=1$ for some $l\in\{i,\ldots,j\}$ and let $k>l$
in $\{i+1,\ldots,j-1\}\cup\{j+1\}$ be the first index such that
$\chi_k\neq 0$, then $\chi_k=-1$.
\item[(ii)] Let $\chi_l=-1$ for some
$l\in\{i+1,\ldots,j-1\}\cup\{i-1,j+1\}$ and let $k>l$ in
$\{i,\ldots,j\}$ be the first index such that $\chi_k\neq 0$, then
$\chi_k=1$. $\Box$
\end{itemize}
\end{cor}
Let $\beta$ be the dimension vector of an indecomposable
representation of $A_n$ and let $\chi=\langle\beta,\cdot\rangle$.
Let $m_1$ be the first vertex such that $\chi(m_1)\neq 0$, in
particular we suppose $\chi(m_1)=1$ and $m_t$ the last vertex such
that $\chi(m_t)\neq 0$, in particular we suppose $\chi(m_t)=1$,
the other case proves in a similar way. Between $m_1$ and $m_t$,
-1 and 1 alternate in correspondence respectively to sinks and to
sources. In this case we have $[\frac t 2]+1=s+1$ occurrences of 1
and $s=[\frac t 2]$ occurrences of -1. We call $i_0=m_1$,
$j_s=m_{t-1}$, $i_1,\ldots,i_{s}$ the sources and
$j_1,\ldots,j_{s-1=p}$ the sinks between $i_0$ and $j_s$. Let $V$
be a representation with $\underline{dim}\,V=\alpha$ such that
$\langle\beta,\alpha\rangle=0$ and
$SL(V)=SL(V_1)\times\cdots\times SL(V_n)$, so we have, by Cauchy
formula
$$
\mathbb{K}[Rep(A_n,\alpha)]^{SL(V)}=SI(A_n,\alpha)=\left(\bigoplus_{\lambda:Q_1\rightarrow\Lambda}\bigotimes_{c\in
Q_1}S_{\lambda(c)}V_{tc}\otimes
S_{\lambda(c)}V_{hc}^*\right)^{SL(V)}
$$
where $\Lambda$ is the set of all partitions.\\
$\chi(k)=0$ for every $k<i_0$ so either
$\lambda(a_{k-1})=\lambda(a_k)$ or
$\lambda(a_{k-1})=0=\lambda(a_k)$ for every $k<i_0$. Since
$\chi(1)=0$ then $\lambda(a_1)=0$ and thus $\lambda(a_k)=0$ for
every $k<i_0$. So we have
$(S_{\lambda(a_{i_0})}V_{i_0})^{SLV_{i_0}}\neq 0$ if and only if
$\lambda(a_{i_0})=(\overbrace{1,\ldots,1}^{\alpha_{a_{i_0}}})$.
Now $\chi(k)=0$ for every $i_0<k<j_1$ and $\chi(j_1)=-1$ then we
have $\lambda(a_{i_0+1})=\lambda(a_{i_0})$ otherwise
$(S_{\lambda(a_{i_0})}V_{i_0+1}^*\otimes
S_{\lambda(a_{i_0+1})}V_{i_0+1})^{SLV_{i_0+1}}$ doesn't have
weight 0. So $\lambda(a_k)=\lambda(a_{i_0})$ for every
$i_0<k<j_1$. For $j_1$ we have $\lambda(a_{j_1})$ and
$\lambda(a_{i_0})$ are complementary with respect to a column of
height $\alpha_{j_1}$ because
$-\lambda(a_{j_1})_h-\lambda(a_{i_0})_{\alpha_{j_1}-h+1}=-1$ for
every $\in\{1,\ldots,\alpha_{j_1}\}$, by proposition \ref{i3}. We
proceed in a similar way with the other vertices until $i_{s}$ for
which $\chi(i_{s})=1$. Since $\chi(k)=0$ for every $k>i_{s}$, we
have either $\lambda(a_{k-1})=\lambda(a_k)$ or
$\lambda(a_{k-1})=0=\lambda(a_k)$ for every $k>i_s$ but because
$\lambda(a_{n-1})=0$, $\lambda(a_k)=0$ for every $k>i_s$. Moreover
$\lambda(a_{i_{s}-1})$ is both a column of height $\alpha_{i_{s}}$
and the complementary of $\lambda(a_{i_{s-1}-1})$ with respect to
a column of height
$\alpha_{i_{s-1}}$.\\
So we proved the following
\begin{lemma}
Let $Q$ be a quiver of type $A_n$, let $\alpha$ be a dimension
vector and $\beta$ be a dimension vector of an indecomposable
representation of $Q$. Let $\chi$ be the weight
$\langle\beta,\cdot\rangle$ and we suppose it is such that
$\chi(i)\neq 0$ for every $i\in I=\{m_j\}_{j\in\{1,\ldots,t\}}$,
where $I$ is a subset of $\{1,\ldots,n\}$. Then the family of
partitions associated to $\chi$ is
$\underline\lambda=(\lambda(a_1),\ldots,\lambda(a_{n-1}))$ such
that $\lambda(a_i)=0$ for every
$i\in\{1,\ldots,m_1-1\}\cup\{m_t,\ldots,n-1\}$, $\lambda(a_{m_1})$
and $\lambda(a_{m_t-1})$ are columns respectively of height
$\alpha_{m_1}$ and $\alpha_{m_t}$ and $\lambda(a_{i})$ is the
complementary of $\lambda(a_{i-1})$ with respect to a column of
height $\alpha_i$ for every $i\in\{m_j\}_{j\in\{2,\ldots,t-1\}}$.
Moreover we have
$\alpha_{m_t}=\alpha_{m_{t-1}}-\alpha_{m_{t-2}}+\ldots\pm
\alpha_{m_{1}}$.
\end{lemma}

 \chapter*{References}\addcontentsline{toc}{chapter}{References}
  \begin{itemize}
  \item[[ARS\!\!\!]] M. Auslander, I. Reiten, S. O. Smal\o,
  \textit{Representation Theory of Artin Algebras}, Cambridge
  Studies in Advanced Mathematics 36, Cambridge University Press,
  1995.

\item[[ASS\!\!\!]] I. Assem, D. Simson, A. Skowronski, \textit{Elements of the
Representation Theory of Associative Algebras}, volume 1, London
  Mathematical Society Students Texts 65, Cambridge University Press, 2006.
\item[[B\!\!\!]] A. Borel, \textit{Linear Algebraic Groups}, 2nd
  ed., Graduate Texts in Mathematics 126, Springer-Verlag, 1991.
  \item[[Bo1\!\!\!]] K. Bongartz, \textit{Degenerations for representations of tame quivers}, Ann. Sci. \'Ecole Normale
Sup. 28 (1995), 647-668.
\item[[Bo2\!\!\!]] K. Bongartz, \textit{On degenerations and extensions of finite dimensional modules}, Advances
Math. 121 (1996), 245-287.
  \item[[BGP\!\!\!]] I. N. Bernstein, I. M. Gelfand, V. A. Ponomarev, \textit{Coxeter functors and Gabriel's theorem}, Uspekhi Mat. Nauk 28, no. 2(170) (1973), 19-33.
\item[[BMRRT\!\!\!]] A. B. Buan, R. Marsh, M. Reineke, I. Reiten, G. Todorov, \textit{Tilting theory and cluster combinatorics},
Adv. Math. 204 (2006), 572-618.
  \item[[D\!\!\!]] C. Di Trapano, \textit{Algebra of
  Semi-invariants
of Euclidean Quivers}, preprint, 2009.
\item[[DP\!\!\!]] C. De Concini, C. Procesi \textit{A
characteristic free approach to invariant theory}, Adv. in Math.
21 (1976), No. 3, 330-354.
\item[[DR\!\!\!]] V. Dlab, C. M. Ringel, \textit{Indecomposable Representations of Graphs and Algebras}, Memoirs
Amer. Math. Soc. 173 (1976).
\item[[DSW\!\!\!]] H. Derksen, A. Schofield, J. Weyman, \textit{On the number of subrepresentations of generic representations of quivers},
math.AG/050739.
  \item[[DW1\!\!\!]] H. Derksen, J. Weyman, \textit{Semi-invariants of quivers and saturation for Littlewood-Richardson coefficients}, J. Amer. Math. Soc. 16 (2000),467-479.
  \item[[DW2\!\!\!]] H. Derksen, J. Weyman, \textit{Generalized
  quivers associated to reductive groups}, Colloq. Math. 94
  (2002), No. 2, 151-173.
  \item[[DW3\!\!\!]] H. Derksen, J. Weyman, \textit{On the
  Littlewood-Richardson polynomials}, J. Algebra 255 (2002),
  247-257.
  \item[[DW4\!\!\!]] H. Derksen, J. Weyman, \textit{On canonical decomposition for quiver representations}, Comp. Math., 133, 245-
265 (2002).
\item[[DW5\!\!\!]] H. Derksen, J. Weyman, \textit{Semi--invariants for quivers with relations}, J. Algebra 258 (2002), 216-227.
\item[[DZ\!\!\!]] M. Domokos, A. N. Zubkov, \textit{Semi-invariants of quivers as determinants},
Transform. Groups 6 (2001), No. 1, 9-24.
  \item[[F\!\!\!]] W. Fulton, \textit{Young tableaux, with
  applications to representation theory and geometry}, London
  Mathematical Society Student Texts 35, Cambridge University
  Press, 1997.
  \item[[FH\!\!\!]] W. Fulton, J. Harris, \textit{Representation Theory; the first course},
Graduate Texts in Mathematics 129, Springer-Verlag, 1991.
\item[[FZ1\!\!\!]] S. Fomin, A. Zelevinsky, \textit{Cluster Algebras I: Foundations}, J. Amer. Math. Soc. 15 (2002),
497-529.
\item[[FZ2\!\!\!]] S. Fomin, A. Zelevinsky, \textit{Cluster Algebras II: Finite type classification}, Invent. Math. 154
(2003), 63-121.
\item[[G\!\!\!]] F. Gavarini, \textit{A Brauer algebra theoretic proof of Littlewood's restriction rules},
J. Algebra 212 (1999), No. 1, 240-271.
\item[[GM\!\!\!]] S. I. Gelfand, Y. I. Manin, \textit{Methods of
Homological Algebra}, Springer Monographs in Mathematics,
Springer-Verlag, 2003.
  \item[[GW\!\!\!]] R. Goodman, N. R. Wallach,
  \textit{Representations and Invariants of the Classical Groups},
  Encyclopedia of Mathematics and its Applications 68, Cambridge
  University Press, 1998.
  \item[[Hu\!\!\!]] J. E. Humphreys, \textit{Linear Algebraic
  Groups},
Graduate Texts in Mathematics 21, Springer-Verlag, 1975.
\item[[IOTW\!\!\!]] K. Igusa, K. Orr, G. Todorov, and J. Weyman, \textit{Cluster complexes via semi-invariants}, Preprint,
arXiv:0708.0798v1, 2007.
\item[[K1\!\!\!]] V. G. Kac, \textit{Infinite root systems~ representations of
graphs and invariant theory}, Invent. Math. 56 (1980), 57-92.
\item[[K2\!\!\!]] V. G. Kac, \textit{Infinite root systems, representations of
graphs and invariant theory II}, J. AIgebra 8 (1982), 141-162.
\item[[KP\!\!\!]] H. Kraft, C. Procesi, \textit{Classical
Invariant Theory. A Primer},\\
http://www.math.unibas.ch/~kraft/Papers/KP-Primer.pdf.
\item[[KR\!\!\!]] H. Kraft, C. Riedtmann, \textit{Geometry of representations of
quivers}, in: Representations of Algebras, London Math. Society
Lecture Notes Series 116, Cambridge Univ. Press, 1986, 109-145.
\item[[La\!\!\!]] S. Lang, \textit{Algebra}, Graduate Texts in
Mathematics 211, Springer-Verlag, 2002.
\item[[LP\!\!\!]] L. Le Bruyn, C. Procesi \textit{Semisimple representations of quivers}, Trans. Amer. Math. Soc. 317 (1990), 585-598.
\item[[Lo\!\!\!]] A. A. Lopatin, \textit{Invariants of quivers under the action of classical groups}, J. Algebra 321 (2009), no. 4, 1079--1106.
 \item[[LoZ\!\!\!]] A. A. Lopatin, A. N. Zubkov, \textit{Semi-invariants of mixed representations of quivers}, Transform. Groups 12 (2007), no. 2, 341-369.
\item[[L\!\!\!]] S. Lovett, \textit{Orbits of orthogonal and symplectic representations of symmetric
quivers}, http://www.enc.edu/~slovett/math/papers/symquiv.pdf.
  \item[[M\!\!\!]] I. G. Macdonald, \textit{Symmetric functions and Hall polynomials},
second edition, with contributions by A. Zelevinsky, Oxford
Mathematical Monographs, Oxford University Press, New York, 1995.
\item[[MWZ1\!\!\!]] P. Magyar, J. Weyman and A. Zelevinsky, \textit{Multiple flag varieties of finite type}, Adv. Math. 141 (1999), no. 1, 97-118.
\item[[MWZ2\!\!\!]] P. Magyar, J. Weyman and A. Zelevinsky, \textit{Symplectic multiple flag varieties of finite type},  J. Algebra 230 (2000), no. 1, 245-265.
\item[[P\!\!\!]] C. Procesi, \textit{Lie Groups. An Approach through Invariants and Representations}, Universitext, Springer, New York, 2007.
  \item[[P1\!\!\!]] C. Procesi, \textit{The invariant theory of n x n matrices}, Adv. in Math. 19 (1976), 306-381.
  \item[[R1\!\!\!]] C. M, Ringel, \textit{Representations of K-species and bimodules}, J. Algebra 41 (1976), 269-302.
  \item[[R2\!\!\!]] C.M., Ringel, \textit{The rational invariants
  of the tame quivers}, Invent. Math. 58 (1980), 217-239.
  \item[[R3\!\!\!]] C. M. Ringel, \textit{Tame algebras and integral quadratic forms}, Lecture
Notes in Mathematics 1099, Springer, 1984.
  \item[[S\!\!\!]] A. Schofield, \textit{Semi-invariants of quivers}, J. London Math. Soc. (3) 65 (1992), 46-64.
  \item[[SK\!\!\!]] M. Sato, T. Kimura, \textit{A classification of irreducible prehomogeneous vector spaces and
their relative invariants}, Nagoya J. Math 65 (1977), 1-155.
\item[[Sp\!\!\!]] T. A. Springer, \textit{Invariant Theory},
Lecture Notes in Mathematics 585, Springer-Verlag, Berlin,
Heidelberg, New York, 1977.
\item[[SV\!\!\!]] A. Schofield , M. Van den Bergh, \textit{Semi-invariants of quivers for arbitrary dimension vectors}, Indag. Math. (N.S.) 12 (2001), 125-138.
  \item[[SW\!\!\!]] A. Skowronski, J. Weyman, \textit{The algebras
  of semi-invariants of quivers}, Trans. Groups 5 (2000), no. 4,
  361-402.
  \item[[W\!\!\!]] J. Weyman, \textit{Cohomology of Vector Bundles
  and Syzygies}, Cambridge Tracts in Mathematics 149, Cambridge
  University Press, 2003.
  \item[[Z\!\!\!]] G. Zwara, \textit{Degenerations for representations of extended Dynkin quivers}, Comment. Math.
Helvetici 73 (1998), 71-88.
  \end{itemize}

\end{document}